\documentclass{amsart}
\usepackage{amsmath,amsthm,amssymb,amsfonts}
\numberwithin{equation}{section}

\newcommand{\supp}{\mathrm{supp}}
\newcommand{\etap}{\tilde{\eta}}

\newcommand{\unif}{\varpi}

\newcommand{\diff}{\mathfrak{d}}
\newcommand{\unift}{\varpi}
\newcommand{\sign}{\mathrm{sign}}

\newcommand{\Tr}{\mathrm{Tr}}
\newcommand{\p}{\mathfrak{p}}

\newcommand{\Ad}{\mathrm{Ad}}

\newcommand{\Zz}{\mathbf{Z}}

\newcommand{\dK}{\mathbf{d}}
\newcommand{\dE}{\mathbf{d}}

\newcommand{\heckeq}{\overline{\mu}_{\q}}
\newcommand{\heckef}{\overline{\mu}_{\fcond}}
\newcommand{\Kmax}{K_{\max}}
\newcommand{\Kmaxg}{K_{\max,{\G}}}
\newcommand{\q}{\mathfrak{q}}
\newcommand{\Lie}{\mathrm{Lie}}
\newcommand{\fcond}{\mathfrak{f}}
\newcommand{\quotg}{\mathbf{X}_{\mathrm{GL}(2)}}

\newcommand{\C}{\mathbb{C}}
\newcommand{\quotG}{\quot_{\G}}
\newcommand{\quotGad}{\quot_{\G,\mathrm{ad}}}
\newcommand{\Siegel}{\Sieg}
\newcommand{\weight}{\mathrm{wt}}
\newcommand{\iso}{\Upsilon}
\newcommand{\heegmap}{\mathscr{H}}

\newcommand{\vol}{\mathrm{vol}}
\newcommand{\ord}{\mathrm{ord}}
\newcommand{\G}{\mathbf{G}}
\newcommand{\lieg}{\mathfrak{g}}

\renewcommand{\l}{\mathfrak{l}}
\newcommand{\Cond}{\mathrm{Cond}}
\newcommand{\set}{\mathcal{S}}
\newcommand{\meas}{\sigma}
\newcommand{\basis}{\mathcal{B}}
\newcommand{\nbar}{\bar{n}}
\newcommand{\Sieg}{\mathfrak{S}}
\newcommand{\height}{\mathrm{ht}}

\renewcommand{\d}{\mathfrak{d}}
\newcommand{\m}{\mathfrak{m}}
\newcommand{\GL}{\mathrm{GL}}
\newcommand{\Supp}{\mathrm{Supp}}
\newcommand{\PGL}{\mathrm{PGL}}
\newcommand{\f}{a([\mathfrak{p}])}

\newcommand{\SL}{\mathrm{SL}}

\newcommand{\kirill}{\mathscr{K}}
\newcommand{\order}{\mathfrak{o}}

\newcommand{\R}{\mathbb{R}}
\renewcommand{\H}{\mathbb{H}}
\newcommand{\n}{\mathfrak{n}}

\newcommand{\Z}{\mathbb{Z}}
\newcommand{\Q}{\mathbb{Q}}

\newcommand{\meast}{\meas^{(2)}}
\newcommand{\adele}{\mathbb{A}}
\newcommand{\Norm}{\mathrm{N}}

\newcommand{\cond}{\mathfrak{p}}

\newcommand{\quot}{\mathbf{X}}
\DeclareFontFamily{OT1}{rsfs}{}
\DeclareFontShape{OT1}{rsfs}{n}{it}{<-> rsfs10}{}
\DeclareMathAlphabet{\mathscr}{OT1}{rsfs}{n}{it}
\newtheorem{ex}{Example}[section]
\newtheorem{lem}{Lemma}[section]
\newtheorem{defn}{Definition}[section]
\newtheorem{prop}{Proposition}[section]
\newtheorem{thm}{Theorem}[section]
\newtheorem{cor}{Corollary}[section]
\newtheorem{hyp}{Hypothesis}[section]
\theoremstyle{remark}
\newtheorem{rem}{Remark}[section]
\begin{document}
\title{Sparse equidistribution problems,
period bounds, and subconvexity}
\author{Akshay Venkatesh}

\begin{abstract}
We introduce a ``geometric'' method to bound
periods of automorphic forms. The key features of this method
are the use of equidistribution results in place
of mean value theorems, and the systematic
use of mixing and the spectral gap. Applications are given
to equidistribution
of sparse subsets of horocycles and to equidistribution of CM points;
to subconvexity 
 of the triple product period in the level aspect over number fields,
 which implies subconvexity for certain standard and Rankin-Selberg $L$-functions;
and to bounding Fourier coefficients of automorphic forms.
\end{abstract}
\maketitle
\tableofcontents

 \section{Introduction.} \label{sec:intro}
\subsection{General introduction.} \label{subsec:intro}
Let $\Gamma \subset G$ be a lattice in an $S$-arithmetic group.
Let $Y \subset \Gamma \backslash G$ be a
subset endowed with a probability measure $\nu$,
and $f$ a function on $\Gamma \backslash G$.
Fixing a basis $\{\psi^{(Y)}_j\}$ for $L^2(Y,\nu)$,
we shall refer to the numbers $\int f \psi^{(Y)}_j d\nu$,
as the {\em periods} of $f$ along $Y$.
Evidently, the periods depend heavily on the choice of basis
for $L^2(Y,\nu)$.
They play a major role in the theory of automorphic forms,
in significant part because they often express
information about $L$-functions.

The present paper is centered around a geometric method yielding upper
bounds for
these periods.
It is applicable, roughly speaking, when considering
the periods of a fixed function $f$ along a sequence of subsets
$(Y_i,\nu_i)$, with the property that the $Y_i$
are becoming equidistributed; that is to say, the $\nu_i$
approach weakly the $G$-invariant measure on $\Gamma \backslash G$.
The key inputs of this method are, firstly,
the {\em equidistribution} of the $\nu_i$, and secondly,
the {\em mixing} properties of certain auxiliary flows. More precisely, we shall need
these properties in a quantitative form; in the cases we consider, this will follow eventually from 
an appropriate spectral gap. 

This situation might seem rather restrictive. However,
it arises often in many natural equidistribution questions
(``sparse equidistribution problems,'' as we discuss below)
as well as in the analytic theory of automorphic forms (especially,
subconvexity results for $L$-functions). There are applications
besides those discussed in the present paper; our aim has not been
to give an exhaustive discussion, but rather just to present
a representative sample of interesting cases. 
We shall explain the method abstractly in Sec. \ref{sec:method}
and will carry out, in the body of the paper, one example of each of the following cases:
$Y_i$ is the orbit
of a unipotent, a semisimple, and a toral
subgroup of $G$. 

In the present paper, we have focused
mostly on the case of $\PGL_2$ and $\GL_2$
over number fields. All our results pertain to this setting,
except for Thm.
\ref{thm:fouriercoefficients}, which applies to a general semisimple group.
The geometric methods of this paper are general
and we hope to analyze further higher rank examples in a future paper.

Throughout the present methods we have tried to use ``soft'' techniques
as a substitute for explicit spectral expansions. However, there still 
seem to be instances where the explicit spectral expansions are important. In a future paper \cite{MV}, joint with P. Michel,
we shall combine ideas drawn from this paper with ideas from Michel's paper \cite{michel};
in that paper, we shall make much more explicit use of spectral decomposition.

We shall use the term ``sparse equidistribution problems'' to describe
questions of the following flavor: Suppose $Z_i \subset Y_i$ is a
subset endowed with a measure $\nu^{Z}_i$, and we would like to
prove that the $\nu^{Z}_i$ are becoming equidistributed. In other
words, we wish to deduce the equidistribution of the ``sparse''
subset $Z_i$ from the known equidistribution of $Y_i$. Examples of
this type of question are Shah's conjecture \cite{Shah} (where the
$Z_i$ are discrete subsets of $Y_i$, a full horocycle orbit) as
well as Michel's results on subsets of Heegner points
\cite{michel} (where the $Z_i$ are subsets of the $Y_i$, the set
of all Heegner points.) The connection to period integrals is as
follows: one can spectrally expand the measure $\nu^{Z}_i$ in
terms of the basis for $L^2(Y_i,\nu_i)$. Using our results for
periods along $Y_i$, it will sometimes be possible to deduce the
equidistribution of $\nu^{Z}_i$.

We now briefly summarize our results.
\begin{enumerate}
\item Sec. \ref{sec:unipotent} considers
where the $Y_i$s are orbits, or pieces of orbits,
of unipotent groups.
The mixing flow is the horocycle flow along $Y_i$.

In Thm. \ref{thm:sparsehoro} (p. \pageref{thm:sparsehoro}) we show
that certain {\em sparse subsets} of horocycles on compact
quotients of $\mathrm{SL}_2(\mathbb{R})$ become equidistributed.
This is progress towards a conjecture of N. Shah. In Thm.
\ref{thm:fouriercoefficients} (p.
\pageref{thm:fouriercoefficients}) we give a fairly general bound
(in the context of an arbitrary semisimple group) on the Fourier
coefficients of automorphic forms. In the case of $G =
\mathrm{SL}_2(\mathbb{R})$ it recovers results of Good \cite{Good} and Sarnak
\cite{SaIMRN}, which resolved a problem of Selberg. The present
proof is more direct, avoiding in particular the triple product
bounds for eigenfunctions.

\item  Sec. \ref{sec:tripleproducts} considers the case when $G =
\PGL_2(F \otimes \mathbb{R})$, where $F$ is a number field, and
$\Gamma$ is a congruence subgroup thereof. The $Y_i$ are a
sequence of closed diagonal $G$-orbits on $\Gamma \backslash G
\times \Gamma \backslash G$. The mixing flow (after lifting to the
adeles) is the diagonal action of  $\PGL_2(\adele_{F,f})$, where
$\adele_{F,f}$ is the ring of finite adeles of $F$.

Prop. \ref{prop:main} and Prop. \ref{prop:mainV} give period
bounds in this context. We refer to Prop. \ref{prop:main} as a {\em subconvex bound for the triple product period}, for the reason that it should be in fact be equivalent to subconvexity for the triple product $L$-function, in the level aspect as one factor varies,
but the necessary computation of $p$-adic integrals (Hypothesis \ref{prop:intrep}) has apparently not yet been done in sufficient generality.  In Thm. \ref{thm:subconvextp} (p. \pageref{thm:subconvextp}) it is shown that these results yield subconvex
bounds, in the level aspect, for standard and Rankin-Selberg $L$-functions attached to  $\PGL_2$.

The results on standard and Rankin-Selberg $L$-functions generalize results of
Duke-Friedlander-Iwaniec and Kowalski-Michel-Vanderkam from the case $F=\Q$.\footnote{We have not attempted to address the issue of varying the
central character. This, in a sense, is the most subtle point, as
is shown by Michel's recent work on Rankin-Selberg convolutions. Our aim
in the present paper has been to show that one can derive a
coherent theory for $\PGL_2$ from the triple product bound of
Prop. \ref{prop:main}. The case of varying central character will be discussed in a future paper with Michel.} The third result, concerning subconvexity of the triple product period in the level aspect,
was not known even over $\Q$; however, Bernstein and Reznikov \cite{BRTriple} have shown
subconvexity for the triple product period in the {\em eigenvalue aspect}. 

%the third
%was not known over $\Q$. \footnote{Note that the last result depends on a
%$p$-adic computation (computing the factors in the Kudla-Harris
%triple integral representation) that, although routine, I have not
% carried out yet in general; the case where the three forms
%are holomorphic has been done by B\"ocherer and Schulze-Pillot, and probably
%one can even deduce the general case from this.}

\item  Sec. \ref{sec:torus1} considers the case when $Y_i$ is a
certain family of noncompact torus orbits on $\Gamma \backslash G$,
where $(\Gamma,G)$ is as in Sec. \ref{sec:tripleproducts}. (In
fact, the $Y_i$ are obtained by taking a fixed noncompact torus
orbit, and translating by a $p$-adic unipotent, where $p$ varies.)
The mixing flow is the action of the {\em adelic} points of the
torus.

We establish in Thm. \ref{thm:subconvexct} (p. \pageref{thm:subconvexct})
subconvexity for character twists of $\GL(2)$
in the level aspect. This was established
for $F=\mathbb{Q}$ by Duke-Friedlander-Iwaniec,
and the special case where $F$ is totally real
and the form holomorphic at all infinite places
was treated by Cogdell, Piatetski-Shapiro and Sarnak.
In particular, (\ref{eq:burgess})
gives a subconvex bound for Gr\"ossencharacter
$L$-functions over $F$, in the level aspect;
this was known over $\Q$ by work of Burgess,
and some special cases were known in the general case.

\item In Sec. \ref{sec:torus2} we consider the case where $Y_i$ is
a (union of) compact torus orbits on $\Gamma \backslash G$, where
$(\Gamma,G)$ are as in Sec. \ref{sec:tripleproducts}.  The
equidistribution of such $Y_i$ will amount to the equidistribution
of Heegner points, and we deduce it from Thm.
\ref{thm:subconvexct} in Thm. \ref{thm:heegner} (p.
\pageref{thm:heegner}). This result generalizes work of Duke over
$\Q$ and was proven, conditionally on GRH, by Zhang, Cohen,
and Ullmo-Clozel (independently). The present work 
makes this result unconditional. 

Applying mixing properties of the adelic torus flow, we obtain in
Thm. \ref{thm:sparseheegner} (p. \pageref{thm:sparseheegner}) we
obtain, under a condition of splitting of enough small primes, the
equidistribution of certain sparse subsets of Heegner points. In
the case $F=\Q$, an {\em unconditional} result of this nature is due to
Michel.\footnote{ Our method is different to Michel's: we do not deduce our
result from results on Rankin-Selberg convolutions, and indeed it
is possible to deduce a subconvexity result from ours. However,
there seem to be some curious parallels between the methods. In fact,
the method of Thm. \ref{thm:sparseheegner} is even more closely related -- as Michel
has pointed out to me -- to the work \cite{DFI} of Duke, Friedlander and Iwaniec. In that paper they amplify
class group $L$-functions but obtain only a conditional result for precisely the same reason
that Thm. \ref{thm:sparseheegner} fails to be unconditional, namely, one cannot guarantee unconditionally the existence of enough small split primes. } 
\end{enumerate}

In that context of $L$-functions, one pleasing feature of the
present method is that it is {\em geometric}: it proceeds not via
Fourier coefficients but via the integral representation. In
practice, this means that there is no difference between Maass or
holomorphic forms, nor between $\Q$ and an arbitrary base field.
Moreover, we do not make use of either the trace formula or the
Kuznetsov formula; indeed, we make no explicit use of families. 

The recent work of Bernstein-Reznikov \cite{BRTriple}
is of a similar flavor. They establish
a ``subconvex'' bound for the triple product when
the eigenvalue of one factor varies, whereas we have
treated the case where the level of one factor varies.
Their method is also geometric in nature,
and moreover their result applies to a nonarithmetic group.
By contrast, the level aspect question is not well-posed
if one leaves the arithmetic setting.

Throughout the paper we have not attempted to optimize the results. 
The input to our method is an equidistribution result. As far as possible
we have tried to establish these results by relatively ``geometric'' methods, 
deriving in the end from the mixing properties of a certain flow. Of course,
it is in many contexts better to use spectral methods, but this
would involve departing from the geometric method that
is intended to be the central theme of this paper. As remarked, we will pursue
such ``spectral'' approaches in a forthcoming paper with P. Michel \cite{MV};
some of the results of this have been discussed in \cite{MV-icm}. 

Finally, implicit in various parts of the paper is ``adelic analysis'', i.e.
the analytic theory of functions on adelic quotients, in the quantitative sense
needed for analytic number theory. There seems to be considerable scope
to develop this theory fully. 

\subsection{Other applications.} \label{subsec:applications}
The method of this paper has other applications not 
elaborated here. We discuss some of them here. 

 There are other subconvexity results that
are naturally approached by the same method: for instance,
a subconvex estimate for $L(\pi, \frac{1}{2}+it)$ where $t$ varies
and $\pi$ is a fixed cuspidal representation of $\GL(2)$
over a number field $F$.  In such a context it is natural to use
the fact that the horocycle flow is (quantifiably) weakly $k$-mixing,
for certain  $k > 1$; the use of this higher order mixing is 
closely related to Weyl's ``successive squaring'' approach
to $\zeta(1/2+it)$. 
%We discuss this example briefly in Sec. \ref{subsec:Lfns}. 
Of course, this particular instance of subconvexity is approachable
by standard methods also; an intriguing question in the subconvexity
context is how to combine the present methods with those such 
as Bernstein-Reznikov.  

 There are certain applications
to effective equidistribution theorems: for instance, it is
also possible to establish some new effective cases of Ratner's
theorem by the same ideas, see Rem. \ref{rem:effratner}. The question
of giving such ``nontrivial'' cases was raised by Margulis
in his talk at the American Institute of Mathematics, June 2004.  Unfortunately,
these new cases are rather artificial. 

One can give certain analytic applications:
let $\Gamma$ be a cocompact subgroup of $\SL(2,\R)$, and let $ \pi \subset L^2(\Gamma\backslash
\SL(2,\R))$ be an irreducible $\SL(2,\R)$-subrepresentation. 
For $m \in \mathbb{Z}$, let $e_m$ be the $m$th weight vector in $\pi$, if defined;
i.e. a vector which transforms under the character $\left( \begin{array}{cc}
\cos(\theta) & \sin(\theta) \\ - \sin(\theta) & \cos(\theta) \end{array}\right) \mapsto
e^{2 \pi i m \theta}$. We normalize it (up to a complex scalar of absolute value $1$)
by requiring that $\|e_m\|_{L^2}  = 1$.  Bernstein and Reznikov proved
the bound $\|e_m\|_{L^{\infty}} \ll (1+|m|)^{1/2}$, and asked \cite[Remark 2.5(4)]{BRSobolev} if any improvement of the
exponent $1/2$ is possible. It is quite easy to deduce from Lem. \ref{lem:ost}
such a bound; indeed, the analytic properties of the $e_m$,
as $|m| \rightarrow\infty$, is connected with the long time behavior of the horocycle flow in the same
fashion that the analytic behavior of Laplacian eigenfunctions are connected to the long time behavior 
of the geodesic flow.  In the time during which this paper was being revised for submission, 
Reznikov has proven independently a result of this type \cite{ReznikovLinfty}. Since the result
he obtains is most likely sharper than that obtained by the technique indicated above, we will not pursue this further, noting only that an advantage of the method we have indicated above is that it is likely to generalize to higher rank.

Moving slightly away from the main subject of the present paper, the idea
of using equidistribution theorems to produce mean value results for $L$-functions
seems capable of application in a variety of settings. In particular,
equidistribution results are readily available on $\mathrm{GL}(n)$, owing to Ratner's work,
whereas trace formulae are extremely unwieldy for $n > 2$. It would be interesting
to see what mean-value statements can be deduced from Ratner-type equidistribution results.

 Historically, one application
of such results has been to nonvanishing results; here the most spectacular
results (e.g. \cite{sound}) have been achieved through the so-called
mollifier technique.  It would be quite interesting to understand
if there is a geometric interpretation of the mollifier technique.

\subsection{Discussion of method: equidistribution, mixing,
and periods.} \label{sec:method}
We now turn to a discussion of the specifics of the
method used in this paper.
This method itself is quite easy to describe.
It consists in essence of two simple steps (see (\ref{eq:equidist})
and (\ref{eq:mixing}) below).

We also remark that the discussion that follows is a relatively 
faithful rendition of the method of the paper. The body
of this paper does not really utilize any new ideas beyond
the ones indicated below.  Most of the bulk consists of
the technical details necessary to connect periods with other objects
of interest (e.g. equidistribution questions or $L$-functions), as
well as setting up the machinery to quantify some standard
equidistribution results. As much as possible, we have tried to
give a self-contained treatment of all these technical details in
Sections \ref{sec:sobolev} -- \ref{sec:rs}.

We hope the
ensuing discussion
serves as a unifying thread for the rest of the paper.
We explain the method first in an abstract setting
(Sec. \ref{sec:abstract}).
We then explain (Sec. \ref{sec:anal1}
and \ref{sec:anal2}) these ideas in a
a more down-to-earth fashion,
 emphasizing the parallel with the analytic
techniques for studying $L$-functions.
Finally, Sec. \ref{subsec:concrete} illustrates these ideas in a simple example -- that of Fourier
coefficients of modular forms.
%and Sec. \ref{subsec:Lfns} motivates (at a fairly vague level) the transition
% example to $L$-functions. 

\subsubsection{Abstract setting.} \label{sec:abstract}

Let $G_2 \subset G_1$ be locally compact groups, $\Gamma \subset
G_1$ a lattice, $X = \Gamma \backslash G_1$. Let $x_i \in X$ and
put $Y_i = x_i  G_2$. We shall suppose that there exists a
$G_2$-invariant probability measure $\nu_i$ on $Y_i$. (This does not precisely cover
all the contexts we consider -- at some points we will consider $Y_i$ which are ``long pieces''
of a $G_2$-orbit rather than a single $G_2$-orbit, but the ideas in that case will be identical
to those discussed here).

 Let $f$ be a
function on $X$ and $\psi_i$ a function on $Y_i$ such that
$\int_{Y_i} |\psi_i|^2 d\nu_i =1$. We will give a bound for the
period $\int_{Y_i} f \psi_i d\nu_i$. 

In words, the idea will be to find certain correlations between the values of
$\psi$ at different points; and then show that the values of $f$
at these same points are ``uncorrelated,'' in some quantifiable sense.
Putting these together will show that the period must be small. 
The ``hard'' ingredient here is some version of the spectral gap,
i.e. quantitative mixing, which is what will show that
the ``uncorrelated-ness'' property of $f$. 

We will suppose that there exists $\meas$, a measure on $G_2$, such that
\begin{equation} \label{eq:correlation}
\psi_i  \star \meas = \lambda_i \psi_i,
\end{equation}
for some $\lambda_i \in \C$.
Here $\star \meas$ denotes the action of $\meas$ by right convolution.
Let $\check{\meas}$ be the image of $\meas$ by the involution
$g \mapsto g^{-1}$ of $G_2$.

Then
\begin{equation}  \label{eq:equidist}\begin{aligned}
\left|\int f \cdot \psi_i d\nu_i \right|^2  = \left| \lambda_i^{-1}
\int_{Y_i} f \cdot (\psi_i \star \meas) d\nu_i \right|^2 =
\left|\lambda_i^{-1} \int_{Y_i} (f \star \check{\meas}) \cdot \psi_i
d\nu_i \right|^2
\\ \leq |\lambda_i|^{-2} \int_{Y_i} |f \star \check{\meas}|^2 d\nu_i,
\end{aligned}\end{equation}
where we have applied Cauchy-Schwarz at the final step.
Now, we are assuming that the $Y_i$ are becoming equidistributed,
and so $\nu_i \rightarrow \nu$, the $G_1$ invariant measure
on $\Gamma \backslash G_1$. Thus
\begin{equation} \label{eq:mixing}\begin{aligned}
\int_{Y_i} |f \star\check{\meas}|^2 d\nu_i \approx
\int_{X} |f \star \check{\meas}|^2 d\nu
\\ = \int_{g,g' \in G_2} \langle g g'^{-1} \cdot f,f \rangle_{L^2(X)}
d\sigma(g) d\sigma(g') , \end{aligned}
\end{equation}
where $g g'^{-1} \cdot f$ denotes the right translate of $f$ by $g
g'^{-1}$.

If the $G_2$-action on $X$ is mixing in a quantifiable way --
i.e., one has strong bounds on the decay of matrix coefficients --
one obtains good upper bounds on the right-hand side of
(\ref{eq:mixing}); in combination with (\ref{eq:equidist}) this
gives an upper bound for the period $|\int_{Y_i} f \psi_i
d\nu_i|$.

The strength of the information required about the mixing varies.
In the cases we study where $G_2$ is amenable, any nontrivial
information will suffice. In the one case where $G_2$ is
semisimple, a strong bound towards Ramanujan is needed. For
instance, in the case of triple products, we need any improvement of the bound that
the $p$th Hecke eigenvalue of a cusp form on $\GL(2)$ is bounded
in absolute value by $p^{1/4}+p^{-1/4}$. (In this normalization, the trivial bound is $p^{1/2} + p^{-1/2}$).

In the rest of this paper, we shall merely apply
this argument many times, with
various different choices for $\Gamma, G_1, G_2$.
The part of the argument which will vary is
quantifying the equidistribution of the $\nu_i$, 
i.e.
 keeping track of the error in the
first approximation of (\ref{eq:mixing}).
Thus we make heavy use of Sobolev norms (Sec. \ref{sec:sobolev}),
which are an efficient method of bounding this error.

In each instance, the proof of the equidistribution
result $\nu_i \rightarrow \nu$ will always be rather straightforward, {\em except} for 
the result of Sec. \ref{sec:torus2}. The equidistribution result needed
for the proof of Thm. \ref{thm:sparseheegner}
is essentially equivalent to the subconvexity result proved
in Sec. \ref{sec:torus1}. A rather striking point is that
a similar logical dependence (although manifested very differently)
is present in the work of Michel.  The meaning of this is unclear
to the author.

In certain specific cases, the above technique
is quite familiar.
When $G_2$ is a one-parameter real group,
the above argument is quite closely related
to standard techniques of analytic number theory.  \footnote{
For example, in certain contexts when $G_2$ is abelian, 
one can push this method further by
squaring multiple times, that is to say,
considering $|\int f \psi_i d\nu_i|^{4},
|\int f \psi_i d\nu_i|^{8}$ and so forth.
In this context, one replaces
the mixing property of the $G_2$
action with results about higher order mixing
of the $G_2$-flow.
Although we will not carry this out in the present paper,
this seems rather closely connected to Weyl's proof of subconvexity for $\zeta(1/2+it)$. }
On the other hand,
when $G_2$ is an adelic group,
and $\meas$ a measure on $G_2$ that corresponds
to the action of Hecke operators (this
is carried out in Sec. \ref{sec:tripleproducts}, for instance),
the above argument will be essentially ``amplification''
in the sense of Friedlander-Iwaniec \cite{fi}.

In the following two sections, we shall attempt to explain more colloquially
the main idea that is at work here, and also discuss how the
method described above fits into the framework of analytic number
theory. Modern proofs of subconvexity, following the path-breaking
work of Friedlander-Iwaniec \cite{fi}, have roughly speaking
consisted of a mean-value theorem and an amplification step. We
shall discuss how the proof indicated above may be viewed as
geometrizing this strategy, where the mean-value 
step is replaced by an equidistribution theorem,
and the amplification step is controlled using mixing. 

Note, in particular, that in the work of Friedlander-Iwaniec,
families of $L$-functions play a central role, whereas the method
above has in a certain sense eliminated the family. Although in
the discussion below we rephrase matters so as to make clear the
connection with the work of Friedlander-Iwaniec, it seems that
from the perspective of the present paper the phrasing in terms of
families is rather artificial.

\subsubsection{Connection with analytic number theory:
Equidistribution, and mean-value theorem for periods.}
\label{sec:anal1}

Follow the notations of the previous section.
We choose an orthonormal basis $\{\psi_{i,j}\}_{j=1}^{\infty}$ for $L^2(Y_i, \nu_i)$
so that $\psi_{i,1} := \psi_i$.

By Plancherel's formula,
$\sum_{j=1}^{\infty} \left|\int f \psi_{i,j} d\nu_i\right|^2
= \int |f|^2 d\nu_i$.
Since $\nu_i \rightarrow \nu$ weakly, and we are holding
$f$ fixed, it follows that:

\begin{equation} \label{eqn:meanvalue}\sum_{j=1}^{\infty} \left| \int f \psi_{i,j} d\nu_i \right|^2  \rightarrow
\int_{\Gamma \backslash G} |f|^2 d\nu,\end{equation}
 as $i \rightarrow \infty$. Thus the equidistribution property of $\nu_i$
underlies a mean-value theorem for the $Y_i$-periods.

In many cases involving automorphic forms,
the periods will essentially be special values of $L$-functions
and (\ref{eqn:meanvalue}) amounts to a mean-value
theorem for $L$-functions. This is fairly well-known;
for example, the mean-value theorem $\int_{-T}^{T} |\zeta(1/2+it)|^4 dt \sim
T \log(T)^4$ is rather closely connected
with the equidistribution properties of the cycle
$\{(1+i/T) x, x\in\mathbb{R}\}$, when projected
to $\SL_2(\Z) \backslash \mathbb{H}$. 
A more striking example is Vatsal's use of equidistribution 
to prove nonvanishing results \cite{vatsal}. In general, it seems
that there are many interesting mean value theorems for $L$-functions
that are connected to equidistribution results. 

In any case, (\ref{eqn:meanvalue})
is not  unrelated to the standard methods of obtaining
such results; however, its primary advantage is that it is often
technically much simpler, for example when working over a number field.

\subsubsection{Connection with analytic number theory (II): Mixing,
and bounds for a single period.} \label{sec:anal2}
We now wish to pass from (\ref{eqn:meanvalue}) to nontrivial upper bounds
for a single period. It is clear that
(\ref{eqn:meanvalue}) implies at once -- by omitting all terms but one --  that
$|\int f \psi_{i,j} d\nu_i| \stackrel{<}{\sim} \|f\|_{L^2(X)}$; we shall refer to
an improvement of this bound as {\em nontrivial}.
It is evident that one must have
some further information about $\{\psi_{i,j}\}$ in order to do this;
otherwise one could simply take $\psi_{i,1}$ to be a multiple
of $f|_{Y_i}$.

In the context of analytic number theory, this is often
carried out by ``shortening the family,'' that is to say:
proving a sharp mean-value theorem of the form of (\ref{eqn:meanvalue}),
but over some subfamily of $\{\psi_{i,j}\}_{j=1}^{\infty}$;
then omitting all terms but $\psi_{i,1} = \psi_i$
will often give a nontrivial upper bound.
In the work of Friedlander-Iwaniec, a weighted mean-value theorem
is derived, which has the same effect as shortening the family.

Such a weighted mean-value theorem is also implicit in our context.
Following the notation of Sec. \ref{sec:abstract},
suppose that there is a {\em fixed} measure $\sigma$ on $G_2$
such that for all $i,j$, we have
$\psi_{i,j} \star\sigma = \lambda_{i,j} \psi_{i,j}$
(some $\lambda_{i,j} \in \C$).
Then, by Plancherel's formula,
and using the fact $\nu_i \rightarrow \nu$, we conclude:
\begin{equation} \label{eq:shorten}
\sum_{j=1}^{\infty} |\lambda_{i,j}|^2 \left|\int f \cdot \psi_{i,j} d
\nu_i\right|^2  \rightarrow  \langle f \star \check{\sigma}, f \star \check{\sigma}
\rangle_{L^2(X)}. \end{equation}
This gives a weighted mean value theorem,
which for appropriate choices of $\sigma$ amounts
to shortening the effective range of summation in (\ref{eqn:meanvalue}).
Moreover, the mixing of the $G_2$-flow bounds the right hand
side of (\ref{eq:shorten}).
In this phrasing, it becomes clear that
the measure $\sigma$ has played the role of an ``amplifier''
and the orthonormal basis for $L^2(Y_i, \nu_i)$ has played the role of the
family.

Having now explained the method
in an abstract context and indicated its equivalence
with other methods,
we now indicate more informally the source of cancellation in periods
that is at the center of our results.

In many natural situations, one obtains a basis for $L^2(Y_i,\nu_i)$
by diagonalizing a geometrically defined algebra of operators on $Y_i$. The
result of this process is that the functions $\{\psi_j\}$
exhibit correlations between their values at different points of $Y$.
For instance (for example when $G_2$ is semisimple), it often will
occur that there is a correspondence
$\mathcal{C}: Y \mapsto Y$ such the value
of each $\psi_j$ at $P \in Y$ and at the collection of points
$\mathcal{C}(P)$ are correlated in some way.
On the other hand (and we shall now speak
quite imprecisely) if the correspondence $\mathcal{C}$
``extends'' to a correspondence $\tilde{\mathcal{C}}: X \mapsto X$
one can often show, using mixing properties of $\tilde{\mathcal{C}}$,
that the values of $f$ at $P$ and $\tilde{\mathcal{C}}(P)$
will be {\em uncorrelated}, at least if $P$ is chosen
at random w.r.t the the uniform measure on $X$.

However, {\em since the $Y_i$ are becoming equidistributed}, it amounts to
almost the same thing to choose $P$ at random w.r.t. $\nu_i$ and
w.r.t. the uniform measure on $X$. Thus, for $\nu_i$-typical $P \in Y_i$,
the values of $f$ at $P$ and $\mathcal{C}(P)$ are uncorrelated,
whereas the values of $\psi_j$ at $P$ and $\mathcal{C}(P)$ are
correlated. One can then play these phenomena against each other
to obtain cancellation in the period integral $\int f \psi_j
d\nu_i$.

%\begin{ex}
%Let $G = \mathrm{SL}_2(\R), \Gamma \subset G$ a Fuchsian group.
%Suppose $u: \mathbb{R}
%\rightarrow G$
%is a one-parameter subgroup, with $U = u(\mathbb{R}))$ and $Y_i$ are a sequence
%of closed $U$-orbits with $U$-invariant measure $\nu_i$:
%we suppose $Y_i = y_i U$, where $y_i u(t_i) = y_i$
%for some $t_i \in \mathbb{R}$.
%In this section we can take the eigenbasis $\psi_{Y_i,j}$
%to be eigenfunctions of the $U$-action, that is,
%$\psi_{Y_i, j}(y_i u(t)) = e(j t t_i^{-1})$.
%
%The values of $\psi_{Y_i,j}$
%at $y \in Y_i$ and $y u(1)$ are correlated,
%as they differ by a multiplicative factor of $e(j/t_i)$.
%On the other hand, the values of $f$ at $x \in \Gamma \backslash G$
%and $x u(1)$ are, on average, somewhat uncorrelated
%if we choose $x$ at random w.r.t. $\mu$; this follows from the quantitative
%mixing of the $U$-flow. Finally, since $\nu_i$ approximates $\mu$ well
%as $i \rightarrow \infty$, we see that $f(x)$ and $f(xu(1))$
%are somewhat uncorrelated for $x$ chosen at random w.r.t. $\nu_i$.
%This, combined with the behavior of $\psi_{Y_i,j}$ is
%enough to obtain a nontrivial bound  for the period integral.  \end{ex}
%
%

\subsubsection{A concrete example.} \label{subsec:concrete}
We shall now discuss 
how to bound Fourier coefficients of a modular form by the methods just described.  Although the material below is essentially redone - with $\SL(2,\R)$ replaced by a general group -- in Sec. \ref{sec:unipotent},
the example below was very important in motivating the author's intuition, and
it seems worthwhile to include it in the introduction.

Let $\Gamma \subset \SL(2,\R)$ be a lattice
containing the element $\left( \begin{array}{cc} 1 & 1 \\0 & 1 \end{array}\right)$. 
Let $f(z)$ be a holomorphic form of weight $2$ w.r.t. $\Gamma$, which 
we write in a Fourier expansion $f(z) = \sum_{n=1}^{\infty}a_n e^{2 \pi i n z}$.
Hecke proved the bound $|a_n| \leq C n$, a bound which was only improved
(for a general -- possibly nonarithmetic -- $\Gamma$) much later, to $|a_n| \leq C n^{5/6}$, by A. Good.
We shall sketch a simple proof of a nontrivial bound $|a_n| \leq C n^{1-\delta}$
along the lines just indicated; for further details, we refer the reader to 
Sec. \ref{sec:fourcoeff}, where the procedure outlined is implemented
for a general semisimple group. 

We note that the ideas that will enter here
are exactly those that will enter into the proof of equidistribution of sparse
subsets of horocycles (see Sec. \ref{subsec:sparsehoro}), or for the
nontrivial bound  for $L^{\infty}$ norms in the weight aspect that is discussed in Sec. \ref{subsec:applications}. The proof below also works for Maass forms
(in that case the result is due to Sarnak).

The Fourier expansion implies that \begin{equation}
\label{eqn:fc} a_n = e^{2 \pi} \int_{x \in \mathbb{R}/\mathbb{Z}}
f(x + \frac{i}{n}) e^{- 2 \pi i n x}. \end{equation} In words, the idea is as follows:
the function $e^{-2 \pi i n x}$ takes the same values at $x, x+\frac{1}{n}, x + \frac{2}{n},
\dots$. On the other hand, the values of the function $f$ at these points are
(in a quantifiable sense) uncorrelated, as we shall deduce from the mixing
properties of the horocycle flow. Playing these two properties against
each other will yield an improvement of the Hecke bound for $|a_n|$. 
\footnote{Underlying this is the usual ``van der Corput'' trick: to bound $\sum_{k=1}^{K} c_k$
it suffices to bound correlations $\sum_{k=1}^{K} c_{k} c_{k+h}$; in effect
we apply this with $c_k = f(\frac{k+i}{n}), K= n$.)}

Let $\tilde{f}$ be the lift of $f$ to $\Gamma \backslash \mathrm{SL}(2, \mathbb{R})$;
that is to say, $$\tilde{f}: \Gamma \left( \begin{array}{cc} a & b \\ c & d \end{array}\right) \mapsto 
f(\frac{ai+b}{ci+d}) (ci+d)^{-2}.$$ Let $x_n = \Gamma \left(\begin{array}{cc} n^{-1/2} & 0 \\ 0 & n^{1/2} \end{array}\right)
$, and put $n(t) = \left( \begin{array}{cc} 1 & t \\ 0 & 1 \end{array}\right)$. 
Then the definitions show that $\tilde{f}(x_n n(t)) = n^{-1} f(\frac{i + t}{n})$; consequently, 
we see that
\begin{equation} \label{eq:fc}a_n = e^{2 \pi}  \int_{t=0}^{n} \tilde{f}(x_n n(t)) e^{-2 \pi i t} dt \end{equation}

(\ref{eq:fc}) expresses the $n$th Fourier coefficient of $f$ as the integral of $\tilde{f}$ over a closed
horocycle of length $n$. Moreover, 
(\ref{eq:fc}) falls into the pattern of Sec. \ref{sec:abstract}, with $G_1 = \SL_2(\mathbb{R})$, 
$G_2 = \{n(t): t \in \mathbb{R}\}$, $Y_n = \{x_n n(t): t\in \mathbb{R}\}$,  and $\psi_n: Y_n \rightarrow \C$
the function given by $x_n n(t) \mapsto e^{-2 \pi i t}$.  The fact that the $Y_n$
are becoming equidistributed amounts to the ``equidistribution of low horocycles,'' cf. \cite{SarnakHoro}. 
 In the language of Sec. \ref{sec:abstract}, we will take
$\sigma$ to be the measure on $G_2 \cong \mathbb{R}$ that is a sum of point masses
$\delta_i$, for integers $i=1, \dots, K$. We now carry out the procedure of Sec. \ref{sec:abstract} in an explicit fashion
in the paragraphs that follow. 

Let $T$ be the operation of right translation by $n(1)$ on $C^{\infty}(\Gamma \backslash \mathrm{SL}_2(\R))$: that is to say, for $F \in C^{\infty}(\Gamma \backslash \SL_2(\R))$,
we put $TF(g) = F (g n(1))$. 

The value of the right-hand side of (\ref{eq:fc}) remains unchanged if we replace $\tilde{f}$
by $T \tilde{f}$; consequently, for any integer $K \geq1$, we have
$$a_n = \frac{e^{2 \pi} }{K} \int_{t=0}^{n} (\sum_{i=0}^{K-1} T^i \tilde{f} (x_n n(t)) e^{- 2 \pi i t} dt.$$
Applying the Cauchy-Schwarz inequality we deduce that 
\begin{equation} \label{eq:Cauchy}|a_n|^2 \leq
\frac{n}{e^{-4 \pi}K^2} \int_{t=0}^{n} \left|\sum_{i=0}^{K-1} T^{i} \tilde{f} (x_n n(t))\right|^2 dt.\end{equation}

We now use come to the equidistribution part of the argument.
The {\em equidistribution of long closed horocycles} asserts that 
the closed horocycle $\{x+iy: 0 \leq x \leq 1\}$ becomes equidistributed in $\Gamma \backslash \mathbb{H}$ as $y \rightarrow \infty$.
Quantitatively, for
any $F \in C^{\infty}(\Gamma \backslash \mathbb{H})$, we have
\begin{equation} \label{eq:closedhoro}\left|\int_{0}^{1} F(x + iy) dx - \int_{\Gamma \backslash \mathbb{H}} F \right| \leq C_F y^{\delta},\end{equation}
for some $C_F$ depending on $F$, and some $\delta$ depending only on $\Gamma$. 
This assertion, originally proved by Sarnak \cite{SarnakHoro} by spectral methods,  can be deduced
quite easily from the mixing properties of the geodesic flow; this is done, in a somewhat
more general context, in Lem. \ref{lem:banana}.

We note that -- a special case of the discussion in Sec. \ref{sec:anal1} -- the equidistribution statement
(\ref{eq:closedhoro}) 
above reflects a mean-value theorem for periods. 
Indeed, if one applies it
to $F =  y^2 |f|^2$, one deduces
the asymptotic for $\sum_{n < X} |a_n|^2$.

In any case, what will be more useful is the version of (\ref{eq:closedhoro}) that is lifted
to $\Gamma \backslash \mathrm{SL}_2(\mathbb{R})$. This asserts that for
any $F \in C^{\infty}(\Gamma \backslash \mathrm{SL}_2(\mathbb{R}))$, we have
\begin{equation} \label{eq:equih} \left|\frac{1}{n} \int_{t=0}^{n} F(x_n n(t)) dt  - \int_{\Gamma \backslash \mathrm{SL_2}(\R)}
F(g) dg\right| \leq C_F n^{-\delta},\end{equation}
where $\delta$ is an explicit constant depending only on $\Gamma$, and $C_F$
is a constant depending on $F$.

From (\ref{eq:Cauchy}) and (\ref{eq:equih}) we conclude that
\begin{equation} \label{eq:gone} |a_n|^2 \ll \frac{e^{4 \pi} n^2}{ K^2} \left(\| \sum_{i=0}^{K-1} T^i \tilde{f} \|^2_{L^2(\Gamma \backslash \SL_2(\R))} 
+ C_{f,K} n^{-\delta}\right).\end{equation}

On the other hand, the explicit derivation of (\ref{eq:equih}) shows that $C_F$
may be bounded by a Sobolev norm of $F$, and consequently the constant $C_{f,K}$
that appears in (\ref{eq:gone}) 
is bounded by $O_f(K^A)$ for some $A > 0$. Thus 
$$|a_n|^2\ll_f n^2 K^{-2}  \left( \| \sum_{i=0}^{K-1} T^i \tilde{f} \|^2_{L^2(\Gamma \backslash \SL_2(\R))} 
+ K^A n^{-\delta}\right).$$

We now use the fact that the horocycle flow is mixing, in a quantifiable way. This amounts
to the assertion that there is an explicit $\delta' > 0$ and constant $C_{f}'$ such that, for $i \in \mathbb{Z}$, 
$|\langle T^{i} \tilde{f}, \tilde{f} \rangle| \ll C_{f}'(1+|i|)^{-\delta'}$. 
It follows easily that
$$\|\sum_{i=0}^{K-1} T^{i} \tilde{f}\|_ {L^2}^2 \ll_f K^{2-\delta'}.$$

 We conclude that
$|a_n| \ll n (K^{-\delta'/2} + K^{A/2-1} n^{-\delta/2})$. Taking  $K$ to be a sufficiently small power of $n$,
we conclude that $a_n$ is bounded by $n^{1-\delta''}$ for some $\delta'' > 0$ depending only on $\Gamma$. 

Clearly $\delta''$ depends only on the spectral gap of $\Gamma \backslash \SL_2(\R)$.
Of course, this dependence does not arise in the ``spectral'' methods. It can be removed
in the above method, but this seems to require some extra input, e.g. the finite-dimensionality
of the space of functionals on an irreducible $\SL_2(\R)$-representation that are invariant
under the subgroup $\{n(t): t \in \R\}$.

\subsubsection{Two other viewpoints on the method of \ref{subsec:concrete}}

There are two other viewpoints which might be helpful in thinking about Section \ref{subsec:concrete}. 
Both of these viewpoints do not literally generalize to the other situations we consider (e.g. triple products)
but may be helpful for intuition. 

\begin{enumerate}
\item The first is based on the following simple principle: suppose that $T$ is a measure-preserving transformation of the probability space $(Y, \nu)$,
and that $T$ is {\em ergodic}. If $\mu_1, \mu_2$ are two $T$-invariant probability measures
with average $\frac{\mu_1 + \mu_2}{2} = \nu$, then $\mu_1 = \mu_2 = \nu$; this follows because
$\nu$ is an {\em extreme point} of the convex set of $T$-invariant probability measures. 
More generally, given any family of probability measures averaging to $\nu$, they must almost all {\em equal} $\nu$. 

We will apply this to $Y = \SL_2(\Z) \backslash \SL_2(\R)$
and $T$ the operation of translation by $n(1)$.

Let $n$ be large; for $t \in \R/\Z$, let
$\mu_t$ be the probability measure that corresponds to normalized counting measure
on $\{x_n n(t+k): k \in \Z, 0 \leq k < n\}$.  Here notation is as prior to \eqref{eq:fc}. 

Then $\int_{0}^{1} \mu_t$ is the measure on the closed horocycle $\{x_n n(t): 0 \leq t \leq n\}$. 
Thus the family of measures $\mu_t$ averages to the measure on a long closed horocycle which,
as we remarked earlier (see \eqref{eq:equih}) approximates the $\SL_2(\R)$-invariant measure $dg$
on $\Gamma \backslash \SL_2(\R)$. But this latter measure $dg$ is ergodic w.r.t. $T$. 
So applying a more quantitative form of the principle discussed above shows that, for almost all $t \in [0,1]$, 
$\mu_t$ must be close to $dg$. It is simple to see that one can use this
to deduce bounds for the Fourier coefficients, via \eqref{eq:fc}. 

\item We will phrase the second rather imprecisely. Consider \eqref{eq:fc}. Our strategy of proof can be rephrased as:

The function $t \mapsto \tilde{f}(x_n n(t))$ is weak-mixing, whereas the function $t \mapsto e^{2 \pi i t}$
is periodic, and a weak-mixing function cannot correlate with a periodic function. 

To explain this statement, we need to explain what it means for a function on the real line
to be {\em weak-mixing}. Consider instead the case of a function $h: \Z \rightarrow \R$. 
Furstenberg's correspondence principle asserts that one can (loosely speaking) associate to this
a dynamical system $(Y, \nu, T)$ in such a way that (again loosely speaking)
$h$ arises by sampling a function on $Y$ along a generic trajectory $y_0, T(y_0), T^2(y_0), \dots$. 
We then say that $h$ is weak-mixing if the system $(Y, \nu, T)$ is so. The fact
that our function $t \mapsto \tilde{f} (x_n n(t))$ (say, when restricted to integer times)
is weak-mixing follows from the equidistribution of long closed horocycles together
with the fact that the horocycle flow is, itself, mixing. 

For more on this point of view, see e.g. \cite[Section 4]{Tao} and \cite[Lemma 5.2]{Tao}
for a version of the statement that weak-mixing functions cannot correlate with periodic ones. 
\end{enumerate} 

\subsection{Connection to existing methods.}
The following comments pertain to the results of the present paper
that concern subconvex bounds for $L$-functions. As we have emphasized above, the
methods presented here are, upon examination, seen to be closely
related to existing methods: in particular, ``Sarnak's spectral
method,'' which gave the only hitherto known instance of
subconvexity over a base other than $\Q$.

Indeed, as we have already indicated, the equidistribution step
of our method can be
seen as the geometric version of a  mean value theorem,
and the rest of the method can be seen as
an amplification step (or ``shortening the family.'')
Nevertheless, the key features of the present method are that it is
essentially
geometric (in that it avoids Fourier coefficients) and adelic (which
allows us to import much from
the modern theory of automorphic forms); it also 
does not use families in any explicit way. 
Once the notation is established -- admittedly
a nontrivial overhead -- the method
allows for very considerable technical simplification.

It is perhaps also noteworthy that the 
method given here does not require any exponential decay information for triple
products. Although such exponential decay information is central
only to subconvexity in the eigenvalue aspect, it has thus far entered
as a technical device even in treatments of the level aspect.

In a sense, the present method bears the same relation to existing
methods as adelic methods do to classical methods in the theory of
automorphic forms. The classical situation has the advantage of
concreteness, and whatever can be carried out in the adelic setting
can be (in principle) carried out in the classical setting.
However there is often a considerable technical and conceptual
advantage in working adelically.

As we have discussed, the connection between
equidistribution results and mean-value theorems for periods -- implicitly
exploited throughout this paper -- appears in the work of V. Vatsal. 

D. Hejhal considered ideas similar to that of Sec. \ref{sec:fourcoeff}
in the context of proving bounds towards Fourier coefficients; see \cite{hejhal}. In the language of this paper, his method
used a measure $\sigma$ (notation of Sec. \ref{sec:method}) with much larger support,
and consequently he was unable to get unconditional results.

Finally, as was remarked in Sec. \ref{subsec:intro}, the main
result of Sec. \ref{subsec:periodbound} is the analogue in the
level aspect of a recent result of Bernstein-Reznikov \cite{BRTriple}: 
they establish a ``subconvex bound'' on triple products as the
eigenvalue of one factor varies. Their methods also are
geometric, avoiding the use of Fourier coefficients.

\subsection{Acknowledgements.}

This paper grew out of my proof of Thm. \ref{thm:sparsehoro}. The
original proof was significantly more complicated, and I am indebted to Elon
Lindenstrauss for his insistence that Thm. \ref{thm:sparsehoro}
should amount to nothing more than equidistribution and mixing. It
was thus his intuition that led to a simplification of the proof
and an important step in my understanding. The idea that the
methods for Thm. \ref{thm:sparsehoro} might be applicable
in a more general setting
arose during conversations with Andreas Str{\"o}mbergsson, who also
made many valuable suggestions about an early version of this paper.   I thank them
both for their significant contributions.

I am very grateful to Gergely Harcos and Philippe Michel for their
encouragement of this project. Philippe
read carefully an early draft of this paper and pointed out many points
where the argument and results could be significantly improved. 
I am also grateful to Peter Sarnak, from whom I learned
much of what I know about this subject.

I have also benefited from several conversations with Joseph Bernstein
and Andr{\'e} Reznikov. I thank for their generosity in sharing
and discussing their elegant ideas. 

I have learned many of the methods that appear here from the work of others.
I mention in particular Peter Sarnak's paper \cite{SaOld},
which uses the idea of changing the test vector;
the Friedlander-Iwaniec idea of amplification and the geometric
version of it that appears in Bourgain-Lindenstrauss \cite{BL};
and the recent work of Bernstein-Reznikov \cite{BR},
in particular their elegant use of Sobolev norms.

This paper suffered a considerable delay before submission.  I would like to thank Philippe Michel
for his encouragement and insistence that it be revised and submitted, without which the delay would have likely been considerably longer. I also thank Nicolas Bergeron and Marina Ratner for
comments that improved the exposition and correctness. 

The ideas of this paper
were worked out during the workshop ``Emerging applications of measure
rigidity,'' AIM, San Francisco and at the Isaac Newton Institute. 

I was supported by the Clay Mathematics Institute
during much of the writing of the paper, and I thank them for their generous support.
I also thank the Institute for Advanced Study for providing excellent working conditions
during the academic year 2005-2006. 
I was also partially supported by NSF grants DMS-0111298 and DMS-0245606. 

\subsection{Structure of paper.}
The logical structure of this paper is as follows: Sec. \ref{sec:notn} introduces all necessary notation. 
The heart of the paper are 
Sec. \ref{sec:unipotent} (unipotent periods),  Sec. \ref{sec:tripleproducts} (the triple product period),
Sections \ref{sec:torus1} and \ref{sec:torus2} (torus periods). The remaining
Sections \ref{sec:sobolev} -- \ref{sec:rs} are of a technical nature, proving various
technical results required in the main text;  at a first reading (or even later) they
should perhaps be referred to only as necessary.

The two examples that best convey the flavor of the paper are Theorem \ref{thm:sparsehoro}
and Prop. \ref{prop:main}. The proofs of these results are relatively self-contained,
and we advise that the reader start with them. 
\section{Notation.}\label{sec:notn}
\subsection{General notation.} \label{subsec:gennotn}
We use the symbol $\ll$ as is standard in analytic number theory: namely, $A \ll B$
means that there exists a constant $c$ such that $A \leq c B$. The notation
$A \ll_{f,g,h} B$ means that the constant $c$ may depend on the quantities $f,g,h$;
the notation $A \ll_{\epsilon} B$ or $A \ll_{\varepsilon} B$ will mean, unless otherwise indicated,
that the stated bound holds for all $\epsilon$ or $\varepsilon >0$. 
In general, we will never explicate the dependence of implicit constants on the number field over which we work; and, by
an abuse of terminology, we will sometimes use the phrase ``absolute constant'' to mean
a constant that depends only on this number field. 

If $Z$ is a space we denote by $\delta_z$ the point
measure at $z \in Z$, i.e. $\delta_z(f) = f(z)$
for $f$ a continuous function on $Z$.

Now let $Z$ be a right $G$-space.
For $f$ a function on $Z$ and $g \in G$,
we write $g \cdot f$ for the right translate
of $f$ by $g$, i.e. $g \cdot f(z) = f(zg)$.
If $\mu$ is a measure on $Z$, we define the translate $g \cdot \mu$
by the rule $g \cdot \mu (g \cdot f) = \mu(f)$. In particular, if $\mu = \delta_{z}$
is the point mass at $z \in Z$, then $g \cdot \mu = \delta_{zg^{-1}}$ is the point mass at $zg^{-1}$.

If $\meas$ is a compactly supported measure on $G$,
we set $f \star \meas
\stackrel{\mathrm{def}}{=}
\int_{g} (g \cdot f) d\meas(g)$,
i.e. $f \star \meas(z) = \int_{g \in G} f(zg) d\meas(g)$.
In particular, if $\delta_{g_0}$ is the point-mass at $g_0$, 
then $f \star \delta_{g_0} = g_0 \cdot f$ is the right translate of $f$ by $g_0$.

If $\meas_1, \meas_2$ are two compactly supported measures on $G$,
we define the convolution $\meas_1 \star \meas_2$ to be
the pushforward to $G$ of $\meas_1 \times \meas_2$ on $G \times G$,
under the multiplication map $(g_1, g_2) \in G \times G \mapsto g_1 g_2$. 
Notations as above, one has the (somewhat unfortunate) compatibility relation
$(f \star \meas_2) \star \meas_1 = f \star (\meas_1 \star \meas_2)$. 

For $\meas$ a measure on a group $G$, we denote by $\check{\meas}$
the image of $\meas$ by the involution $g \mapsto g^{-1}$, and by
$\|\meas\|$ the total variation of $\meas$.

If $G$ is a Lie group, we denote by $\Ad(g)$ the endomorphism ``$X \rightarrow g X g^{-1}$''
of its Lie algebra. 

If $B \subset A$ is a finite index subgroup of the group $A$, then
we denote by $[A:B]$ the index of $B$ in $A$. 

If $h$ is an entire function, the notation $\int_{\Re(s) = \sigma} h(s) ds$ 
denotes the line integral along the line $\Re(s) = \sigma$ from $\sigma - i \infty$ to $\sigma + i \infty$.
The notation $\int_{\Re(s) \gg 1} h(s) ds$ denotes $\int_{\Re(s) = \sigma} h(s) ds$ for sufficiently large $\sigma$; in the contexts where we use this notation, the answer will be constant when $\sigma$ is sufficiently large. 
\subsection{Classical modular forms.}
As usual 
$\mathbb{H}$ denotes the upper half plane, i.e. $\{z \in \C: \mathrm{Im}(z) > 0\}$.
It admits the usual action of $\mathrm{SL}(2,\R)$ by fractional linear transformations. 

\subsection{Number fields and associated notations.} \label{subsec:NF}
Let $F$ be a number field. Throughout the paper we shall regard $F$ as fixed: that is to say,
we allow implicit constants in $\ll, \gg$ may depend on $F$ without explicit statement. 

We set $F_{\infty} = F \otimes \R$, $\adele_F$ 
the ring of adeles of $F$, $\adele_{F,f}$ the ring of finite
adeles. Thus $\adele_F = F_{\infty} \times \adele_{F,f}$. We will
fix {\em once and for all} an additive character $e_F: \adele_F /F \rightarrow
\mathbb{C}$, and denote by $e_{F_v}$ the induced
additive character of $F_v$. 

For each place $v$ we have a canonical ``absolute value'' $x \mapsto |x|_v$
on $F_v^{\times}$, namely, $|x|_v = \mathrm{meas}(x S)/\mathrm{meas}(S)$
for any Haar measure, $\mathrm{meas}$, on $F_v^{\times}$, and any subset $S$ of positive measure. 

 The same definition defines a character $\adele_F^{\times}/F^{\times}
\rightarrow \R_{>0}$, which we denote by $a \mapsto |a|_{\adele}$,
or simply by $a \mapsto |a|$ if it is clear from context.  We
denote by $\adele_F^{1}$ the subgroup of $\adele_F^{\times}$
consisting of adeles of norm $1$; then the quotient
$\adele_F^{1}/F^{\times}$ is compact. For a finite place $v$
of $F$, we denote by $\order_{F_v}$ the maximal compact subring of
the completion $F_v$, by $\mathfrak{q}_v$
the maximal ideal of $\order_{F_v}$, and by $q_v$ the cardinality of the residue
field.

We shall generally denote ideals of $\order_F$
by gothic letters $\mathfrak{l}, \mathfrak{q}, \mathfrak{n}$,
etc.
If $\fcond$ is an integral ideal of $\order_F$,
we set $\Norm(\fcond) := |\order_F/\fcond|$ to be its norm.
Moreover, we shall denote $\order_{\fcond} := \prod_{\q|\fcond}
\order_{\q}$.
Here $\order_{\q}$ denotes the completed ring, not the localized ring,
i.e. $\order_{\fcond}$
is the inverse limit of the rings $\order_F/ \fcond^N$.

We denote by $\diff$ the different of the character $e_F$, i.e.
$\diff$ is a fractional ideal so that $\diff_v^{-1}$
is, for every finite place $v$,
the largest $\order_{F_v}$-submodule of $F_v$ upon which $e_F$ is
trivial.

\subsection{Adele groups and their function spaces} \label{sec:adelicfunctionspace}

Let $\G$ be a connected reductive algebraic group over a number field $F$,
and let $\Zz$ be its center.
Denote by $\adele_{F,f}$ the ring of finite adeles,
and fix for each finite place $v$ a maximal open compact subgroup
$K_{v,\G} \subset \G(F_v)$ with the property that
$\Kmaxg := \prod_{v \, \mathrm{finite}} K_{v,\G}$ is a
maximal open compact subgroup of $\G(\adele_{F,f})$.
Put $\quot_{\G} = \G(F) \backslash \G(\adele_F)$,
$\quotGad = \Zz(\adele_F) \G(F) \backslash \G(\adele_F)$.
Then $\quotGad$ has finite volume with respect to any $\G(\adele_F)$-invariant 
measure. 

Let $\omega:\Zz(\adele_F) \rightarrow \C^{\times}$ be a unitary character.
We define the space $C^{\infty}_{\omega}(\quotG)$
to be the space of functions on $\quotG$
whose stabilizer in $\Kmaxg$ has finite index,
which transform under $\Zz(\adele_F)$ by $\omega$,
and so that the function $g \mapsto f(xg)$
is a $C^{\infty}$ function of $g \in \G(F_{\infty})$,
for each $x \in \quotG$.
Similarly one defines an $L^2$-space $L^2_{\omega}(\quotG)$, or simply $L^2$
if the central character $\omega$ is clear from context, by completing the space of compactly supported functions in $C^{\infty}_{\omega}(\quotG)$
with respect to the Hilbert norm
$\|f\|_2 := \left( \int_{\quotGad} |f(g)|^2 dg\right)^{1/2}$. 

For $\psi \in C^{\infty}_{\omega}(\quotG)$,
we denote by $K_{v,\psi}$ the stabilizer of $\psi$ in $K_{v,\G}$,
and put 
\begin{equation} \label{eq:kpsidef} K_{\psi} =\prod_{v \, \mathrm{finite}} K_{v,\psi}.\end{equation} We note that
$K_{\psi}$ is, in general, a proper subgroup of the stabilizer
of $\psi$ in $\Kmaxg$.
%we denote by $K_{\psi} \subset \Kmaxg$ the stabilizer of $\psi$
%in $\Kmaxg$.

For $\psi \in C^{\infty}_{\omega}(\quotG)$, we define the finite set
of places $\Supp(\psi)$ to be those finite $v$ for which $K_{v,\G}$ does not fix $\psi$, i.e.
\begin{equation} \label{eq:Cond}
\Supp(\psi) \stackrel{\mathrm{def}}{=}
\{v: K_{v,\G}  \neq K_{v,\psi}\}.
\end{equation}

It is convenient to introduce some
notions of ``size'' on $\G(\adele_F)$.
Let $\mathfrak{g}$ be the Lie algebra
of $\G(F_{\infty})$. It is a finite dimensional real vector space; fix an arbitrary norm on it. 
For $g_{\infty} \in \G(F_{\infty})$,
we denote by $\|g_{\infty}\|$ the operator norm
of the adjoint endomorphism
$\Ad(g_{\infty}^{-1}):  \mathfrak{g} \rightarrow \mathfrak{g}$.
If $v$ is a finite place of $F$ and $g_v \in \G(F_v)$,
we set $\|g_v\| = [K_{v,\G} g_v K_{v,\G} : K_{v,\G}]$, i.e.
the number of right- $K_{v,\G}$ cosets in $K_{v,\G} g_v K_{v,\G}$. 
For $g_f = (g_v)_{v \, \mathrm{finite}} \in \G(\adele_{F,f})$
we put $\|g_f\| = \prod_{v} \|g_v\|$.
Finally for $g_{\adele} = (g_{\infty}, g_f) \in \G(F_{\infty}) \times
\G(\adele_{F,f})$,
set $\|g_{\adele}\| = \|g_{\infty}\| \cdot \|g_f\|$.

We remark that $\|g_{\infty}\|, \|g_f\|, \|g_{\adele}\|$ are all invariant
by the center of $\G$.

\subsection{The groups $\G= \GL(2)$ and $\G=\PGL(2)$ and some of their subgroups. } \label{sec:gl2pgl2}
We will deal most often with the cases of
$\G = \GL(2)$ (resp. $\G = \PGL(2)$). 
In that setting we shall write $\quotg$ (resp. $\quot$) for $\quot_{\G}$.

We will make use of the following algebraic subgroups of $\GL_2$,
which we will often also regard as algebraic subgroups of $\PGL_2$ in the
obvious way:
$$ N = \left(\begin{array}{cc} 1 & * \\ 0 & 1 \end{array}\right), \ \
B = \left(\begin{array}{cc} * & * \\ 0  & * \end{array}\right), \ \
A = \left(\begin{array}{cc} * & 0 \\ 0 & * \end{array}\right), \ \
Z = \left(\begin{array}{cc} x & 0 \\ 0 & x \end{array}\right).$$

If $R$ is any ring and $x \in R, y \in R^{\times}$,
we denote\footnote{
(In Section \ref{subsec:sparsehoro} alone, we will use slightly different notation for $a(y)$
to accomodate the fact that we deal with $\mathrm{SL}_2$ rather than $\GL_2$. We make
the relevant notation clear in that section.}
\begin{multline} \label{eq:notn}
n(x) = \left(\begin{array}{cc} 1 & x \\ 0 & 1 \end{array}\right), \ \ \
\overline{n}(x) = \left(\begin{array}{cc} 1 & 0 \\ x & 1
\end{array}\right),
\ \ a(y) =
\left(\begin{array}{cc} y & 0 \\0 & 1
\end{array}\right), \\
a'(y) =
 \left(\begin{array}{cc} 1 & 0 \\0 & y
 \end{array}\right),
 \ \ \
 z(y) = \left(\begin{array}{cc} y & 0 \\0 & y
 \end{array}\right),
 w = \left(\begin{array}{cc} 0 & 1 \\ -1 & 0 \end{array}\right).
 \end{multline}
 all elements of $\GL_2(R)$.

If $v$ is a place of $F$ and $x \in F_v, y \in F_v^{\times}$,
we denote by $n_v(x)$ (resp. $a_v(x)$) the element $n(x)$ (resp. $a(y)$)
considered as an element of $\GL_2(\adele_F)$
via the natural inclusion $\GL_2(F_v) \hookrightarrow \GL_2(\adele_F)$.

For each place $v$, we let $K_v$ be the standard maximal compact subgroup
of
$\GL_2(F_v)$, i.e. $K_v$ is the stabilizer of the norm on $F_v^2$ given by
$\sqrt{|x|_v^{2/\deg(v)}+|y|_v^{2/\deg(v)}}$ if $v$ is archimedean, where $\deg(v)$
is the degree\footnote{Recall that $
|x|_v$, for a complex place $v$ and $x \in F_v$, is the {\em square } of the usual absolute value on $\mathbb{C}$!} of $F_v$ over $\mathbb{R}$; and $\max(|x_v|, |y_v|)$ if
$v$ is nonarchimedean.  Thus, in particular, $K_v
=\GL_2(\order_{F,v})$ if $v$ is nonarchimedean. We put $\Kmax = \prod_{v \, \mathrm{finite}} K_v$. 
$K_v$ (respectively $\Kmax$) is a maximal compact subgroup of 
$\GL_2(F_v)$ (respectively $\GL_2(\adele_{F,f})$), and
(by projection) can also be regarded as a maximal compact subgroup of 
$\PGL_2(F_v)$ (respectively $\PGL_2(\adele_{F,f})$). 
Similarly $\Kmax \times K_{\infty}$ is a maximal compact subgroup of $\GL_2(\adele_F)$,
and may also be regarded as a maximal compact subgroup of $\PGL_2(\adele_F)$. 
%\footnote{We remark on a slight change of notation from the setting of $\G$ a general reductive group:
%in that setting $\Kmaxg$ denotes a maximal compact subgroup of $\G(\adele_{F,f})$, 
%whereas here $\Kmax$ is a maximal compact subgroup of $\GL_2(\adele_F)$.}. We put moreover
%$K_{\infty} = \prod_{v \, \mathrm{infinite}} K_v$

For $\q$ a finite prime of $F$,
we denote by $\unif_{\q} \in F_{\q}$ a uniformizer,
and by $[\unif_{\q}]$ the element
of $\adele_F^{\times}$ that is the image
of $\unif_{\q}$ under the natural inclusion $F_{\q}^{\times} \hookrightarrow \adele_F^{\times}$.

Let $\q$ be a finite prime of $F$. It will
be convenient to define certain open compact
subgroups of $K_{\q}$.  For each $e_{\q} > 0 $, 
we define $K[\q^{e_{\q}}] \subset K_{\q}$ (resp. $K_0[\q^{e_{\q}}] \subset K_{\q}$)
to be the 
be the kernel of $\GL_2(\order_{\q}) \rightarrow \GL_2(\order_{\q}/\unif_{\q}^m)$
(resp. the preimage, under this map, of the upper triangular matrices). 
Thus 
$$K_{0}[\q^{e_{\q}}] = \{\left(\begin{array}{cc} a & b \\ c & d \end{array}\right):
a,b,d \in \order_{\q}, c \in \q^{e_{\q}}, ad-bc \in \order_{\q}^{\times}\}.$$

Now let $\fcond$ be a fractional ideal, not necessarily prime, of $F$.
Factorize $\fcond = \prod_{\q} \q^{e_{\q}}$ into prime ideals. We
define elements $[\fcond] \in \adele_F^{\times}, a([\fcond]),
n([\fcond]) \in \GL_2(\adele_F)$ via:
\begin{equation} \label{eq:fconddef}
[\fcond] = \prod_{\q | \fcond} [\unif_{\q}]^{-e_{\q}}, \ \ \ \
n([\fcond]) := \prod_{\q|\fcond}
n_{\q}(\unif_{\q}^{-e_{\q}}),  \ \ \ \
a([\fcond]) = \prod_{\q | \fcond}
a_{\q}(\unif_{\q}^{-e_{\q}}).\end{equation}

Suppose $\chi: \adele_F^{\times}/F^{\times} \rightarrow \C^{\times}$
is a character.
We define
$$\chi(\fcond) = \begin{cases} 0, \chi \mbox{ ramified at any place
dividing }\fcond, \\
 \prod_{\q|\fcond} \chi(\unif_{\q})^{e_{\q}}, \mbox{ else.} \end{cases}$$

\subsection{Measures.} \label{subsec:measures}

The choice of measure is not especially important, as we are only
interested in upper bounds; thus, so long as we are consistent,
the precise selection does not matter. We choose a ``standard''
set of measures here; at times in the text, especially when carrying
out equidistribution arguments, it will be more convenient to use probability measures, and
we will indicate when this is the case. 

We denote by $\mu_{\quot}$ the $\PGL_2(\adele_F)$-invariant probability measure on $\quot$. 
We shall sometimes simply denote it by $dx$. 

Let $v$ be a finite place of $F$. Unless explicitly stated otherwise,
the measures on $\GL_2(F_v)$, $\PGL_2(F_v)$,
$F_v$ and $F_v^{\times}$ are the Haar measure
which assigns $\GL_2(\order_{F_v})$ (resp. $\PGL_2(\order_{F_v})$, $\order_{F_v}$, $\order_{F_v}^{\times}$)
the total mass $1$.

For $v$ archimedean, endow $F_v$ with a multiple of Lebesgue  measure $c_v dx$,
where the constants $c_v$ are fixed arbitrarily in such a way that the induced product measure on $F_{\infty}$ satisfies  $\vol(F_{\infty}/\order_F) = 1$;
equivalently, the product measure on $\adele_F$ satisfies $\vol(\adele_F/F) = 1$.   In particular, this product measure on $\adele_F$ is self-dual
with respect to $e_F$. 
We endow  $F_v^{\times}$ with the measure $d^{\times}x = \frac{dx}{|x|_v}$, where $dx$ is Lebesgue measure.

These choices induce a Haar measure on $N(F_v)$, by means
of the identification $x \mapsto n(x)$; similarly, the identifications
$(y,y') \mapsto a(y) a'(y')$ and $y \mapsto z(y)$ induce Haar
measures on $A(F_v)$ and $Z(F_v)$. Equip $K_v$
with the measure of mass $1$, and give $\GL_2(F_v)$
the measure arising from the Iwasawa decomposition $N(F_v) \times A(F_v) \times K_v$.
Equip $\PGL_2(F_v) = \GL_2(F_v)/Z(F_v)$ with the ``quotient'' measure. 

We then take the measures 
on $\GL_2(\adele_F), \PGL_2(\adele_F), \adele_F, \adele_F^{\times}$
to be the corresponding product measures.

The measure on any discrete group (e.g. $\PGL_2(F)$, considered
as a subgroup of $\PGL_2(\adele_F)$) will be counting measure.

Usually (indeed, unless otherwise specified) we shall use the $\PGL_2(\adele_F)$-invariant probability measure on 
$\quot = \PGL_2(F) \backslash \PGL_2(\adele_F)$. This does not coincide
with the quotient measure induced from $\PGL_2(\adele_F)$, but they
differ by some constant depending only on $F$. On the few occasions
we shall have occasion to use the latter measure, we will indicate this. 
%and the measure on the quotient space $\PGL_2(F) \backslash \PGL_2(\adle_F)$
%will be the induced (quotient) measure unless specified. Similarly we obtain a measure on $\GL_2(F) \backslash \GL_2(\adele_F)$.  

\subsection{Projection onto locally constant functions.} \label{subsec:proj}
For equidistribution questions it is usually convenient
to deal with the constant function and its orthogonal complement
separately. Some minor complications arise in our case
since the ambient spaces are not connected. In fact:
The space $C^{\infty}(\quotG)$
is a direct limit of function spaces $C^{\infty}(\quotG/K)$
where $K \subset \Kmaxg$ has finite index.
Unless $\G$ is simply connected, the manifolds
$\quotG/K$ need not be connected. 

Of course, to deal with this,
one can (if $\G$ is semisimple) simply replaces the notion of {\em constant} function
by {\em locally constant} function.  However, in the general case of $\G$ reductive, 
matters
are slightly complicated by the necessity of dealing with central characters.

Since we will only use this definition when
$\G$ is a product of $\GL(2)$s,
we restrict ourselves to that setting. First suppose
that $\G = \GL(2)$.
We define a projection $\mathscr{P}: C^{\infty}_{\omega}(\quotg)
\rightarrow C^{\infty}_{\omega}(\quotg)$ via
\begin{equation} \label{eq:projdef} \mathscr{P}f(x) = \int_{h \in SL_2(F) \backslash \SL_2(\adele_F)}
f(hx) dh = \sum_{\chi^2 = \omega} \chi(x) \int_{\quot} f(y)
\overline{\chi(y)}
dy,\end{equation}
where $dh$ is the $\SL_2(\adele_F)$-invariant probability measure,
$dy$ the $\GL_2(\adele_F)$-invariant probability measure on $\quot$, 
$\chi$ ranges over characters of $\adele_F^{\times}/F^{\times}$ with square $\omega$, 
$\chi(y)$ the function on $\quotg$ defined by $g \mapsto \chi(\det(g))$,
and the second equality is easily verified.
We note, in particular, that the $\chi$ sum is finite
(any $\chi$ for which the corresponding term is nonvanishing must
be unramified outside $\Supp(f)$).

Then $\|\mathscr{P}f\|_{L^{\infty}} \leq \|f\|_{L^{\infty}}$, as is clear
from the first equality of (\ref{eq:projdef}),
and $\mathscr{P}$ is a self-adjoint projection w.r.t. $L^2$, as is clear
from the second equality. 
We say a function $f$ is {\em totally nondegenerate}
if $\mathscr{P} f=0$.

If $\G = \GL(2) \times \GL(2)$, 
and $\omega = (\omega_1, \omega_2)$ is a character of the center $\mathbf{Z}(\adele_F) = \adele_F^{\times} \times \adele_F^{\times}$, we denote by $\mathscr{P}_1$ the operator on $C^{\infty}_{\omega}(\quotG)$ given by 
 $$\mathscr{P}_1f(x_1, x_2)
 = \int_{h \in SL_2(F) \backslash \SL_2(\adele_F)}
 f(hx_1, x_2) dh = \sum_{\chi^2 = \omega_1} \chi(x_1)
 \int_{\quot} f(y, x_2) \overline{\chi(y)} dy.$$
 We define $\mathscr{P}_2$ similarly, interchanging
 the role of the first and second coordinate.
 The operators $\mathscr{P}_{j}$ for $j=1,2$ commute,
 satisfy $\|\mathscr{P}_j f\|_{L^{\infty}} \leq \|f\|_{L^{\infty}}$
 and are commuting self-adjoint projections on $L^2$.
 We say that a function $f$ is
 {\em totally nondegenerate} if $\mathscr{P}_1 f = \mathscr{P}_2 f = 0$.

\begin{lem} \label{rem:whyproj}
Let $v$ be a place of $F$. 
The projection $\mathscr{P}$ acts by the identity on 
the subspace $W \subset L^2_{\omega}(\quotG)$ spanned
by one-dimensional representations of $\GL_2(F_v)$ occurring
in $L^2_{\omega}(\quotG)$. 

Similarly,  $\mathscr{P}_1$ (resp. $\mathscr{P}_2$) acts by the identity on
the space $W_1$ (resp. $W_2$) spanned by one-dimensional
representations of $\GL_2(F_v)$
occurring in $L^2(\quot \times \quot)$
for the action on the first (resp. second) factor. 
\end{lem}

\proof This follows from the spectral decomposition for $\GL(2)$.
For instance, it is known that the space $W$ is precisely
the span of functions of the form $g \mapsto \chi(\det (g))$, 
where $\chi$ ranges over characters of $\adele_F^{\times}/F^{\times}$
satisfying $\chi^2= \omega$. 
 \qed

\subsection{Hecke operators and bounds towards the Ramanujan conjecture.}\label{subsec:Hecke}
Let $\mathfrak{l}$ be a prime ideal of $\order_F$ and $r$
an integer $\geq 1$.
Let $F_{\mathfrak{l}}$ be the completion of $F$ at the prime
$\mathfrak{l}$.
Take the Haar measure on $\GL_2(F_{\mathfrak{l}})$ so that
it assigns mass $1$ to $\GL_2(\order_{F_{\mathfrak{l}}})$.
Define the measure $\mu^{*}_{\mathfrak{l}^r}$ on $\GL_2(F_{\mathfrak{l}})$
to the restriction of Haar measure to the set
$\GL_2(\order_{F_{\mathfrak{l}}}) \cdot \left(\begin{array}{cc}
\unif_{\l}^r & 0 \\ 0 & 1
\end{array} \right) \GL_2(\order_{F_{\mathfrak{l}}})$, so
that the total mass of $\mu^{*}_{\mathfrak{l}^r}$ is $\Norm(\mathfrak{l})^{r-1}
(\Norm(\mathfrak{l})+1)$. 

Moreover, set
\begin{equation} \label{eq:heckeopdef}
\mu_{\mathfrak{l}^r} = \frac{1}{\Norm(\mathfrak{l})^{r/2}}\sum_{k \leq
\frac{r}{2}} \mu^{*}_{r-2k}, \ \
\overline{\mu}_{\l^r} := \frac{\mu_{\l^r}}{\|\mu_{\l^r}\|},
\end{equation}
where $\| \cdot \|$ denotes total variation. Thus $\overline{\mu}_{\l^r}$
is
a probability measure.
Via the natural inclusion of $\GL_2(F_{\mathfrak{l}})$
in $\GL_2(\adele_{F,f})$, we may regard $\mu_{\mathfrak{l}^r}$ as a
compactly
supported measure on $\GL_2(\adele_{F,f})$;
by abuse of notation, we will not introduce a different symbol for this
measure.
If $\mathfrak{n}$ is an integral ideal of $\order_F$,
factorize $\mathfrak{n} = \prod_i \mathfrak{l}_i^{r_i}$
and put $\mu_{\mathfrak{n}} = \prod \mu_{\mathfrak{l}_i^{r_i}},
\overline{\mu}_{\mathfrak{n}} = \prod
\overline{\mu}_{\mathfrak{l}_i^{r_i}}$.
Here $\prod$ is taken to mean convolution of measures on
$\GL_2(\adele_{F,f})$.

Convolution by $\mu_{\mathfrak{n}}$ on $L^2(\quot)$ corresponds to the $\mathfrak{n}$th
Hecke operator; in this normalization
the Ramanujan conjecture corresponds to it
having eigenvalues $\leq 2$ in absolute value.
%More generally, \mu_{\mathfrak{l}^r}$ acts by the trace of the $r$th symmetric power of the
%Hecke matrix, normalized so that temepred representations correspond to unitary Hecke matrix. 

The adelic measures $\mu_{\n}$ satisfy the usual multiplication laws, appropriately interpreted:
if $\mathfrak{n}$ and $\mathfrak{m}$ are ideals, then
%One would like to say that
%$\mu_{\mathfrak{n}} \star \mu_{\mathfrak{m}} = \sum_{\d |
%(\mathfrak{m},\mathfrak{n})}
%\mu_{\n\m\d^{-2}}$; this is not precisely true as an equality
%of measures on $G$.
%However, the two measures have the same action on any function invariant
%by
%$\prod_{\mathfrak{l} |\n \m} \PGL_2(\order_{F_{\mathfrak{l}}})$, i.e.
\begin{equation} \label{eq:conv}
\int_{\PGL_2(\adele_F)} h(x) d(\mu_{\mathfrak{n}} \star \mu_{\mathfrak{m}})(x)   =
\sum_{\d |
(\mathfrak{m},\mathfrak{n})}
\int_{\PGL_2(\adele_F)} h(x) d\mu_{\n\m\d^{-2}}(x), \end{equation}
whenever $h$ is a function on $\PGL_2(\adele_F)$
that is invariant under $\PGL_2(\order_v)$ for all $v|\n \m$.

\begin{defn} \label{def:ramanujan}
Set $\alpha$ be a bound towards Ramanujan for $\GL_2$
over $F$, i.e. $\alpha$ is so that $\mu_{\l}$ acts on any $\PGL_2(\order_{F_{\l}})$-invariant
cuspidal eigenfunction by an eigenvalue  $\leq  \Norm(\l)^{\alpha}+\Norm(\l)^{-\alpha}$ in absolute value. 
\end{defn}

 Thus $\alpha = 0$ corresponds to the Ramanujan conjecture,
$\alpha =1/2$ the trivial bound.
By work of Kim and Kim-Shahidi, we can take
$\alpha = 3/26$. For our applications, any value of $\alpha$
less than $1/4$ would suffice.

Note that we shall slightly vary this notation (but in a reasonably compatiable way) in Section \ref{subsec:sparsehoro} and Section \ref{subsec:weirdnotn}. In those parts,  
we shall deal with a (not necessarily arithmetic) quotient $\Gamma \backslash \SL_2(\R)$,
and $\alpha$ will denote a number so that $L^2(\Gamma \backslash \SL_2(\R))$
does not contain any complementary series with parameter $\geq \alpha$. 
(Here the complementary series is understood to be parameterized by $(0,1/2)$). This is compatible with the above notation, however; e.g. 
if $\Gamma$ were a congruence subgroup, $\alpha = 3/26$ would again be admissible.

\subsection{Sobolev-type norms on real and adelic quotients.}
\subsubsection{General comments.}\label{subsec:generalcomments}
Let $M$ be a real manifold. Recall that the Sobolev norm on
$C^{\infty}(M)$ controls, roughly speaking,
the $L^p$ norm of a function together with the $L^p$ norms of certain
derivatives. These norms will be tremendously useful
throughout the paper to control equidistribution rates.
We shall use both the relatively simple definition
when $M = \Gamma \backslash G$ and an adelic variant.

First let us remark on the use of $L^p$-Sobolev norms for $p > 2$.
This is solely to do with noncompactness. 
If we were to deal only with {\em compact} quotients, then
the $L^2$-Sobolev theory would always suffice. However,
in the noncompact case, the $L^2$-Sobolev norms do not
(e.g.) give good bounds on the size of a function high in a cusp.
There are, of course, various ways to rectify this; for example
we could include weights that measure the height into the cusp. 
We have chosen instead to use $L^p$-norms with $p > 2$, which
is technically very simple, but has some disadvantages (e.g. it does not induce
a Hilbert space structure). 

Note that we will allow our
seminorms and norms to take the value $\infty$.
Thus a seminorm
on a complex vector space $V$ will be a function
from $V$ to $\mathbb{R}_{\geq 0} \cup \{\infty\}$
satisfying \begin{enumerate}
\item $\|\lambda v \| = |\lambda| \|v\|$, for any $v \in V$
such that $\|v\| < \infty$;
%\item $\|\lambda v\| = \infty$ if $\lambda \neq 0$
%and $\|v\| =\infty$
\item $\|v_1 + v_2\| \leq \|v_1\|+\|v_2\|$
if both $\|v_1\|$ and $\|v_2\|$ are not infinite.
\end{enumerate}
It is a norm if additionally $\|v\|=0$ implies $v = 0$.
Note that giving such a seminorm on $V$
is equivalent to giving a subspace $V_f \subset V$
together with a finite-valued seminorm on $V_f$.
Indeed take $V_f =\{v \in V: \|v\| <\infty\}$,
equipped with the restriction of $\|\cdot\|$.

We remark that we do not require that our norms be complete. 

\subsubsection{Non-adelic setting.} \label{sec:sobreal}
Suppose $\Gamma \subset G$ is a lattice in a connected semisimple Lie
group.
Fix for all time a basis $\basis$ for the Lie algebra $\lieg$ of $G$
and a norm $\|\cdot\|$ on $\lieg$.
For $g \in G$, we denote by $\|g\|$ the operator
norm of $\mathrm{Ad}(g^{-1}): \lieg \rightarrow \lieg$, i.e. the map
$X \mapsto g^{-1} X g$. 

For $f \in C^{\infty}(\Gamma \backslash G)$,
and $1 \leq p \leq \infty$,
we put 
\begin{equation} \label{eq:spddef} S_{p,d} = \sum_{\ord(\mathcal{D})\leq d}
\|\mathcal{D} f(g)\|_{L^p(\Gamma \backslash G)}.\end{equation}
Here $\mathcal{D}$ ranges over all monomials
in $\basis$ of order $\leq d$, and $\mathcal{D}$ acts on $f$ by right differentiation. (For example, $X \in \mathfrak{g}$
acts on $f$ via $Xf(g) = \frac{d}{dt} f(g e^{t X})$.) 

Changing $\basis$
only distorts $S_{p,d}$ by a bounded factor.
(That is to say, if $S_{p,d}'$ is the norm
obtained by replacing $\basis$ by another basis, then
there are positive reals $c_1, c_2$, possibly depending on $d$, such that $c_1 S_{p,d} \leq S_{p,d}' \leq c_2 S_{p,d}$. )

We will often use the following simple remark: Fix a Riemannian metric $d(\cdot, \cdot)$ on $G$
and suppose $g \in G$ belongs to some fixed compact set.
Then, for $f \in C^{\infty}(\Gamma \backslash G), x \in \Gamma \backslash G$, we have $|f(xg) - f(x)| \ll S_{\infty,1}(f) d(g,1)$. 
Indeed, we may assume that $g$ is close to the identity and write $g=\exp(X)$, with $X \in \lieg$;
now apply the mean value theorem to $t \mapsto f(x e^{tX})$. 

Moreover, the following elementary properties are easily verified (we only need
them in the case $p = \infty$).
\begin{lem} \label{lem:realsob}
Let $f_1, f_2 \in C^{\infty}(\Gamma \backslash G)$
and $g \in G$.
Then
\begin{equation} \begin{aligned}
S_{\infty,d}(f_1  f_2) \ll_{d} S_{\infty,d}(f_1) S_{\infty,d}(f_2) \\
S_{\infty,d}(g \cdot f_1) \ll_{d} \|g\|^{d} S_{\infty,d}(f_1)
\end{aligned}\end{equation} \end{lem}

We remark that, in the case $p=2$,  
the rule (\ref{eq:spddef}) also defines a system of Sobolev norms on any unitary $G$-representation;
the case discussed above corresponds to the unitary representation $L^2(\Gamma \backslash G)$. 

\subsubsection{Adelic Sobolev norms} \label{subsec:adelicsobolev}
Let's first describe what the point is intended to be (evidently
there are many ways of implementing it, cf. Rem. \ref{rem:sobolevnorm}).  We would like
to put a norm on the adelic function space, suitable
for controlling e.g. period integrals. 
Consider $C^{\infty}_{\omega}(\quotG)$
in the case of $\G= \SL_2$, $ F= \Q$, $\omega = 1$
as a direct limit of spaces $C^{\infty}(\Gamma_i \backslash \SL_2(\R))$,
where $\Gamma_i$ ranges over some class of congruence subgroups
of $\Gamma_0 := \SL_2(\Z)$. 
We equip each quotient $\Gamma_i \backslash \SL_2(\R)$ with the $\SL_2(\R)$
invariant probability measure. 
Then, on each space $C^{\infty}(\Gamma_i \backslash \SL_2(\mathbb{R}))$
we have the norm $S_{p,d}$ defined in the previous section. On the other hand,
typical bounds on automorphic forms have an implicit dependence on the ``level'',
i.e. the index $[\Gamma: \Gamma_i]$, so one would like
to have a norm that increases with the level. 
The most naive candidate is, fixing
a real number $\beta >0$, to define
the ``norm'' of $f \in C^{\infty}(\Gamma_i \backslash \SL_2(\R))$ to be $[\Gamma:\Gamma_i]^{\beta}  S_{p,d}(f)$. This unfortunately does not quite make sense when we pass to the direct limit:
however, we can ``force it to make sense'' by considering the {\em maximal norm} on the direct limit
whose restriction to each $C^{\infty}(\Gamma_i \backslash \SL_2(\R))$ is bounded
above by $[\Gamma: \Gamma_i]^{\beta} S_{p,d}(f)$. This will suffice for our purposes.

Let us formalize these ideas. 
In what follows we return to the setting of $\G$ a reductive group over
$F$.
The adelic Sobolev norms will be a family of norms on $C^{\infty}_{\omega}(\quotG)$
indexed by a triple $(p,d,\beta)$. The $d$ and $\beta$ indicate,
approximately speaking, how stringently one should
``penalize'' rapid variation at the infinite and finite places
respectively.

Let $p \geq 1, k \in \mathbb{N}, \beta \geq 0$.
Fix a basis $\basis = \{X_i\}$ for the real Lie group $\mathrm{Lie}(\G(F_{\infty}))$.
Recalling the definition of $K_{\psi}$ from (\ref{eq:kpsidef}), 
we define
the pre-Sobolev functions $PS_{p,d, \beta}$ on $
C^{\infty}_{\omega}(\quotG)$ via:
\begin{multline} PS_{p,d,\beta}(\psi) =
[\Kmaxg:K_{\psi}]^{\beta}
\sum_{\ord(\mathcal{D})  \leq d}
\|\mathcal{D} \psi\|_{L^p(\quotGad)} \\ =
\prod_{v \, \mathrm{finite}} [K_{v,\G} : K_{v,\psi}]^{\beta}
\sum_{\ord(\mathcal{D}) \leq d} \|\mathcal{D} \psi \|_{L^p(\quotGad)},
\end{multline}
where the sum ranges over $\mathcal{D}$ that are monomials
in $\basis$ of order $\leq d$.

The function $PS_{p,d,\beta}$ does not satisfy the triangle inequality.
We define the $(p,d,\beta)$-Sobolev norm
$S_{p,d,\beta}$ to be the maximal seminorm on $C^{\infty}_{\omega}(\quotG)$
satisfying $S_{p,d,\beta}(\psi) \leq PS_{p,d,\beta}(\psi)$.
Explicitly,
\begin{equation} \label{eq:concretenorm}S_{p,d,\beta}(\psi) = \inf
\{\sum_{i=1}^{n} PS_{p,d,\beta}(\psi_j):
\sum_{i=1}^{n} \psi_j = \psi,\, \psi_j \in
C^{\infty}_{\omega}(\quotG).\}\end{equation}
In fact, it is clear that the right-hand side
of (\ref{eq:concretenorm}) defines a seminorm
that is dominated by $PS_{p,d,\beta}$ (take the collection
$\{\psi_i\}$ to consist of $\{\psi\}$ alone); moreover,
it is evidently maximal in the class of such seminorms.
Finally, as $PS_{p,d,\beta}(\psi) \geq \|\psi\|_{L^p}$,
the $S_{p,d,\beta}$ are in fact norms on $C^{\infty}_{\omega}(\quotG)$.

It will often be useful to omit the argument $\beta$
and set it to a ``default'' value of $1/p$.
We therefore define $S_{p,d} := S_{p,d,1/p}$, for $p \neq 0$,
and $S_{\infty, d} := S_{\infty,d,0}$.

{\bf Notational convention:}
We will very often have cause
to bound linear functionals $L$ on $C^{\infty}_{\omega}(\quotG)$
by Sobolev norms. In writing statements of the form
$|L(f)| \ll S_{p,d,\beta}(f)$, we will always allow
the implicit constant to depend on $p,d$ and $\beta$ without
explicitly saying so.

\subsection{Adelic Sobolev norms -- a slight generalization.}\label{subsec:ocean11}
The notations of this section will only be required in Sec. \ref{sec:torus2}.
We recommend it be omitted at a first reading.

In the discussion at the start of Sec. \ref{subsec:adelicsobolev}, 
we did not address what class of subgroups $\Gamma_i$ to consider (should we take {\em all} finite index subgroups of a fixed $\Gamma_0$ or some subclass?)  Implicitly,
such a choice was made in defining the Sobolev norms of the previous section. 
The Sobolev norms introduced in the previous section are good for most of our purposes.
However, roughly speaking, they have the following defect: they only measure
the index of a stabilizer of a function $f \in C^{\infty}_{\omega}(\quotG)$.

That this definition might
lead to some peculiar results can be already seen in
the case $\G = \PGL(2)$, $F=\Q$. Let $\chi_p$ be the character of $\adele_{\Q}^{\times}/\Q^{\times}$
that corresponds to the quadratic Dirichlet character of $\Q$ with conductor $p$, a prime number. 
Then the function $g \mapsto \chi_p(\det(g))$ descends to a function $f$ on $\quot$,
and it is easy to check that $[\Kmax:K_f]=2$, for any $p$. 
Thus the index of this stabilizer does not reflect the conductor of the underlying representation
(which, by any reasonable definition of conductor, should grow as $p$ increases). 
 In this section we shall introduce
a slight modification of the definitions which avoids this problem. (This problem would not occur for $\SL(2)$). 

This is a purely technical matter, and it seems there is much scope for giving better
and more natural definitions.  We restrict ourselves to the case $\G = \GL(n)$.  For a finite place $\q$
and $m \geq 0$, we put  $K[\q^m] \stackrel{\mathrm{def}}{=} \mathrm{ker}(\GL(n, \order_{\q})
\rightarrow \GL(n, \order_{\q}/\varpi_{\q}^m \order_{F_{\q}})$, where $\varpi_{\q}$ is a uniformizer in $F_{\q}$. 
Now, for $\psi \in C^{\infty}_{\omega}(\quotG)$, put $K_{\q,\psi}^{*}$ to be the largest
subgroup $K[\q^m]$ which stabilizes $\psi$, and put $K_{\psi}^{*} = \prod_{\q} K_{\q,\psi}^{*}$,.

We define the $\star$-pre-Sobolev norm $PS_{p,d,\beta}^{*}$ by the rule
$$PS^{\star}_{p,d,\beta}(\psi) =
[\Kmaxg:K^*_{\psi}]^{\beta}
\sum_{\ord(\mathcal{D})  \leq d}
\|\mathcal{D} \psi\|_{L^p(\quotGad)} 
$$
and we define the $\star$-Sobolev norm $S_{p,d,\beta}^{*}$ to be the
maximal seminorm dominated by $PS_{p,d,\beta}^{*}$. 
Clearly $S_{p,d,\beta}^{*} \geq S_{p,d,\beta}$. 

The eventual purpose of this is that $S_{p,d,\beta}^{*}$ (unlike $S_{p,d,\beta}$) will never be ``too small''
on an automorphic representation whose conductor is large. This can be quantified, although we do not do so in the present document. 

%
% This is quantified
%in Lem. \ref{lem:compare}. [[change if you delete last section]] 

\begin{rem} \label{rem:sobolevnorm}
Evidently the definitions of this section and the previous are not the only
``sensible'' way of defining a notion of adelic Sobolev norms. The results of this
paper do not require any more sophisticated definition, although
this would certainly be of help in optimizing the results.

However,
  it would
be interesting
 to impose a system of Sobolev-type norms in a less
{\em ad hoc} fashion. Moreover, it would be pleasant if the system
of norms had nice interpolation properties (this often is very
helpful for getting sharp results). For example it would be nice if as one varied
$\beta$ one got a family of interpolation spaces.

We remark on a simple way of defining Hilbertian norms
which seems (more) appropriate to the adelic context. 
Let $K[\q^m]$ be as above, and let $E_{\q^m}$
be the averaging projection onto the $K[\q^m]$-fixed vectors, i.e.
$E_{\q^m}(v) = \int_{k \in K[\q^m]} k \cdot v dk$, where
the measure is the Haar probability measure. 
Then $e_{\q^m} := E_{\q^m} - E_{\q^{m-1}}$ 
is a projection.  If $\mathfrak{f} = \prod_{i} \q_i^{m_i}$ is
an arbitrary integral ideal, put $e_{\mathfrak{f}} := \prod_{i}
e_{\q_i^{m_i}}$. Now put $P(s) = \sum_{\mathfrak{f}} e_{\mathfrak{f}} \Norm(\mathfrak{f})^s$.
Then $f \mapsto \sum_{\ord(\mathcal{D}) \leq d}  \| \mathcal{D} \cdot P(s) \cdot  f\|_{L^2}$
defines a Hilbert norm which seems to have reasonably pleasant formal properties.  
In fact it is majorized (up to constants) by a norm of the type described above,.

J. Bernstein has a more canonical notion of norms on 
representation spaces of $p$-adic groups, and he
has informed me that these norms have adelic analogues. I do 
not know the relation.  The norms
arising from his constructions are Hilbertian.
\end{rem}

\subsection{Some properties and uses of the Sobolev norms.}
We briefly summarize certain results that will be used in the text. Detailed proofs are given
in Sec. \ref{sec:sobolev}. 

For general $\G, \omega$ we have:
\begin{eqnarray}
\label{eq:sobproducteq}S_{p,d,\beta}(F_1 F_2) \ll_{d} S_{2p,d,\beta}(F_1)
S_{2p,d,\beta}(F_2).  \\
\label{eq:sobgroupeq}S_{p,d,\beta}(g \cdot F) \ll \|g_{\infty}\|^{d}
\|g_f\|^{\beta} S_{p,d,\beta}(F). \end{eqnarray}

\eqref{eq:sobproducteq}, proved in Lem. \ref{lem:sobproduct}, and
\eqref{eq:sobgroupeq}, proved in Lem. \ref{lem:sobgroup},
give some basic stability properties of Sobolev norms.

Now we specialize
to some results for $\GL(2)$
and $\PGL(2)$. Let $F \in C^{\infty}(\quot \times \quot)$, let $\q$ be a prime ideal of $\order_F$, 
and suppose $F$ is invariant by $\PGL_2(\order_{F_{\q}}) \times
\PGL_2(\order_{F_{\q}})$.  Then:
\begin{multline} \label{eq:heckesoboleveq} 
\left|
\int_{\quot} F(x, xa([\q])) dx  -
\sum_{\chi^2=1}\chi([\q]) \int_{\quot} F(x,y) \chi(x) \chi(y)
d\mu_{\quot}(x)
d\mu_{\quot}(y)
\right|
\\ \ll_{\epsilon} \Norm(\q)^{\frac{2\alpha - 1}{p}+\epsilon} S_{p,d}(F).
\end{multline}

{\large \eqref{eq:heckesoboleveq}, proved in 
}Lem. \ref{lem:heckesobolev},  quantifies Hecke equidistribution. 
To understand the relation, take $F$ to be a pure tensor:
$F(x,y) = f_1(x) f_2(y)$. Then \eqref{eq:heckesoboleveq}
in effect bounds the inner product $\langle T_{\q} f_1, f_2 \rangle$,
where $T_{\q}$ is the Hecke operator corresponding to $\q$.

\subsection{Cusp forms, $L$-functions and the analytic conductor.} 

As a general remark on notation -- and a mild abuse of notation-- by {\em cuspidal representation} we
shall always mean {\em unitary cuspidal representation}.  This is automatic for $\PGL(2)$
but not for $\GL(2)$. 

\subsubsection{$L$-functions.}
Let $\pi = \otimes_v \pi_v$ be an automorphic cuspidal representation of $\GL(n)$
over $F$.  We denote by $L_v(s,\pi_v)$ the local $L$-factor
of the  representation $\pi_v$; when it causes no confusion, we will sometimes
abbreviate this to $L(s, \pi_v)$.

We write $L(s, \pi) := \prod_{v \, \mathrm{finite}} L_v(s,\pi_v)$ for the (finite part of)
the global $L$-function attached to $\pi$, and $\Lambda(s,\pi) :=
\prod_{v} L_v(s,\pi_v)$ for the (completed) $L$-function attached to $\pi$. 
\subsubsection{The analytic conductor of Iwaniec-Sarnak.} \label{iscond}
We recall the definition in the context where it will arise.
Let $\pi = \otimes \pi_v $ be a cuspidal representation of $\GL(n)$ over $F$.

For each finite place $v$ we denote by $\Cond_v(\pi)$
the conductor, in the sense of Jacquet, Piatetski-Shapiro, and Shalika,
of $\pi_v$; thus $\Cond_v(\pi) = q_v^{m_v}$, where $m_v$ is the
smallest non-negative integer such that  $\pi_v$
possesses a fixed vector under the subgroup of $\GL_n(\order_{F_v})$
consisting of matrices whose bottom row is congruent to $(0,0,\dots, 0,1)$
modulo $\varpi_v^m$.

 For each infinite place $v$, let $\Gamma_v(s) = \pi^{-s/2} \Gamma(s/2)$
 or $(2 \pi)^{-s} \Gamma(s)$ according to whether $v$ is real or complex respectively,
 and put $\deg(v) = [F_v:\mathbb{R}]$. 
Let $\mu_{j,v} \in \mathbb{C}$ satisfy
$L(s, \pi_v) =  \prod  \Gamma_v(s+\mu_{j,v})$, and
put $\Cond_v(\pi) = \prod_{v} (1 + |\mu_{j,v}|)^{\deg(v)}$.
We then put $\Cond(\pi) = \prod_{v} \Cond_v(\pi)$
(this is within a constant factor of the Iwaniec-Sarnak definition).
Moreover, we put $\Cond_{\infty}(\pi) = \prod_{v \, \mathrm{infinite}} \Cond_v(\pi)$
and $\Cond_f(\pi) = \prod_{v \, \mathrm{finite}} \Cond_v(\pi)$ (the  ``infinite'' 
and ``finite'' parts of the conductor). 

We will occasionally refer to the ``finite conductor'' of $\pi$
as the {\em ideal} $\prod_{v} \mathfrak{q}_v^{m_v}$, where $\mathfrak{q}_v$
is the prime ideal corresponding to the finite place $v$; then $\Cond_f(\pi)$
is the norm of this ideal. Hopefully the distinction between the two usages will be 
clear from context.  

%In fact, the analytic conductor $\Cond(\pi)$ has a pleasing connection with the adelic Sobolev norms:
%the logarithm of the analytic conductor of $\pi$ is equal, up to a multiplicative constant,
%to the logarithm of the minimal Sobolev norm of any vector in the space of $\pi$. 
%We will later use a special case of this result.  

\begin{rem} \label{rem:explication}
(Explication for $\GL(1)$ in the archimedean case)
Let us be slightly more explicit in the case of a unitary character $\omega$ of $\adele_F^{\times}/F^{\times}$. If $v$ is real, then there is $t \in \mathbb{R}$ such that $\omega(x) = |x|^{it}$ for $x >0$;
then $\Cond_v(\omega) \asymp (1+|t|)$. If $v$ is complex, then there is $t \in \mathbb{R}, N \in \mathbb{Z}$ such that $\omega(r e^{i \theta}) = |r|^{it} e^{i N \theta}$; then 
$\Cond_v(\omega) \asymp (1 + |t| + N)^2$.

We can heuristically summarize this:  in the real case,
$\omega$ is approximately constant in a neighbourhood of the identity of 
size $\Cond_v(\omega)^{-1}$; in the complex, case $\omega$ is approximately
constant in a disc around the identity of area $\Cond_v(\omega)^{-1}$. 
\end{rem}

\subsubsection{Cusp forms} \label{subsec:cuspformsbdd}
If $\pi$ is a cuspidal representation of $\GL_2(\adele_F)$ or $\PGL_2(\adele_F)$, 
it will be convenient to denote by $\pi_{\infty}$ the archimedean representation
(of $\GL_2(F_{\infty})$ or $\PGL_2(F_{\infty})$) that corresponds to $\pi$. 

By the dual $\widehat{\GL_2(F_{\infty})}$ or $\widehat{\PGL_2}(F_{\infty})$, 
we shall mean the space of irreducible, admissible representations.
We say a subset of this dual is {\em bounded} if the corresponding set of Langlands
parameters is bounded. We may define, in an evident way, the conductor
$\Cond(\pi_{\infty})$ for $\pi_{\infty} \in \widehat{\GL_2(F_{\infty})}$; with this
definition, a subset is bounded exactly when $\Cond$ takes bounded values on it. 

In a similar fashion, we define the notion of a bounded subset of $\widehat{GL_2(F_v)}$
or $\widehat{PGL_2(F_v)}$ for any place $v$, where, again
$\widehat{\GL_2(F_v)}$ denotes the set of irreducible, admissible representations.

\section{Unipotent periods.} \label{sec:unipotent}

In this section, we will make systematic use of the
(nonadelic) Sobolev norms $S_{\infty,d}$ on homogeneous
spaces $\Gamma \backslash G$. In rough terms, $S_{\infty,d}$ controls the $L^{\infty}$ norm of the
first $d$ derivatives. See Sec. \ref{sec:sobreal}.
Note in particular that $S_{\infty,0}$ is just the $L^{\infty}$ norm.
\subsection{Equidistribution of  sparse subsets of horocycles.}
\label{subsec:sparsehoro}Let $\Gamma \subset \SL_2(\R)$ be a cocompact lattice.
For this section alone, we will use mildly different notation to that of Sec. \ref{sec:gl2pgl2},
to accommodate the fact we deal with $\SL_2$ and not with $\PGL_2$. 
For $x \in \R$, put 

\begin{equation} \label{eqn:nan} n(x) =\left( \begin{array}{cc}1 & x \\ 0 &
1
\end{array}\right), a(x) = \left(\begin{array}{cc} x^{1/2} & 0 \\ 0 &
x^{-1/2}\end{array}\right), \nbar(x) = \left(\begin{array}{cc}
1 & 0 \\ x & 1 \end{array}\right).\end{equation}

We denote by $C(\Gamma \backslash \SL_2(\R))$
(resp. $C^{\infty}(\Gamma \backslash \SL_2(\R))$)
the space of continuous (resp. smooth) functions on the compact real
manifold $\Gamma \backslash \SL_2(\R)$.
We denote by $dg$ the measure on $\SL_2(\R)$
that descends to a probability measure on the quotient
$\Gamma \backslash \SL_2(\R)$.
Finally, we denote by $\mathbb{H}$
the usual upper half plane
$\{z: \Im(z) > 0\}$
with the standard action of $\SL_2(\R)$.

\begin{thm} \label{thm:sparsehoro}
There exists $\gamma_{\mathrm{max}} > 0$, depending on $\Gamma$,
such that $\{x_0 n(j^{1+\gamma}): j \in \mathbb{N}\}$
is equidistributed, for
any $x_0 \in \Gamma \backslash\SL_2(\R)$ and any $0
\leq \gamma < \gamma_{\mathrm{max}}$.  In other words,
for any $f \in C(\Gamma \backslash \SL_2(\R))$,

$$\lim_{N \rightarrow \infty} \frac{\sum_{j=1}^{N} f(x_0 n(j^{1+\gamma}))}{N}
= \int_{\Gamma \backslash \SL_2(\R)} f(g) dg.$$

If $\lambda_1$ is the smallest nonzero
eigenvalue of the Laplacian on $\Gamma \backslash \mathbb{H}$,
put
 $\alpha = \begin{cases} 0,  & \lambda_1 \geq 1/4 \\
\sqrt{1/4 - \lambda_1}, & \mbox{ else.} \end{cases}.$ Then we can
take $\gamma_{\max} = \frac{(1-2\alpha)^2}{16(3 - 2 \alpha)}$.
\end{thm}

This result represents (extremely modest) progress towards a
conjecture of N.Shah, which asserts that the statement should
remain valid for {\em any} $\gamma > 0$. The method is not restricted to
sequences of the specific type in Thm. \ref{thm:sparsehoro}, and
we have also not optimized the maximal value for
$\gamma_{\mathrm{max}}$. Nevertheless the method is fundamentally
limited. As it presently stands, it does not seem capable of achieving even $\gamma = 1$. See 
also Remark \ref{rem:effratner} (page \pageref{rem:effratner}). 

The dependence of $\gamma_{\max}$ on $\Gamma$ can likely be removed, but this seems to require
using further input (cf. last paragraph of Section \ref{subsec:concrete}).

The proof follows the line of Sec. \ref{sec:abstract},
with $G_1 = \mathrm{SL}_2(\R)$, $G_2 = \{n(x): x  \in \mathbb{R}\}$.
The $Y_i$ are not quite closed $G_2$-orbits, but rather
long pieces of general $G_2$-orbits. The
basis $\{\psi_{i,j}\}$ for $Y_i$ will correspond
to additive characters of $G_2 \cong \mathbb{R}$.

Let $f \in C^{\infty}(\Gamma \backslash \mathrm{SL}_2(\mathbb{R}))$,
and let $\alpha$ be as in the statement of Thm. \ref{thm:sparsehoro}.
Let $T \geq 1$.
Let $\psi$ be a fixed nontrivial character of the additive group of $\mathbb{R}$.
Let $g$ be a fixed smooth function of compact support on $\mathbb{R}$
satisfying $\int_{-\infty}^{\infty}g(x) dx = 1$. 
We denote by $\langle \cdot, \cdot \rangle_{L^2(\Gamma \backslash G)}$
the inner product in the Hilbert space $L^2(\Gamma \backslash G)$.

We set:
\begin{equation} \label{protomunu}
\nu_T(f) = \frac{1}{T}\int_{0}^{T} f(x_0 n(t)) dt,
\ \mu_{T,\psi}(f) = \frac{1}{T}\int_{0}^{T} \psi(t) f(x_0 n(t)) dt.\end{equation}

Remark first that the measures $\nu_T$ are equidistributed as $T \rightarrow \infty$,
in the following quantitative sense: for $f \in C^{\infty}(\Gamma \backslash
\SL_2(\R))$,
\begin{equation} \label{eq:protoass}\left|\nu_T(f) - \int_{\Gamma \backslash
\SL_2(\R)}f(g)dg\right| \ll T^{-\kappa_1} S_{\infty,
1}(f),\end{equation} for any $\kappa_1 < \frac{1/2-\alpha}{2}$.
This is proven in Lem. \ref{lem:longhoroequi}, without taking any pains to optimize the  exponent.
(We prove it to keep the paper self-contained. However, we emphasize that neither result nor proof is new;  see \cite{Ratner1} and \cite{Ratner2}.  A precise analysis of the equidistribution of long horocycles
is carried out in \cite{FlamForni}.)
\begin{lem} \label{lem:ost}
Suppose $\int_{\Gamma \backslash \SL_2(\R)} f(g) dg = 0$. Then:
\begin{equation}  \label{eq:muboundone}
|\mu_{T,\psi}(f)| \ll T^{-b} S_{\infty,1}(f),\end{equation} whenever
$b < \frac{ (1-2\alpha)^2}{8 (3- 2\alpha)}$ and the implicit
constant is independent of $\psi$.
\end{lem}
Remark that if $\psi$ is wildly oscillatory,
cancellation in $\mu_{T,\psi}$ can
be proved directly by integration by parts;
on the other hand, if $\psi$ is almost constant,
the cancellation in $\mu_{T,\psi}$
arises from the equidistribution of the
horocycle $x_0 G_2$. It is therefore the intermediate
case in which (\ref{eq:muboundone}) is of interest.

\proof Let $H \geq 1$, and let $\sigma_H$  be the measure on $N(\mathbb{R})$ defined
by $\sigma_H(g) = \frac{1}{H}\int_{0}^{H} \psi(x) g(n(x))
dx$, for $g$ a function on $N(\mathbb{R})$.

For $f \in C^{\infty}(\Gamma \backslash G)$, we denote by $f \star
\sigma_H$ the right convolution of $f$ by $\sigma_H$. Then it is
easy to verify that $\left|\mu_{T,\psi}(f) - \mu_{T,\psi}(f \star
\meas_H)\right| \ll \frac{H}{T}S_{\infty,0}(f).$  Here $\star
\meas_H$ denotes right convolution by $\meas_H$. On the other
hand, by Cauchy-Schwarz, $|\mu_{T,\psi}(f \star \meas_H)|^2 \leq
\nu_T(|f \star \meas_H|^2 )$. Thus, expanding $|f \star
\meas_H|^2$ and applying (\ref{eq:protoass}), we conclude
\begin{equation}\begin{aligned}
|\mu_{T,\psi}(f)| \ll \frac{H}{T} S_{\infty,0}(f) +
\left(\frac{1}{H^2} \int_{(h_1, h_2) \in [0,H]^2}
\left|
\nu_T(n(h_1)f \cdot \overline{n(h_2) f})\right| dh_1 dh_2 \right)^{1/2}
\\ \ll  \frac{H}{T} S_{\infty,0}(f)
+ \left(\frac{1}{H^2} \int_{(h_1, h_2) \in [0,H]^2} \left|\langle
n(h_1-h_2) f,f \rangle_{L^2(\Gamma \backslash G) } \right| dh_1
dh_2 \right)^{1/2}
\\ + \left(T^{-\kappa_1} \sup_{(h_1, h_2) \in [0,H]^2} S_{\infty,1}(n(h_1) f
\cdot \overline{n(h_2) f})\right)^{1/2}
\end{aligned}\end{equation}
Utilising bounds towards matrix coefficients -- see Sec. \ref{subsec:sl2mc}, esp. (\ref{eq:sl2mc}) -- and
basic properties of Sobolev norms (see 
\footnote{Lem. \ref{lem:realsob}
would actually give the exponent $(1+|h_1|+|h_2|)^4$ in the latter inequality, but it is
easy to see directly the stronger result.} Lem. \ref{lem:realsob}), we note:
\begin{multline}\langle n(h) f, f \rangle \ll_{\epsilon} (1+|h|)^{2 \alpha-1+\epsilon}
S_{\infty,1}(f)^2,
\\ |S_{\infty,1}(n(h_1)f \cdot \overline{n(h_2) f})|
\ll (1+|h_1|+|h_2|)^2 S_{\infty,1}(f)^2.\end{multline} Thus,
$|\mu_{T,\psi}(f)| \ll_{\epsilon} \left(\frac{H}{T} + H^{\alpha-1/2+\epsilon} +
T^{-\kappa_1/2} H\right) S_{\infty,1}(f).$ Choose $H$ so that
$H^{\alpha - 1/2} = H T^{-\kappa_1/2}$ to 
obtain the claimed result. \qed

\proof (of Thm. \ref{thm:sparsehoro}). 
Given Lem. \ref{lem:ost}, the Theorem follows quite readily
by Fourier-expanding the measure on $\mathbb{R}$
that is a sum of point masses at $j^{\gamma}$, for $j \in \mathbb{N}$.
The argument that follows formalizes a minor variant of this argument (we first
consider instead a sum of point masses along arithmetic progressions
which approximate $\{j^{\gamma}: j \in \mathbb{N}\}$). 

Let $x_0 \in \Gamma
\backslash \SL_2(\R)$, $f \in C^{\infty}(\Gamma \backslash
\SL_2(\R))$. We first claim that, if $b$ is as in the previous
Lemma, $f$ so that $\int_{\Gamma \backslash \SL_2(\R)} f(g) dg
=0$, and $K \geq 1$, then:
\begin{equation} \label{eq:quickequi}
\frac{\sum_{0 \leq j < K^{1/b-1}} f(x_0 n(K
j))}{K^{1/b-1}} \rightarrow
0,\end{equation}
as $K\rightarrow \infty$.  In other words, $K^{1/b-1}$ points,
distributed along a horocycle with spacing $K$, become equidistributed.

This follows from Lem. \ref{lem:ost}: put $g_{\delta}(x) =
\max(\delta^{-2}( \delta - |x|),0)$, a function on $\mathbb{R}$.
For $\lambda \in \mathbb{R}$, write $a_{\lambda} =
K^{-1} \int_{\mathbb{R}} \exp(-2 \pi i K^{-1} \lambda t) g_{\delta}(t)
dt$. Then $\sum_{j \in \Z} g_{\delta}(t + K j) = \sum_{k \in \Z}
\exp(2 \pi i K^{-1} k t) a_k$. Moreover, a simple computation
shows that $\sum_{k \in \Z} |a_k| \ll  \delta^{-1}$. 

Choose $\varepsilon > 0$ so that $b+ \varepsilon$ still satisfies
the inequality of Lem. \ref{lem:ost}. 
By Lem.
\ref{lem:ost},
\begin{equation} \label{eq:raghu}\left| \int_{t=0}^{T} dt \left(\sum_{j \in \Z}
g_{\delta}(t + Kj)\right) f(x_0 n(t)) \right| \ll_f T^{1-b-\varepsilon} \sum_{k} |a_k|
\ll T^{1-b-\varepsilon} \delta^{-1}.\end{equation}

$g_{\delta}$ has integral $1$ and is supported in
a $\delta$-neighbourhood of $0$; in particular, 
the left-hand side of (\ref{eq:raghu})
differs from $\sum_{j \in \Z, 0 \leq  K j \leq T} f(x_0 n(K j))$
by an error that is $\ll_{f} (1 + T K^{-1} \delta)$. 
Thus $ \left| \sum_{j \in \Z, 0 \leq K j \leq T} f(x_0 n(K j)) \right|
\ll_{f} (1 + T K^{-1} \delta + T^{1-b-\varepsilon} \delta^{-1})$ from which
 (\ref{eq:quickequi}) readily follows.

We now deduce the Theorem from (\ref{eq:quickequi}). 
Let $T_0 \in \mathbb{N}$ be large.
Then, for $t$ small,
we have $(T_0 + t)^{1+\gamma} = T_0^{1+\gamma} + (1+\gamma) T_0^{\gamma} t
 + O(t^2 T_0^{\gamma-1})$. In particular,
$(T_0+t)^{1+\gamma}$ is well-approximated by the linear function
$T_0^{1+\gamma} + (1+\gamma) T_0^{\gamma}t$ in the range
where $|t| \ll T_0^{\frac{1-\gamma}{2}}$.

The claim of Thm. \ref{thm:sparsehoro} follows from
(\ref{eq:quickequi}) as long as $\frac{1-\gamma}{2\gamma}> 1/b-1$;
in particular, any $\gamma < b/2$ will do.  \qed

\begin{rem} \label{rem:effratner}
The applicability of this method is not, of course, restricted to
sequences of the form $t^{1.01}, \ t \in \mathbb{N}$. In certain
(very artificial) settings, one can prove some effective instances
of Ratner's theorem in non-horospherical cases by the same
technique. The question of proving such results, even in very
special cases, was raised by Margulis in his talk in the workshop
``Emerging applications of measure rigidity.''

 For instance, let $r > 1$ be an integer, and
let $\Gamma$ be a cocompact lattice in
$\mathrm{SL}_2(\mathbb{R})^r$. For concreteness, let us suppose
that $\Gamma$ is a congruence lattice. Let $\mathbf{n}(x_1, \dots,
x_r) = (n(x_1), n(x_2), \dots, n(x_r)) \in \SL_2(\R)^r$. Then it
is well-known that one can give an effective version of the
equidistribution of $\mathbf{n}(\R^r)$, since $\mathbf{n}(\R^r)$
is horospherical. Using the method described above, one can extend
this slightly: let $V \subset \R^r$ be a linear subspace. One may
show that $\mathbf{n}(V)$-orbits are effectively equidistributed
if $\mathrm{dim}(V)/r$ is sufficiently close to $1$. However, the
small codimension condition is crucial for this method and it does
not seem that it would extend beyond this case.

\end{rem}
\subsection{Fourier coefficients of automorphic forms.}
\label{sec:fourcoeff}

In the notations of Sec. \ref{sec:intro},
if we take for $Y_i$ the closed orbits
of a unipotent group, the resulting
periods are so-called ``Fourier coefficients
of automorphic forms.'' We shall
give a general nontrivial bound in that context.
(The word ``nontrivial'' must be interpreted with care;
see the discussion at the end of this section.)
Our methods are restricted
to the case of horospherical unipotent subgroups.

We have made no effort to optimize the exponents
of the results, nor even to state
a result of maximal generality. In fact, one can
considerably increase the scope of Thm. \ref{thm:fouriercoefficients},
since we deal in the present
section only with {\em closed orbits of horospherical subgroups},
one can profitably apply spectral theory. We do not carry this out,
instead using \cite{KM} to give equidistribution statements in a fairly soft
fashion.

Let $G$ be a connected semisimple real Lie group,
$\Gamma \subset G$ a lattice, $K \subset G$ the maximal
compact subgroup,
$\lieg$ the Lie algebra of $G$,
and $H \in \lieg$ a semisimple element. Fix a norm $\|\cdot\|$ on the real vector space $\lieg$. 
Let $\exp: \lieg \rightarrow G$ be the exponential map.
Fix a Haar measure on $G$ so that $\Gamma \backslash G$ has volume $1$. 
Let $\mathfrak{u}$ be the sum
of all negative root spaces for $H$ 
and let $U = \exp(\mathfrak{u}) \subset G$.
Let $x_0 \in \Gamma \backslash G$
be so that $x_0 U$ is compact.
Let $x_t = x_0 \exp(t H)$, and let $\Delta_t$
be the stabilizer of $x_t$ in $U$.
We denote by $\langle \cdot, \cdot \rangle_{L^2(\Gamma \backslash G)}$
the inner product in the Hilbert space $L^2(\Gamma \backslash G)$.

We shall analyze periods of a fixed function
along $x_t U$ as $t$ varies.
The proofs follow Sec. \ref{sec:abstract}
with $G_1 = G, G_2 = U$, $Y_i = x_t U$,
and $\psi_{i,j}$ corresponding to characters of $U$.

Let $T > 0$ and let $\psi$ be any character
of $U$ trivial on $\Delta_T$ (:= $\exp(-TH) \Delta_0 \exp(TH)$).
We define $\nu_T, \mu_{T,\psi}$ in a closely analogous fashion to
(\ref{protomunu}):
\begin{equation} \label{eq:lipatti}\nu_T(f) =
\frac{\int_{\Delta_T \backslash U} f(x_T u)
du}{\mathrm{vol}(\Delta_T\backslash U)},
\mu_{T,\psi}(f) = \frac{\int_{\Delta_T \backslash U} f(x_T u) \psi(u) du}
{\vol(\Delta_T \backslash U)} .\end{equation}

Let $f, g \in C^{\infty}(\Gamma \backslash G)$. 
Let $E \in \mathfrak{u}$ have unit length w.r.t. the fixed norm $\|\cdot\|$ on $\lieg$. 
It is proven by Kleinbock-Margulis in \cite{KM} -- see also Lem. \ref{lem:assmix} and Lem. \ref{lem:banana} -- that
there are $\kappa_1, \kappa_2 > 0$ such that:
\begin{equation} \begin{aligned} \label{eq:assmix}
%\left|\langle \exp(tH) \cdot f, g \rangle
%- \int_{\Gamma \backslash G} f \int_{\Gamma \backslash G} g\right|
%\leq \exp(-\kappa_1|t|)  S_{\infty,\dim(K)}(f) S_{\infty ,\dim(K)}(g)\\
\left| \langle \exp(s E) \cdot f, g \rangle_{L^2(\Gamma \backslash G)} - \int_{\Gamma
\backslash G} f \int_{\Gamma \backslash G} g \right|
\ll (1+|s|)^{-\kappa_1} S_{\infty,\dim(K)}(f) S_{\infty,\dim(K)}(g)
\end{aligned}\end{equation}
\begin{equation} \label{eq:banana}|\nu_T(f) - \int_{\Gamma \backslash G} f| \ll e^{-\kappa_2T} S_{\infty,\dim(K)}(f).\end{equation}
(\ref{eq:assmix}) and (\ref{eq:banana}) assert, respectively, quantitative mixing of the flow generated by $U$ on $\Gamma \backslash G$,
and the equidistribution of the orbit $x_T U$ as $T \rightarrow \infty$. 
We may assume that $\kappa_1 < 1$ (since making $\kappa_1$ smaller
does not change the truth of (\ref{eq:assmix})). This will ease the notation
 in the proof of the Theorem.

\begin{thm}\label{thm:fouriercoefficients}

There exists $\kappa_3 > 0$ such that,
for any $f \in C^{\infty}(\Gamma \backslash G)$
satisfying $\int_{\Gamma \backslash G} f =0$, we have:
\begin{equation}\label{eq:fourcoeff}
|\mu_{T,\psi}(f)| \ll \exp(-T \kappa_3) S_{\infty,\dim(K)}(f)
\end{equation}
for all $T \geq 0$, and for all characters $\psi$ of $U$
trivial on $U_T$.

Indeed, if $o$ is the order of the polynomial map
$\mathbb{R} \rightarrow \mathrm{End}(\lieg)$ defined by
$s \mapsto \mathrm{Ad}(\exp(s H))$, then any $\kappa_3 < \frac{\kappa_1 \kappa_2}{2(2o \dim(K)+\kappa_1)}$ 
is admissible,
$\kappa_{1,2}$ being as in (\ref{eq:assmix}) and (\ref{eq:banana}). 
\end{thm}

The relevance to this to ``Fourier coefficients'' in the classical sense may not be immediately
clear; after the proof, we give the example of $\mathrm{SL}_2(\mathbb{R})$
to illustrate. 

Also, observe that the estimate (\ref{eq:fourcoeff}) is uniform in $\psi$. In fact, just as in
Lem. \ref{lem:ost}, the case when $\psi$ is  constant amounts to (\ref{eq:banana}), whereas
the case where $\psi$ is highly oscillatory could be handled by integration by parts. 
It is, again, the intermediate case where (\ref{eq:fourcoeff}) has content. 

\proof
We first remark that, other than being
in a slightly more general setting,
the proof
is almost exactly the same
as the proof of Lem. \ref{lem:ost}. 

The signed measure $\mu_{T,\psi}$
satisfies $\mu_{T,\psi}( u \cdot f) = \overline{\psi(u)} \mu_{T,\psi}(f)$,
for $u \in U$.

Take $E \in \mathfrak{u}$ of unit length w.r.t. the norm $\|\cdot\|$ on $\lieg$. 
Fix $H \geq 1$. Let $\meas$ be the measure
on $U$ defined via the rule
$$\meas(h) = \frac{1}{H} \int_{s=0}^{H} \psi(\exp(sE)) h(\exp(sE))ds,$$
for $h$ any continuous compactly supported function on $U$.
Then $\mu_{T,\psi}(f)  =\mu_{T,\psi}(f \star \meas)$; thus:
\begin{multline}|\mu_{T,\psi}(f)|^2 = |\mu_{T,\psi}(f \star \meas)|^2 \leq \nu_T(|f \star
\meas|^2) \\ \ll 
\int_{\Gamma \backslash G} |f \star \meas|^2
+ e^{-\kappa_2 T} S_{\infty,\dim(K)}(|f \star \meas|^2) \\ \leq
\int_{\Gamma \backslash G} |f \star \meas|^2
+ H^{2o \dim(K)}  \exp(-\kappa_2 T) S_{\infty,\dim(K)}(f)^2.\end{multline}
where we have applied Cauchy-Schwarz followed
by (\ref{eq:banana}),
noting that by Lem. \ref{lem:realsob}, we have
that $S_{\infty,\dim(K)} (|f \star \meas|^2) \ll
S_{\infty,\dim(K)} (f\star \meas)^2  \ll H^{2o \dim(K)} S_{\infty,\dim(K)}(f)^2$.
Here $o$ is chosen as in the statement of the Theorem.

By (\ref{eq:assmix}), we see that \begin{multline} \int_{\Gamma \backslash G} |f \star \meas |^2  \ll 
\left( \frac{1}{H} \int_{0}^{H} (1+|t|)^{-\kappa_1} dt \right) S_{\infty,\dim(K)}(f)^2\\  \ll H^{-\kappa_1}
S_{\infty,\dim(K)}(f)^2.\end{multline} Indeed, this follows simply by expanding the leftmost expression. Thus
$$|\mu_{T,\psi}(f)|^2 \ll  (H^{-\kappa_1} + H^{2o \dim(K)}\exp(-\kappa_2 T)) S_{\infty,\dim(K)}(f)^2.$$
We choose $H$ so that $H^{2o \dim(K)+\kappa_1} = \exp(\kappa_2 T)$ to conclude.
\qed

\begin{rem}
We now explain, when we specialize $G =\mathrm{SL}_2(\R)$,
why this recovers Sarnak's result \cite{SaIMRN},
which was the first improvement of the Hecke bound for nonarithmetic groups.
Take $H = \left(\begin{array}{cc} -1 & 0 \\ 0 & 1 \end{array}\right)
\in \mathfrak{sl}_2$, so that $U =
 \left(\begin{array}{cc} 1 & * \\ 0 & 1 \end{array}\right)$.
Let $\Gamma \subset G$  be a nonuniform lattice so that
$\Gamma \cap U = \{\left(\begin{array}{cc} 1 & n \\ 0 &1 \end{array}\right),
n \in \Z\}$, and take
$x_0 = \left(\begin{array}{cc} 1 & 0 \\ 0 & 1 \end{array}\right)$.
Let $f(x+iy) = \sum_{n \neq 0}
a_n\sqrt{y} K_{i \nu}(2 \pi n y)
e^{2 \pi i n x}$ be a Maass cusp form
of eigenvalue $1/4 +\nu^2$ on $\Gamma \backslash \mathbb{H}$,
where $\mathbb{H}$ denotes the upper half-plane; it lifts
to a function on $\Gamma \backslash \SL_2(\R)$,
viz. $g \mapsto f(g.i)$. 

Then the Theorem implies (in concrete language) that
there exists $\delta > 0$ such that, for any $y \leq 1$
and any $n \in \mathbb{Z}$,
\begin{equation} \label{eq:horo}
\int_{0 \leq x \leq 1} f(x+iy) e(nx) dx \ll y^{\delta}.\end{equation}
Taking $y \asymp n^{-1}$ in (\ref{eq:horo}),
one easily deduces that the Fourier coefficients $a_n$
satisfy the ``nontrivial'' bound $|a_n| \leq n^{1/2 - \delta}$.
\end{rem}

\begin{rem}
The bound Thm. \ref{thm:fouriercoefficients}
is nontrivial in that
it improves, as $t \rightarrow \infty$ on the trivial bound:
$$\left|\frac{\int_{\Delta_t \backslash U} f(x u) \psi(u) du}
{\vol(\Delta_t \backslash U)} \right| \ll_{f} 1,$$
which follows from Cauchy-Schwarz.

However, if there is another
interpretation for the Fourier coefficients,
it is not always the case that Thm. \ref{thm:fouriercoefficients}
improves on ``trivial'' bounds arising from that interpretation.

For instance, Fourier coefficients of cusp forms
on $\mathrm{GL}_n$ admit a spectral interpretation, that is to say, they are connected
to the eigenvalues of Hecke operators. 
In that case, Thm. \ref{thm:fouriercoefficients}
does not give anything even approaching the bounds of Jacquet-Piatetski-Shalika.
Another example is when $G = \widetilde{\SL}_2(\R)$,
and one takes for $f$ the Shimura lift
of a cusp form of integral weight. In that case
the (absolute values of the squares of) square-free Fourier coefficients of $f$
are given by special values of a twisted $L$-function;
but the estimate above does not even recover the convexity bound
(in fact, the method as indicated {\em cannot} recover this bound, even under
optimal assumptions.)

It seems as though, in these cases,
there is extra cancellation in the unipotent integrals for
subtle arithmetic reasons. The crude methods indicated
above do not detect this.
\end{rem}
\begin{rem}
We remark that, in the proof just given, the constant $\kappa_3$
depends on the spectral gap for $\Gamma \backslash G$. This
dependence can very likely be removed in many cases, including the
case of $G =\mathrm{SL}_2(\mathbb{R})$, but we do not carry this
out; again, cf. the last paragraph of Section \ref{subsec:concrete}. In the higher rank case, if $G$ has property (T), one has in any case a uniform spectral gap
and this point becomes irrelevant. 

It seems worthwhile to remark that, whereas the proof above
is clearly not unrelated to that of Sarnak, it does not
require any information on the decay of triple products (in
particular, the deep ``exponential decay'' results proved by
Sarnak and Bernstein-Reznikov).

We also remark that the proof indicated above,
although it can be optimized in various ways, probably does
not lead to as good an exponent as the work of Good,
and the later refinement of Sarnak's result due to Bernstein-Reznikov.
Its advantage lies, rather,  in its robustness and general applicability. 
\end{rem}
\section{Semisimple periods: triple products in the level aspect.}
\label{sec:tripleproducts}

In this section, we will give bounds for the triple product period
on $\PGL_2$ over a number field $F$. We will use the notation
of Sec. \ref{sec:notn}; in particular, $\adele_F$ is the adele ring of $F$ and
and $\quot = \PGL_2(F) \backslash \PGL_2(\adele_F)$. 

In Sec. \ref{subsec:periodbound}, we will give a special case of the triple product bound (Prop. \ref{prop:main}) 
which does not require Sobolev norms to state. For some applications
we will require a slight generalization, which will require the Sobolev norms of
Sec. \ref{subsec:adelicsobolev}. This will be carried out in Sec. \ref{subsec:technical}
(see Prop. \ref{prop:mainV}).  

\subsection{Period bound for triple products.}

\label{subsec:periodbound}
 We now give a period bound for triple
products on $\PGL_2$. Although it is unfortunately somewhat
disguised in the adelic language, the situation and method
corresponds to that of Sec. \ref{sec:abstract} with $G_1 =
\PGL_2(F_S) \times \PGL_2(F_S)$, $G_2 =\PGL_2(F_S)$ embedded
diagonally. Here $S$ is a set of places of $F$ containing all
infinite places, and $F_S = \prod_{v \in S} F_v$.

For $1 \leq p \leq \infty$ we will write $L^p$
for $L^p(\quot)$. Thus, e.g., $\|f_1\|_{L^4}$
denotes $\left(\int_{\quot} |f_1(x)|^4 dx\right)^{1/4}$.

%Let $\pi$
% be an automorphic cuspidal representation of $\PGL_2(\adele_F)$
%The conductor $\cond := \cond(\pi)$ of $\pi$
% is an integral ideal of $\order_F$; we will assume $\cond$
%prime for simplicity of exposition.

\begin{prop} \label{prop:main} (``Subconvexity for the triple product period.'') 
Let $\pi$ be an automorphic cuspidal representation of $\PGL_2(\adele_F)$
with prime finite conductor $\cond$. 
Let $f_1, f_2 \in C^{\infty}(\quot)$
be totally nondegenerate\footnote{See Sec. \ref{subsec:proj}; equivalent to ``orthogonal to locally constant functions'' in this case.} and such that 
$f_1, f_2$ are $\PGL_2(\order_{F_{\cond}})$-invariant. 
Let $\varphi \in \pi$ and
suppose\footnote{Recall that a finite place $v$ belongs to the support of $f \in C^{\infty}(\quot)$
exactly when $\PGL_2(\order_v)$ does not fix $f$.} that there exists $b \in \mathbb{R}$ such that\footnote{This assumption 
 (\ref{eq:cond}) is purely technical and the reader
may safely assume that $\varphi$ is spherical at all places away from
$\cond$ and $f_1, f_2$ are everywhere spherical without losing the gist of the argument. It is not used in the present document, but will probably
be of use in establishing polynomial dependence of subconvex bounds.
 (\ref{eq:cond}) ensures, among other things, that there are many places when all of $f_1, f_2, \varphi$ are unramified, so that we can use
the Hecke operators at those places.}
\begin{equation} \label{eq:cond}
\prod_{\mathfrak{q} \in 
\Supp(\varphi) \cup \Supp(f_1) \cup \Supp(f_2) } \Norm(\q)
\leq \Norm(\cond)^b\end{equation}
Put $I(\varphi) = \int_{\quot} f_1(g)
f_2(g \f) \varphi(g) dg,$ where $\f$ is as in (\ref{eq:fconddef}) 
and $dg$ is the $\PGL_2(\adele_F)$-invariant probability measure. 
Then \begin{equation} \label{eqn:est}|I(\varphi)| \ll_{b,\epsilon, F}
\|f_1\|_{L^{4}} \|f_2\|_{L^{4}} \|\varphi\|_{L^2}
\Norm(\cond)^{\epsilon - \frac{(1 - 4\alpha) (1-2 \alpha)}{4 (3 - 4 \alpha)}} \end{equation}\end{prop}

We refer to Prop. \ref{prop:main} as {\em subconvexity for the triple product period}, cf. first assertion
of Theorem \ref{thm:subconvextp}. 
We note that, with $\alpha = 3/26$ (Kim's bound) we have $\frac{(1-4 \alpha) (1- 2 \alpha)}{4(3-4 \alpha)}
> 1/26$. 
As we will see, Prop. \ref{prop:main} is a very strong result that implies many subconvexity results on $\PGL(2)$. 

First let us explain the content of Prop. \ref{prop:main} in a classical setting,
and how it can be regarded as the 
type of period bound discussed in the introduction. Suppose $F= \Q$;
let $p \geq 1$ and
let $$\Gamma_0(p) = \{ \left( \begin{array}{cc} a& b\\c& d \end{array} \right) \in \GL_2(\mathbb{Z}):
p|c \},$$ and put $Y(p) = \Gamma_0(p) \backslash \PGL_2(\mathbb{R})$.
Then there is an embedding $Y(p) \rightarrow Y(1) \times Y(1)$
which corresponds to the graph of the $p$th Hecke correspondence on $Y(1)$;
the image is a certain closed orbit of the diagonal $\PGL_2(\R)$. 
Let $f_1, f_2$ be fixed functions on $Y(1)$ and $\varphi$ a Maass form
on $Y_0(p)$. Then the function $f_1 \times f_2: (x_1,x_2) \mapsto f_1(x_1) f_2(x_2)$
is a function on $Y(1) \times Y(1)$, and we can construct
its restriction $f_1 \times f_2|_{Y(p)}$ by means of the embedding indicated above. 
 Then, translating from adelic to classical, one finds that (\ref{eqn:est}) furnishes precisely an estimate
for $\int_{Y(p)} (f_1 \times f_2)|_{Y(p)} \varphi$. When we vary $p, \varphi$ and hold $(f_1, f_2)$ fixed, such an estimate
falls precisely into the pattern described in the introduction: we are computing
the periods of the fixed function $f_1 \times f_2$ along the varying sequence of sets $Y(p)$.
The fact that the $Y(p)$ become equidistributed in $Y(1) \times Y(1)$ is precisely
equivalent to the equidistribution of  $p$-Hecke orbits on $Y(1)$. 
Moreover, the key property of $\varphi$ that is used
is the fact that $\varphi$ is an eigenfunction
of many Hecke operators; this is used to construct the measure $\sigma$, in the notation of 
Sec. \ref{sec:method}. 

Prior to beginning the proof, we make some comments about
applications and generalizations; for details, we refer to Sec. \ref{sec:app}. 
The implicit constant of (\ref{eqn:est}) is independent of $f_1$, $f_2$.
Taking $f_1, f_2$ to be a pair of cusp forms 
Prop. \ref{prop:main}
implies (conditional on some computations of $p$-adic integrals that we state 
as Hypothesis Prop. \ref{prop:intrep}) subconvexity for certain triple product $L$-functions.
Similarly, taking $f_1, f_2$ to be
a cusp form and an Eisenstein series, resp. a pair of Eisenstein series, Prop. \ref{prop:main}
implies subconvexity for Rankin-Selberg convolutions and standard $L$-functions. 

The latter applications are rather delicate because Eisenstein
series are not in $L^4$. To get around this we will eventually replace the Eisenstein series
by an appropriate wave-packet (cf. proof of Thm. \ref{thm:subconvextp}).

%It is not critical that $\pi$ be cuspidal:
%if $\pi$ is an Eisenstein representation,
%the fact that vectors in $\pi$ lie almost in $L^2$ suffices.

\proof
Clearly we may assume that $\|\varphi\|_{L^2} = \|f_1\|_{L^4}
= \|f_2\|_{L^4} = 1$. It follows that
$\|f_1\|_{L^2} \leq 1$ and $\|f_2\|_{L^2} \leq 1$.

We shall moreover assume, for simplicity, that
$\varphi$ is spherical at all finite places $v \neq \cond$.
The reader may verify that the proof carries through
to the more general situation of Prop. \ref{prop:main} without modification.

Put $q = \Norm(\cond)$.
Let $\meas$ be a (signed real) measure on $\PGL_2(\adele_F)$
such that $\varphi \star \check{\meas} = \lambda \varphi$,
for some $\lambda \in \C$. We shall assume
that $\mathrm{supp}(\meas)$ commutes with $\PGL_2(\Q_{\cond})$; we will choose $\meas$ later. 
Set further $\Psi(x) = f_1(x) f_2(x  \f) \in C^{\infty}(\quot)$.
Then

\begin{equation} \label{eq:john1}\begin{aligned}\lambda \cdot I(\varphi) =
\int_{\quot} \Psi(x) \cdot (\varphi \star \check{\meas})(x) dx =
\int_{\quot} (\Psi \star \meas)(x) \cdot \varphi(x) dx
 \leq \left(\int_{\quot} |\Psi \star \meas|^2 dx\right)^{1/2}
\\ = \left(\int_{\quot} \int_{(g,g') \in \PGL_2(\adele_F)^2} \left(g \cdot \Psi\right)
\overline{\left(g' \cdot \Psi\right)}
d\meas(g) d\meas(g') dx \right)^{1/2}
\\ = \left(\int_{\quot} \int_{(g,g') \in \PGL_2(\adele_F)^2}f_1(xg) f_2(x \f g)
\overline{f_1(xg') f_2(x \f g')} d\meas(g) d\meas(g') dx
\right)^{1/2}  \\ =
\left(\int_{\quot} \int_{(g,g') \in \PGL_2(\adele_F)^2}f_1(xg)f_2(x g\f)
\overline{f_1(xg')f_2(x  g'\f)} d\meas(g) d\meas(g') dx
\right)^{1/2}
\end{aligned} \end{equation}
In the last step, we have used the fact that $\PGL_2(\Q_\cond)$,
and thus $\f$, commutes with $\mathrm{supp}(\meas)$.

For any two functions $h_1, h_2$ on $\quot$, both
right invariant by $\PGL_2(\order_{F_{\cond}})$,
the assumed bound on Ramanujan (Def. \ref{def:ramanujan}) implies:
\begin{equation} \label{eq:hecke}\begin{aligned}
\left|\int_{\quot} h_1(x) h_2(x \f)dx -
\sum_{\chi^2=1}  \chi(\cond) \int_{\quot}h_1(x) \chi(x) dx
\int_{\quot} h_2(x) \chi(x) dx \right| \\ \leq 2 q^{\alpha-1/2} ||h_1||_{L^2}
||h_2||_{L^2}\end{aligned} \end{equation}

Here $\chi$ ranges over all characters of $\adele_F^{\times}/F^{\times}$
such that $\chi^2=1$, and $\chi(x)$ denotes the function
$g \mapsto \chi(\det(g))$ on $\PGL_2(F) \backslash \PGL_2(\adele_F)$.
Indeed, to see (\ref{eq:hecke}), we note that the quantity inside the absolute value
on the left
hand side of (\ref{eq:hecke}) equals $\langle h_1 - \mathscr{P} h_1,  \f \cdot (h_2 - \mathscr{P} h_2 )\rangle_{L^2}$,
where $\mathscr{P}$ is as in Sec. \ref{subsec:proj} (in this case, the orthogonal projection onto the  locally constant functions).  But the $L^2$ orthogonal projection
$\mathrm{Id} - \mathscr{P}$ kills all one-dimensional $\PGL_2(F_{\cond})$ representations, from
which the result follows easily.

The functions $h_j(x) = f_j(xg) \overline{f_j(xg')} \,
(j=1,2)$ are $\PGL_2(\order_{F_{\cond}})$-invariant for $g,g' \in
\mathrm{supp}(\meas)$
since $\mathrm{supp}(\meas)$ commutes with $\PGL_2(F_{\cond})$
and $\cond \notin \Supp(f_1) \cup \Supp(f_2)$.
Moreover, $\|h_j\|_{L^2} \leq \|f_j\|_{L^4}^2$.
Apply (\ref{eq:hecke}) to these $h_j$, and substitute
in (\ref{eq:john1}). It results:
\begin{multline} \left| \lambda \cdot I(\varphi)\right|^2   \ll
     q^{\alpha-1/2}  \|\meas\|^2  \\ +
\sum_{\chi^2 = 1} \int_{(g,g')}
|\langle g^{-1} g' \cdot f_1, f_1 \otimes \chi \rangle
|\langle g^{-1}  g'\cdot f_2, f_2\otimes\chi \rangle d|\meas|(g) d|\meas|(g')
    \end{multline}
    where $|\meas|$ is the total variation measure associated to $\meas$,
  $\|\meas\| = |\meas|(\quot)$ is the total variation of $\meas$, 
  $f_i \otimes \chi$ is the function $x \mapsto f_i(x) \chi(\det(x))$, and
 brackets $\langle \cdot, \cdot \rangle$ denote inner product
in the Hilbert space $L^2(\quot)$; we will suppress the reference
to $L^2(\quot)$ here and in the rest of the argument.

Put $\meast = \check{|\meas|} \star |\meas|$. We may rewrite the previous result as:
\begin{multline} \label{eq:periodnew}
|\lambda|^2 |I(\varphi)|^2
\ll  q^{\alpha-1/2}\|\meas\|^2  \\ +
\left(
\int_{g}
\sum_{\chi^2 = 1}|\langle g \cdot f_1, f_1  \otimes
\chi\rangle| \cdot
|\langle g \cdot  f_2,f_2 \otimes \chi\rangle| d\meast(g)
\right )\end{multline}

We shall take $\meas$ in (\ref{eq:periodnew}) to be a linear combination of
Hecke operators. We follow the notations introduced
in Sec. \ref{subsec:Hecke}.
For $\mathfrak{n} \notin \Supp(\varphi)$, 
we denote by
$\lambda(\mathfrak{n})$ be the $\mathfrak{n}$th Hecke eigenvalue of
$\varphi$, i.e. $\varphi \star \mu_{\mathfrak{n}} = \lambda(\mathfrak{n}) \varphi$.
 With our normalizations, the Ramanujan
conjecture amounts to $|\lambda(\mathfrak{l})| \leq 2$ for $\mathfrak{l}$
prime.

Let $a_{\mathfrak{n}}$ be a sequence of complex numbers indexed
by integral ideals of $\order_F$.
Assume moreover that $a_{\mathfrak{n}} = 0$ whenever $\mathfrak{n}$
is divisible by any place in $\Supp(\varphi) \cup \Supp(f_1) \cup
\Supp(f_2)$.
  Let $\meas$ be the measure
on $\PGL_2(\adele_{F,f})$ defined by $\sum_{\n} a_\n \mu_{\n}$.
Then $\meas$ is symmetric under $g \mapsto g^{-1}$,
and $|\meas| = \sum_{\n} |a_{\n}| \mu_{\n}$. 
Moreover, $\varphi \star \meas = \lambda \varphi_3$,
where $\lambda = \sum_n a_\n \lambda(\n)$.
From the assumed bound on Ramanujan (see Sec. \ref{subsec:decay}, esp.
equation (\ref{eq:mcbound}))\footnote{
A small caution here is that the vectors $f_1 \otimes \chi$
and $f_2 \otimes \chi$ need not
be invariant by $\GL_2(\order_{F_{\cond}})$, if $\cond|\mathfrak{n}$, because $\chi$ may be
ramified. 
However, $\GL_2(\order_{F_{\cond}})$ always
fixes the line spanned by either of these vectors,
and the bound of (\ref{eq:mcbound}) depends only on the dimension
of the $\GL_2(\order_{F_{\cond}})$-span of the vectors in question.}, an elementary computation shows:

\begin{equation} \label{eq:useful}
\left|\int_{g \in \PGL_2(\adele_F)} \left|\langle g \cdot f_1,
f_1 \otimes \chi \rangle \langle g \cdot f_2, f_2  \otimes \chi\rangle
\right|
d \mu_{\n}(g) \right| \ll_{\epsilon}
\Norm(\n)^{2 \alpha - 1/2+\epsilon}
\end{equation}
Moreover, for fixed $g \in \Supp(\mu_{\n})$,
the inner product $\langle g f_1, f_1 \otimes \chi\rangle$
is nonvanishing only if $\chi$ is unramified at those places
at all places not in $\Supp(f_1)$ and not dividing $\n$. 
The number of such quadratic characters is $O_{\epsilon}(\Norm(\n)^{\epsilon} q^{\epsilon})$, where
the implicit constant (as always) is allowed to depend on the base field $F$. 
Thus:
\begin{equation} \label{eq:useful2}
\left|
\sum_{\chi^2 = 1}\int_{g} \left|\langle g \cdot f_1,
f_1 \otimes \chi \rangle \langle g \cdot f_2, f_2  \otimes \chi\rangle
\right|
d \mu_{\n}(g) \right| \ll_{\epsilon}q^{\epsilon} \Norm(\n)^{2 \alpha - 1/2+\epsilon}
\end{equation}

The total variation of $\meas$ may be computed:
\begin{equation} \label{eq:mutot}
\|\meas\| \ll_{\epsilon}
\sum_{\n} \Norm(\n)^{1/2 +\epsilon} |a_\n|
\end{equation}
Using (\ref{eq:conv}) to compute $\meas^{(2)}$, and
combining (\ref{eq:periodnew}), (\ref{eq:useful}) and
(\ref{eq:mutot}), we conclude:
\begin{equation} \label{eq:subconvex}\left|I(\varphi)
\right|
\ll_{\epsilon} q^{\epsilon}
\frac{
\left(\left(\sum_{\n} \Norm(\n)^{1/2 +\epsilon} |a_\n|\right)^2 q^{\alpha-1/2}
+
\sum_{\n,\m}
\sum_{\d | (\n,\m)} \left(\Norm\left(\frac{\n\m}{\d^2}\right)\right)^{2 \alpha -1/2}
%d^{-1} \lambda_1(\frac{\n\m}{\d^2}) \lambda_2(\frac{\n\m}{\d^2})
|a_\n| |a_\m |
\right)^{1/2}
}{
\left|\sum_{n}
a_\n \lambda(\n)\right|}.\end{equation}

The choice of $a_{\mathfrak{n}}$ follows an idea of Iwaniec; we slightly
modify the standard choice so that we do not need to appeal to Ramanujan on average. \footnote{The argument that follows was improved by a suggestion of P. Michel.}
Fix $K$ with $q^{1/1000} \leq K \leq q^{1000}$. 
 Let $S$ be the set of
prime ideals $\mathfrak{l}$ such that $\Norm(\mathfrak{l}) \in [K,2K]$
and $\mathfrak{l} \notin \Supp(f_1) \cup \Supp(f_2) \cup \Supp(\varphi)$.
In view of the assumptions, $|S| \gg_{\epsilon} K^{1-\epsilon}$. 

For $z \in \mathbb{C}$ we put $\sign(z) =z/|z|$ for $z \neq 0$ and 
$\sign(0) = 1$. Put
\begin{equation} \label{eq:achoice} a_{\mathfrak{n}}= \begin{cases}
\overline{\sign(\lambda(\n))}, \n \in S \\
\overline{\sign(\lambda(\n^2))},  \n = \mathfrak{l}^2, \mathfrak{l} \in S \\
0, \mbox{ else}.\end{cases} \end{equation}

Then $\left|\sum_{n}
a_\n \lambda(\n) \right| \gg_{\epsilon, F} K^{1-\epsilon}$,
$(\sum_{\n} \Norm(\n)^{1/2 +\epsilon} |a_\n|)  \ll_{\epsilon}
K^{2+\epsilon}$,
and
\begin{equation} \label{eq:pm}
\sum_{\n,\m}
\sum_{\d | (\n,\m)} \left(\Norm\left(\frac{\n\m}{\d^2}\right)\right)^{2
\alpha -1/2} |a_\n| |a_\m |
%\sum_{(\n, \m ) = 1} + \sum_{(\n,\m) \neq
%1, \n \neq \m} + \sum_{\n = \m} \\
%\ll K^2 \cdot \left(K^{2 \alpha} K^{4 \alpha - 1} +
%K^{\alpha} K^{6 \alpha - 3/2} + K^{8 \alpha -2} \right)
%+ K \cdot K^{\alpha + \epsilon} +  K(K^{2 \alpha +  \epsilon} + K^{\epsilon})
%\\
\ll K^{4 \alpha+1}.\end{equation}

We deduce from (\ref{eq:subconvex}) that
\begin{equation}|I(\varphi)|\ll_{\epsilon, F} (qK)^{\epsilon}
\frac{\left(K^{4}  q^{\alpha-1/2}   + K^{1+ 4
\alpha}\right)^{1/2}}{K} \end{equation}
Taking $K = q^{\frac{1/2-\alpha}{3 - 4 \alpha}}$, we obtain
$|I(\varphi)| \ll_{\epsilon,F}  q^{\epsilon + \frac{ (2 \alpha - 1/2) (1/2-\alpha)}{3 - 4 \alpha}}$.  \qed

\subsection{A technical generalization.} \label{subsec:technical}
For certain applications,
we shall require a slight generalization of Prop. \ref{prop:main}
in which the role of  $g \mapsto f_1(g) f_2(g \f)$
is replaced by $g \mapsto F(g,g\f)$, where $F$
is a function on $\quot \times \quot$ that is not necessarily of product
type. Although the method of proof is identical to Prop. \ref{prop:main}
the details are slightly more technical; in particular, to
state the result we will have need of the adelic Sobolev norms
discussed in Sec. \ref{subsec:adelicsobolev}.
We shall also use the notion of totally nondegenerate for functions
on $\quot \times \quot$: see Sec. \ref{sec:adelicfunctionspace}.

\begin{prop} \label{prop:mainV}
Suppose $F \in C^{\infty}(\quot \times \quot)$
is totally nondegenerate. Suppose moreover that there is $b \in \mathbb{R}$ with 
\begin{equation} \label{eq:cond2}
\prod_{\mathfrak{q} \in \Supp(F) \cup \Supp(\varphi)} \Norm(\q)
\leq \Norm(\cond)^b\end{equation}

%\begin{equation}\label{eq:sideortho} \int_{\quot} F(x_0, x) dx =
%\int_{\quot} F(x,x_0) dx = 0.\end{equation}
Let $\pi$ be a cuspidal representation of $\PGL_2$
over $F$, with conductor $\p$, and put
$I(\varphi) = \int_{\quot} F(x, x\f) \varphi(x) dx$,
for $\varphi \in \pi$.
Then, for any $p>4, d \gg 1$,
$$|I(\varphi)| \ll_{b,\epsilon}  \Norm(\cond)^{-\beta+\epsilon}
\|\varphi\|_{L^2} S_{p, d,2/p}(F),$$
where $\beta = \frac{(1-2\alpha) (1-4\alpha)}{p(7-4
\alpha)}$.
\end{prop}

With Kim's bound $\alpha = 3/26$, we obtain $\frac{(1-2 \alpha) (1-4\alpha)}{7-4\alpha} > 1/17$. 
\proof
The proof follows closely the proof of Prop. \ref{prop:main};
the only difference is that we apply \eqref{eq:heckesoboleveq} (proved in
Lem. \ref{lem:heckesobolev}) in place of
(\ref{eq:hecke}).
Again, we may freely assume that $\|\varphi\|_{L^2} = 1$; again
we put $q = \Norm(\cond)$. 

Let notations be as in the proof of
Prop. \ref{prop:main};
in particular, $\meas$ is a signed real measure on
$\PGL_2(\adele_{F,f})$ whose
support commutes with $\PGL_2(F_{\cond})$, and $\lambda \in \C$
satisfies $\varphi \star \check{\meas} = \lambda \varphi$.
Proceeding  as in that proof, and
in particular as in (\ref{eq:john1}), we obtain:
\begin{multline} |\lambda I(\varphi)|^2
\leq \int_{\quot}
\int_{(g,g') \in \PGL_2(\adele_{F,f})^2}  F((x,x) (g,g) (1,\f)) \cdot \\
\overline{F((x,x) (g',g') (1,\f))} d\meas(g)
d\meas(g')\end{multline}

For any $g,g' \in \PGL_2(\adele_{F})$, set $F_{g,g'}(x_1, x_2) = F((x_1, x_2) (g,g))
\overline{F((x_1,x_2) (g',g'))}$.
Then $F_{g,g'}$ is invariant by $\PGL_2(\order_{F_{\cond}}) \times
\PGL_2(\order_{F_{\cond}})$ for $g,g' \in \mathrm{supp}(\meas)$.
By Hecke equidistribution in the form of
\eqref{eq:heckesoboleveq} (proved in 
Lem. \ref{lem:heckesobolev}) we see that for $p>2, d \gg 1$: 
\begin{multline}\left|\int_{\quot} F_{g,g'}((x,x)(1,\f)) dx -
\sum_{\chi^2=1} \chi(\cond) \int_{\quot \times
\quot }
F_{g,g'}(x_1,x_2 )\chi(x_1) \chi(x_2) dx_1 dx_2\right| \\ \ll_{\epsilon} q^{(2 \alpha -1)/p+\epsilon}
S_{p,d}(F_{g,g'}).\end{multline}

Here, as before, $\chi(x)$ denotes the function on $\quot$ defined by
$g \mapsto \chi(\det(g))$.
By definition $|\int_{\quot \times \quot} F_{g,g'}((x_1, x_2))   \chi(x_1) \chi(x_2)|
 = |\langle ( g'^{-1} g, g'^{-1} g) F, F \otimes (\chi,\chi)\rangle_{L^2(\quot \times \quot)}|$.

By the basic properties of adelic Sobolev norms
(\eqref{eq:sobproducteq} and
\eqref{eq:sobgroupeq}, proofs in Lem. \ref{lem:sobproduct} and Lem. \ref{lem:sobgroup})
\begin{equation} \begin{aligned}S_{p,d} (F_{g,g'} ) :=
S_{p,d,1/p} (F_{g,g'})
\ll S_{2p,d,1/p}((g,g) \cdot F)
S_{2p,d,1/p}((g',g')\cdot F) \\ \leq \|g\|^{2/p}
\|g'\|^{2/p} S_{2p,d,1/p}(F)^2.
\end{aligned}\end{equation}

Let us remark that the factors $\|g\|^{2/p}$ and $\|g'\|^{2/p}$ arises in the following way:
Lem. \ref{lem:sobgroup} actually gives a factor $\|(g,g)\|^{1/p}$, where the norm $\|\cdot\|$ (as in
Sec. \ref{sec:adelicfunctionspace})
is computed in $\PGL_2(\adele_{F,f})^2$; this equals $\|g\|^{1/p}$ where the norm
is computed in $\PGL_2(\adele_{F,f})$. 

Choose $\meas$ as in the proof of Prop. \ref{prop:main} (see esp.
paragraph before  (\ref{eq:useful}))
and choose the coefficients $a_{\n}$ as in that proof (see (\ref{eq:achoice})). 
In particular, $\|\meas\| \ll K^{2+\epsilon}$.
Since $F$ is totally nondegenerate, the
matrix coefficients $\langle (g'^{-1} g, g'^{-1} g) F, F \rangle$
satisfy bounds that are of the same quality
as in the proof of Prop. \ref{prop:main}; in particular, as in
(\ref{eq:useful2}):
$$\sum_{\chi^2 = 1} \int_{g} \left|\langle (g,g) F, F \otimes (\chi,\chi)\rangle
\right| d\mu_{\n}(g) \ll   q^{\epsilon}\Norm(\n)^{2 \alpha - 1/2+\epsilon} \|F\|^2$$

Finally $\|g\| \ll_{\epsilon} K^{2+\epsilon}$ for all $g \in \mathrm{supp}(\meas)$.
Proceeding just as in the previous proof,
$$|I(\varphi)| \ll_{\epsilon} (qK)^{\epsilon} \frac{\left(K^4 q^{(2\alpha -1)/p}
K^{8/p} S_{2p,d,1/p}(F)^2
+ K^{1+4\alpha} \|F\|_{L^2(\quot \times\quot)}^2\right)^{1/2}
}{K}.$$

Consequently, for any $p >2$,
$$|I(\varphi)| \ll_{\epsilon} (qK)^{\epsilon}\frac{
\left(K^8 q^{(2\alpha-1)/p} S_{2p,d,1/p}(F)^2+ K^{1+4 \alpha}
\|F\|^2_{L^2}\right)^{1/2}}{K} $$

To conclude, choose $K = q^{\frac{1-2\alpha}{p(7-4\alpha)}}$
and replace $p$ by $p/2$ (thus, e.g., $p >2$ becomes $p>4$). \qed
\section{Application to $L$-functions.} \label{sec:app}

We now present the first applications to subconvexity. The rough idea is simply
that certain $L$-functions are expressed as period integrals
of the type that are bounded by Prop. \ref{prop:main} and Prop. \ref{prop:mainV}.
There is one significant issue in implementing this (rather evident) idea:
namely, the integral representation that we use for Rankin-Selberg and the standard
$L$-functions involve Eisenstein series, which are not in $L^2$; this causes
problems in applying Prop. \ref{prop:main}!

Thus we need to {\em regularize}. Two natural ways of doing this are to replace
an Eisenstein series by a ``wave-packet''; or to use a suitable form of truncation in the defining integrals. 
In the present paper we will use the wave-packet technique; in the paper \cite{MV} we shall
also use truncation. 

 Let us briefly describe the wave packet technique
in a classical language. Roughly speaking we can express
the Rankin-Selberg $L$-function of two classical forms $f,g$ 
via an integral of the form $L(s, f \times g) = \int_{z} f(z) g(z) E(s,z) $,
for some Eisenstein series $E(s)$. We now regularize, replacing $E(s,z)$ by a wave packet. 
Let $h(s)$ be any holomorphic function: then 
\begin{equation} \label{eq:regularize} \int_{Re(s) = 1/2} h(s) L(s, f \times g) ds = \int_{z} f(z) g(z) \int_{\Re(s) = 1/2} h(s) E(s,z).\end{equation}
We wish to eventually recover an upper bound for $L(1/2, f \times g)$ (say) 
from the left-hand side, so we take $h(s) = \overline{ L(1-\overline{s}, f\times g)}$.
Then $h(s) L(s, f \times g)$ is positive along $\Re(s) = 1/2$.  
To apply Prop. \ref{prop:main} to the right-hand side of (\ref{eq:regularize}), we shall moreover need to control the behavior of the regularized Eisenstein series
$E_h = \int_{\Re(s) = 1/2} h(s) E(s,z)$; this type of analysis is carried out in 
Sec. \ref{sec:eisreg2} and Sec. \ref{theartoffugue}, the main point being that
the divergence of the Eisenstein series comes entirely from the constant term.

 It is worth remarking that
Iwaniec's bounds for the $L$-function near $1$ enter rather crucially in this analysis: in effect,
we bound $E_h$ by an easy argument involving shift of contours; this 
necessitates that $h$ be estimated on a line $\Re(s) = -\varepsilon$, which
amounts to estimating $L(s, f \times g)$ for $\Re(s) = 1+\varepsilon$. 

 In what follows
we have not attempted to obtain polynomial dependence in all parameters.
This is not hard to do --- and,
at its essence, a statement that one can find {\em analytically} suitable
test vectors in a Rankin-Selberg integral; but we have not
done so here. On the other hand, we give full details of this procedure
in the proof of Thm. \ref{thm:subconvexct} (in which the polynomial dependence
is particularly useful for applications).

\begin{thm} \label{thm:subconvextp}
Let $\pi_1, \pi_2$ be fixed automorphic cuspidal representations
of $\PGL_2$ over $F$; fix $t \in \R$. Let $\pi$
be an automorphic cuspidal representation
with conductor $\cond$, a prime ideal that is prime to the conductors
of $\pi_1$ and $\pi_2$.  

Then, assuming Hypothesis \ref{prop:intrep}:

\begin{equation} \label{eq:triple}
L(\frac{1}{2}, \pi_1 \otimes \pi_2 \otimes \pi) \ll_{\pi_{\infty}}
\Norm(\cond)^{1 - \frac{1}{13}} \end{equation}

and, unconditionally:

\begin{equation} \label{eq:double}
|L(\frac{1}{2}+it, \pi_1 \otimes \pi)|^2 \ll_{ \pi_{\infty}}
\Norm(\cond)^{1-\frac{1}{100}} \end{equation}

\begin{equation} \label{eq:single}
|L(\frac{1}{2} + it, \pi)|^4 \ll_{\pi_{\infty}} \Norm(\cond)^{1-\frac{1}{600}}
\end{equation}

In these statements, the notation $\ll_{\pi_{\infty}}$
indicates an implicit constant that depends continuously on the local archimedean representation $\pi_{\infty}$ of $\GL_2(F_{\infty})$ underlying $\pi$. 

%\begin{equation} \label{eq:symsquare}
%L(\frac{1}{2} + it, \mathrm{Sym}^2 \pi_3)^2 \ll \Norm(\cond)^{1-?}
%\end{equation}
\end{thm}

Note we make no claim about the dependency of the implicit constant on $t, \pi_1, \pi_2$;
as remarked above, this dependence could be made polynomial in the conductors, but this
would require more careful analysis of the archimedean integrals. \footnote{ It is important to note, however, that this is 
an entirely local problem; it is intended that this will be carried out in a more general context in
\cite{MV}.  Both for applications and to illustrate procedure, we have carried out this type of analysis for the results
on subconvexity of character twists in Section \ref{sec:torus1}. Those results are proved with polynomial dependence on all parameters.}

We remark that we have used H.Kim's bound $\alpha = 3/26$; any
value of $\alpha$ less than $1/4$ would give subconvexity and
under Ramanujan one obtains for (\ref{eq:triple}) the exponent
$5/6$. The exponents for (\ref{eq:double}) and
(\ref{eq:single}) can be improved, e.g. the present proof does not
take into account the fact that unitary Eisenstein series satisfy
Ramanujan! 

\subsection{Results relating periods and integral representations.}
For the convenience of the reader, we summarize here the results
that relate periods and integral representations (proved in later sections).
Roughly speaking, any integral representation for an $L$-function
expresses it as a period integral
with certain test vectors belonging to the space of an automorphic cuspidal representation.

A delicate point, which is quite relevant to issues of polynomial dependence
in auxiliary parameters, is precisely {\em which} test vectors. In principle,
the proofs of results about integral representation give explicit test vectors.
In practice, it is tedious to extract these explicit test vectors. Our policy
throughout this paper is the lazy one: to deduce results, as far as possible, by formal arguments and
without choosing explicit test vectors. The price of this
is that we will not obtain not quite the $L$-function, but rather
a holomorphic function that differs from the $L$-function by some
harmless factors. 

More precisely, the content of the Proposition (Prop. \ref{prop:rs}) that follows is that one can write down an integral representation $I(s)$
for the $L$-functions of interest, so that:
\begin{enumerate} \item $I(1/2)$ is not too much smaller than $L(1/2)$ -- or with $1/2$ replaced
by the point of interest -- so that a bound for $I(1/2)$ gives a bound for $L(1/2)$. 
\item  $I(s)$ is not too much bigger than $L(s)$ for {\em any} $s$. This type of  control will be useful in shifting contours. 
\end{enumerate}
One might prefer to get $I(s) = L(s)$ but we don't need this stronger statement.

As is discussed at length in Sec. \ref{sec:eisenstein},
to a Schwarz function $\Psi$ on $\adele_F^2$ is associated
a family of Eisenstein series $E_{\Psi}(s, g)$ on $\quot$, 
which varies meromorphically in the parameter $s \in \C$.

\begin{prop}
\label{prop:rs}
Let $s_0,t_0, t_0' \in \C$. 
Let $\pi_1$ be a fixed automorphic cuspidal representation 
of $\PGL_2(\adele_F)$ and $\pi$ an automorphic cuspidal
representation of prime conductor $\cond$; assume that the finite conductors of $\pi, \pi_1$ are coprime. 

Denoting by $\pi_{\infty}$ the representation of $\PGL_2(F_{\infty})$ corresponding to $\pi$, 
suppose that $\Cond(\pi_{\infty})$ is bounded above; equivalently,
 $\pi_{\infty}$ belongs to a bounded subset\footnote{See Section \ref{subsec:cuspformsbdd} for definition} of the dual $\widehat{\PGL_2}(F_{\infty})$
 (in what follows the implicit
constants may depend on these bounds).

There exists a fixed finite set $\mathcal{F}$
of Schwarz Bruhat functions on $\adele_F^2$
and a real number $C > 0$ so that:

There exist vectors $\varphi_1 \in \pi_1, \varphi \in \pi$ and $\Psi \in \mathcal{F}$ so that
$$\Phi(s) := \Norm(\cond)^{1-s}
\frac{\int_{\quot}\varphi(g) \varphi_1(g \f ) E_{\Psi}(s, g) dg}{ \Lambda(s, \pi_1 \otimes \pi)}$$
is holomorphic and satisfies:
\begin{enumerate}
\item  $|\Phi(s_0)| \gg 1$ and $|\Phi(s)| \ll C^{|\Re(s)|} (1+|s|)^C$;
\item At any nonarchimedean place $v$ such that $\pi_1$ is unramified, $\Psi_v$
is invariant by $\PGL_2(\order_{F_v})$; 
at any nonarchimedean place $v$ such that
$\pi_{1}$ and $\pi$ are both unramified,
 $\varphi$,  $\varphi_1$ are both invariant by $\PGL_2(\order_{F_v})$. 
 \item $\|\varphi_1\|_{L^{\infty}} \ll 1$ and  $ \|\varphi\|_{L^2(\quot)} \ll_{\epsilon}
\Norm(\cond)^{\epsilon}$.
\end{enumerate}

Moreover, there exist vectors $\varphi \in \pi, \Psi_1, \Psi_2 \in \mathcal{F}$ so that:
\begin{equation}\label{eq:losing1} \Phi(t,t') = \Norm(\cond)^{1/2-t} 
\frac{
\int_{\quot} \varphi(g) E_{\Psi_1}(g, \frac{1}{2}+t) E_{\Psi_2}(g \f, \frac{1}{2}+t') dg}
{ \Lambda(\frac{1}{2} + t + t', \pi)
\Lambda(\frac{1}{2} + t -t', \pi)}
\end{equation}
is holomorphic and satisfies:
\begin{enumerate}
\item $|\Phi(t_0, t_0')| \gg 1$ and $|\Phi(t, t')| \ll   C^{ |\Re(t)| + C|\Re(t')|}(1+|t|+|t'|)^C $. 
\item For any nonarchimedean place $v$, each $\Psi_1$ and $\Psi_2$
is invariant by $\PGL_2(\order_{F_v})$; for each place at which $\pi$ is unramified,
the same is true of $\varphi$. 
\item$\|\varphi\|_{L^2(\quot)} \ll_{ \epsilon}
\Norm(\cond)^{\epsilon}$.
\end{enumerate}
\end{prop}
\proof
Lem. \ref{lem:rsone} and Lem. \ref{lem:rsoneeis}. 
%\footnote{It should be noted that in Section \ref{sec:rs}
%we use the measure on $\quot$ which is the quotient measure from $\PGL_2(\adele_F)$, rather
%than the probability measure on $\quot$. The probability measure is used by default everywhere else in this paper, including the Proposition, but because the two measures differ by a constant depending only on $F$, the switch of measure does not affect the truth of any of the statements.}
\qed 

In effect, we could achieve ``$\Phi = 1$'' in Prop. \ref{prop:rs} by a more careful choice of local data;
but this is irrelevant for the purpose of global estimation. 

\subsection{Proof of Thm. \ref{thm:subconvextp}.}

\proof  (of (\ref{eq:triple})). It follows from Hypothesis
\ref{prop:intrep} and Proposition \ref{prop:main}. 

\proof (of (\ref{eq:double})).
The basic idea is that the Rankin-Selberg convolution
is a triple product, with one factor being an Eisenstein series.
However, one cannot naively apply Prop. \ref{prop:main}
since Eisenstein series do not belong to $L^4(\quot)$.
To avoid this, we will use a wave-packet of Eisenstein series.

First, we can assume from the start that $\pi_{\infty}$
belongs to a bounded subset of the dual $\widehat{\PGL_2(F_{\infty})}$.  
The implicit constants in the proof that follow depend on this subset. 
We denote by $\Lambda$ the completed $L$-function.
We begin by remarking that since $\Norm(\cond) \rightarrow \infty$
we may assume that $\pi_1$ is not isomorphic to $\pi$, or to
any quadratic twist of $\pi$. In particular, 
we are free to assume that $\Lambda(s, \pi_1 \otimes \pi)$
has no poles. Moreover, the finite conductor
of
$\Lambda(s,\pi_1 \otimes \pi)$ differs from $\Norm(\cond)^2$
by an absolutely bounded constant.

Fixing $t_0 \in \mathbb{R}$, 
let $\Psi, \varphi, \varphi_1, E_{\Psi}(g,s),\Phi$ be as in Prop. \ref{prop:rs} with $s_0 = 1/2+it_0$, so that
$|\Phi(1/2+i t_0)| \gg 1$. 
For simplicity we write simply $E(g,s)$ for $E_{\Psi}(g,s)$. 
Fix $\kappa > 0$. In the rest of the proof we omit the subscript
$\kappa, \epsilon$ from $\ll$, with the understanding
that all implicit constants depend on $\kappa$ and $\epsilon$.
Put $$I(s) = \Norm(\cond)^{s-1} \Lambda(s,\pi_1 \otimes \pi) \Phi(s) =
\int_{\quot} \varphi_1(g \f)  E(s,g)
\varphi(g)dg$$
From Iwaniec's upper bounds for $L$-functions \cite[Chapter 8]{iwaniec}, the functional equation for
$\Lambda$, and the bounds on $\Phi$ furnished by Prop. \ref{prop:rs}, 
%A result of Molteni implies
%$|\Lambda(1 + \kappa + it)|  \ll_{\epsilon} (1+|t|)^{-4}
%\Norm(\cond)^{\epsilon}.$
%From this and the functional equation, we deduce
%$|\Lambda(-\kappa + it) | \ll
%\Norm(\cond)^{1 + 2 \kappa+ \epsilon}(1+|t|)^{-4}$.
\begin{equation} \label{eq:ibound} |I(1+\kappa + it)| \ll
(1+|t|)^{-4}
\Norm(\cond)^{\kappa + \epsilon},
|I(-\kappa+it)| \ll (1+|t|)^{-4}
\Norm(\cond)^{ \kappa + \epsilon}.\end{equation}

Put $h(s) = s(1-s)(s-\frac{1}{2})^2 \overline{I(1-\overline{s})}$. Then $h(s)$ is holomorphic
in $-\kappa \leq \Re(s)  \leq 1+\kappa$ and $h(\frac{1}{2}) =
0$.\footnote{\label{extrasing} The fact that we impose $h(1/2) = 0$ has a very concrete
meaning in classical terms. Fix, for example, a form $f$
and $t \in \mathbb{R}$.
Consider $\sum_{g} |L(\frac{1}{2} + it, f \otimes g)|^2$
where the sum is taken over a basis of holomorphic Hecke eigenforms
of level $N$ and trivial Nebentypus. If $t =0$,
this has the asymptotic behaviour $N \log(N)^3$.
On the other hand, if $t \neq 0$, it behaves
like $N \log(N)$. Forcing $h(1/2) = 0$ ``counteracts'' this extra
singularity.}
Moreover $h(s)$ has rapid decay as $\Im(s) \rightarrow \infty$,
in view of the $\Gamma$-factors of the completed $L$-function.
Put $E_h(g) = \int_{\Re(s) = 1+\kappa} h(s) E(s,g)$.
It is proved in Lem. \ref{eislinfty} that, for such $h$,
$\|E_h(g )\|_{L^{\infty}}
\ll \int_{-\infty}^{\infty} \left(
|h(-\kappa + it)| + |h(1+\kappa+it)|\right) dt.$
Applying (\ref{eq:ibound}), we conclude that
$\|E_{h}(g) \|_{L^{\infty}} \ll \Norm(\cond)^{\kappa + \epsilon}$.

On the other hand,  we see from the definition of $I(s)$ that 
\begin{equation} \label{eq:m}\int_{\Re(s) = 1+\kappa} h(s) I(s) =
\int_{\Re(s) = 1+\kappa} h(s) ds \int_{\quot} \varphi_1(g \f)  E(s,g)
\varphi(g)dg.\end{equation}
The double integral on the right hand side of (\ref{eq:m})
is absolutely convergent and orders may be switched;
thus $$\int_{\Re(s) = 1+\kappa} h(s) I(s) =
\int_{\quot} \varphi_1(g \f) E_h(g) \varphi(g) dg =
\int_{\quot} \varphi_1(g \f) E_h^{0}(g) \varphi(g)dg, $$
where $E_h^{0} := \mathscr{P}(E_h)$ is 
totally nondegenerate (see Section \ref{subsec:proj}) and satisfies
$\|E_h^{0}\|_{L^{\infty}} \ll_{\epsilon} \Norm(\cond)^{\kappa+\epsilon}$.  

We then deduce from Prop. \ref{prop:main} that:
\begin{equation} \label{eq:bound}\left|\int_{\Re(s) = 1+\kappa} h(s) I(s)\right|
\ll
\|\varphi_1\|_{L^4(\quot)}
\Norm(\cond)^{-\frac{1}{26} + \kappa+ \epsilon}\end{equation}

$I(s)$ and $h(s)$ both decay exponentially rapidly as $\Im(s) \rightarrow \infty$.
It is therefore simple
to justify shifting the line of integration in (\ref{eq:bound})
to $\Re(s) = 1/2$.
We deduce thereby that
\begin{equation}\label{eqn:one}\left|\int_{\Re(s) = \frac{1}{2}}t^2 |I(\frac{1}{2} + it)|^2 dt \right|
\ll \|\varphi_1\|_{L^4(\quot)} \Norm(\cond)^{
\kappa   - \frac{1}{26} +
\epsilon}.\end{equation}

From (\ref{eq:ibound})
we deduce bounds on $I$ and $I'$
inside the strip $0 \leq \Re(s) \leq 1$
 by the maximal modulus principle. In particular,
\begin{equation} \label{eqn:two}
|I'(\frac{1}{2} + it)|
\ll_{t} \Norm(\cond)^{\kappa+\epsilon}.\end{equation}

Combining (\ref{eqn:one}) and (\ref{eqn:two}),
and recalling that $\kappa$ is arbitrary,
we obtain $|I(1/2+it)| \ll_t \Norm(\cond)^{\frac{1}{5 \cdot 26}+\epsilon}$.
%On the other hand,
%for any real-valued function $f(x)$ on $\mathbb{R}$,
%an easy computation with the mean-value theorem shows $$|f(0)|^4
%\leq 12 \max_{x \in [-1,+1]} |f'(x)|^3
%\int_{-1}^{1} x^2 |f(x)| dx.$$
%Applying this to $f(x)=I(\frac{1}{2} + ix) I(\frac{1}{2} - ix)$,
%and using (\ref{eqn:one}) and (\ref{eqn:two}),
%we conclude
%$$I(\frac{1}{2})^8  \ll
% \Norm(\cond)^{7 \kappa - \frac{1}{12} +  \alpha/2+  \epsilon}.$$
Thus $|\Lambda(\frac{1}{2}+it_0, \pi_1 \otimes \pi)|
\ll_{\epsilon,t} \Norm(\cond)^{\frac{1}{2} -\frac{1}{130} +
\epsilon}$.
\qed

\proof (of (\ref{eq:single}).) The proof is similar
to that of (\ref{eq:double}), but a slightly
more elaborate regularization is required, since
we shall proceed from the expression (\ref{eq:losing1}) of $L(s,\pi)$
as a triple product against {\em two} Eisenstein series. Again
we may assume from the start that $\pi_{\infty}$ is confined
to a bounded subset of $\widehat{\GL_2(F_{\infty})}$; the implicit constants will, again,
depend on this subset. 

 Let $\Lambda(s,\pi)$ be the completed $L$-function attached
to $\pi$. Fixing $t_0, t_0' \in  i \mathbb{R}$,  Prop. \ref{prop:rs} gives the existence
of Eisenstein series $E_{\Psi_1}(g,s) = E_1(g,s), E_{\Psi_2}(g,s) = E_2(g,s)$ on $\quot$,
and $\varphi \in \pi$ so that 
$$\Phi(t, t') : = \Norm(\cond)^{1/2-t}
\frac{\int_{\quot} \varphi(g) E_1(g, \frac{1}{2}+t) E_2(g \f, \frac{1}{2}+t') dg}{ \Lambda(\frac{1}{2} + t + t', \pi)
\Lambda(\frac{1}{2} + t -t', \pi)}$$
satisfies $|\Phi(t_0, t_0')| \gg 1$ and $\Phi(t,t') \ll C^{|\Re(t)| + |\Re(t')|} (1+|t|+|t|'|)^C$. 

We put
\begin{equation}\label{eq:idef}\begin{aligned}I(z_1,z_2) = \Phi(z_1, z_2) \Norm(\cond)^{z_1-1/2}
\Lambda(\frac{1}{2}+z_1+z_2,\pi) \Lambda(\frac{1}{2}+z_1-z_2) \\ = \Phi(z_1, z_2)
\Norm(\cond)^{\frac{z_1+z_2}{2}-1/4} \Lambda(\frac{1}{2} + z_1 +
z_2, \pi) \Norm(\cond)^{\frac{z_1-z_2}{2}-1/4} \Lambda(\frac{1}{2}
+ z_1 - z_2,\pi)\end{aligned}\end{equation} Then $I(z_1,z_2)$ is a
holomorphic function of $(z_1,z_2) \in \C^2$. $I(z_1,z_2)$ has
rapid decay along ``vertical lines'', that is, for $\sigma,
\sigma'$ in a fixed compact set and $(t,t') \in \mathbb{R}$ we
have $I(\sigma + it, \sigma' + it') \ll_N (1+|t|+|t'|)^{-N}$.

%\begin{equation}\label{eq:iwan}\left|I(t,t')\right|
%\ll \Norm(\cond)^{\kappa}  \mbox{  when }(|\Re(t)|,|\Re(t')|) \in \{
%(1/2+\kappa,0), (0,1/2+\kappa)\}
%\end{equation}
Let $\kappa > 0$ be fixed. From (\ref{eq:idef}), Iwaniec's bounds
for $L$-functions near $1$, and the rapid decay of $I$ along
``vertical lines,'' we obtain by the maximal modulus principle:
$$(1+|z_1|+|z_2|)^N \max( |I(z_1,z_2)|, |\partial_1 I(z_1,z_2)|,
|\partial_{2} I(z_1,z_2)|) \ll_N \Norm(\cond)^{\kappa}, \ \
|\Re(z_1)| + |\Re(z_2)|  \leq 1/2+\kappa,$$ where $\partial_1$
(resp. $\partial_2$) is the operator of differentiation w.r.t.
$z_1$ (resp. $z_2$). Put \begin{equation} \label{eq:pauline} h(z_1,z_2) = z_1^2 z_2^2 (1/4-z_1^2)
(1/4-z_2^2)\overline{I(-\overline{z_1},-\overline{z_2})}. \end{equation} Then, in the notation of Sec.
\ref{theartoffugue} (esp. Def. \ref{def:hdef}) , $h$ belongs to the space
$\mathcal{H}^{(2)}(\kappa)$ and satisfies $\|h\|_N \ll_N
\Norm(\cond)^{\kappa}$.

Put $I = \int_{\Re(z_1) = \Re(z_2) = 0} h(z_1, z_2) I(z_1, z_2)
dz_1 dz_2$. Then: \begin{equation} \begin{aligned}
I
= \int_{\Re(z_1) = 0, \Re(z_2)=0}
h(z_1, z_2) dz_1 dz_2 \int_{\quot}
\varphi(g) E_1(g, 1/2+z_1) E_2(g \f,1/2+z_2) dg
\\  = \int_{\quot} \varphi(x) E_{h}((x,x) (1,\f)) dx
\end{aligned} \end{equation}
where the function $E_h$ on $\quot \times \quot$
is defined by $$E_h(g_1, g_2)  = \int_{\Re(z_1) = 0,
\Re(z_2) = 0} h(z_1, z_2) E_1(g_1,1/2+z_1) E_2(g_2,1/2+z_2),$$
and the interchange of orders is justified by the
(easily verified) absolute convergence of the double integral defining $I$.
Note that $E_h(g_1, g_2)$ is totally nondegenerate (see Section \ref{subsec:proj} for definition). 
We now apply Prop. \ref{prop:mainV} to conclude that
$|I| \ll_{p,d} S_{p,d,2/p}(E_h) \|\varphi\|_{L^2} \Norm(\cond)^{-\frac{1}{17 p}}$ for any $p>4, d
\gg 1$.
We note at this point that the requirement $p>4$ makes
it critical that the regularized Eisenstein series $E_h$ belong to $L^4$;
the trivial fact that Eisenstein series belong to $L^{2-\epsilon}$ is far
from sufficient.

On the other hand, by
Lem. \ref{eislinfty3}, $S_{p,d,2/p} (E_h) \ll \|h\|_N$ for
some $N$, possibly depending on $p,d$. By Prop. \ref{prop:rs}, 
$\|\varphi\|_{L^2} \ll_{\epsilon} \Norm(\cond)^{\epsilon}$. 

 Thus
$|I| \ll_{\epsilon} \Norm(\cond)^{\epsilon-1/68}$.
Now, by the definition of $h$ (\ref{eq:pauline}) we have
$ I = 
\Norm(\cond)^{-1} \int_{(t,t') \in \mathbb{R}^2}
(1/4+t_1^2) (1/4+t_2^2) t_1^2 t_2^2 | I(it_1, i t_2)|^2 dt_1 dt_2$. 
%\prod_{\epsilon,\epsilon' \in \{\pm 1\}^2}
%\Lambda(1/2+i \epsilon t + i \epsilon' t') \\ =
%\\ = \Norm(\cond)^{-1} \int_{(t,t') \in \R^2}(t_1 t_2)^2 (1/4+t_1^2) (1/4+t_2^2)
%\left|\Lambda(1/2 + i (t+t')) \Lambda(1/2+i(t-t'))\right|^2 dt dt'
Thus we obtain:
\begin{equation} \label{eq:whatisI} \int_{(t,t') \in \mathbb{R}^2} |I(it_1, i t_2)|^2 t_1^2 t_2^2 dt_1 dt_2 \ll_{\epsilon} \Norm(\cond)^{1-1/68+\epsilon}.\end{equation}

Using (\ref{eq:whatisI}), and the given properties of $\Phi$,  we deduce that
$|\Lambda(\frac{1}{2}+t_0+t_0') \Lambda(\frac{1}{2} +  t_0 - t_0')|^2 \ll_{t_0, t_0'} \Norm(\cond)^{1-1/600}$
in a similar fashion to the conclusion of the proof
of (\ref{eq:double}). We take $t_0' = 0$ to conclude. 
\qed

\section{Torus periods (I): subconvex bounds
for character twists over a number field.} \label{sec:torus1}

In this section we shall work in considerable generality; we shall
derive subconvex bounds without any assumptions of prime or squarefree conductor,
and obtaining polynomial dependence in all auxiliary parameters. This
is useful for applications, but will involve some notational overhead. As a result,
we have sacrificed good exponents for simplicity at many steps.

\begin{thm}  \label{thm:subconvexct}
Let $\pi$ be a (unitary) cuspidal representation of $\GL_2(\adele_F)$, and
$\chi$ a unitary character of $\adele_F^{\times}/F^{\times}$, with finite conductor $\fcond$. Then there is $N
> 0$ such that
\begin{eqnarray}
\label{eq:cusp}L(\frac{1}{2}, \pi \times \chi) \ll
\mathrm{Cond}(\pi)^{N} \Cond_{\infty}(\chi)^N \Norm(\fcond)^{1/2-\frac{1}{24}}, \\
\label{eq:burgess} L(\frac{1}{2}, \chi) \ll \Cond_{\infty}(\chi)^N
\Norm(\fcond )^{1/4 -\frac{1}{200}}\end{eqnarray}
\end{thm}

Note that the result also implies a corresponding statement for the $L$-functions evaluated 
at $\frac{1}{2} + it$, since one may replace $\chi$ by $\chi |\cdot|^{it}$. 

Since it is perhaps hidden in the proof where the polynomial
dependence on conductor arises, we would like to explicate it now.
If $\pi$ is an automorphic cuspidal representation
with analytic conductor $\Cond(\pi)$, there exists
a vector $\psi \in \pi$ with Sobolev norms
$S_{2,d,\beta}(\psi) \ll \Cond(\pi)^{\mathrm{const}\max(\beta,d)}$.
Moreover, one can choose such a $\psi$ to be a
``good'' test vector w.r.t. certain toral periods.
Thus the analytic conductor enters precisely through the minimal
Sobolev norm of a suitable vector belonging to the space of $\pi$.
We note that the test vectors we choose are smooth but not $K$-finite at infinite places; this
idea has been heavily exploited in the previous work of Bernstein and Reznikov. 

Note that some cases of Thm. \ref{thm:subconvexct} --
where $\pi$ has trivial central character and $\pi$ is quadratic -- 
are subsumed by the previous result Thm. \ref{thm:subconvextp}. Nevertheless,
we have chosen to
give a distinct presentation since
the method is entirely different, it is simpler in the present method to deal
with the case of noncuspidal $\pi$. Also, we shall consistently
deal in the present section with $\GL(2)$, rather than $\PGL_2$.
Thus $\omega$ will be a unitary character of $\adele_F^{\times}/F^{\times}$,
and $C^{\infty}_{\omega}(\quotg)$ the space
of functions on $\quotg = \GL_2(F) \backslash \GL_2(\adele_F)$
with central character $\omega$.

In the case $F=\Q$, the subconvexity result  (\ref{eq:burgess})
for characters is due to Burgess.
Burgess' method gives a much better exponent; of course there
is considerable scope for improvement in the present technique also.

For the ease of the reader, we briefly explain in advance the
points of our proof in classical language. The discussion that
follows is not a completely faithful rendition of the proof, but
it hopefully conveys the main ideas. While it follows the pattern
of all the proofs in this paper, one minor complication is that we
deal with integrals w.r.t. certain measures of infinite mass.

\begin{enumerate} \label{enum:strategy}
\item  \label{firststep}If $f$ is a Maass form on $\SL_2(\Z)
\backslash \mathbb{H}$, the integral \begin{equation} \label{classical} 
\frac{1}{q}\int_{y=0}^{\infty}\sum_{x=1}^{q} \chi(x) f(\frac{x}{q}
+ i y) \frac{dy}{y}\end{equation} equals, up to some $\Gamma$-factors,
$\frac{1}{\sqrt{q}} L(\frac{1}{2}, f \times \chi)$. This is an exercise
in Hecke-Jacquet-Langlands theory. The version of
this equality that we shall use is proved in  Lem.
\ref{lem:callas} when $f$ is a cusp form and Lem. \ref{lem:hjleis}
when $f$ is Eisenstein.

\item  \label{secondstep} It will then suffice to bound
$\sum_{x=1}^{q} \chi(x) f(\frac{x}{q}+iy)$ for each fixed value of
$y$. As it turns out, the crucial range of $y$ is around $y =
q^{-1}$; the contribution of other $y$s are small for relatively
trivial reasons (use the Fourier expansion). This
is roughly a geometric form of the approximate functional equation: it says
that the Fourier coefficients $a_n(f)$ with $n \asymp q$ are most
important to determining the $L$-function. 
The general version
of this fact is proven in Lem. \ref{lem:callas2}.

 \item \label{thirdstep}
 In the
range when $y \asymp q^{-1}$, the set $\{\frac{x}{q} + iy\}_{\{1
\leq x \leq q-1\}}$ is roughly equidistributed, because it is
(with the exception of two points) the orbit of $i q y \in \mathbb{H}$ by the $q$th
Hecke operator. This is easy to quantify and
actually can be regarded as a statement about equidistribution of
{\em $p$-adic} horocycles.The general version of
this is proved in Lem. \ref{lem:measboundvar}.

\item  \label{fourthstep} We are now in a situation where we are
trying to bound the period of $f$ along the roughly
equidistributed set $\{\frac{x}{q}+iy\}_{\{1 \leq x \leq q-1\}}$.
To do this, we apply mixing properties of the adelic torus flow,
in the same fashion as the previous proofs of this paper. This
shows that $\sum_{x=1}^{q} \chi(x) f(\frac{x}{q} + iy)$ is small.
%The final result is Prop. \ref{prop:measprop}, proven on p.
%\pageref{finallyproof}.
\end{enumerate}

The computations that underlie steps (\ref{firststep}),
(\ref{secondstep}) and (\ref{thirdstep}) are fairly routine but
technically complicated. We have therefore carried them out in 
Sec. \ref{subsec:hjl}. In the sections that follow, we merely quote the results
and carry out what amounts to step (\ref{fourthstep}).

\subsection{Relating integral representations and periods.}
Let $z \in \mathbb{R}$ and let $\mu_z, \nu_z, \mu,\nu$ be the measures on $
\quotg$ defined by
\begin{multline}\label{eq:nudef}
\mu_z(f) = \int_{|y| =z} f(a(y)
n([\fcond])) \chi(y) d^{\times} y , \ \ \mu = \int_{z \in \mathbb{R}^{\times}}
\mu_z d^{\times}z, \\
\nu_z(f) =  \int_{|y|=z} f(a(y) n([\fcond])) d^{\times}y, \ \ \nu
= \int_{z \in \mathbb{R}^{\times}} \nu_z d^{\times}z
.\end{multline} In both cases, the measure $d^{\times}y$ is the
probability measure invariant by $\adele_F^{1}/F^{\times}$ and the
measure $d^{\times}z$ is a Haar measure on $\mathbb{R}^{\times}$.
Thus $\mu_z, \nu_z$ are probability measures, whereas $\mu,\nu$
have infinite mass. It is simple to see that the integrals
defining $\mu(f), \nu(f)$ converge absolutely if $f$ is a function decaying
rapidly enough at the cusps, e.g.
satisfying $|f(x)| \ll \height(x)^{-\varepsilon}$ (notation of Sec. \ref{sec:redtheory}), for any
$\varepsilon > 0$. Note also the analogy between these measures
and those used in the analysis of unipotent periods (cf.
(\ref{protomunu}).)
Classically, $\nu_z(f)$ should be thought of the measure
on $\SL_2(\Z) \backslash \mathbb{H}$ defined by
$\sum_{0 \leq x \leq q-1} f(\frac{x}{q} + iz)$,
and $\mu_z(f)$ the measure on $\SL_2(\Z) \backslash \mathbb{H}$
defined by $\sum_{0 \leq x \leq q-1} f(\frac{x}{q}+iz) \chi(x)$.  (These
statements are not to be interpreted precisely; they are for intuition only).

Here is the Proposition that formalizes (\ref{firststep}) and (\ref{secondstep}) of the discussion 
above, 
in the cuspidal case. 

\begin{prop} \label{prop:hjlcusp}
Let $\pi$ be a cuspidal representation on $\GL_2(\adele_F)$,
$\chi$ a character of $\adele_F^{\times}/F^{\times}$ of finite conductor $\mathfrak{f}$. 
Write $L_{unr}(s, \pi \times\chi) = \prod_{\chi_v \, \mathrm{unram.}} L_v(s,\pi\times \chi)$, where
the product is taken over all finite places at which $\chi$ is not ramified.

Let $d,\beta \geq 0$.
Let $g_{+}, g_{-}$ be positive smooth functions on $\mathbb{R}_{\geq 0}$
such that $g_{+} + g_{-} = 1$, $g_{+} (t)= 1$ for $t \geq 2$
and $g_{-}(t) = 1$ for all $t \leq 1/2$.  

 Then there exists $\varphi \in \pi$ such that, with
\begin{equation} \label{eq:tomyknees1} \Phi(s) = \Norm(\fcond)^{1/2} \frac{\int_z \mu_z(\varphi) |z|^{s-1/2} d^{\times}z}{L_{unr}(s, \pi \times \chi)}\end{equation}
then $\Phi(s)$ is holomorphic and satisfies:
\begin{enumerate}
\item $ |\Phi(s)| \ll_{\Re(s),\epsilon} \Norm(\fcond)^{\epsilon}$ and $|\Phi(\frac{1}{2})| \gg_{\epsilon} \Norm(\fcond)^{-\epsilon}$. 
%\begin{equation} \left|\int_z \mu_z(\varphi) |z|^{s-1/2} d^{\times}z \right| \ll_{\Re(s)} 
%\Norm(\fcond)^{-1/2+\epsilon}  L_{unr}(s) 
%\end{equation}
%\item 
%\begin{equation} \label{eq:john2}
% \left|\int_{z} \mu_z(\varphi) d^{\times}z \right| \gg_{\epsilon}
%\Norm(\fcond)^{-1/2-\epsilon} |L(\frac{1}{2}, \pi \otimes \chi)|\end{equation}
\item
$\varphi$ is new at every finite place (i.e.,
for each finite prime $\q$ it is invariant by $K_0[\q^{s_{\q}}]$, 
where $s_{\q}$ is the local conductor of $\pi$).
\item The Sobolev norms of $\varphi$ satisfy the bounds 
\begin{equation} \label{eq:john171}
S_{2,d,\beta}(\varphi) \ll_{\epsilon} \Cond_{\infty}(\pi)^{2 d+\epsilon}
\Cond_{f}(\pi)^{\beta+\epsilon}
 \Cond_{\infty}(\chi)^{1/2 + 2d}\end{equation}
\item  The integration of (\ref{eq:tomyknees1}) may be ``truncated without significant change'' to the region
$z$ around $\Norm(\fcond)^{-1}$; more formally:

$$\left| \int_{z} \mu_z(\varphi) g_{+}(z/T)d^{\times}z   \right| \ll (\Norm(\fcond) T)^{-1/2}
(T \Cond(\pi) \Cond(\chi))^{\epsilon}$$
$$ \left| \int_{z} \mu_z(\varphi) g_{-}(z/T) d^{\times} z  \right| \ll (\Norm(\fcond) T)^{1/2}
(T \Cond(\chi))^{\epsilon}(\Cond_{\infty}(\chi) \Cond(\pi))^{1+\epsilon} $$
\end{enumerate}
\end{prop}

\proof Lem. \ref{lem:callas} and Lem. \ref{lem:callas2}. \qed

We next give the corresponding result for the ``noncuspidal case.''
We recall that the Eisenstein series $E_{\Psi}(s,g)$ associated
to a Schwarz function $\Psi$ on $\adele_F^2$ are discussed in Sec. \ref{sec:eisenstein}. 
The normalization is so that the functional equation interchanges $s$ and $1-s$. 
$\bar{E}(s,g)$ denotes, as explained in that section
(cf. \eqref{eq:baredef}) the truncated Eisenstein series obtained by subtracting the constant term; 
it is a function on $B(F) \backslash \GL_2(\adele_F)$. 

\begin{prop} \label{prop:hjleis}
Let $s_0, s_0' \in \C$, and suppose that $\chi$ is ramified at at least one finite place. 
Let $g_{\pm}$ be as in Prop. \ref{prop:hjlcusp}.
There is an absolute $C > 0$ (i.e. depending only on $F$)
and a choice of $\Kmax$-invariant Schwarz function $\Psi$ (depending on $\chi$) so that if we put
$$ \Phi(s,s') := \Norm(\fcond)^{1/2} \frac{ \int_{y \in
\adele_F^{\times}/F^{\times}} \bar{E}_{\Psi}(s,a(y)n ([\fcond])) \chi(y) |y|^{s'}
d^{\times}y}{ L(\chi,s+s') L(\chi,1-s+s')}$$
where $\bar{E}$ is defined as in (\ref{eq:baredef}),  then the integral
defining $\Phi$ is absolutely convergent when $\Re(s),\Re(s') \gg 1$. Moreover, 
$\Phi$ extends from $\Re(s), \Re(s') \gg 1 $ to a holomorphic
function on $\C^2$,  satisfying
\begin{enumerate}
\item  $|\Phi(s_0, s_0')| \gg 1$ and $|\Phi(s,s')| \ll C^{1+|\Re(s)| + |\Re(s')|} (1+|s|+|s'|)^C$.

Moreover, given $N > 0$ we have that 
\begin{equation} \label{Bergeron} |\Phi(s,s')|
(1+|s|+|s'|)^{N} \ll_{\Re(s), \Re(s'), N} \Cond_{\infty}(\chi)^{N'}\end{equation} where $N'$ and the implicit constant may be taken
to depend continuously on $N, \Re(s), \Re(s')$. 
\item $\Psi$, and so also
$E_{\Psi}(s,g)$ is invariant by $\Kmax$;
\item Let $h \in \mathcal{H}(\kappa)$ be as in (\ref{eq:normdef}), 
and put $E_h := \int_{\Re(s) \gg 1} h(s) E_{\Psi}(s, g) dg$. 
Then, for each $d,\beta$, there is $N > 0$ such that
 $S_{\infty,d,\beta}(E_h) \ll_{\kappa} \|h\|_0 \Cond_{\infty}(\chi)^{N}$,
where the norm $\|h\|_0$ is defined in (\ref{eq:normdef}). \item 
We have $\mu_z(E_h) \ll_{K,\Psi,h} \min(|z|^K, |z|^{-K})$ for each\footnote{The implicit constants here are totally unimportant; this estimate will be used only to verify that certain integrals converge.} $K \geq 1$.  Moreover, 
there is $N > 0$ such that
$$ \left| \int_{z} \mu_z(E_h) g_{+}(z/T)d^{\times}z \right|  \ll (\Norm(\fcond) T)^{-1/2}
(T \Cond(\chi) )^{\epsilon} \|h\|_N $$
$$  \left| \int_{z} \mu_z(E_h) g_{-}(z/T) d^{\times} z \right| \ll  (\Norm(\fcond) T)^{1/2}
(T \Cond(\chi) )^{\epsilon} \Cond_{\infty}(\chi)^{1+\epsilon}  \|h\|_N$$
\end{enumerate}
\end{prop}
\proof Lem. \ref{lem:hjleis} and Lem. \ref{lem:callas3}.  \qed

\subsection{Proof of theorem -- cuspidal case.} \label{sec:proofcusp}

Let $\chi$ be a character of $\adele_F^{\times}/F^{\times}$, of
varying conductor $\fcond$. Put $q = \Norm(\fcond)$.

We shall need the following estimate, proved in Lem. \ref{lem:measboundvar}.
It amounts in essence to a statement about the equidistribution of 
$p$-adic horocycles (classically, these roughly correspond to a statement
about the equidistribution of $\{\frac{x}{q} + i z\}_{0 \leq x \leq q-1}$,
if $z \asymp q^{-1}$).

For any function $f$ that is invariant by
$\prod_{\q} K_0[\q^{s_{\q}}]$, we have (with $\m = \prod_{\q} \q^{s_{\q}}$)
\begin{equation}\label{eq:mdr1}
\left|\nu_z(f) - \int_{\quot} f\right| \ll_{\epsilon}
q^{\alpha-1/2 + \epsilon} \Norm(\m)^{3/2+\epsilon} \max(qz, \frac{1}{qz})^{1/2}
S_{2,d}(f),  \,\,f \in C^{\infty}(\quot)\end{equation}

\proof \label{finallyproof}(of Thm. \ref{thm:subconvexct} -- cuspidal case.)
Choose $f \in \pi$ to be the ``$\varphi$'' of Prop. \ref{prop:hjlcusp}, so that
$|\mu(f)| \gg_{\epsilon} \Norm(\fcond)^{-1/2-\epsilon} |L_{unr}(1/2, \pi \times \chi)|$. 
For each ramified prime $\q$ of $\pi$, let $\q^{s_{\q}}$ be the local conductor
of the local representation $\pi_{\q}$.  Set $\mathfrak{m} := \prod_{\q} \q^{s_{\q}}$,
the finite conductor of $\pi$. 

Let $K\geq 1$ be an integer satisfying $K \leq \Norm(\fcond)$.
Let $\set$ be the set of prime ideals of $\order_F$,
with norm lying in $[K,2K]$, and satisfying $(\n , \fcond) = 1$
and $(\n, \Supp(f)) = 1$.  Fix $\n_0 \in \mathcal{S}$. 
For each prime ideal $\n \in \set$, let $\unif_{\n} \in F_{\n}$
be a uniformizer.
We define a measure $\meas$ on $\GL_2(\adele_F)$
so that
\begin{equation} \label{sigmameasuredefn}\meas =
|\set|^{-1}\sum_{\n \in \mathcal{S}} \chi(\unif_{\n}) \chi(\unif_{\n_0})^{-1} 
\delta_{a_{\n}(\unif_{\n}) a_{\n_0}(\unif_{\n_0})^{-1}}.\end{equation}
 Clearly $\meas$ has total mass $1$.
Moreover, $\mu(f) = \mu(f \star \meas)$ and $f \star \meas$ is
invariant by $K_0[\q^{s_{\q}}]$ for
each $\q|\m$.

 Choose $\kappa$ ``slightly smaller than $1$,'' to be
specified later. Our aim is now to cut the $z$
integration in $\mu = \int_{z} \mu_z d^{\times}z$ into three ranges,
the crucial range of which will be $q^{-2+\kappa} \ll z \ll q^{-\kappa}$;
this avoids the pain of dealing with the infinite mass measure $\mu$. 
Let $g_{+}, g_{-}$ be as in Prop. \ref{prop:hjlcusp}.
  Define $h(t)$
by the rule $g_{-}(\frac{z}{q^{-2+\kappa}}) + h(t) + g_{+}(\frac{t}{q^{-\kappa}})= 1$. 
\begin{multline} \label{kd}
|\mu(f)|^2 = |\mu(f \star\meas)|^2
= \left|\int_{z} \mu_z(f \star \meas) d^{\times}z\right|^2
\\ \ll \left|\int  
g_{-}(\frac{z}{q^{-2+\kappa}}) \mu_z(f \star \meas) \right|^2
 + \left|\int h(z) \mu_z(f \star \meas)
 d^{\times} z \right|^2
+  \left|\int g_{+}(\frac{z}{q^{-\kappa}}) \mu_z(f\star \meas) d^{\times}z \right|^2
\end{multline}
%\begin{multline}\mu_z(f) = |\set|^{-1}\int_{|y| = z} f(a(y) n([\fcond]))
%\sum_{\n \in \set} \chi(\omega_n) \chi(a(y) \omega_n)
%By (\ref{eq:mdr1})

By Prop. \ref{prop:hjlcusp},  the first and last term
 (without the square, e.g. $\left|\int g_{+}(\frac{z}{q^{-\kappa}}) \mu_z(f\star \meas) d^{\times}z \right|$)
 are $ \ll 
\Cond(\pi)^{1+\epsilon} \Cond_{\infty}(\chi)^{1+\epsilon}  q^{\frac{\kappa-1}{2}+\epsilon}$. More explicitly, we note that
\begin{equation} \label{eq:jay} \mu_z (f \star \meas) = |S|^{-1} \sum_{\n \in \mathcal{S}} \mu_{\Norm(\n)^{-1}\Norm(\n_0) z}(f) \end{equation}
Now $1/2 \leq \Norm(\n) \Norm(\n_0)^{-1} \leq 2$ for all $\n$ -- this was the purpose of the factors
involving $\n_0$ in (\ref{sigmameasuredefn}) -- and so one easily deduces a bound
on $\int g_{-}(\frac{z}{q^{-2+\kappa}}) \mu_z(f \star \meas) d^{\times}z$ from the final assertion of
Prop. \ref{prop:hjlcusp}.  Similarly for the term involving $g_{+}$.

That the first and last
terms of \eqref{kd} should be less significant may be seen in the classical setting from the Fourier expansion;
it should be regarded as a geometric version of the approximate functional equation). 

As for the intermediate term, we note $$\int_{z} h(z) \mu_z (f) 
= \int_{y \in \adele_F^{\times}/F^{\times}} h(|y|) \chi(y) f(a(y) n([\fcond])) d^{\times}y.$$ Applying
Cauchy-Schwarz, and the fact that $\int_{\adele_F^{\times}/F^{\times}} h(|y|)d^{\times}y 
\ll \log(q)$, we get:
\begin{multline}\left|\int_{z} h(z) \mu_z(f \star \sigma)
\right|^2 \ll_{\epsilon} q^{\epsilon} \int_{z}
h(z) \nu_z(|f \star \sigma|^2) d^{\times}z
\\ \ll_{\epsilon} q^{\epsilon} \int_{\quot} |f \star \sigma|^2 d\mu_{\quot}
+ q^{\alpha - \kappa/2+\epsilon} \Norm(\m)^{3/2+\epsilon} S_{2,d}(|f \star
\meas|^2),\end{multline}
where we have applied (\ref{eq:mdr1}).
%Put $I(f) = \int_{|y| \geq q^{-1}} f(a(y) n([\fcond]))
%\chi(y) d^{\times}y$ and fix $\kappa > 0$.
%Then
%we may write $I(f) = I_{1} + I_{2}$, where
%\begin{equation}\begin{aligned} I_1 =
%\int_{q^{-1}}^{q^{-\kappa}} \mu_z(f) d^{\times}z
%= \int_{q^{-1}}^{q^{-\kappa}} \mu_z(f \star \mu) d^{\times}z,  \\
%I_2 =
%\int_{q^{-\kappa}}^{\infty} \mu_z(f) d^{\times}z \end{aligned}\end{equation}
%
%Then, applying
%Cauchy-Schwarz and (\ref{eq:mdr1}), $$|\mu_z(f)|^2 =| \mu_z(f \star \meas)|^2 \leq \nu_z(|f \star
%\meas|^2) \ll \int_{\quot} |f \star \meas|^2 d\mu_{\quot}+
%|z|^{1/2} q^{\alpha} S_{\infty,d,\beta}(|f \star \meas|^2)$$
By Lemmas \ref{lem:sobproduct} and \ref{lem:sobgroup}, \begin{multline}
\label{eq:chembai}S_{2,d,\beta}(|f \star \meas|^2) \ll
S_{4,d,\beta}(f \star \meas)^2  \\ \ll
 (\mathrm{sup}_{g \in \mathrm{supp}(\meas)}\|g\|)^{ 2\beta}
S_{4,d,\beta}(f)^2 \ll K^{4 \beta}
S_{2,d',\beta+3/2}(f)^2.\end{multline}
where the last line holds for $d' \gg d$, and we have used Lem. \ref{lem:cuspfunctions} (which bounds the $L^{\infty}$ norm of a cusp form in terms of $L^2$ norms), together with the easily
verified fact that $\mathrm{sup}_{g \in \mathrm{supp}(\meas)} \|g\| \ll K^2$. 

By bounds towards Ramanujan,
\begin{equation} \label{eq:boundstowardsR}
 \|f \star \meas\|_{L^2}^2 =
\int_{g,g'} \langle g^{-1} g' f, f  \rangle d\meas(g) d\meas(g') \ll
K^{2 \alpha - 1} \|f\|_{L^2}^2.\end{equation}
%and $\|g\|_{L^2} = \left(\int_{\quot}
%|g|^2 \right)^{1/2}
%\leq \|\psi \|_{L^{\infty}} \|\psi \star \meas\|_{L^2}$.
%If $\psi$ is orthogonal to the constants,
%$\|\psi \star \meas\|_{L^2} \ll K^{2 \alpha - 1}$.

%Thus
%$|\mu_z(f)|^2 \ll_{\epsilon}
%K^{2\alpha-1} + K^{1+\epsilon} |z|^{1/2} q^{\alpha},$
%and it follows that
%$|I_1| \ll (K^{ \alpha - 1/2} + q^{-\kappa/4 +  \alpha/2}
%K^{1/2+\epsilon}) S_{\infty,d,1/2}(f).$
Thus:
\begin{multline} 
|\mu(f)| \ll_{\epsilon}
\Cond(\pi) \Cond_{\infty}(\chi) ( \Cond(\chi) \Cond(\pi))^{\epsilon} q^{\frac{\kappa-1}{2}} \\ 
+ \left( K^{\alpha - 1/2} q^{\epsilon} +\Norm(\m)^{3/4+\epsilon}  q^{-\kappa/4+\alpha/2+\epsilon} K 
\right)  S_{2,d,2}(f) \\
\ll q^{\epsilon} (q^{(\kappa-1)/2} + K^{\alpha-1/2} +  q^{\alpha/2-\kappa/4} K)
\Cond_{\infty}(\chi)^N \Cond(\pi)^N,\end{multline} 
%We take $\kappa =
%\frac{2 (3-\alpha-2\alpha^2)}{7-6\alpha}$,
%$K=q^{\frac{1-\kappa}{1-2\alpha}}$, and conclude that $|I|
%\ll_{\epsilon} q^{-\frac{(1-2\alpha)^2}{2 (7-6\alpha)}}
%S_{2,d,1+\epsilon}(f)$.  
for appropriate $N > 0$.  We have used Prop. \ref{prop:hjlcusp}, (3) at the last step. 

Prop. \ref{prop:hjlcusp} guarantees that
$|L_{unr}(1/2, \pi \times \chi)| \ll_{\epsilon} q^{1/2+\epsilon} |\mu(f)|$. 
From this, optimizing $\kappa,K$, and applying trivial bounds at ramified places, we obtain the conclusion, taking $\alpha =
3/26$.   \qed

%here's the justification at least of $1/25$. we want basically to minimize $q^{(\kappa - 1)/2} + K^{\alpha-1/2} + 
%K q^{\alpha/2-\kappa/4}$. Set them all equal. Then $K = q^{(1- \kappa)/(1-2\alpha)}$.
%Also, $ K  = q^{3 \kappa/4 - 1/2 - \alpha/2}$. So $(1-\kappa)/(1-2 \alpha) = 3 \kappa/4 - 1/2 -\alpha/2$. 
%We get $\kappa = \frac{2 (3- \alpha - 2 \alpha^2)}{(7- 6 \alpha)}$ So with $\alpha = 3/26$, 
%just take $\kappa = 9/10$ and $K = q^{1/8}$. Then $q^{(\kappa-1)/2} = q^{-1/20},
%K^{\alpha - 1/2} = q^{-5/104)}$ and $K q^{\alpha/2 - \kappa/4} = q^{1/8 - 9/40 + 3/52} = q^{3/52-1/10}$. 
%All are better than $q^{-1/25}$. 
\subsection{Proof of theorem -- noncuspidal case.}
We turn to the proof of 
 (\ref{eq:burgess}).
This is very similar, but we implement a mild regularization procedure
to deal with the Eisenstein series, just as in the case of Rankin-Selberg $L$-functions. 

\proof (of (\ref{eq:burgess}).) 
We may assume that $\chi$ ramifies at least at one finite place. 

Let $\Psi$ be a Schwarz
function on $\adele_F^2$, $E(g,s) := E_{\Psi}(g,s)$ the corresponding
Eisenstein series, chosen as in Prop. \ref{prop:hjleis} with $s_0 = 1/2, s_0' = 0$. 

Let $\kappa' > 0$, let
$h$ be holomorphic in an open neighbourhood of the vertical strip
$-\kappa' \leq \Re(s) \leq 1+ \kappa'$ and put $E_h(s) =
\int_{\Re(s) = 1+\kappa'} h(s) E(g,s) ds$. Then if $h(0) =
h(\frac{1}{2}) = h(1) =0$, it follows from the third assertion of Prop. \ref{prop:hjleis}
that $$S_{\infty,d,\beta}(E_h) \ll_{\epsilon,d}  \Cond_{\infty}(\chi)^N \|h\|_0$$
 for appropriate $N =N(d,\beta)> 0$. 
Here, as in (\ref{eq:normdef}) with $\kappa$ replaced by $\kappa'$, the norm $\|h\|_N $ is defined to be $ \int_{-\infty}^{\infty} \left(|h(1+\kappa'+it)|+ |h(-\kappa' +
it)|\right) (1+|t|)^{N} dt$. 

Put, in the  notation of Prop. \ref{prop:hjleis}, $I(s) = \Phi(s,0) L(\chi, s) L(\chi, 1-s)$.  
Then: \begin{equation} \label{eq:jeff} 
\int_{z } \mu_z(E_h) d^{\times}z
% \int_{y \in \adele_F^{\times}} E_{h}(a(y) n([\fcond]))
%\chi(y)d^{\times} y 
%= \int_{\Re(s) = 1+\kappa} h(s) \int_{y \in \adele_F^{\times}/F^{\times}}
%E(s,a(y) n[\fcond]) \chi(y) d^{\times} y \\ 
%= \int_{\Re(s) = 1+\kappa} h(s) \int_{y \in \adele_F^{\times}/F^{\times}}
%\bar{E}(s,a(y) n[\fcond]) \chi(y) d^{\times} y 
 = \Norm(\fcond)^{-1/2} \int_{\Re(s) = 1/2}h(s) I(s) ds.\end{equation}
 
 This is established in \eqref{wiw}; for now, we remark that that this is ``almost''
 obvious from Proposition \ref{prop:hjleis},  the only additional point being that one can replace
 $E$ by $\bar{E}$, and this is exactly where the fact that $\chi$ is ramified at a finite place comes in -- to kill
 the constant term of the Eisenstein series.\footnote{The classical version of this fact -- see \eqref{classical} --
 it is clear that the $\chi$-sum will kill any constant term of $f$, as long as $\chi$ is not trivial.}
% 
% Here the replacement of $E$ by $\bar{E}$ is permissible by the assumption
% that $\chi$ is ramified at one finite place, since $E - \bar{E}$ is left invariant by $N(\adele_F)$ 
% (cf. discussion after  \eqref{bathsheba}). At the last stage, we have shifted contours,
% which is justifiable in view of the holomorphy of $I(s)$. 
 
Take $h = (s-1/2)^2 s (1-s) \overline{I(1-\overline{s})}$. The ``good''  analytic properties of $h$, e.g.
rapid decay along vertical lines,  follow\footnote{This point was not clear in a previous version; thanks to N. Bergeron for pointing this out.} from 
\eqref{Bergeron}.  In particular, $h$ belongs to the function spaces $\mathcal{H}(\kappa)$
defined in (\ref{eq:normdef}) for any $\kappa > 0$, and the norms $\|h\|_N$ are all bounded by suitable
powers of $\Cond_{\infty}(\chi).q$. 

Then (\ref{eq:jeff}) becomes
\begin{equation} \label{eq:mmi1} \int_{t = -\infty}^{\infty} t^2 (1/4+t^2) |I(\frac12 + it)|^2 = q^{1/2} \int_{z} \mu_z(E_{h}) d^{\times} z \end{equation}

To bound the right-hand side, we proceed as in Sec. \ref{sec:proofcusp}, but with $f$ replaced by $E_h$. 
We use notation as in that Section, except replacing the ``$h$'' defined before \eqref{kd}
by $1-g_{-}-g_{+}$ to avoid clashing with its alternate usage here. 
 
One proves as in that Section, that for $d \gg 1$:

\begin{multline}  \left| \int_{z} (1-g_{-}-g_{+}) \mu_z (E_h \star \sigma) \right| \ll_{\epsilon}  q^{\epsilon} \left(K^{\alpha-1/2}  + q^{\alpha/2-\kappa/4} 
 (\mathrm{sup}_{g \in \mathrm{supp}(\meas)}\|g\|)^{1/2} \right)
S_{4,d,1/2}(E_h) \\  \ll
q^{\epsilon}   \left(K^{\alpha-1/2}  + q^{\alpha/2-\kappa/4} K \right) \|h\|_0 \Cond_{\infty}(\chi)^N
\end{multline}
for some appropriate $N > 0$. At the last stage we have applied Prop. \ref{prop:hjleis} to control the Sobolev norm.

 Prop. \ref{prop:hjleis} also guarantees that, for appropriate $N > 0$, we have:
\begin{multline} \label{eq:mmi2}  \left|\int  
g_{-}(\frac{z}{q^{-2+\kappa}}) \mu_z(E_h \star \meas) \right|
+  \left|\int g_{+}(\frac{z}{q^{-\kappa}}) \mu_z(E_h \star \meas) d^{\times}z \right| \\  \ll 
\Cond_{\infty}(\chi)^{1+\epsilon}  q^{\frac{\kappa-1}{2}+\epsilon} \|h\|_N,
\end{multline}

Combining \eqref{eq:mmi1} and \eqref{eq:mmi2}, we obtain as in the previous Section the bound,
for sufficiently large $N$: 
\begin{multline} \int_{t=-\infty}^{\infty} t^2 (\frac{1}{4}+t^2)
|L(\frac{1}{2}+it,\chi) L(\frac{1}{2}-it, \chi) |^2 |\Phi(1/2+it,0)|^2  \\ \ll_{\epsilon} \left(
 q^{\frac{\kappa-1}{2}}  
+   K^{\alpha-1/2}  + q^{\alpha/2-\kappa/4} K\right) \|h\|_N  q^{1/2+\epsilon}\Cond_{\infty}(\chi)^N
\end{multline} 
One applies the convexity bound to bound $\|h\|_N$, obtaining  
$$\int_{-\infty}^{\infty} t^2 |L(\frac{1}{2} + it,\chi)| L(\frac{1}{2} - it, \chi)|^2 |\Phi(1/2+it,0)|^2 \ll \Cond_{\infty}(\chi)^N q^{24/25},$$ where we have increased $N$ as necessary. From this we get
$L(\frac{1}{2}, \chi) \ll \Cond_{\infty}(\chi)^N q^{1/4-1/200}$. 
\qed

\section{Torus periods (II): equidistribution of compact torus orbits.}
\label{sec:torus2}

It has been independently shown by Zhang \cite{Zhang}, Clozel-Ullmo \cite{ClozelUllmo} and P. Cohen
\cite{Cohen}
that the subconvexity result Thm. \ref{thm:subconvexct} implies 
the equidistribution of Heegner points over totally real fields; in particular,
they pointed out that GRH implies this equidistribution. Thm. \ref{thm:subconvexct}
makes this result unconditional.

The main aim of this section is to explain how one can
obtain certain conditional results about equidistribution
of subsets of Heegner points, and how this
fits into the general framework of ``sparse equidistribution questions.'' In particular, this approach 
does not rely on reducing questions about subsets of Heegner points to subconvexity, but 
rather approaches the equidistribution question directly. 

The proofs of the results (and various supporting Lemmas) will only be sketched, and we will
confine ourselves for simplicity to the case of narrow class number $1$; we will in any case present an unconditional
approach, based on combining the ideas of this paper with the ideas of Michel, 
in the paper \cite{MV} (joint with P. Michel).  We nevertheless
feel that the ideas presented here may be of use in other contexts. Indeed,
this section is of a different flavor to the other Sections; it uses ``adelic analysis'' more
genuinely.

In fact, we shall need a mild refinement of the results of \cite{ClozelUllmo}, which
allow better control of the dependence on the test vectors. 
We state this refinement without proof in Thm. \ref{thm:heegner}; the proof
is an exercise in explicating some of the proofs in \cite{ClozelUllmo}. 
%We sketch this refinement in Sec. \ref{subsec:curef}.  [[change if you delete]] 
\subsection{Equidistribution of Heegner points.} \label{subsec:equidistHeeg}
We recall the definition of Heegner points.
Let $F$ be a totally real number field of degree $d$ over $\Q$. 
For simplicity we shall confine ourselves to the case where the ring of integers
of $F$ has narrow class number $1$.
 This assumption does not change any of the technical details, which are in any
 case carried out adelically; it  simply allows us to be a little more explicit
 about the torus orbits we consider. 
Let $E = F(\sqrt{-\dE})$ be a totally imaginary quadratic
extension of $F$, where $\dE \in \order_F$ is totally positive and
squarefree. Here ``squarefree'' means that it is of valuation $\leq 1$ at all finite places.

 Let $T_E$ be the torus
$\mathrm{Res}_{E/F}(\mathbb{G}_m)/\mathbb{G}_m$; we embed $T_E$ in
$\PGL_2$ via (in obvious notation):
\begin{equation} \label{eq:janani}\iota_E: x + y \sqrt{-\dE} \mapsto \left(\begin{array}{cc} x & y \\ -y\dK & x \end{array}\right).\end{equation}

Regard $\dE$ as an element of $F \otimes \R$ via the inclusion
$F \hookrightarrow F \otimes \R$. Since it is totally positive,
it possesses a unique totally positive square root,
$\sqrt{\dE} \in F\otimes \R$.
Set $[\dE]_{\infty} = \left(\begin{array}{cc} 1 & 0 \\ 0  & \sqrt{\dE} \end{array}\right) \in \PGL_2(F \otimes \R)$.
We define a map $\heegmap:  T_E(\adele_F)/T_E(F)
 \rightarrow
\quot$
via
\begin{equation} \label{eq:heegmap}
\heegmap: x \mapsto \iota_E(x) [\dE]_{\infty},\end{equation}
where we regard $[\dE]_{\infty} \subset \PGL_2(F \otimes\R) \subset
\PGL_2(\adele_F)$ acting by right translation on $\quot$.
Denote by $\Norm(\dE)$ the absolute norm of $\dE$, 
i.e. $\Norm(\dE) = |\order_F/\dE \order_F|$.

The $F$-torus $T_E$ is anisotropic, and
there is a unique $T_E(\adele_F)$-invariant
probability measure on $T_E(\adele_F)/T_E(F)$.
Let $\nu_E$ be its image by the map $\heegmap$.

\begin{thm} \label{thm:heegner}
Set $E = F(\sqrt{-\dE})$, where $\dE \in \order_F$ is totally
positive and squarefree. The measures $\nu_E$ become
equidistributed as $\Norm(\dE) \rightarrow \infty$. Indeed, there
exist $\delta > 0, d,\beta$ such that for $f \in
C^{\infty}(\quot)$ we have
$$ \left|\int f d\nu_E
- \int_{\quot} f(x) dx \right|
\ll \Norm(\dE)^{-\delta} S^{*}_{\infty,d,\beta}(f).$$
\end{thm}

Recall the definition of $S^{*}$ from Sec. \ref{subsec:ocean11}.
We do not give the proof; as we have remarked it can be obtained by following
the computations of \cite{ClozelUllmo} a little more explicitly. 
%\proof 
%See Sec. \ref{subsec:curef}. (We note that the proof is only sketched).  [[change if you delete]]
%\qed
%\footnote{As presently written,
%Thm. \ref{thm:subconvexct} only implies it
%assuming, e.g., that the smallest prime divisor
%of $\dE$ is getting large, owing
%to the coprimality assumption of Thm. \ref{thm:subconvexct}.
%This is inessential; I plan to remove it in a later version.}

One recovers from Thm. \ref{thm:heegner} the equidistribution of
certain Heegner points associated to $E = F(\sqrt{-\dE})$ as $\dE$
varies. 
%\footnote{Actually, with our definitions, I think the Heegner point
%orbit defined above may have conductor involving primes dividing $2$.}
%\footnote{An entertaining consequence of this
%equidistribution is the following: denote by $\mathbb{P}^1$ the
%projective line, and define as usual the hight of a point in
%$\mathbb{P}^1(\overline{\Q})$. For any number field $E$, let
%$\min(E)$ be the lowest height of a point of $\mathbb{P}^1(E)$.
%Then the above result can be used to show that $\min(E)$
%approaches a constant as $E$ varies through a sequence of
%quadratic extensions of a fixed field. The question of the
%lowest-height point on the projective line was raised to me by J.
%Ellenberg.}
 Thm. \ref{thm:heegner} also gives an effective rate of
equidistribution for Heegner points with polynomial dependence on
the level and the eigenvalue of a test function. This rather
innocuous polynomial dependence (in the level aspect, at least)
will in fact play a crucial role in our deduction of the
equidistribution of sparse subsets in the following section.

\subsection{Equidistribution of subsets of Heegner points.}
We turn to certain conditional
%\footnote{
%Note that it is plausible that these could
%be made unconditional; I am not sure. It maybe
%that the condition on small split primes
%is weak enough that it is either trivially true
%or within the range of subconvexity results like
%(\ref{eq:burgess}). But I don't know, this requires
%being careful about $d,\beta$ in Thm. \ref{thm:heegner}.}
results
on equidistribution of sparse subsets. $F$ being as in 
Sec. \ref{subsec:equidistHeeg},
let $E_i = F(\sqrt{\dE_i})$ be a sequence of distinct quadratic,
totally imaginary, extensions of $F$.
For each $E_i$, let $S_i \subset T_{E_i}(\adele_F)/T_{E_i}(F)$
be a subgroup of finite index $m_i$. Let $\mu_{E_i}^{S_i}$
be the image of the Haar probability measure on $S_i$
by the map $\heegmap$.

The import of the next theorem is that, {\em if $E_i$
has enough small split primes}, one can obtain the equidistribution
of the measures $\mu_{E_i}^{S_i}$ as $i \rightarrow \infty$. 
This result is quite similar to the results of Duke-Friedlander-Iwaniec \cite{DFI} in the case $F=\Q$,
although the method is at least superficially rather different. 
One can also contrast with
Michel's striking result, for $ F = \Q$,
that gives a comparable result but without the condition
on enough small split primes.  Our method
is different to Michel, who deduces
the result from his subconvexity bound for Rankin-Selberg $L$-functions.
\footnote{We note that this bound of Michel is considerably
deeper than (\ref{eq:double}),
since it deals with varying central character. For some speculative discussion on the ``reason''
that Michel's method can avoid this condition, see the last paragraph of \cite{MV-icm}.}
In the present approach, we prove
the equidistribution theorem
directly.
% \footnote{
%We remark again on the surprising point
%that, although the approach of Michel and the approach of the present paper
%are apparently very different,
%they both show a similar inductive structure
%in that they both reduce the problem to the subconvexity
%of character twists.} 
 In a sequel to this paper, the author and P. Michel
combine the methods here with some methods developed by Michel
to make the results of this section unconditional.

To quantify the existence of enough small split primes, one might impose the condition (as does Linnik
\cite{linnik})
that the $E_i$ vary through a sequence of quadratic
extensions that split at a fixed prime of $F$.
We will prefer to take a more quantitative
approach, which will yield a stronger result at the price of a stronger
assumption.
In that regard we introduce the following notation:
For $\delta > 0$, we put
$$\weight(E_i, \delta) = \#\{\mathfrak{q} \subset \order_F
\mbox{ prime and split in }E_i,
\Norm(\dE_i)^{\delta} \leq \Norm(\mathfrak{q}) \leq
2\Norm(\dE_i)^{\delta}\}.$$

\begin{thm} \label{thm:sparseheegner}
There exists $\delta_1 > 0$ such that, if
$\frac{m_i}{\min(\Norm(\dE_i)^{\delta_1(1/2-\alpha)},
\weight(E_i,\delta_1)^{1/2})} \rightarrow 0$, the
sequence $\mu^{S_i}_{E_i}$ converges, as $i \rightarrow \infty$,
to the invariant measure on $\quot$.
\end{thm}
\proof
This is deduced from Thm. \ref{thm:heegner}
by using the mixing properties of the $T_E(\adele_F)$-flow.
Indeed, we fix an index $i$ and a corresponding field $E_i$.
Let $\delta_1 > 0$ be fixed.
Let $\set$ be the set of prime
ideals of $F$ which split in $E_i$ and with norm in $[\Norm(\dE_i)^{\delta_1}, 2 \Norm(\dE_i)^{\delta_1}]$.
For each $\q \in \set$, the torus
$T_{E_i}(F_{\q})$ is isomorphic to $F_{\q}^{\times}$.
Fix an isomorphism $\iso_{\q}: T_{E_i}(F_{\q}) \rightarrow F_{\q}^{\times}$,
and let $\unif_{\q}$ be an element in $T_{E_i}(F_{\q})$
such that $\iso_{\q}(\unif_{\q})$ has valuation $\pm 1$ in $F_{\q}^{\times}$.

Let $\chi$ be a character of $T_{E_i}(\adele_F)/T_{E_i}(F)$,
trivial on $S_i$.
Let $\nu_{E_i}$ be as defined prior to Thm. \ref{thm:heegner},
and define
$$\mu_{E_i}(f) = \int_{t \in T_{E_i}(\adele_F)/T_{E_i}(F)}
f(\heegmap(t)) \chi(t) dt,$$
where $dt$ is the Haar probability measure on
$T_{E_i}(\adele_F)/T_{E_i}(F)$. Let $\meas$ be the
probability measure $\frac{1}{|\set|}\sum_{\q \in \set}
\chi(\unif_{\q}) \delta_{\unif_{\q}}$ on
$T_{E_i}(\adele_F)$. Then 
$\mu_{E_i}(f) = \mu_{E_i}(f \star \heegmap_{*}\meas)$, where
$\heegmap_{*} \meas$ denotes the image of $\meas$ by the map
$\heegmap$.

By Cauchy-Schwarz, and Thm. \ref{thm:heegner},
\begin{multline}|\mu_{E_i}(f \star \heegmap_{*}\meas)|^2 \leq \nu_{E_i}(|f \star
\heegmap_{*}\meas|^2)
\\ \leq \|f \star  \heegmap_{*} \meas\|_{L^2}^2 +
O\left(\Norm(\dE_i)^{-\delta}S^{*}_{\infty,d,\beta}(|f \star
\heegmap_{*} \meas|^2)\right),\end{multline}
where $\delta,d,\beta$ are as in Thm. \ref{thm:heegner}.
Now,  appropriate variants of Lem. \ref{lem:sobproduct}
and \ref{lem:sobgroup} (for $S^{*}$ instead of $S$) show that
\begin{multline} \label{eq:ariyakudi}
S^{*}_{\infty,d,\beta}(|f \star \heegmap_{*}\meas|^2)
\ll S^{*}_{\infty,d,\beta}(f \star \heegmap_{*}\meas)^2
\\ \ll \sup_{g \in \mathrm{supp}\heegmap_{*}\meas} \|g\|^{6 \beta}
S^{*}_{\infty,d,\beta}(f)^2 \ll \Norm(\dE_i)^{6 \delta_1 \beta}
S^{*}_{\infty,d,\beta}(f)^2\end{multline}
 and bounds towards Ramanujan show that
\begin{equation} \label{eq:Rambounds}
\|f \star\heegmap_{*}\meas\|_{L^2}^2 \ll \left(\Norm(\dE_i)^{\delta_1 (2
\alpha-1)} +|\set|^{-1}\right) \|f\|_{L^2}^2. \end{equation}

We note that
(\ref{eq:ariyakudi}) and (\ref{eq:Rambounds})
are very closely analogous to
(\ref{eq:chembai}) and (\ref{eq:boundstowardsR}),
with $K$ replaced by $\Norm(\dE_i)^{\delta_1}$.
In the context of (\ref{eq:boundstowardsR}), the set
$\set$ has size $K^{1-\epsilon}$;
thus the term
$|\set|^{-1}$ that appears in (\ref{eq:Rambounds}) could
be neglected.

Recalling the definition of $\mu_{E_i}$, we conclude
\begin{multline} \label{eq:mubound}
\left| \int_{t \in T_{E_i}(\adele_F)/T_{E_i}(F)}
f(\heegmap(t)) \chi(t) dt \right|
\\ \ll
\left(\Norm(\dE_i)^{ 3 \delta_1 \beta - \delta/2} +\Norm(\dE_i)^{\delta_1
(\alpha-1/2)} + |\set|^{-1/2}
\right) S^*_{\infty,d,\beta}(f).\end{multline}

Summing the left-hand side of (\ref{eq:mubound})
over all $m_i$ characters $\chi$ of $T_{E_i}(\adele_F)/T_{E_i}(F)$ that are trivial on $S_i$,
and substituting $|\set| =\weight(E_i,\delta_1)$, we obtain:
$$\left|\mu_{E_i}^{S_i}(f)  \right|
\ll  m_i
\left(\Norm(\dE_i)^{3 \delta_1 \beta - \delta/2} +\Norm(\dE_i)^{\delta_1
(\alpha-1/2)} + \weight(E_i,\delta_1)^{-1/2}
\right) S^{*}_{\infty,d,\beta}(f).$$

Choosing $\delta_1$ sufficiently small
(the exact value will depend on the value
of $\beta,\delta$ from Thm. \ref{thm:heegner})
we obtain the claimed conclusion.  \qed

%$$|\mu_{K_i}(f)| \ll \left(\Norm(\dK_i)^{-\min(\frac{\delta}{4},
%\frac{\delta (\alpha-1/2)}{4 \beta}} + \weight(K_i, \delta_1)^{-1/2}\right)
%S_{\infty,d,\beta}(f).$$
%\qed

\section{Background on 
Sobolev norms and reduction theory.}
\label{sec:sobolev}

The rest of the paper consists of technical Lemmas. The sections that follow
are arranged to be used as a reference, rather than to be read through. 

\subsection{Formal properties of the Sobolev norms.}
We begin by explicating certain formal properties of the Sobolev norms defined
in Sec. \ref{subsec:adelicsobolev}.
 
\begin{rem} \label{rem:def}
The following properties of this definition are formal and will be
repeatedly used:

\begin{enumerate}  
\item \label{remdef1}
Translations by $\Kmaxg$ preserve $S_{p,d,\beta}$,
i.e. $S_{p,d,\beta}(k \cdot f) = S_{p,d,\beta}(f)$ for
$k \in \Kmaxg$.
\item  \label{remdef2}
If $L : C^{\infty}_{\omega}(\quotG) \rightarrow \C$ is a linear functional
and $|L(\psi)| \leq PS_{p,d,\beta}(\psi)$,
then also $|L(\psi)| \leq S_{p,d,\beta}(\psi)$.
Indeed $\psi \mapsto |L(\psi)|$ is itself a seminorm on
$C^{\infty}_{\omega}(\quotG)$.
%\item \label{variant} {\em Variant:
%Suppose $K_0 \subset \Kmax$ is a
%compact subgroup, and $L: C^{\infty}_{\omega}(\quotg) \rightarrow
%\C$ a linear functional. Suppose that $|L(f)| \leq
%PS_{p,d,\beta}(\psi)$ for all $\psi$ that are $K_0$-invariant. Then,
%for all $K_0$-invariant $\psi$, we have in fact $|L(f)| \leq
%S_{p,d,\beta}(\psi)$. To see this, apply the preceding remarks to the
%averaged functional $\int_{k \in K_0} (k \cdot L) dk$.

\item \label{end}
Suppose that $E : C^{\infty}_{\omega}(\quotG) \rightarrow
C^{\infty}_{\omega}(\quotG)$
is a linear endomorphism satisfying $PS_{p,d,\beta}(Ef)
\leq  A \cdot PS_{p,d,\beta}(f)$, for some $A \in \mathbb{R}$.
Then also $S_{p,d,\beta}(Ef) \leq A S_{p,d,\beta}(f)$.
Indeed, $f \mapsto A^{-1} S_{p,d,\beta}(Ef)$
is a seminorm dominated by $PS_{p,d,\beta}$.
\item  \label{endvar}
We shall need a slight variant of (\ref{end}) in the case where
we are studying only the space of $f$ with some invariance property. 

Suppose that $E : C^{\infty}_{\omega}(\quotG) \rightarrow C^{\infty}_{\omega}(\quotG)$
is a linear endomorphism, $M$ is a finite set of finite places, and for each
$v \in M$ we are given an open compact $K_{1,v} \subset K_v$.
Suppose moreover
that $PS_{p,d,\beta}(Ef) \leq A \cdot PS_{p,d,\beta}(f)$ for some $A \in \mathbb{R}$
and for all $f$ which are $\prod_{v \in M} K_{1,v}$-fixed. 
Then, for all $f$ which are $\prod_{v \in M} K_{1,v}$-fixed, we have in fact
 $$S_{p,d,\beta} (Ef) \leq A \prod_{v \in M} [K_v: K_{1,v}]^{\beta} S_{p,d,\beta}(f).$$
Indeed, put $K_{1,M} = \prod_{v \in M} K_{1,v}$ and let $\Pi$ be the averaging operator
$\int_{k \in K_{1,M}} \pi(k) dk$, where $K_{1,M}$ is endowed with the Haar probability measure. 
Then apply (\ref{end}) above to the operator $f \mapsto E(\Pi f)$. 
\end{enumerate}
\end{rem}

\begin{lem} \label{lem:sobproduct}
Let $F_1 \in C^{\infty}_{\omega_1}(\quotG)$,
$F_2 \in C^{\infty}_{\omega_2}(\quotG)$.
Then $$S_{p,d,\beta}(F_1 F_2) \ll_{d} S_{2p,d,\beta}(F_1)
S_{2p,d,\beta}(F_2).$$
\end{lem}

Note that $F_1 F_2 \in C^{\infty}_{\omega_1 \omega_2}(\quotG)$.
\proof
Put $F = F_1 F_2$.
For any monomial $\mathcal{D}$ of degree $d$ in $\basis$, we can write
$\mathcal{D}(F_1 F_2) = \sum_{\alpha \in \mathcal{I}}( \mathcal{D}_{\alpha, 1} F_1)(
\mathcal{D}_{\alpha,2} F_2)$,
where $\alpha$ range over an index set $\mathcal{I}$ whose
size is bounded by a constant depending only on $d$,
and the $\mathcal{D}_{\alpha,\star}$ are certain monomials
in $\mathcal{B}$ satisfying
$\ord(\mathcal{D}_{\alpha,1}) + \ord(\mathcal{D}_{\alpha,2}) = d$.
It follows that
$$\|\mathcal{D} F\|_{L^p(\quotGad)}
\leq \sum_{\alpha \in \mathcal{I}} \left(\int_{\quotGad}
|\mathcal{D}_{\alpha,1}F_1|^p |\mathcal{D}_{\alpha,2}F_2|^p\right)^{1/p}.$$
Applying Cauchy-Schwarz, we conclude
\begin{equation} \begin{aligned}
\|\mathcal{D} F\|_{L^p(\quotGad)}
 \leq \sum_{\alpha \in \mathcal{I}}
\|\mathcal{D}_{\alpha,1}F_1\|_{L^{2p}(\quotGad)}
 \|\mathcal{D}_{\alpha,2}F_2\|_{L^{2p}(\quotGad)}\end{aligned}\end{equation}

Clearly, for each finite place $v$, we have $K_{v,F} \supset K_{v,F_1} \cap
K_{v,F_2}$;
in particular $[\Kmaxg: K_F] \leq [\Kmaxg:K_{F_1}] [\Kmaxg:K_{F_2}]$.
It follows that

\begin{multline}[\Kmaxg:K_F]^{\beta} \sum_{\mathcal{D}} \|\mathcal{D}
F\|_{L^p(\quotGad)}
\\ \ll \left([\Kmaxg:K_{F_1}]^{\beta} \sum_{\mathcal{D}}\|\mathcal{D}
F_1\|_{L^{2p}}\right)\left(
[\Kmaxg: K_{F_2}]^{\beta} \sum_{\mathcal{D}} \|\mathcal{D}
F_2\|_{L^{2p}}\right),\end{multline}
where the implicit constant depends only on $d$, and in all
three instances $\mathcal{D}$ varies over the set of monomials
in $\mathcal{B}$ of degree $\leq d$. 

That is to say, there is a constant
$C = C(d)$ such that $$PS_{p,d,\beta}(F_1 F_2) \leq
C \cdot PS_{2p,d,\beta}(F_1) PS_{2p,d,\beta}(F_2).$$
From (\ref{eq:concretenorm}) we deduce
$$S_{p,d,\beta}(F_1 F_2) \leq C\cdot S_{2p,d,\beta}(F_1)
S_{2p,d,\beta}(F_2),$$
as required.
\qed

We recall the definition of $\|g\|$ for $g\in \G(F_{\infty}), \G(\adele_{F})$
etc. from Sec. \ref{sec:adelicfunctionspace}. 

\begin{lem} \label{lem:sobgroup}
Let $F \in C^{\infty}(\quotG)$ and $g = (g_{\infty}, g_f) \in \G(\adele_F)$.
$$S_{p,d,\beta}(g \cdot F) \ll \|g_{\infty}\|^{d}
\|g_f\|^{\beta} S_{p,d,\beta}(F).$$
\end{lem}
\proof
Put $F' = (g_{\infty}, g_f) \cdot F$,
where $g_f = (g_v)_{v \, \mathrm{finite}}$.
For each finite place $v$, we note that
$K_{v,F'} \supseteq g_v K_{v,F} g_v^{-1} \cap K_{v,\G}$.
The index $[K_{v,\G}:K_{v,F'}]$ is therefore
bounded above by the number of cosets $x g_v K_{v,F}$ in
$K_{v,\G} g_v K_{v,F}$. Clearly this is bounded
above by the number of left $K_{v,F}$
cosets in $K_{v,\G} g_v K_{v,\G}$; but the number of such cosets
is precisely $\|g_v\| \cdot [K_{v,\G}:K_{v,F}]$.
It now follows easily from the definitions
that $PS_{d,\beta,f}(F') \ll \|g_{\infty}\|^d\|g_f\|^{\beta}
PS_{d,\beta,f}(F)$. Applying Rem. \ref{rem:def} to the endomorphism 
$F \mapsto (g_{\infty}, g_f) \cdot F$,  we obtain the claim.
\qed

The following crude Lemma is as much of interpolation as we need.
It will be applied, in practice, where $E$ is
a composite of a Hecke operator and a certain $L^2$-projection.

\begin{lem}
\label{lem:interp}
Let $E$ be a linear endomorphism of $C^{\infty}_{\omega}(\quotG)$
which commutes with $\G(F_{\infty}) \times \Kmaxg$.
Suppose there are real numbers $A,B > 0$ 
such that for any
$f \in C^{\infty}_{\omega}(\quotG)$, we have
$\|Ef||_{L^2} \leq A \|f\|_{L^2}, \|E f\|_{L^{\infty}} \leq B
\|f\|_{L^{\infty}}$. Then for $2 \leq p \leq \infty$:
%$$S_{2,0,0}(Ev) \leq A S_{,0,0}(v), \ \
%S_{\infty,0,0}(Ev) \ll B S_{\infty,0,\beta_0}(v)$$
%Then for any $2 \leq p \leq \infty$:
$$S_{p,d,\beta}(Ev) \leq A^{2/p} B^{1-\frac{2}{p}} S_{p,d,\beta}(v).$$

(We admit also $B = \infty$, in which case the $L^{\infty}$ hypothesis
should be seen as void, and the result becomes
$S_{2,d,\beta}(Ev) \leq A S_{2,d,\beta}(v)$.) 
\end{lem}
\proof
%For every subgroup $K_1 \subset \Kmaxg$,
%the endomorphism $E$ descends to
%an endomorphism of $V:=C^{\infty}_{\omega}(\quotG)^{K_1}$,
%the space of $K_1$-fixed vectors.
%The operator norm of $E$ on $V$,
%w.r.t the $L^2$ (resp. $L^{\infty}$ norms) is $\leq A$
%(resp. $\leq B$.)
By interpolation, the operator norm
of $E$ w.r.t. the $L^p$ norm on $C^{\infty}_{\omega}(\quotG)$ is
$\leq A^{2/p} B^{1-2/p}$. Moreover, the assumption on $E$ shows that
$K_{E f} \supset K_f$.

It follows that for $f \in C^{\infty}_{\omega}(\quotG)$ we have the inequality
$$PS_{p,d,\beta}(E f) \leq
A^{2/p} B^{1-2/p} PS_{p,d,\beta}(f).$$
Rem. \ref{rem:def} implies the conclusion.
\qed

%The following Lemma will be applied, in practice,
%where $\meas$ is the measure corresponding to a Hecke operator.
%\begin{lem} \label{lem:sobhecke}
%Let $F \in C^{\infty}(\quotG)$.
%Let $\meas$ be a probability measure on $\G(\adele_{F,f})$
%that satisfies $k \cdot \meas \cdot k^{-1} = \meas$,
%for all $k \in \Kmaxg$.
%Let $\|\meas\|_{\mathrm{Op},2}$ be the
%operator norm of convolution by $\meas$
%on $L^2_0(\quotG)$.
%Then, for $p \geq 2$,
%$$S_{p,d,\beta}(F \star \meas - \int_{\quot} F d\mu_{\quotG}) \ll
%\|\meas\|_{\mathrm{Op},2}^{2/p}
%S_{p,d,\beta}(F).$$
%\end{lem}
%\proof
%Consider
%the linear endomorphism $E$
%of $C^{\infty}(\quot)$ defined by
%$E (F) = F \star \meas - \int_{\quotG} F d\mu_{\quotG}$.
%
%Since (by assumption) the operation of right convolution by $\meas$
%commutes with $\Kmaxg$, it follows that
%$K_{E(F)} \supset K_F$.
%The operator norm of the linear endomorphism $E$,
%is  $ \leq 2$ w.r.t the $L^{\infty}$ norm,
%and $\|\sigma\|_{2,\mathrm{Op}}$ w.r.t the $L^2$ norm.
%By Riesz-Thorin interpolation,
%its operator norm w.r.t. the $L^p$ norm is $\leq
%\|\meas\|^{2/p}_{\mathrm{Op},2}$.
%An easy computation with the definition shows that $PS_{p,d,\beta}(E(F)) \leq
%\|\meas\|^{2/p}_{\mathrm{Op},2} PS_{p,d,\beta}(F)$.
%Now apply Rem. \ref{rem:def}. \qed
%

\subsubsection{Computing Sobolev norms in the Kirillov model.}
In the present section, let $v$ be an archimedean place of $F$. 

Let $\pi_v$ be a generic unitary irreducible representation of $\GL_2(F_v)$.
Recall that this means that $\pi_v$ is realized in a space of functions 
$\kirill$
(the Kirillov model, consisting of restrictions of functions in the Whittaker
model to the diagonal torus) on $F_v^{\times}$. 
Recall also the definition of the local conductor $\Cond_v(\pi_v)$ from 
Sec. \ref{iscond}. 

In this model, the diagonal torus acts by translation and upper triangular matrices act through multiplication by characters:
that is to say, for $f \in \kirill$,  $y_1, y_2 \in F_v^{\times}$, $z \in F_v$ we have the rules
\begin{equation} \label{eq:kirillov}
\pi(a(y_1)) f: y_2 \mapsto f(y_1 y_2), \pi(n(z)) f: y_2 \mapsto f(y_2) e_{F_v}(z y_2).\end{equation}
From these facts it is easy to verify that  the space of smooth vectors in $\pi_v$ contains all compactly supported smooth functions on $F_v^{\times}$.   
Moreover, \begin{equation} \label{eq:innerproduct} \|f\|_2^2 = \int_{F_v^{\times}} |f(y)|^2 d^{\times}y \end{equation} defines a $\GL_2(F_v)$-invariant inner product on $\kirill$. 

We will eventually have occasion to choose test vectors in $\pi_v$ in this model,
and wish to evaluate the ``Sobolev norms'' of the resulting vectors.

\begin{lem} \label{lem:sobnormk}
Suppose $F_v \cong \mathbb{R}$. 
 Let $f \in \kirill$ be $C^{\infty}$ and compactly supported. 
Then $$\sum_{\ord(\mathcal{D}) \leq k} 
\|\mathcal{D} f\|_2 \ll \Cond_v(\pi_v)^{2 k} \left( \sum_{j=0}^{2 k} \int_{\mathbb{R}^{\times}}
(|y|+|y|^{-1})^{2 k} \left| \frac{d^j f}{d^j y}\right|^2 d^{\times} y\right)^{1/2},$$
where the $\mathcal{D}$ sum ranges over all monomials in a fixed basis
for $\Lie(\GL_2(F_v))$ of degree $\leq k$. 

Suppose $F_v \cong \mathbb{C}$, and suppose $f \in \kirill$ is $C^{\infty}$ and compactly supported. 
Then
$$\sum_{\ord(\mathcal{D}) \leq k} 
\|\mathcal{D} f\|_2\ll \Cond_v(\pi_v)^{k} \left( \sum_{0 \leq i+j \leq 2k} \int_{\mathbb{C}^{\times}}
(|z|+|z|^{-1})^{2k} \left| \frac{\partial^{i+j} f}{\partial^{i} z \partial^{j} \bar{z}}\right|^2 d^{\times}z \right)^{1/2}$$
\end{lem}
\proof
We prove only the case with $F_v \cong \mathbb{R}$, the complex
case being similar.

Let $h,e$, $f$, $z$ be nonzero elements of the (real) Lie algebras of $\GL_2(\mathbb{R})$,
defined via $$h = \left( \begin{array}{cc} 1 & 0 \\ 0 & -1 \end{array} \right), e =\left( \begin{array}{cc}
0 & 1 \\0 & 0 \end{array} \right), f =\left( \begin{array}{cc}
0 & 0 \\ 1  & 0 \end{array} \right) , z= \left( \begin{array}{cc} 1 & 0 \\ 0 & 1 \end{array} \right)$$
These satisfy the usual commutation relations $[h,e] = 2e, [h,f] = 2f, [e,f] = h$. 
Let $\lambda$ be the scalar by which the Casimir operator $\frac{1}{2} h^2 + ef+ fe$
acts, and $\nu$ the scalar by which $z$ acts; then $1+|\lambda| +|\nu|^2 \ll \Cond_v(\pi_v)^2$. 

It is easy to see how $h,e$ act on $\kirill$: $h$ acts 
by a multiple of the differential operator $c_1 \nu + y \frac{d}{dy}$ and $e$ acts by multiplication by $c_3 y$,
for some constants $c_1, c_2, c_3$. 
The Casimir operator
$\frac{1}{2} h^2 + e f + f e =  \frac{1}{2}h^2 + 2 e f - h$ acts by the scalar
$\lambda$;  so it follows that for $v \in \kirill$ we have
$e f v = \frac{1}{2} (\lambda + h - h^2) v$. 
In particular, $e$ acts on any compactly supported function via the differential operator $c_1' y^{-1}+
c_2' \frac{d}{dy} + c_3' y \frac{d^2}{dy^2}$, for certain constants $c_1', c_2', c_3'$,
satisfying $|c_1'|, |c_2'|, |c_3'| \ll \Cond_v(\pi_v)^2$.  (In fact,  $|c_i'| \ll \Cond_v(\pi_v)^{3-i}$.)

Any monomial of degree $k$ in $h,e,f,z$ is therefore a sum of terms
$c_{\gamma \delta} y^{\gamma} \partial_y^{\delta}$,
where $|c_{\gamma \delta}| \ll \Cond_v(\pi_v)^{2k},  |\gamma| \leq k, \delta \leq 2k$. 
The claimed result follows in the case $F_v \cong \mathbb{R}$. 

 A similar
proof holds for $F_v  \cong \mathbb{C}$.  \qed

\subsection{Reduction theory.} \label{sec:redtheory}

Recall that $F_{\infty} := F \otimes_{\Q} \mathbb{R}$.
Let $K_{\infty}, K_v, \Kmax$ be as in Sec. \ref{sec:gl2pgl2}.
Then $K_{\infty} \times \Kmax$ is a maximal compact subgroup
of $\GL_2(\adele_F)$.
Given $g \in \GL_2(\adele_F)$ we may always
write $g =  \left(\begin{array}{cc} 1 & t \\ 0 & 1 \end{array}\right)\left(
        \begin{array}{cc}  x &0 \\ 0 & y \end{array}\right) k$,
with $t \in \adele_F, x,y \in \adele_F^{\times}, k
\in K_{\infty} \times \Kmax$.
We set $\height(g) = |x y^{-1}|_{\adele}$; this is well-defined,
although $x,y$ are not unique.

Then $\height$ descends to a function $B(F) \backslash \GL_2(\adele_F)
\rightarrow \mathbb{R}_{>0}$.  Explicitly,
\begin{equation}\label{eq:heightformula}\height \left(\begin{array}{cc} a & b \\ c & d \end{array}\right)
= \frac{|ad-bc|_{\adele}}{\prod_{v} \|(c_v, d_v)\|^2},\end{equation}
where one defines $\|(c_v, d_v)\| =\max(|c_v|_v,|d_v|_v)$
for $v$ finite, and \begin{equation} \label{eq:inftynormdef} \|(c_v, d_v)\| = \left( |c_v|_v^{2/\deg(v)} + |d_v|_v^{2/\deg(v)} \right)^{\deg(v)/2}\end{equation}
for $v$ infinite, where $\deg(v) = [F_v: \mathbb{R}]$. 

Define $\Siegel(T) \subset B(F) \backslash \GL_2(\adele_F)$
to be $\Siegel(T) := \{ g : \height(g) \geq T\}$.
Then, for all $T > 0$
the natural projection $\Pi: \Siegel(T) \rightarrow
\GL_2(F) \backslash \GL_2(\adele_F)$ has finite fibers;
for sufficiently large $T$, it is injective,
and for sufficiently small $T$ it is surjective.
This is the content of reduction theory for $\GL_2$.
As a consequence, the complement of $\Pi(\Siegel(T))$ has
compact closure, modulo the center, for each $T$. 

%Moreover, for all $T  >0$, there exists $\Omega := \Omega(T)$
%compact in $\GL_2(\adele_F)$ such that
%\begin{equation} \label{eq:compactcomplement}Z(\adele_F) \Omega \cup
%\Pi(\Siegel(T)) = \GL_2(F) \backslash
%\GL_2(\adele_F).\end{equation}

Fix $T_0$ such that $\Pi: \Siegel(T_0) \rightarrow \GL_2(F) \backslash
\GL_2(\adele_F)$ is injective.
%, and let $\Omega = \Omega(T_0)$ be as above.
Then we define a function $\height: \GL_2(F) \backslash \GL_2(\adele_F)
\rightarrow \mathbb{R} $ via
the rule
$$\height(g) = \begin{cases} \height(g'), \mbox{if } g = \Pi(g')
\mbox{ for some }g' \in \Siegel(T_0), \\
T_0, \mbox{ else.} \end{cases}$$

In fact, it is clear that $\height$ descends to a function $\quot = \PGL_2(F) \backslash
\PGL_2(\adele_F) \rightarrow \mathbb{R}_{\geq T_0}$.

\begin{lem} \label{lem:redtheoryone}
Let $U \subset \GL_2(F_{\infty})$ be compact
and $x \in \quotg$. The fibers of the map
$U \times \Kmax \rightarrow \quotg$ defined by
$(u,k) \mapsto xuk$ have size bounded by $O(\height(x))$,
where the implicit constant depends on $U$.
\end{lem}
\proof

Suppose $g \in \GL_2(\adele_F) $ is a lift of $x \in \quotg$. 
Consider the map $U \times \Kmax \rightarrow \quotg$
given by $(u,k) \mapsto guk$, as above.
Let $(u,k)$ belong to a fiber of maximal size. Call this size $M$.
Then
\begin{multline}M = \#\{\gamma \in \GL_2(F): g uk  = \gamma g u' k',  \exists
\, \,  u' \in U, k' \in \Kmax\}
\\ \leq \#\{\gamma: g u'' k'' = \gamma g, \exists \, \,  u'' \in U\cdot U^{-1}, k'' \in
\Kmax\}.\end{multline}

Set $V = U \cdot U^{-1}$, a compact subset of $\GL_2(F_{\infty})$.
The definition of $\Sieg(T)$ shows that 
there exists a constant $c < 1$, depending on $V$, such that
$\Siegel(T) \cdot V \cdot \Kmax \subset \Siegel(c T)$.
Choose $T$ so large that the projection $\Siegel(c T)
\rightarrow \quotg$ is injective.
It will suffice to show, for each $g \in \Siegel(T)$,
that
\begin{equation} \label{eq:card}\#\{\gamma \in \GL_2(F): \gamma g \in g V \Kmax\}
\ll \height(g).\end{equation}

Both $g$ and $g V K$ belong entirely to $\Siegel(cT)$. 
By the choice of $T$, $\gamma g \in g V \Kmax$ implies
$\gamma$ in $B(F)$.
Write $\gamma = a_{\gamma} n_{\gamma}$,
with $a_{\gamma} \in A(F)$ and $n_{\gamma} \in N(F)$;
also, write $g = n_g a_g k_g$ with $n_g \in N(\adele_F),
a_g \in A(\adele_F), k_g \in K_{\infty} \times \Kmax$.
We are free to adjust $g$ on the left by
an element of $N(F)$, since doing so will not 
affect the cardinality of the set $\{\gamma \in \GL_2(F): \gamma g \in g V \Kmax \}$. 
 We may thereby assume that $n_g$ lies in a fixed
compact subset of $N(\adele_F)$.
Thus we can write $g = a_g k_g'$,
where $k_g' := a_g^{-1} n_g a_g k_g$ lies
in a certain fixed compact subset $\Omega$ of $\GL_2(\adele_F)$.

Now, $\gamma g \in g V K$ implies
that $a_g^{-1} a_{\gamma} n_{\gamma} a_g \in
\Omega V \Kmax \Omega^{-1}$. Noting
that $a_g^{-1} a_{\gamma} n_{\gamma} a_g = a_{\gamma} 
a_g^{-1} n_{\gamma} a_g$, we deduce that
$a_{\gamma}$ lies in a fixed compact
subset of $\GL_2(\adele_F)$, depending only on $U$;
thus the number of possibilities
for $a_{\gamma}$ are $\ll_U 1$.
Moreover, it now follows that $a_g^{-1} n_{\gamma} a_{g}$
lies in a compact subset of $\adele_F$ depending only on $U$.

Thus, if we write $a_g =\left(\begin{array}{cc} x & 0 \\ 0 & y
\end{array}\right), n_{\gamma} = \left(\begin{array}{cc} 1 & \beta \\ 0 & 1
\end{array}\right)$, then $ \beta \in xy^{-1}  \Omega'$,
where $\Omega' \subset \adele_F$ is a compact subset
that depends only on $U$. It is easy
to see that the number of possibilities for $\beta$ is $\ll_U 1+|xy^{-1}|_{\adele_F}$.
But $|xy^{-1}|_{\adele_F} = \height(g)$, which is a function that is bounded away from zero, 
and we are done.
\qed

\begin{lem} \label{lem:redtheoryoneb}
Let notations be as in the previous Lemma \ref{lem:redtheoryone}. 
Consider the composite map $U \cdot \Kmax \stackrel{\Pi}{\rightarrow} \quotg \rightarrow \quot$. 
Each fiber of this map may be written as the union of at most $O(\height(x))$ sets
each of the form $y Z(\adele_F) \cap U \Kmax$, where $Z$ is the center of $\GL_2$
and $y \in \GL_2(\adele_F)$. 
\end{lem}
\proof
Let $\bar{x}$ be the image of $x$ in $\quot$. Let $u,u' \in U, k, k' \in \Kmax$. Suppose that $\bar{x} u k = \bar{x} u' k'$ in 
$\quot$. Then there is $z \in \adele_F^{\times}$ and $\gamma \in \GL_2(F)$ such that
\begin{equation} \label{whysolong} x u k = \gamma x u' k' a(z,z), 
\mbox{   equality in }\GL_2(\adele_F)\end{equation}
For fixed $u,k$ and $\gamma$, the set of $u'k'$ satisfying (\ref{whysolong})
is visibly the intersection of $U \Kmax$ with a fixed $Z(\adele_F)$-coset. 
This coset depends only on the class of $\gamma$ in $\PGL_2(F)$,
so it suffices to show that those $\gamma \in \GL_2(F)$ that occur in equalities
such as (\ref{whysolong}) for varying $u,k,u',k' $ represent at most $O(\height(x))$
distinct cosets $\gamma Z(F)$ in $\PGL_2(F)$. 

 Taking determinant followed by the norm
$\adele_F^{\times}/F^{\times} \rightarrow \mathbb{R}$, we conclude that 
$|z|_{\adele}$ belongs to a compact subset of $\mathbb{R}^{\times}$ that depends only on $U$.
The norm map $\adele_F^{\times}/F^{\times} \rightarrow \mathbb{R}^{\times}$
being proper, it follows that $z$ itself belongs to a compact subset $
\adele_F^{\times}/F^{\times}$ that depends only on $U$.  

In particular, there is a compact subset $\Omega \subset F_{\infty}^{\times}$,
depending only on $U$, and a finite subset $P \subset \adele_F^{\times}$, containing $1$ and
also depending only on $U$, such that $z \in F^{\times} \Omega . P . \prod_{v \, \mathrm{finite}}
\order_{F,v}^{\times}$. 
Let $\tilde{U} = U \cdot \{a(z_{\infty},z_{\infty}): z_{\infty} \in \Omega\}$.  Given a solution to (\ref{whysolong}), write $z = \delta z_{\infty} p o$, with $\delta \in F^{\times}, z_{\infty} \in \Omega,
p \in P, o \in \prod_{v \, \mathrm{finite}} \order_{F,v}^{\times}$. 
Then $$x u k = \gamma a(\delta,\delta)  x a(p,p) u' a(z_{\infty}, z_{\infty}) k' a(o,o),$$
in particular, taking $\tilde{u} = u' a(z_{\infty}, z_{\infty}) \in \tilde{U}, k'' = k' a(o,o) \in \Kmax$,
the image of $x a(p,p) \tilde{u} k''$ in $\quotG$ coincides with $x u k$.  So 
the number of possibilities for the $Z(F)$-coset of $\gamma$ is bounded
above by the fibers of the map $P \times \tilde{U} \times \Kmax \rightarrow \quotG$
given by $(p, \tilde{u}, k) \rightarrow x a(p,p) \tilde{u} k$. 
The result follows from Lem. \ref{lem:redtheoryone}.

\qed

We shall now need a quantitative version of certain statements in reduction
theory.
The subsequent Lemma is a fancier version of the following statement:
the number of $\gamma \in \mathrm{SL}(2, \mathbb{Z})$ 
that map a fixed $z \in \mathbb{H}$ to the Siegel set $$\{x+iy: 0 \leq x \leq 1, y \geq T\}$$
is $\ll 1+T^{-1}.$

\begin{lem} \label{lem:redtheorytwo}
Let $g \in \GL_2(\adele_F)$ and $Y > 0$ a positive real number. 
Then
\begin{equation} \label{eq:nottoo}
\#\{\gamma \in B(F) \backslash \GL_2(F):
\height(\gamma g) \geq Y\}  \ll_{\epsilon} 1+Y^{-1-\epsilon}.\end{equation}
Here the implicit constant is independent of $g$.
Moreover, suppose $g \in \Siegel(T)$ with $T \geq 1$. Then:
\begin{equation} \label{eq:bigone}
\mathrm{sup} \{\height(\gamma g): \gamma \notin B(F)\} \leq
T^{-1}.\end{equation}
\end{lem}

\proof
The proof of (\ref{eq:nottoo}) is not difficult,
generalizing in a straightforward way the proof with $F=\Q$.
However, it is somewhat notationally tedious; the (hypothetical) reader may wish 
to simply work out the proof for $F=\Q$, where it 
is equivalent to the following fact: the number of primitive vectors
in a unimodular sublattice of $\mathbb{R}^2$ that are contained in 
an $R$-ball is $\ll (1+R^2)$, uniformly in the lattice. 
(The result also be deduced if one admits some basic facts from the theory of Eisenstein
series over $F$, but we wish to rather deduce these basic facts from the present Lemma). 
We also remark that the entire content of (\ref{eq:nottoo}) lies in the uniformity in $g$.

Without loss of generality, we take $g \in 
%= \left( \begin{array}{cc} a &  b \\ c & d \end{array} \right) \in
 \Siegel(T_0)$,
where $T_0$ is sufficiently small that the map $\Siegel(T_0) \rightarrow \quotg$ is surjective. 
So  $g =  \left(\begin{array}{cc} 1 & t \\ 0 & 1 \end{array}\right)\left(
        \begin{array}{cc}  x &0 \\ 0 & y \end{array}\right) k$ with $|x y^{-1}|_{\adele} \geq T_0$. 
        Moreover, replacing $g$ by $ gz$, for any $z \in Z(\adele_F)$
         does not affect the problem,
        so we may take $y=1$.
Then, for $\gamma 
=  \left( \begin{array}{cc} a &  b \\ c & d \end{array} \right)  \in \GL_2(F)$, we have \begin{equation} \label{eq:height} \height(\gamma g) = \frac{|x|_{\adele}}{\prod_{v} \|(x_v c, ct_v + d)\|_v^2},\end{equation} 

%The class of $\gamma$
%in $B(F) \backslash \GL_2(F)$, and also the right hand side of (\ref{eq:height}), depend only on 
%the ratio $d/c$. Thus
%the left hand side of (\ref{eq:nottoo}) is bounded above by
%$1 + \#\{d \in F: \prod_{v} \|(1, d y_v) \|_v^2 \leq |xy|_{\adele} Y^{-1}$. 

The equivalence class of $\gamma$ in $B(F) \backslash \GL_2(F)$
depends only the pair $(c,d) \in F^2$, considered up to $F^{\times}$
equivalence (i.e., it depends only on $c/d \in F \cup \{\infty\}$.)
It suffices, then, to estimate the number
\begin{equation} \label{eq:no2} \#\{ [c:d] \in \mathbb{P}^1(F), \prod_{v} \|x_v c, ct_v +d \|_v^2  \leq Y^{-1} |x|_{\adele}\},\end{equation}

If $\Omega$ is any fixed compact subset of $\adele_F^{\times}$, 
then for $\omega \in \Omega$ we have $\prod_{v} \|x_v \omega c, ct_v+d\|_v \asymp_{\Omega} \prod_{v} \|x_v c, ct_v+d\|_v.$
Consider $\mathbb{R}_{>0}$ as embedded in $\adele_F^{\times}$ via
$\mathbb{R}_{>0} \hookrightarrow \adele_{\Q}^{\times} \hookrightarrow \adele_F^{\times}$. 
Then there is a compact subset $\Omega \in \adele_F^{\times}$ such that
$\adele_F^{\times} = F^{\times} \cdot \Omega \cdot \mathbb{R}_{>0}$. 

The size of (\ref{eq:no2}) is unaffected by the substitution $(x,t) \mapsto (x \tau, t \tau)$,
for any $\tau \in F^{\times}$. 
In view of the above remarks we may assume -- decreasing $Y$ by
a constant that depends only on $F$ -- that $x \in \mathbb{R}_{> 0}$.
Moreover, the size  of (\ref{eq:no2}) is also unaffected by the substitution 
$t \mapsto t + \tau$, for $\tau \in F$.
We may therefore assume that $|t|_v \leq 1$
for all finite places $v$.

Fix a set of representatives
$\mathfrak{J}_1, \dots, \mathfrak{J}_h$ for the class group of $\order_F$;
we will assume each $\mathfrak{J}_i$ is integral.
 For any $[c:d] \in \mathbb{P}^1(F)$, we may find a representative $(c,d)$
so that the ideal $c \order_F + d \order_F$ is one of the $\mathfrak{J}_i$;
moreover, replacing $(c,d) \in \mathfrak{J}_i^2$ by $(\epsilon c, \epsilon d)$
for $\epsilon \in \order_F^{\times}$ does not change the class $[c:d]$.

The restrictions on $x,t$ imply that $\|(x_v c, ct_v +d)\|_v = \|(c,d)\|_v$ 
for all finite $v$. 
Then $\prod_{v \, \mathrm{finite}}\|(c,d)\|_v = \Norm(\mathfrak{J}_i)^{-1}$, the inverse of the norm
of $\mathfrak{J}_i = c \order_F + d \order_F$. So it will suffice to bound,
for each $1 \leq i \leq h$, the quantity
$$\#\{(c,d) \in \mathfrak{J}_i^2 / \order_F^{\times}: c\order_F + d \order_F =\mathfrak{J}_i : 
\prod_{\infty|v} (|x_v|^2 |c|_v^2 + |ct_v+d|_v^2) \leq Y^{-1} |x|_{\adele} \Norm(\mathfrak{J})_i^{2}\}.$$

Since $\mathfrak{J}_i$ belongs to a finite set, the quantity $\Norm(\mathfrak{J})$ is bounded;
thus, decreasing $Y$ again as necessary, it suffices to estimate for each $1 \leq i \leq h$
$$\#\{ (c,d) \in \mathfrak{J}_i^2/\order_F^{\times}: 
c\order_F + d \order_F = \mathfrak{J}_i,
\prod_{\infty|v}\|(x_v c, ct_v +d)\|_v^2 \leq Y^{-1} |x|_{\adele}\}$$

There is only one term corresponding to $c=0$. 
Otherwise, $(c)$ is a principal ideal divisible by $\mathfrak{J}_i$;
let $\mathcal{P}$ be the set of integral principal ideals. 
Then the size of the set above is precisely
\begin{equation} \label{eq:notaword} \sum_{ (c) \in \mathcal{P}} \#\{d \in \mathfrak{J}_i : (c) + d \order_F = \mathfrak{J}_i,
\prod_{\infty|v} \|(x_v c, ct_v +d)\|_v^2 \leq Y^{-1} |x|_{\adele}\}\end{equation}
We note that the size of the inner set is independent of the choice
of generator for the principal ideal $(c)$.  Moreover, the inequality
of (\ref{eq:notaword})
implies that the norm $\Norm((c))$ of the principal ideal $(c)$ satisfies $\Norm((c))^2 \leq Y^{-1} |x|_{\adele}^{-1}$.

Let us estimate the number of $d$ that can correspond to a fixed principal ideal $(c)$
in (\ref{eq:notaword}). 
Recall that $|x|_{\adele} \geq T_0$ and that $x_{\adele}$
is in the image of the embedding $\mathbb{R}_{>0} \hookrightarrow \adele_Q^{\times}
\hookrightarrow \adele_F^{\times}$. In particular, $|x|_v$
is bounded below at each infinite place. Moreover, 
since $\order_F^{\times}$ is a cocompact subgroup
of the elements of $F_{\infty}^{\times}$ with norm $1$, we can choose
a representative for the principal ideal $(c)$ so the same is true of $|c|_v$. 
Note that (cf. \ref{eq:inftynormdef}) that $\|(x_v c, ct_v +d)\|_v \asymp
(|x_v c|_v + |c t_v+d_v|_v)$. 

So in fact, again decreasing $Y$ as necessary, it will suffice to estimate
\begin{equation} \label{eq:aintnojoke}
\sum_{(c) \in \mathcal{P}: \Norm(c) \leq Y^{-1/2} |x|_{\adele}^{-1/2}}
\#\{ d \in \mathfrak{J}_i : \prod_{\infty|v} (1 + |ct+d|_v)^2 \leq Y^{-1} |x|_{\adele} \}\end{equation}

To estimate the right-hand side, first observe that if $\{M_v\}_{\infty|v}$ is any set of positive real numbers indexed
by the infinite places of $F$, then $\#\{d \in \mathfrak{J}_i: |ct+d|_v \leq M_v
\mbox{ for } \infty|v\} \ll \prod_{\infty|v} (1+M_v)$.  Indeed, by subtraction,
it will suffice to estimate $\#\{d \in \mathfrak{J}_i: |d|_v \leq 2 M_v \mbox{ for } \infty|v\}$;
this amounts to counting points in the lattice $\mathfrak{J}_i \subset F_{\infty}$ in a region that is
the product of a box and a disc; the result is then clear. 

Next, if $T \geq 1$, then the subset $\{(y_1, \dots, y_d): \prod_{i} (1+y_i) \leq T\}$
in $\mathbb{R}_{>0}^{d}$ is contained in the union of $O_{\epsilon}(T^{\epsilon})$
boxes $\{(y_1, \dots, y_d): y_i \leq M_i\}$, where $\prod_{i} (1+M_i) \ll T$. 
We may assume $Y^{-1} |x|_{\adele} \geq 1$, else (\ref{eq:aintnojoke})
has no solutions. 
We conclude that the number of $d$ attached to each principal ideal $(c)$
in (\ref{eq:aintnojoke}) is $\ll_{\epsilon} (Y^{-1/2} |x|_{\adele}^{1/2})^{1+\epsilon}$.

The number of possibilities for $(c)$ is bounded
by the number of integral ideals with norm $\leq Y^{-1/2} |x|_{\adele}^{-1/2}$,
which is $\ll_{\epsilon}( Y^{-1/2} |x|_{\adele}^{-1/2})^{1+\epsilon}$. Finally
there is  one class with $c=0$. We conclude
that the number of pairs $(c,d)$ up to equivalence
is $\ll Y^{-1-\epsilon}+1$. This proves (\ref{eq:nottoo}).

As for (\ref{eq:bigone}), suppose $g \in\Siegel(T)$,
so we may write $g =
\left(\begin{array}{cc} x & z \\ 0 & y \end{array}\right) k$
with $k \in K_{\infty} \times \Kmax$, and $|xy^{-1}|_{\adele} \geq T$.
Suppose $\gamma = \left(\begin{array}{cc} \alpha & \beta \\ \alpha'
& \beta'
\end{array}\right)$. If $\gamma \notin B(F)$, then $\alpha' \neq 0$.
In that case,
following the notation of (\ref{eq:heightformula}), we have:
$$\prod_{v} \|(\alpha'_v ,\beta'_v) g_v \|
\geq  \prod_{v} |\alpha'_v x_v|_v = |x|_{\adele}.$$
and therefore, by (\ref{eq:heightformula}),
$\height(\gamma g) \leq |\det(g) x^{-2}|_{\adele}
\leq T^{-1}$.  \qed

\section{Background on quantitative equidistribution results.}
The aim of this section is to quantify various standard
equidistribution results (equidistribution of long horocycles,
Hecke points, etc.), using the adelic Sobolev norms.  As such
neither the results nor the methods are new; we just collect together
those results we need and provide brief proofs. 

As regards the origin of the ideas used here, we have drawn in particular from the work
of Clozel-Ullmo, Linnik, Oh, Margulis, Ratner and Sarnak. 

\subsection{Decay of matrix coefficients.}
\label{subsec:decay}
\subsubsection{Local setting.}
Our fundamental tool in establishing all these results
is the spectral gap, i.e., quantitative mixing
properties of real and $p$-adic flow. As such,
we begin by recalling the basic
relevant bound on matrix coefficients. 

Let $0 \leq \alpha \leq 1/2$. 
Let $v$ be a place of $F$, and suppose that
$(V, \pi)$ is a unitary representation of $\GL_2(F_v)$
which does not contain, in its spectral decomposition,
any complementary series with parameter
$\geq \alpha$. (More formally: $V$ does not {\em weakly contain}
such a representation). Thus $\alpha = 0$ corresponds to $V$ being tempered,
and $\alpha = 1/2$ corresponds to $V$ having
no almost invariant vectors.

Then for $w_1, w_2$ any two
$K_v$-finite elements of $V$, satisfying $\langle w_1, w_1 \rangle= \langle w_2, w_2 \rangle = 1$, and any $x \in F_v$
we have the bound on matrix coefficients given by
\begin{equation}\label{eq:mcbound}\langle \pi(a(x) w_1, w_2 \rangle \ll_{\epsilon,F} \dim(K_v w_1)^{1/2}
\dim(K_v w_2)^{1/2} (1+|x|_v)^{\alpha - 1/2 + \epsilon}.\end{equation}

The implicit constant of (\ref{eq:mcbound}) depends only on $\epsilon$.
Since we do not know of an available reference, we briefly sketch
an argument for (\ref{eq:mcbound}). In the case where $\alpha = 0$,
i.e. $V$ is {\em tempered}, then (\ref{eq:mcbound}) is
proven in \cite{CHH}. In the general case, let $(\sigma_{1/2-\alpha},W)$
be the complementary series with parameter $1/2-\alpha$;
let $v^{0} \in W$ be a unitary spherical vector. 
Then the representation $V \otimes W$ is tempered. Indeed
it suffices -- again by \cite{CHH} -- to verify that a dense set of matrix coefficients are in $L^{2+\epsilon}$,
which follows by direct computation.
Now one may estimate the matrix coefficient
$\langle a(x) w_1 \otimes v^{0} , w_2 \otimes v^{0} \rangle$ by appealing
again to \cite{CHH}.  On the other hand
$\langle a(x) w_1 \otimes v^{0}, w_2 \otimes v^{0}  \rangle
= \langle a(x) w_1, w_2 \rangle \langle a(x) v^{0}, v^{0}\rangle$,
and an easy computation shows that $\langle a(x) v^0, v^0 \rangle
\gg_{\epsilon} (1+|x|_v)^{-\alpha -\epsilon}$. Thus (\ref{eq:mcbound}) follows. 
(This argument is a variant of an argument that appears at the end of \cite{CHH}.)

Let us record a useful further variant. 
Suppose $v$ is finite. Let $K_1, K_2 \subset K_v$
be subgroups and let $\meas$ be the $(K_1, K_2)$-bi-invariant
probability measure supported on $K_1 a(x) K_2$.  

Then 
\begin{equation} \label{eq:measopnormbound} \|v \star \meas\|_2 \ll
 [K_v:K_1]^{1/2} [K_v:K_2]^{1/2} (1+|x|_v)^{\alpha-1/2+\epsilon} \|v\|_2.\end{equation}
 Indeed, for $i=1,2$ let $\Pi_{K_i}$ be the projection operator $w \mapsto \int_{K_i} k w$
on $V$, where $K_i$ is endowed with the Haar probability measure.  Then:
\begin{multline}  \|v \star \meas\|_{2}
= \sup_{w \in V} \frac{\langle v \star \meas, w \rangle}{\|w\|_2}
\\ = \sup_{w \in V} \frac{\langle a(x) \Pi_{K_1} v , \Pi_{K_2} w\rangle}{\|w\|_2}
\leq [K_v:K_1]^{1/2} [K_v:K_2]^{1/2} (1+|x|_v)^{\alpha-1/2+\epsilon} \|v\|_2 \end{multline}

\subsubsection{Variant for $\Gamma \backslash \SL_2(\mathbb{R})$.} \label{subsec:sl2mc}
Let $0 \leq \alpha \leq 1/2$, suppose $G = \SL_2(\mathbb{R})$, and let $V$ be a unitary representation 
of $G$ such that $V$ does not weakly contain any complementary series
with parameter $\geq \alpha$.  The normalization is again so that $\alpha =0$ corresponds
to tempered and $\alpha = 1/2$ corresponds to $V$ not having almost invariant vectors. 

Then one has the following
variant of (\ref{eq:mcbound}), proved by the same method:
\begin{equation}\label{eq:mcbound2}\langle \pi(\left(\begin{array}{cc} y^{1/2} & 0 \\ 0 & y^{-1/2} \end{array} \right) w_1, w_2 \rangle \ll_{\epsilon} \dim(\mathrm{SO}(2) \cdot w_1)^{1/2}
\dim(\mathrm{SO}(2) \cdot w_2)^{1/2} (1+|y|)^{\alpha - 1/2 + \epsilon}.\end{equation}

It is convenient to extend the validity of (\ref{eq:mcbound2}) beyond the $K$-finite
space by replacing $\dim(\mathrm{SO}(2) w_i)$ by appropriate Sobolev norms.  We confine
ourselves to the case of main interest, where $V$ is the orthogonal
complements of the constants in $L^2(\Gamma \backslash \SL_2(\mathbb{R}))$, 
where $\Gamma$ is a lattice in $\SL_2(\mathbb{R})$.  The estimates we are about to describe
are, again, not new; estimates for effective mixing of geodesic and horocycle flows
in this setting are contained in \cite{Ratner3}. 

For our purposes it would be optimal to use fractional Sobolev norms; since we have not defined these, we
shall use a rather crude form of interpolation instead. 

Thus let $f_1,f_2 \in C^{\infty}(\Gamma \backslash \SL_2(\mathbb{R}))$. 
One expands both $f_1$ and $f_2$ into a sum of $\mathrm{SO}(2)$-types
and applies (\ref{eq:mcbound2}). 
Indeed, write for $i \in \{1,2\}$ an expansion $f_i = \sum_{n=-\infty}^{\infty} f_i^{(n)}$, where $f_i^{(n)}$
transforms under the character $\left( \begin{array}{cc} \cos(\theta) & \sin(\theta) \\ -\sin(\theta) & \cos(\theta) \end{array} \right) \mapsto e^{i n \theta}$. 
Expanding:
\begin{multline} \langle \left(\begin{array} {cc} y^{1/2} & 0 \\ 0 & y^{-1/2} \end{array}\right) f_1, f_2 \rangle
= \sum_{n,m \in \mathbb{Z}} \langle \left(\begin{array} {cc} y^{1/2} & 0 \\ 0 & y^{-1/2} \end{array}\right) f_1^{(n)}, f_2^{(m)} 
\rangle \\
 \ll_{\epsilon} (1+|y|)^{\alpha - 1/2 +\epsilon} \sum_{n,m}
\|f_1^{(n)}\|_2 \|f^{(m)} \|_2
= (1+|y|)^{\alpha - 1/2+\epsilon} \left( \sum_{n} \|f_1^{(n)}\|_2 \right)
\left( \sum_m \|f_2^{(m)} \| _2\right) \end{multline}
Our definitions of the Sobolev norms (Sec. \ref{sec:sobreal}) are so that
$S_{2,1}(f_1)^2 \gg \sum_{n} (1+|n|)^2 \|f_1^{(n)}\|_2^2$, and similarly
for $f_2$. On the other hand, it is an elementary estimate that 
$$\left( \sum_{n} \|f_1^{(n)} \|_2 \right)^2 \ll_{\epsilon}
\left( \sum_n \|f_1^{(n)}\|_2^2 (1+|n|)^{2}\right)^{1/2+\epsilon} \left( \sum_{n} \|f_1^{(n)}\|_2^2\right)^{1/2-\epsilon}$$

It follows from this that for any $k,k' \in \mathrm{SO}(2)$ we have the matrix coefficient bound:
\begin{multline} \label{eq:sl2mc}|\langle k \left(\begin{array} {cc} y^{1/2} & 0 \\ 0 & y^{-1/2} \end{array}\right) k' f_1, f_2 \rangle| \\
 \ll (1+|y|)^{\alpha-1/2 + \epsilon} (S_{2,1}(f_1) S_{2,1}(f_2))^{1/2+\epsilon}
\|f_1\|^{1/2-\epsilon} \|f_2\|^{1/2-\epsilon},\end{multline}
at least for $f_1,f_2$ which are $\mathrm{SO}(2)$-finite. But the general case of smooth
$f_1, f_2$ follows from density. 

Note that in (\ref{eq:sl2mc}) that the factor $\|f_1\|^{1/2-\epsilon} S_{2,1}(f_1)^{1/2+\epsilon}$
is a crude substitute for the fractional ($1/2+\epsilon$-) Sobolev norm of $f_1$.

\subsection{Pointwise bounds.}
In this section, we make free use of the adelic Sobolev norms introduced in Sec. \ref{subsec:adelicsobolev}. 
We recall the definition $S_{p,d} :=S_{p,d,1/p}$. We also
recall that in statements of the form $|L(f)| \ll S_{p,d}(f)$, for certain
linear functionals $L$, we shall allow the implicit constant of $\ll$
to depend on $p$ and $d$ without explicit mention. 
\begin{lem} \label{lem:adelicsobolev}
Let $f \in C^{\infty}_{\omega}(\quotg)$ and let $x \in \quotg$.
Then, for any $p \geq 2$ and $d \gg 1$,
\begin{equation} \label{eq:pointsob}|f(x)| \ll  \height(x)^{1/p}
S_{p,d}(f)\end{equation}

Moreover, if $F \in C^{\infty}(\quot \times \quot)$,
and $p> 2, d \gg 1$,
\begin{equation}\int_{\quot} |F(x,x)| dx \ll  S_{p,d}(F)
  \label{eq:diagsob} \end{equation}
\end{lem}

\proof
As in (\ref{eq:kpsidef}), set $K_f = \prod_{v \, \mathrm{finite}} K_{v,f}$, where $K_{v,f}$ is the stabilizer of $f$ in $K_v$. 
Fix an open neighbourhood of the identity $U \subset \GL_2(F_{\infty})$. Consider
the map $\Pi: U \cdot K_{f} \rightarrow \quot$ defined by
$(u,k) \mapsto xuk$.
By Lem. \ref{lem:redtheoryoneb},
the fibers are unions of at most $O(\height(x))$ sets, each of the form $y Z(\adele_F) \cap
U K_f$. 
Moreover, for any $y \in U \cdot K_f$, the measure
of $\{z \in \adele_F^{\times}:  y a(z,z) \in U \cdot K_f\}$ is bounded
above by a constant depending only on $U$. Indeed, the set
of such $z$ is contained in a fixed compact subset of $\adele_F^{\times}$
that depends only on $U$.

Equip $U \cdot K_f$
with the restriction of Haar measure from $\GL_2(\adele_F)$.
From the preceding paragraph, one easily deduces that the push-forward of this measure to $\quot$, 
 under $(u,k) \mapsto x uk$,
is bounded above by $C \cdot \height(x)$ times the measure on $\quot$,
where the constant $C$ depends only on $U$.

Then:
\begin{multline} \label{eq:homer}
\int_{u \in U} |f(xu)|^p =
\vol(K_f)^{-1} \int_{u \in U, k \in K_f} |f(xuk)|^p du dk 
\\ \ll
\height(x) [\Kmax:K_{f}] \int_{\quot} |f(x)|^p d\mu_{\quot}(x)
\end{multline}
(\ref{eq:homer}) holds with $f$ replaced
by $\mathcal{D} f$, for $\mathcal{D}$ any fixed
monomial in $\Lie(\GL_2(F_{\infty}))$.
The standard Sobolev estimate, applied
to the function $u \mapsto f(xu)$ on the real
manifold $U$, implies that $|f(x)| \ll \height(x)^{1/p} PS_{p,d,1/p}(f)$ for sufficiently large $d$. 
 (Indeed, it suffices
to take any $d > \mathrm{dim}(U)/2 = 2[F:\Q]$.) Then Remark \ref{rem:def}, (\ref{remdef2}) implies the conclusion. 

As for the second conclusion, we proceed
in a similar fashion as above (with $\quot$
replaced by $\quot \times \quot$) to obtain the 
estimate $|F(x,y)| \ll \height(x)^{1/p} \height(y)^{1/p} S_{p,d}(F)$. 
It is easy to see that $\int_{\quot} \height(x)^{2/p} dx < \infty$
for $p >2$, and the conclusion follows.  \qed

The next lemma quantifies the rapid decay of a cuspidal function,
or more generally a truncated automorphic function, in the cusp. 
Recall that for $T_0 > 0$ we have defined the Siegel domain $\Sieg(T_0)$
in Section \ref{sec:redtheory}. 

\begin{lem} \label{lem:constterm}
Let $f \in C^{\infty}_{\omega}(\quotg)$.
Put $f^{N}(g) = \int_{N(F) \backslash N(\adele_F)}f(ng) dn$,
where the measure on $N(F) \backslash N(\adele_F)$ is the $N(\adele_F)$-invariant probability measure. 
Then for $x \in \Sieg(T_0)$, $p \geq 2, k \geq 0$ and $d \gg 1$, 
\begin{equation} \label{eq:constterm}|f(x)- f^{N}(x)| \ll_{T_0}  \height(x)^{1/p-k}
S_{p,d,1/p + k}(f)\end{equation}
\end{lem}
\proof
We may assume that $x \in N(\adele_F) A(\mathbb{R}) \Omega (K_{\infty} \times \Kmax)$,
for some fixed compact set $\Omega \subset A(\adele_{F})$.
Here $A(\R)$ is regarded as a subset of $A(F_{\infty})$ via the natural
inclusion $\R \hookrightarrow F_{\infty}$. 
 Write accordingly
$x = n a \omega k$, where $\omega \in \Omega$.

Consider the function on $F_{\infty}$ defined by $g(t)=  f(n(t) x) - f^N(x)$.
It is invariant by the lattice $\Lambda = \{t \in \order_F: n(t) \in \GL_2(F_{\infty}) 
\omega K_f \omega^{-1} \}$, , where $K_f$ is again as in $(\ref{eq:kpsidef})$.  One sees that, since $\omega$ belongs to the fixed compact $\Omega$,
the covolume bound
$\vol(F_{\infty}/\Lambda) \ll [\Kmax:K_f]$.  Moreover, since $\Lambda$ may be
regarded as a fractional
ideal of $\order_F$, 
 the homothety class of $\Lambda$ lies in a fixed compact
set in the space of homothety classes of lattices in $F_{\infty}$. 
Also,  $g(t)$ defines a function on $F_{\infty}/\Lambda$, with integral $0$.

Suppose now that $G$ is a smooth function on $\R^d/L$, for some $d >1$ and some lattice $L \subset \R^d$,
with integral $0$. Let $\|G\|_{(i)} = \sup_{z \in \R^d/L} |\mathcal{D} G|$, where 
$\mathcal{D}$ varies over all monomials in $\partial_1, \dots, \partial_d$ of exact order $i$. 
Then an elementary argument shows that $\|G\|_{(0)} \ll \vol(\R^d/L)^{i/d}
\|G\|_{(i)}$, and the implicit constant may be taken to vary continuously with the
homothety class of $L$. 

Apply this lemma to the function $g$ on $F_{\infty}/\Lambda$, with $i = k [F:\Q]$
for some $k \geq 1$. The norm $\|g\|_{(k [F:\Q])}$, in the sense of the above paragraph, 
is bounded, by Lem. \ref{lem:adelicsobolev} and an elementary computation, by
$\height(x)^{-k} (\height(x)^{1/p} PS_{p,d'}(f))$, for some $d' \gg 1$.  It follows 
that $$\sup|g(t)| \ll \height(x)^{1/p-k} PS_{p,d'}(f) [\Kmax:K_f]^{k} = \height(x)^{1/p-k}
PS_{p,d',1/p+k}(f).$$
Applying Rem. \ref{rem:def}, (\ref{remdef2}), we conclude $|f(x) - f^N(x)| \ll \height(x)^{1/p - k} S_{p,d',1/p+k} (f)$. 
\qed

\begin{lem} \label{lem:cuspfunctions}
Suppose $f \in C^{\infty}_{\omega}(\quotg)$ is cuspidal. 
Then $S_{\infty, d, \beta}(f) \ll S_{2,d',\beta+3/2}(f)$, for
sufficiently large $d'$.
\end{lem}
\proof
By Lem. \ref{lem:constterm} for $f$ cuspidal, applied with $p=2, k=1$, we see that $|f(x)| \ll S_{2,d,3/2}(f)$ for $d \gg 1$.  Applying
this inequality to $\mathcal{D} f$, for $\mathcal{D}$ in the
universal enveloping algebra of $\GL_2(F_{\infty})$, we see
that $PS_{\infty,d,\beta}(f) \ll PS_{2,d', \beta+3/2}(f)$, for
$d'$ sufficiently large.
This equality holds for cuspidal $f$. 

 Let $\Pi$ be the $L^2$-orthogonal
projection onto the space of cuspidal functions;
then $\Pi$ commutes with $\GL_2(\adele_F)$, and it follows 
that $$ PS_{\infty,d,\beta}(\Pi f) \ll PS_{2,d',\beta+3/2}(\Pi f) \leq PS_{2,d',\beta+3/2}(f)$$
for arbitrary $f \in C^{\infty}_{\omega}(\quotg)$. 
Now Rem. \ref{rem:def}, (\ref{end}) (or, more precisely, a trivial modification thereof) 
 implies the conclusion. 
\qed

\subsection{Equidistribution of long horocycles and closed horospheres.}
\label{subsec:longhorohoro}
Let $G$ be a semisimple group, $\Gamma \subset G$ a lattice, $U$
a unipotent subgroup of $G$. It is well-known that one can prove,
in a quantitative fashion, the equidistribution of $U$-orbits on $\Gamma \backslash G$
if $U$ is a {\em horospherical} subgroup, i.e. the unipotent radical of a proper
parabolic subgroup. We shall quantify two instances of this that
will be of interest to us.  

We emphasize that neither the results nor the techniques of this section are new; we have included
proofs {\em only} to keep the present paper as self-contained as possible. 

Effective estimates for equidistribution of long horocycles on quotients
of $\SL_2(\R)$ are already implicit in the work of Ratner \cite{Ratner1} and \cite{Ratner2}, where
the effective mixing of the horocycle flow is used. 
We will proceed in a closely related fashion, using the mixing property of the
Cartan action; again, this is definitely not new and appears already, although in a different context,
in the doctoral thesis of Margulis (reprinted in \cite{Margulisthesis}).  

\subsubsection{Equidistribution of long horocycles in hyperbolic $2$-space.} \label{subsec:weirdnotn}
Let $\Gamma \subset \SL(2, \R)$ be a lattice 
such that $L^2(\Gamma \backslash \SL(2,\R))$ does not contain
any complementary series representation with parameter $\geq \alpha$,
for any $0 \leq \alpha < 1/2$. (That is: $\alpha \in [0,1/2)$ is such that all nonzero eigenvalues
of the hyperbolic Laplacian $-y^2 (\partial_{xx} + \partial_{yy})$ on $\Gamma \backslash \mathbb{H}^2$
are bounded below by $1/4-\alpha^2$).

We define $n, a, \bar{n}$ as in \eqref{eqn:nan}.  

%\begin{lem}
%Suppose $\left(\begin{array}{cc} 1 & 1 \\ 0 & 1\end{array}\right) \in \Gamma$.
%Let $f \in C^{\infty}(\Gamma \backslash \SL_2(\R))$. Then
%\begin{equation}
%\left|\int_{t \in \R/N \Z} f(a(N^{-1/2}) n(t) dt - \int_{\Gamma \backslash \SL_2(\mathbb{R})}
%f(g) dg\right| \ll N^{-1/2} S_{2,1+\epsilon}(f)
%\end{equation}  
%\end{lem}
%\proof
%Here it is convenient to use explicit spectral expansion, which gives somewhat better
%estimates. For simplicity, assume there is only one cusp for the action of $\Gamma$ on $\SL_2(\R)$;
%the general case differs only in adding an indexing of the cusps. 
%For $s \in \mathbb{C}$, let $\pi_s$ be the unitary induction to $\mathrm{SL}_2(\mathbb{R})$
%of the character $a(y) \mapsto |y|^s$. One can realize $\pi_s$ as a meromorphically varying
%family of representations on the fixed Hilbert space $L^2(\mathrm{SO}_2)$; the theory
%of Eisenstein series supplies a meromorphic family of intertwining operators
%$I(s): \pi_s \rightarrow C^{\infty})\Gamma \backslash \mathrm{SL}_2(\mathbb{R})$. 
%\qed

The following Lemma quantifies the equidistribution of long horocycles. 
Results of this type are already implicit in \cite{Ratner1} and \cite{Ratner2}. 
 This problem is analyzed in much more detail than we go into, in \cite{Strom} and \cite{FlamForni}. 
 
\begin{lem} \label{lem:longhoroequi}
Assume $\Gamma$ is cocompact, and let $x_0 \in \Gamma \backslash \SL_2(\R)$.
% Let
%$g$ be a smooth function of compact support on the real line satisfying
%$\int_{-\infty}^{\infty}g(t) = 1$. 
\begin{equation} \label{eq:wanted}\left|\frac{1}{T}  \int_{t =0}^{T} f(x_0 n(t)) dt - \int_{\Gamma \backslash \SL_2(\R)}
f(g) dg \right| \ll_{\epsilon} T^{\frac{\alpha-1/2}{2}+\epsilon} S_{\infty,1}(f). \end{equation}
\end{lem}

\proof
The idea (which is certainly not new -- cf. remarks at start of Section \ref{subsec:longhorohoro})  is that, upon flowing a small ball
in $\Gamma \backslash G$ for a long time by the geodesic flow, it turns into a narrow neighbourhood
of a long horocycle. One thereby can deduce the equidistribution of the long horocycle
from the mixing properties of the geodesic flow.

Let $N, A, \bar{N}$ be the images of $n, a, \bar{n}$ respectively. 
Let $g_1$ be a smooth function of compact support on the real line,
with integral $\int_{-\infty}^{\infty}g_1(x) dx = 1$. 
It will remain fixed for all time throughout our arguments.
Fix $1 > \delta > 0$ and let $g_{\delta}: \mathbb{R} \rightarrow \mathbb{R}$
be the convolution of the characteristic function of $[0,1]$ 
with $g_1(x/\delta) \delta^{-1}$; that is to say
$$g_{\delta}(x) =\delta^{-1}  \int_{t=0}^{1} g_1(\frac{x-t}{\delta}) dt.$$
 Then $g_{\delta}$ is a smooth
function of integral $1$, which is supported in a small interval
around $[0,1]$. 
 
%Define a function $g_{\delta}: A. \bar{N} \rightarrow \C$ via the rule 
%$g_2(a(y) \bar{n}(z)) = g_1(y/\delta) g_1(z)$.
  Define
a probability measure $\mu_{\delta}$ on $\Gamma \backslash \SL_2(\mathbb{R})$
via the rule
$$\mu_{\delta} (f) = \delta^{-1}\int_{x,y,z \in \mathbb{R}}
 f(x_0 n(x) a( e^y) \bar{n}(z)) g_{\delta}(x) g_1(y/\delta) g_1(z) 
dxdydz.$$ 
In words, $\mu_{\delta}$ is a measure supported on a small box around $x_0$; this box
has width $O(1)$ in the $N$ and $\bar{N}$ directions, and $O(\delta)$ in the $A$ direction.
When we flow this by $A$, it will become a measure supported along a box that
closely approximates an $N$-orbit.

We observe that
$$\mu_{\delta}(a(T^{-1}) f) = \frac{1}{\delta } \int_{x,y,z} f(x_0 a(T)^{-1} n(x) a(e^y) \bar{n}(z))
g_{\delta}(x/T) g_1(y/\delta) g_1(T z) dx dy dz.$$ 
On the other hand, for any fixed $x_1 \in 
\Gamma \backslash \mathrm{SL}_2(\mathbb{R})$, we note that
\begin{equation} \label{eq:bananaman}\left|\delta^{-1}T \int f(x_1 a(e^y) \bar{n}(z)) g_1(y/\delta) g_1(Tz)  dy dz - f(x_1)\right|
\ll \max(T^{-1}, \delta) S_{\infty,1}(f).\end{equation}
Indeed, (\ref{eq:bananaman}) merely quantifies the fact that the right-hand side integral
is against a probability measure supported in a very small ball (of size $\min(\delta, T^{-1})$) around $x_1$.
Since $\Gamma \backslash \SL_2(\R)$ is assumed compact, the implicit
constant of (\ref{eq:bananaman}) may be taken independent of $x_1$. 

Consequently,
\begin{equation} \label{eqn:ga} \left| \frac{1}{T} \int_{t \in \mathbb{R}} g_{\delta}(t/T) f(x_0 a(T)^{-1} n(t)) dt 
 - \mu_{\delta}(a(T)^{-1}  f) \right| \ll \max(T^{-1}, \delta)
 S_{\infty,1}(f).\end{equation}
 
 On the other hand, the measure $\mu_{\delta}$
 has a continuous distribution function $h_{\delta}$, i.e.
 $\mu_{\delta}(f) = \int_{\Gamma \backslash \SL_2(\R)} f(g) \cdot h_{\delta}(g) dg$,  and 
 $\mu_{\delta}(a(T)^{-1} f)$ may be estimated using
 (\ref{eq:sl2mc}), i.e. the decay of matrix coefficients.
 
  A routine computation   
  shows that
 $\|h_{\delta}\|_{L^2} \ll \delta^{-1/2}$ and $S_{2,1}(h_{\delta}) \ll \delta^{-3/2}$; on account
 of the cocompactness of $\Gamma \backslash \SL_2(\R)$, both these are estimates
 are uniform in $x_0$. 
 
   Using (\ref{eq:sl2mc}) now yields:
 \begin{equation} \label{eqn:sa} |\mu_{\delta}(a(T)^{-1} f) -\int_{\Gamma \backslash \SL_2(\mathbb{R})}
 f(g) dg| \ll_{\epsilon} T^{\alpha-1/2} 
 S_{2,1}(f) \delta^{-1-\epsilon}.\end{equation}

Finally, note that if $\chi_{[0,T]}$ denotes the characteristic function of $[0,T]$ in the
 real line, then $\frac{1}{T} \int_{t \in \mathbb{R}} |g_{\delta}(t/T) - \chi_{[0,T]}(t)| dt \ll \delta$. 
 It follows that 
 \begin{equation} \label{eqn:re} \frac{1}{T} \left| \int_{0}^{T} f(x_0 a(T)^{-1} n(t)) dt - \int_{t} g_{\delta}(t/T)  f(x_0 a(T)^{-1} n(t)) dt \right|
 \ll \delta \cdot S_{\infty,0}(f).\end{equation}
  
  Combining (\ref{eqn:ga}), (\ref{eqn:sa}) and (\ref{eqn:re}), and replacing
  $x_0$ by $x_0 a(T)$, we conclude that the left hand side of (\ref{eq:wanted})  is bounded by
 $$O_{\epsilon}\left(S_{\infty,1}(f) ( \max(T^{-1}, \delta) +  T^{\alpha-1/2} \delta^{-1-\epsilon} + \delta)\right).$$
 We choose $\delta^2 = T^{\alpha-1/2}$ to obtain the claimed conclusion. 
\qed
\subsubsection{Equidistribution of large horospheres on higher rank groups.}
We now prove quantitative equidistribution of large closed horospheres.
This result is well-known and generalizes the result of Sarnak, 
that the closed horocycle $\{x+iy\}_{0 \leq x \leq 1}$ is equidistributed
in $\SL_2(\Z) \backslash \mathbb{H}$, as $y \rightarrow 0$.

We shall follow the notation of Sec. \ref{sec:fourcoeff}, which we briefly reprise. 
Let $G$ be a connected semisimple (real) Lie group,
$\Gamma \subset G$ a lattice, $K \subset G$ the maximal
compact subgroup,
$\lieg$ the Lie algebra of $G$,
and $H \in \lieg$ a semisimple element. Fix arbitrarily a norm $\|\cdot\|$ on $\lieg$. We equip $G$ with the Haar measure in which $\Gamma \backslash G$ has volume $1$. 
Let $\exp: \lieg \rightarrow G$ be the exponential map.
Let $\mathfrak{u}$ be the sum
of all negative root spaces for $H$,
and let $U = \exp(\mathfrak{u}) \subset G$.
Let $x_0 \in \Gamma \backslash G$
be so that $x_0 U$ is compact; note that the existence of such $x_0$
implies that $\Gamma \backslash G$ is noncompact.

Let $x_t = x_0 \exp(t H)$, and let $\Delta_t$
be the stabilizer of $x_t$ in $U$.
We denote by $\langle \cdot, \cdot \rangle_{L^2(\Gamma \backslash G)}$
the inner product in the Hilbert space $L^2(\Gamma \backslash G)$.

\begin{lem} \label{lem:assmix}
There is $\kappa_1  > 0$ such that,
for any $f,g \in C^{\infty}(\Gamma \backslash G)$
and for any $U \in \mathfrak{u}$ with unit length (w.r.t. the fixed norm $\|\cdot\|$ on $\lieg$)
we have:
\begin{equation} \begin{aligned}
\left|\langle \exp(tH) \cdot f, g \rangle
- \int_{\Gamma \backslash G} f \int_{\Gamma \backslash G} g\right|
\ll \exp(-\kappa_1|t|)  S_{\infty,\dim(K)}(f) S_{\infty ,\dim(K)}(g)\\
\left| \langle \exp(sU) \cdot f, g \rangle - \int_{\Gamma
\backslash G} f \int_{\Gamma \backslash G} g \right|
\ll (1+|s|)^{-\kappa_1} S_{\infty,\dim(K)}(f) S_{\infty,\dim(K)}(g)
\end{aligned}\end{equation}
\end{lem}

Of course the constant $\kappa_1$ will depend on the choice of the norm $\|\cdot\|$. 

\proof
This follows from a nice result of Kleinbock and Margulis: see \cite{KM}.
(The orthogonal complement $L_0^2$ of the identity representation in $L^2(\Gamma \backslash G)$
is isolated, by \cite[Thm 1.12]{KM}, from the trivial representation in the unitary dual of $\widehat{G}$. A sufficiently
high tensor power of $L_0^2$ is therefore tempered, whereupon one applies
the bounds of \cite{CHH}.)
Note that \cite{KM} only
claims the result (in effect) with $S_{\infty,d}$ for some $d$;
the fact that we can take $d = \dim(K)$ follows by explicating the argument just sketched. 
\qed

Recall the definition of $\nu_T$ from Sec. \ref{sec:fourcoeff}, that is to say:
$\nu_T(f) =
\frac{\int_{\Delta_T \backslash U} f(x_T u)
du}{\mathrm{vol}(\Delta_T\backslash U)}.$ Thus $\nu_T$ is the measure
supported on a closed horosphere, and this horosphere expands as $T \rightarrow \infty$.
One deduces from Lem. \ref{lem:assmix} that the measures $\nu_T$ are equidistributed
as $T \rightarrow \infty$: 
\begin{lem} \label{lem:banana}
Set $\kappa_2= \frac{\kappa_1}{\dim(G) + \dim(K)+ 1}$,
$\kappa_1$ being as in the previous Lemma.
Then, for $T \geq 0$ and $f \in C^{\infty}(\Gamma \backslash G)$,
$$|\nu_T(f) - \int_{\Gamma \backslash G} f| \ll e^{-\kappa_2T} S_{\infty,d}(f).$$
\end{lem}
\proof
The idea is identical to Lem. \ref{lem:longhoroequi} and we refer to the first paragraph
of that proof for a description of it. 

 Fix a left-invariant Riemannian metric on $G$. This
descends to a metric on $\Gamma \backslash G$. We first choose
some ``smoothing kernels'' on $G$. For each $\epsilon > 0$, choose
a function $k_{\epsilon} \in C^{\infty}(G)$ such that
$k_{\epsilon}$ is positive, supported in an
$\epsilon$-neighbourhood of the identity, $\int_{G} k_{\epsilon} =
1$, and so that for any $X_1, X_2, \dots, X_l \in \lieg$ we have:
\begin{equation} \label{eq:kder}\sup_{g \in G} |X_1 \dots X_l k_{\epsilon}| \ll_{X_1, \dots, X_l}
\epsilon^{-l - \dim(G)}.\end{equation}
It is easy to see this is possible (for example:
choose an appropriate sequence of functions on $\mathfrak{g}$ and transport
to $G$ via the exponential map.)

The measure $\nu_0$ is a $U$-invariant probability measure
supported on the closed orbit $x_0 U$.
$\nu_0 \star k_{\epsilon}$
is supported in an $\epsilon$-neighbourhood
of $x_0 U$ and is given by integration
against a $C^{\infty}$ density function $g_{\epsilon}$, that is:
$\nu_0 \star k_{\epsilon}(f) = \int_{\Gamma \backslash G}
f g_{\epsilon}.$

Moreover, it follows from (\ref{eq:kder})
that $g_{\epsilon}$ satisfies the bounds $S_{\infty,l}(g_{\epsilon})
\ll \epsilon^{-l-\dim(G)}$, for any $l \geq 0$.

The translate of $\nu_0 \star k_{\epsilon}$
by $\exp(- T H)$ is supported in
an $\epsilon$-neighbourhood of $x_T U$;
note it is essential that $T \geq 0$ for this.
(Recall -- Sec. \ref{subsec:gennotn} -- our conventions are such that the right translate
of the point mass at $x$ by $g \in G$ is the point mass at $xg^{-1}$.) 

In fact, one verifies that
\begin{equation} \label{eq:pallavi}\left|\nu_T(f)  - \nu_0 \star k_{\epsilon}
(\exp(TH) \cdot f) \right|
\ll \epsilon S_{\infty,1}(f).\end{equation}

(Indeed, let $g \in \supp(k_{\epsilon})$ and
let $\delta_g$ be the point mass at $g$.
It suffices to check that the identical bound holds
for $\left| \nu_T(f) - \nu_0 \star \delta_g( \exp(TH) \cdot f)\right|$,
which equals $
\left|\nu_T(f) - \nu_0(g \exp(TH) \cdot f)\right|$. 
Let $\mathfrak{b}$ the sum of non-negative root spaces
for $H$ on $\mathfrak{g}$. 
If $\epsilon$ is sufficiently small, 
we may write $g  = u m$, with $u \in \exp(\mathfrak{u})$
and $m \in \exp(\mathfrak{b})$. 
Moreover, again if $\epsilon$ is sufficiently small, $u,m$
lie in a $C \epsilon$-neighbourhood of the identity,
for some fixed constant $C$. 
Then $\exp(-TH) g \exp(TH) = u' m'$, with $u' \in \exp( \mathfrak{u})$
and where $m' \in \exp(\mathfrak{b})$ is in a $C' \epsilon$-neighbourhood
of the identity, for some absolute $C'$. 
Also, $\nu_0(g \exp(TH) \cdot f) = \nu_T(m' f)$. 
Thus
it suffices to bound $\left| \nu_T(f) - \nu_T( m' 
f)\right|$.  But the $L^{\infty}$ norm of
$f - m' \cdot f$ is $\ll \epsilon S_{\infty,1}(f)$.) 

On the other hand, by Lem. \ref{lem:assmix}, for $T \geq 0$:
we have
\begin{multline}\left|\nu_0 \star k_{\epsilon} (\exp(TH) \cdot f) - \int_{\Gamma \backslash G} f\right|
= \left|\langle g_{\epsilon}, \exp(T H) f \rangle_{L^2(\Gamma \backslash
G)}  - \int_{\Gamma \backslash G} f\right| \\
\ll \exp(-\kappa_1 T) S_{\infty,\dim(K)}(f) \epsilon^{-\dim(G) -
\dim(K)}.
\end{multline}
It follows from this and (\ref{eq:pallavi}) that
$$|\nu_T(f) - \int_{\Gamma \backslash G} f| \ll (\epsilon + \exp(- \kappa_1
T)\epsilon^{-\dim(G)-\dim(K)})
S_{\infty,\dim(K)}(f).$$
 To conclude, take $\epsilon =
\exp(-\frac{\kappa_1 T}{\dim(G)+\dim(K)+1})$.
\qed

\subsection{The equidistribution of Hecke orbits and $p$-adic
horocycles.}

In this section, we prove some ``$p$-adic'' equidistribution statements, 
pertaining to the equidistribution of Hecke points and $p$-adic horocycles.

In the Lemmas that follow,
$\fcond$ will be a prime
ideal of $F$, $\overline{\mu}_{\fcond}$,
the normalized Hecke measure
defined subsequent to (\ref{eq:heckeopdef}),
and $[\fcond]$ as defined in Sec. \ref{sec:gl2pgl2}.

The first Lemma is an adelic version of the fact
that the Hecke orbit $T_q(z)$ of a point $z \in \SL(2,\Z) \backslash \mathbb{H}$
is equidistributed, as $z \rightarrow \infty$. 

\begin{lem} \label{lem:heckeequiprim}
Let $f \in C^{\infty}_{\omega}(\quotg)$
and $\fcond$ an ideal of $F$.
%with $(\fcond,\Supp(f))=1$.
Then, for $x_0 \in \quotg$, $d\gg 1$,
\begin{equation} \label{eq:phil} \left|f \star \overline{\mu}_{\fcond}(x_0)
- \sum_{\chi^2=\omega}\chi([\fcond] )\chi(x_0) \int_{x \in \quot} f(x) \chi(x) d
  \mu_{\quot}(x)\right| \ll 
\Norm(\fcond)^{\alpha-1/2+\epsilon} \height(x)^{1/2} S_{2,d}(f).\end{equation}
Here $\chi(x)$ denotes the function $g \mapsto \chi(\det(g))$ on $\quot$.
\end{lem}
\proof
Let $\mathscr{P}$ be the projection defined in Sec. \ref{subsec:proj}.
Let $E$ be the endomorphism $f \mapsto (f- \mathscr{P}f) \star \heckef$ of $C^{\infty}_{\omega}(\quotg)$. The operator $E$ has
norm $\ll_{\epsilon} \Norm(\fcond)^{\alpha-1/2+\epsilon}$ w.r.t. the $L^2$ norm 
(this follows from Lem. \ref{rem:whyproj} and the bounds of Sec. \ref{subsec:decay}). 
By Lem. \ref{lem:interp} it follows that  
the operator norm of $E$ w.r.t $S_{2,d,\beta}$ is also $\ll_{\epsilon} \Norm(\fcond)^{\alpha-1/2+\epsilon}$.

%has operator norm $\leq 2$ w.r.t the $L^{\infty}$ norm.
% (Note that, 
% $\heckef$ being a probability measure, convolution by it
 %also does not increase $L^{\infty}$ norms).  
The left hand side of (\ref{eq:phil}) is exactly $Ef(x_0)$. Now apply Lem. \ref{lem:adelicsobolev},
with $p=2$, to conclude. 
\qed

%\begin{lem}
%Let $A,B > 0$ and suppose $E: C^{\infty}(\quot) \rightarrow C^{\infty}(\quot)$
%is a linear endomorphism satisfying, for all $v \in C^{\infty}(\quot)$
%\begin{equation}
%S_{2,d}(Ev) \leq A S_{2,d}(v), \ \ S_{\infty,d}(E v) \leq
%BS_{\infty,d}(v)\end{equation}
%Then also $S_{p,d}(Ev) \leq C S_{p,d}(v)$ for all
%$v \in C^{\infty}(\quot)$, where $C = B^{2/p} C^{1-2/p}$.
%\end{lem}
%\proof Interpolation. \qed

The next Lemma is an adelic version of the following (again closely connected to equidistribution of Hecke points). 
Let $Y(p)$ be embedded in $Y(1) \times Y(1)$ (notation of discussion
after Prop. \ref{prop:main}). Then $Y(p)$ is equidistributed as $p \rightarrow \infty$. 
The quantification of this is slightly complicated by noncompactness;
in particular, we must use Sobolev norms $S_{p,d}$ for $p > 2$. (Cf. discussion
in Sec. \ref{subsec:generalcomments}). 

\begin{lem} \label{lem:heckesobolev}
Let $\q$ be a prime ideal of $\order_F$.
Let $F \in C^{\infty}(\quot \times \quot)$
be $\PGL_2(\order_{F_{\q}}) \times
\PGL_2(\order_{F_{\q}})$ invariant.Then, for any $d \gg 1, p > 2$,
\begin{multline}
\label{eq:heckesobolev}\left|
\int_{\quot} F(x, xa([\q])) dx  -
\sum_{\chi^2=1}\chi([\q]) \int_{\quot} F(x,y) \chi(x) \chi(y) d\mu_{\quot}(x)
d\mu_{\quot}(y)
\right|\\ \ll_{\epsilon} \Norm(\q)^{\frac{2\alpha - 1}{p}+\epsilon} S_{p,d}(F).\end{multline}
\end{lem}
\proof
Let $\sigma$ be the measure $\delta_1 \times \heckeq$
on $\PGL_2(F_{\q}) \times \PGL_2(F_{\q})$, where $\delta_1$ is the  measure
consisting of a point mass at the identity. 
Recalling (see (\ref{eq:heckeopdef}) in the case of a prime ideal, and Section \ref{sec:gl2pgl2} for the definition of $K_{\q}$) that
$\heckeq$ is the $K_{\q}$-bi-invariant
probability measure supported on $K_{\q} a([\q]) K_{\q}$, we note that
$$(F \star \sigma)(x,x) = \int_{k_1, k_2 \in K_{\q}} 
F(x, x k_1 a([\q]) k_2) dk_1 dk_2 = \int_{K_{\q}} F(xk, x k a([\q])) dk,$$
where we equip $K_{\q}$ with the Haar measure of mass $1$,
and we use the $\PGL_2(\order_{F_{\q}})$-invariance of $F$ at
the second step. It follows that 
\begin{equation} \label{eq:duck1} \int_{\quot} F(x,x
a([\q])) dx = \int_{\quot} (F\star \sigma)(x,x) dx.\end{equation}

Let $\mathscr{P}_2$ be as in Sec. \ref{subsec:proj}.
Let $E$ be the endomorphism of $C^{\infty}(\quot \times \quot)$ defined by
$E(F) = (F - \mathscr{P}_2 F) \star \sigma$.
Combining (\ref{eq:duck1}) and the easily verified equality $$\int_{\quot} (\mathscr{P}_2 F \star \sigma) (x,x)
dx  = \sum_{\chi^2 =1} \chi([\q]) \int_{\quot} F(x,y) \chi(x) \chi(y) d\mu_{\quot}(x)
d\mu_{\quot}(y),$$ we see that the left hand side of (\ref{eq:heckesobolev})
is precisely $\int_{\quot} EF(x,x) dx$.

Since $\mathscr{P}_2$ does not increase $L^{\infty}$ norms,
and $\sigma$ is a probability measure, it follows that
the operator norm of $E$ w.r.t the $L^{\infty}$ norm is $\leq 2$.
Moreover, the operator norm of $E$ w.r.t the $L^2$ norm is $\ll \Norm(\q)^{\alpha-1/2}$,
as follows from Lem. \ref{rem:whyproj}. 

Lem. \ref{lem:interp} now implies that for $2 \leq p \leq \infty$
we have the majorization $S_{p,d}(EF) \ll \Norm(\q)^{\frac{2\alpha-1}{p} + \epsilon} S_{p,d}(F)$. 
Now Lem. \ref{lem:adelicsobolev} shows that 
 $$|\int_{\quot} EF(x,x) dx| \ll S_{p,d}(EF)
\ll \Norm(\q)^{\frac{2\alpha-1}{p}+\epsilon} S_{p,d}(F)$$ for $p>2, d \gg 1$; whence the conclusion of  the Lemma. 
\qed

The next Lemma shows the equidistribution of certain $p$-adic horocycle orbits,
as $p$ varies. 
The idea will be as follows: (speaking very loosely, in the case of $\SL_2$) a typical
$p$-adic horocycle orbit, when projected to $\SL(2,\Z) \backslash \mathbb{H}$, looks like
$\{z+\frac{i}{p}\}_{0 \leq i \leq p-1}$. This set looks very much
like the image, under the $p$-Hecke operator, of the point $pz$.
Thus one can deduce distribution properties of the $p$-adic horocycle
orbit from some standard facts about Hecke operators.  

This is a rather {\em ad hoc} argument. Let us say a few words about why
this problem does not quite fit into the usual setup of such questions. We are proving
statements about the distribution of e.g. $p$-adic horocycles {\em when $p$ varies}. This does
not fit easily into the usual context of such matters, where one considers e.g. a fixed unipotent
flow on an $S$-arithmetic homogeneous space.  It would be interesting to have a more conceptual
and natural way of treating such questions, in the aspect where ``$p$ varies.''

\begin{lem} \label{lem:padicunipotent}
Let $f \in C^{\infty}(\quot)$ and let $\fcond$ be an integral
ideal of $\order_F$, factorizing as $\fcond = \prod_{\q}
\q^{e_{\q}}$.  For each $\q|\fcond$, let $s_{\q} \geq 0$
be a non-negative integer, and suppose $f$ is invariant by
$\prod_{\q|\fcond} K_0[\q^{s_{\q}}]$.
Put $\m = \prod_{\q|\fcond} \q^{s_{\q}}$. 
Let $\eta_{\fcond}$ be the Haar probability measure on $\prod_{\q|\fcond}
N(\q^{-e_{\q}} \order_{\q})$ and $dh$ the Haar probability measure on $\SL_2(F) \backslash \SL_2(\adele_F)$.

Then, for $y \in F^{\times} \backslash \adele_F^{\times}$,  \begin{multline} \left|f \star \eta_{\fcond} (a(y)) - \int_{h
\in SL_2(F) \backslash \SL_2(\adele_F)} f(ha(y)) dh \right|
\\ \ll_{\epsilon} \Norm(\fcond)^{\alpha - 1/2+\epsilon}
\max(\Norm(\fcond) |y|, \frac{1}{\Norm(\fcond)|y|}) ^{1/2} 
\Norm(\m)^{3/2+\epsilon} S_{2,d}(f)
\end{multline}
\end{lem}

\proof
As usual let $K_{v,f}$ be the stabilizer of $f$ in $K_v =
\GL_2(\order_{v})$, so that $K_{\q,f} $ contains $K_0[\q^{s_{\q}}]$
for each $\q|\fcond$. 

We now define a measure 
$\etap_{\q}$ on $\PGL_2(F_{\q})$ for each
$\q|\fcond$.  It will ``approximate'' $\eta_{\fcond}$ but will be composed of Hecke operators.

For those $\q$ such that $s_{\q} = 0$,
put \begin{equation} \label{eq:etapdef} \etap_{\q} = \Norm(\q)^{-e_{\q}/2} \delta_{a(\varpi^{-e_{\q}})} \star \mu_{\q^{e_{\q}}} -
\Norm(\q)^{\frac{-e_{\q}-1}{2}} \delta_{a(\varpi^{-e_{\q}-1})} \star \mu_{\q^{e_{\q}-1}}.\end{equation}
(We refer to Sec. \ref{subsec:Hecke} for definitions of $\mu_?$ appearing above.) 
For $\q$ such that $s_{\q} \geq 1$, we set
$\sigma_{\q}$ to be the
unique bi-$K_0[\q^{s_{\q}}]$-invariant probability measure on $K_0[\q^{s_{\q}}] 
a(\varpi^{e_{\q}}) K_0[\q^{s_{\q}}]$, normalized to have mass $1$, and we put
$\etap_{\q} = \delta_{a(\varpi^{-e_{\q}})} \star \sigma_{\q}$. 
Finally, set
$\etap_{\fcond} = \prod_{\q|\fcond} \etap_{\q}$.

One then verifies by a direct computation that \begin{equation} \label{eq:hecketheory} f \star
\eta_{\fcond}  =  f \star \etap_{\fcond}\end{equation}
The intuition for this statement, in the classical setting, as as follows:
let $z \in \SL_2(\Z) \backslash \mathbb{H}$. Then (for a prime number $p$) the set $\{z+i/p\}_{0 \leq i \leq p-1}$
is the $p$-Hecke orbit of $pz$, with the point $p^2 z$ removed. 
In the case $e_{\q}= 1$, the first  term on the right hand side of (\ref{eq:etapdef}) corresponds to the $p$-Hecke orbit of $pz$,
and the second term corresponds to removing the point $p^2 z$. 

More formally, to verify (\ref{eq:hecketheory}),
the unramified computation, at those places where  $s_{\q} = 0$, is easy; the ramified computation is just Hecke theory at ramified primes,
see e.g. \cite[Prop 3.33]{Sh}. \footnote{For the unramified assertion, let $\mathcal{B} =\PGL_2(F_{\q})/K_{\q}$
and let $x_0 \in \mathcal{B}$ be the identity coset. The set $\mathcal{B}$
has the structure of the vertices of a $q_v+1$-valent tree.
Let $S_1$ be the set of all vertices at distance
$ e_{\q} - 2i $ (some $i \geq 0$) from $a(\varpi^{-e_{\q}}) x_0$. Let $S_2$ be the
set of all vertices at even distance $\leq e_{\q}-1 - 2i$ (some $i \geq 0$) from 
$a(\varpi^{-e_{\q}-1}) x_0$.  Then $S_2 \subset S_1$
and $S_1 - S_2$ is precisely the $n(\q^{-e_{\q}})$-orbit of $x_0$. 
As for the ramified case:
one notes that, if $s_{\q} > 0$, the Haar measure on $K_0[\q^{s_{\q}}]$
is just the pushforward of the Haar measure on $ n(\order_{\q}) \times a(\order_{\q}^{\times})
\times \bar{n}(\q^{s_{\q}}) $ by the product map $(n,a,\bar{n}) \mapsto n a \bar{n}$.}

The projection $\mathscr{P}$ of Sec. \ref{subsec:proj} commutes with the action of $\GL_2(\adele_F)$,
and so (\ref{eq:hecketheory}) holds also with $f$ replaced by $\mathscr{P}f$
or $f - \mathscr{P} f$. 
Moreover, $\mathscr{P}f \star \eta_{\fcond} = \mathscr{P}f$. 
It
follows that
\begin{equation} \label{eq:raymond}f \star \eta_{\fcond}(x) = \int_{h \in SL_2(F)
\backslash \SL_2(\adele_F)} f(hx) dh + (f - \mathscr{P}f) \star
\etap_{\fcond}(x).\end{equation}

Set $\bar{f} = f - \mathscr{P}f$. Then, expanding the term $\bar{f} \star \etap_{\fcond}$:
\begin{multline} \label{eq:hagrid}\bar{f} \star \etap_{\fcond}  =
\sum_{S \subset \{\q|\fcond, s_{\q} = 0\}}
\bar{f} \star
\prod_{\q|\fcond: s_{\q} \geq 1} \delta_{a(\unif^{-e_{\q}})} \star \sigma_{\q}
\star \\  \prod_{\q|\fcond: s_{\q} = 0, \q \notin S }
\left( \Norm(\q)^{-e_{\q}/2} \delta_{a(\unif^{-e_{\q}})} \star 
\mu_{\q^{e_{\q}}} \right)  \star \prod_{\q \in S} 
\left(- \Norm(\q)^{-\frac{e_{\q}+1}{2}} \delta_{a(\unif^{-e_{\q}-1})}\star \mu_{\q^{e_{\q}-1}} \right) 
\end{multline}

We now specialize to the case under consideration where $x = a(y)$
for some $y \in \adele_F^{\times}$. 
For $S \subset \{\q|\fcond, s_{\q} =0\}$ set $$\sigma_S = 
\prod_{\q|\fcond: s_{\q} \geq 1} \sigma_{\q}
 \prod_{\q|\fcond: s_{\q} = 0, \q \notin S }
\Norm(\q)^{-e_{\q}/2} \mu_{\q^{e_{\q}}} \star \prod_{s_{\q} = 0, \q \in S} -\Norm(\q)^{-\frac{e_{\q}+1}{2}} \mu_{\q^{e_{\q}-1}}.$$

With this notation, we have:
\begin{equation} \label{eq:joey}\bar{f} \star \etap_{\fcond}(a(y)) =
\sum_{S \subset \{\q|\fcond, s_{\q} = 0\}}
 \bar{f} \star  \sigma_S \left( a(y  [\fcond]) \prod_{s_{\q} = 0, \q \in S} a([\q])\right)
\end{equation}

%Fix $S \subset \{\q|\fcond: s_{\q} = 0\}$ and put
%$$\sigma_S = 
%\prod_{\q|\fcond: s_{\q} \geq 1} \sigma_{\q}
%\star \prod_{\q|\fcond: s_{\q} = 0, \q \notin S }
%\mu_{\q^{e_{\q}}} \star \prod_{\q \in S} (-\Norm(\q))^{-1}
%\mu_{\q^{e_{\q}-1}} \delta_{a(\unif_{\q}^{-1})}.$$

Now apply Lem. \ref{lem:adelicsobolev}
to see that, for any $z \in \adele_F^{\times}$ and $d \gg 1$,
we have \begin{equation} \label{eq:krisna}\bar{f} \star \sigma_S (a(z)) \ll 
\max(|z|, |z|^{-1})^{1/2} PS_{2,d}(\bar{f} \star \sigma_S),\end{equation}
where we have used the easily verified fact
that $\height(a(z)) \asymp \max(|z|,|z|^{-1})$.

Now, for any $f \in C^{\infty}(\quot)$, we have
$$[\Kmax:K_{\bar{f} \star \sigma_S}] \leq \prod_{s_{\q}\geq 1}
[K_{\q}: K_0[\q^{s_{\q}}]] [\Kmax: K_f] \ll_{\epsilon} \Norm(\m)^{1+\epsilon}
[\Kmax:K_f].$$
 By the bounds
on matrix coefficients (\ref{eq:measopnormbound}), and recalling that $\m = \prod_{\q|\fcond} \q^{s_{\q}}$,
we compute that 

%\footnote{For example, let $\pi$
%be a representation of $\PGL_2(F_{\q})$ and $v$
%in the space of $\pi$. To bound $\|v \star \sigma_{\q}\|_2 =
%\sup_{\|w\|_2 =1} \langle v \star \sigma_{\q}, w\rangle$. Moreover,
%one can restrict to $w$ that are $K_0[\q^{s_{\q}}]$-invariant.
%So $\|v \star \sigma_{\q}\| \leq [K_{\q}: K_0[\q^{s_{\q}}]] \Norm(\q)^{\alpha-1/2+\epsilon}$.}
%
%$$\|\bar{f} \star \sigma_S\|_{L^2} 
%\ll_{\epsilon} \|f\|_{L^2}  \prod_{\q|\fcond} \Norm(\q)^{s_{\q}+\epsilon} (\Norm(\fcond) \prod_{\q \in S}
%\Norm(\q)^{-1})^{\alpha - 1/2+\epsilon} \prod_{\q \in S} \Norm(\q)^{-1}.
%$$
%A similar bound holds if we replace $f$ by $\mathcal{D} f$,
%for any $\mathcal{D}$ belonging to the universal enveloping algebra
%of $\GL_2(F_{\infty})$. Thus
$$PS_{2,d}(\bar{f} \star \sigma_S) \ll_{\epsilon}
\Norm(\m)^{3/2+\epsilon}
(\Norm(\fcond) \prod_{\q \in S}
\Norm(\q)^{-1})^{\alpha - 1/2+\epsilon} \prod_{\q \in S} \Norm(\q)^{-1} PS_{2,d}(f).$$

Combining this with (\ref{eq:joey}) and (\ref{eq:krisna}), we find
that for each $S \subset \{\q : s_{\q} \geq 1\}$
\begin{multline} \label{eq:deb} \left| \bar{f} \star \etap_{\fcond} (a(y)) \right| \ll_{\epsilon}
\Norm(\fcond)^{\alpha-1/2+\epsilon} \Norm(\m)^{3/2+\epsilon}
 \max(\Norm(\fcond) |y|, \Norm(\fcond)^{-1}|y|^{-1})^{1/2}
 PS_{2,d}(f)
%\sum_{S
%\subset \{\q|\fcond, s_{\q} =0\}} 
%\max(\Norm(\fcond) |y| \prod_{\q \in S} \Norm(\q),
%\Norm(\fcond)^{-1} |y|^{-1} \prod_{\q \in S}
%\Norm(\q)^{-1})^{1/2}  \\
%(\Norm(\fcond) \prod_{\q \in S}
%\Norm(\q)^{-1})^{\alpha - 1/2+\epsilon} \prod_{\q \in S} \Norm(\q)^{-1}
%\\
%\ll \Norm(\fcond)^{\alpha-1/2+\epsilon} PS_{2,d}(f) 
%\max(\Norm(\fcond) |y|, \Norm(\fcond)^{-1} |y|^{-1})\sum_{S \subset \{\q|\fcond, s_{\q} =0\}}
%\prod_{\q notin S} \Norm(\q)^{-\alpha}
%\\ \ll \Norm(\fcond)^{\alpha - 1/2+\epsilon} PS_{2,d}(f)  
\end{multline}
%From the definitions, we now deduce that
%$$PS_{2,d} (\bar{f} \star \sigma_S) 
%\ll  PS_{2,d}(f) $$
%
%To estimate
%$\bar{f} \star \sigma_{\fcond}(a(y [\fcond]))$ we
%combine the expansion (\ref{eq:hagrid}) with Lem.
%\ref{lem:adelicsobolev}. We make use of the fact that
% the height of $a(y) \in \quot$
%is $\asymp \max(|y|, |y|^{-1})$, and that
%for each finite place $v$ we have $K_{v, \bar{f} \star \sigma_{\fcond}}
%= K_{v,f}$.  After some routine computations one obtains:
%\begin{multline}\left| f \star \eta_{\fcond}(a(y)) - \int_{h \in
%SL_2(F) \backslash \SL_2(\adele_F)} f(h a(y)) dh \right| \\
%\ll_{\epsilon}
%%\Norm(\fcond)^{\alpha - 1/2+\epsilon}
%%\max(\Norm(\fcond)y,\frac{1}{\Norm(\fcond)y})^{1/2} [\Kmax: K_f]^{3/2}
%%S_{2,d}(f)  \\
% \Norm(\fcond)^{\alpha - 1/2 +\epsilon} \max(\Norm(\fcond) |y|,
%\frac{1}{\Norm(\fcond)|y|})^{1/2}   PS_{2,d,1}(f).\end{multline}
This bound is valid for all $f \in C^{\infty}(\quot)$, not merely those $f$ that
are invariant by $\prod_{\q|\fcond} K_0[\q^{s_{\q}}]$. 
Apply Rem. \ref{rem:def}, (\ref{end}) to the endomorphism $f \mapsto \bar{f} \star \etap_{\fcond}$; this shows
that (\ref{eq:deb}) remains valid, for any $f \in C^{\infty}(\quot)$, if we replace $PS_{2,d}$ by $S_{2,d}$ on the right hand side. 
Now, specialize to the case where $f \in C^{\infty}(\quot)$ is actually
$\prod_{\q|\fcond} K_0[\q^{s_{\q}}]$-invariant and apply (\ref{eq:raymond}) to obtain the conclusion of the Lemma. 
\qed

The Lemma that follows states an adelic version of the following fact:
the measure on $\SL_2(\Z) \backslash \H$ defined by $\nu_y := q^{-1} \sum_{0 \leq x \leq q-1} \delta_{\frac{x}{q} + iy}$, approximates the uniform measure if $y \asymp q^{-1}$;
more precisely we have an inequality that $\left|\nu_y(f) - \int_{\SL_2(\Z) \backslash \H}f \right|$
is bounded by $\max(qy, \frac{1}{qy})^{1/2} q^{-\delta} S(f)$, where $S$ is an appropriate Sobolev norm
and $\delta > 0$.

\begin{lem} \label{lem:measboundvar}
Let $f \in C^{\infty}(\quot)$ and let notations be as in Sec.
\ref{sec:torus1}(see esp. \eqref{eq:nudef}).  In particular, $\fcond$ is an integral ideal of
$\order_F$, $q =\Norm(\fcond)$, $[\fcond]$ is as in
(\ref{eq:fconddef}) and $$\nu_z(f) = \int_{|y| = z, y \in
\adele_F^{\times}/F^{\times}} f(a(y) n([\fcond])) d^{\times}y.$$

Suppose $f$ is invariant by $K_0[\q^{s_{\q}}]$, for each $\q | \fcond$,
and put $\m = \prod_{\q|\fcond} \q^{s_{\q}}$. Then
\begin{multline}
\label{eq:mss}\left|\nu_z(f)-\int_{\quot} f(x)
d\mu_{\quot}(x)\right| \\
 \ll_{\epsilon} \Norm(\fcond)^{\alpha-1/2+\epsilon} \Norm(\m)^{3/2+\epsilon} \max(\Norm(\fcond)z,\frac{1}{\Norm(\fcond) z})^{1/2}
S_{2,d}(f).\end{multline}
\end{lem}
\proof For each $\q | \fcond$ and integer $0 \leq e \in \mathbb{Z}$, let
$\eta_{\q^{e}}$ be the Haar probability measure on the group
$N(\q^{-e} \order_{\q})$. 
Then, since
the assumption implies that $f$ is right invariant by $a_{\q}(\order_{\q}^{\times})$, for each $\q$ dividing $\fcond$, we see that for any $x \in \quot$:
\begin{multline}\int_{y \in \order_{F_{\q}}^{\times}}
 f(x a(y) n_{\q}(\unif_{\q}^{-e})) d^{\times}y
 = \int_{y \in
\order_{F_{\q}}^{\times}} f(x n_{\q}(y \unif_{\q}^{-e}))
d^{\times}y \\ = \frac{ f \star (\eta_{\q^e} - \Norm(\q)^{-1}
\eta_{\q^{e-1}})(x)} {1-\Norm(\q)^{-1}}
\end{multline}

It follows that $$\nu_z(f) = \prod_{\q|\fcond}(1-\Norm(\q)^{-1})^{-1}
\cdot  \int_{y \in \adele_F^{\times}/F^{\times}, |y| = z} f \star
\prod_{\q|\fcond} {(\eta_{\q^e_{\q}} - \Norm(\q)^{-1}
\eta_{\q^{e_{\q}-1}})} (a(y)).$$ We conclude by applying the
previous Lemma. 
\qed

%
%
%
%
%\begin{lem} \label{lem:abheri}
%Let $f \in  C^{\infty}_{\omega}(\quotg)$ and put
%$$I_T(f) = \int_{x \in F \backslash \adele_F}
%\int_{y \in F^{\times} \backslash \adele_F^{\times},
%T \leq y \leq 2T}
%f(n(x) a(y)) dx d^{\times}y.$$
%Then for $d \gg 1$,
%$\left|\int_{z = T}^{2T} I_z(f) d^{\times}z \right|
%\ll_{\epsilon} T^{\epsilon} S_{2,d}(f) $.
%\end{lem}
%\proof
%Let $K_f \subset \Kmaxg$ be the stabilizer of $f$.
%Consider the function $L$ on $K_{\infty}$ defined by
%$$L(k)= I_T(k \cdot f).$$
%Then $L(k)$ is a positive function and satisfies
%$$\int_{k \in K_{\infty}} |L(k)|^2
%\ll \frac{[\Kmaxg: K_f]}{T} \int_{g \in B(F) \PGL_2(\adele_F):
%T \leq \height(g) \leq 2T} |f(g)|^2 dg.$$
%
%Applying (\ref{eq:nottoo}),
%%we conclude that for $d \gg 1$ ation $\|L(k)\|_{L^2(K_{\infty})}
%\ll  T^{\epsilon} S_{2,d}(f)$.
%A similar estimate holds with $f$ replaced by any right-invariant
%derivative.
%Applying the usual Sobolev estimate on the real manifold $K_{\infty}$,
%we obtain the conclusion.
%\qed

\section{Background on Eisenstein series.} \label{sec:eisenstein}

This section essentially develops the theory of Eisenstein series on $\PGL_2$ over a number field.
This is needed for the Rankin-Selberg method that we reprise in the next section, which
in turn is used in the text to relate a period integral with an $L$-function. 

Let $Z$ be a topological space. In this section, we will often speak -- in various contexts, 
often with $Z = \quot$ or $\GL_2(\adele_F)$ -- of a function $F(s,z)$ on $\C \times Z$
being ``holomorphic'' or ``holomorphic in $s$.'' For the purposes of this document, this can 
be assumed to mean that the function is jointly continuous and holomorphic for each $z$ individually. 

Note that $s \mapsto \int_{Z} F(s,z) dz$, if absolutely convergent and uniformly so in $s$, defines a holomorphic function. Indeed, it suffices to verify that its integral over a closed curve in the $s$-variable
is zero, which follows by Fubini's theorem. 

Similarly, we will say that $F(s,z)$ is meromorphic if there exists a holomorphic function
$h(s)$ so that $h(s) F(s,z)$ is holomorphic.

\subsection{Construction and basic properties of the Eisenstein series.}
We recall the Eisenstein series that we shall have need of and its basic
properties, following Jacquet \cite[\S19]{jacquet}. We will need Eisenstein series
only on $\PGL_2$. 

\subsubsection{Schwarz functions.}
Let $\Psi$ be a Schwarz-Bruhat function on $\adele_F^2$,
i.e. $\Psi$ is a finite linear combination
of functions $\prod_{v} \Psi_v$,
where each $\Psi_v$ is locally constant
of compact support, for $v$ finite,
 $\Psi_v$ is a Schwarz function on $F_v^2$ for $v$ infinite, and $\Psi_v$
 is the characteristic function of $\order_v^2$ for almost all $v$. 

If $v$ is a real place, choose $a_v \in F_v$ so that
$e_{F_v}(x) = e^{2 \pi i a_v x}$, and
say a Schwarz function $\Psi_v$ on $F_v^2$ is {\em standard}
if it is the product of a polynomial and $e^{-\pi |a_v|_v (|x|_v^2+|y|_v^2)}$. If $v$ is a complex place, choose
$a_v \in F_v$ so that $e_{F_v}(x) = e^{2 \pi i \Tr_{\C/\R}(a_v x)}$; 
we say that a Schwarz function $\Psi_v$ is {\em standard}
if it is the product of a polynomial and $e^{-2 \pi |a|_v^{1/2} (|x|_v+|y|_v)}$. 
The significance of this normalization is twofold: a standard function is automatically $K_v$-finite and also
the class of standard functions is self-dual under the Fourier transform corresponding to the 
character $e_{F_v}$.

If $V$ is a real vector space, then by a {\em Schwarz norm}
on the space of Schwarz functions on $V$, we 
shall mean a norm $\mathcal{S}$ of the form
\begin{equation} \label{schwarznormdef} \mathcal{S}(\Psi) =  \sup_{\mathcal{D}} \sup_{x} |(1+\|x\|)^M\mathcal{D} \Psi|, \end{equation} for 
some finite collection of constant-coefficients differential operators $\mathcal{D}$ on $V$
and some norm $\|x\|$ on $V$. 

Put, for $g \in \GL_2(\adele_F)$, $$f_{\Psi}(s,g) =
|\det(g)|^{s}
\int_{t \in
\adele_F^{\times}}\Psi((0,t)g)|t|^{2s}  d^{\times} t.$$
The integral converges absolutely for $\Re(s) > 1/2$
and extends to a meromorphic function of $s$ with possible
poles at most at $s=0,1/2$. Moreover, for all $s$,
\begin{equation} \label{eq:ftrans}
f_{\Psi}(\left(\begin{array}{cc} a & x \\ 0 & b \end{array}\right) g)
= |a/b|^{s}f(g).\end{equation}
Put $E_{\Psi}(s,g) = \sum_{\gamma \in B(F) \backslash \GL_2(F)}
f(s,\gamma g)$. This converges when $\mathrm{Re}(s) > 1$,
extends to a meromorphic function of $s$ with a simple pole at $s=0,1$
and satisfies the functional equation $$E_{\Psi}(s,g)
= E_{\widehat{\Psi}}(1-s,g),$$ where $\widehat{\Psi}$ is the Fourier transform 
\begin{equation} \label{fourierpsidef} \widehat{\Psi}(x_1, y_1) = \int_{\adele_F^2} \Psi(x,y) e_F(x_1 y - y_1 x) dx dy.\end{equation}
Moreover, the pole at $s=1$ is the constant
function with value $c_1 \int_{\adele_F^2} \Psi(x,y) dxdy$,
and the pole at $s=0$ is the constant function with value
$c_2 \Psi(0)$, where $c_1, c_2$ are constants (depending only on the choice of measure).  Finally, for any fixed $g$ the function $s \mapsto s (1-s) E_{\Psi}(s,g)$
decays rapidly in vertical strips, i.e. $(1+|s|)^N |s(1-s) E_{\Psi}(s,g)|$ is bounded
in any strip $A \leq \Re(s) \leq B$. The proof of all these properties
follows from ``Poisson summation'' for $F^2 \subset \adele_F^2$, and we omit them. 

Moreover, the association $\Psi \mapsto E_{\Psi}$ is twisted-equivariant for the
natural $\GL_2(\adele_F)$-action on the space of Schwarz functions
and on $C^{\infty}(\quot)$: that is to say, 
\begin{equation} \label{eq:twistedequivariant} E_{h . \Psi}(s,g) = |\det(h)|^{-s} (h \cdot E_{\Psi}(s,g)), \end{equation} where $h \cdot $ denotes right translation by $h$. 

We give an example with $F= \mathbb{Q}$ (cf. \cite[(3.29)]{iwaniec}). 
\begin{ex}\label{qexam}
Suppose $F = \mathbb{Q}$, $\Psi = \prod_{v} \Psi_v$ where,
for each finite $v$, $\Psi_v$ is the characteristic function
of the maximal compact of $F_v$, and $\Psi_{\infty}(x,y)=
e^{-\pi(x^2+y^2)}$.
Then $E_{\Psi}(g)$ is determined by its restriction to $\SL_2(\R)$.

%invariant by the standard maximal compact
%of $\GL_2(\adele_\Q)$, and also by the center of $\GL_2(\adele)$.
%Consequently, it is determined
%by its restriction to $\SL_2(\R)$.
%
%For $g \in \SL_2(\R)$, we have thereby
%\begin{equation} \label{epsi}E_{\Psi}(s,g) = \zeta(2s) \sum_{(c,d)=1} \int_{t \in \R^{\times}}
%e^{-t^2 ||(c,d) g||^2} |t|^{2s} d^{\times} t \end{equation}

Moreover, $E_{\Psi}(s,g)$ descends from a function of $g \in \SL_2(\R)$
 to a function $E^*(s,z)$ on $\mathbb{H} =
\SL_2(\R)/\mathrm{SO}_2$, where the identification is $g \mapsto g \cdot i$.  In fact,
\begin{equation} \label{eisid}E^*(s,z) = \pi^{-s} \Gamma(s) \zeta(2s) \sum_{[c:d] \in \mathbb{P}^1(\mathbb{Q})} \frac{y^s}
{|cz+d|^{2s}}.\end{equation}
If we put $\xi(s) = \pi^{-s/2} \Gamma(s/2) \zeta(s)$,
then $E^*(s,z)$ has the Fourier expansion
\begin{equation} \label{eq:four}E^*(s,z) = \xi(2s) y^{s} + \xi(2-2s) y^{1-s}
+4 \sqrt{y}\sum_{n \in \mathbb{N}} K_{s-1/2}(2 \pi n y)
\cos(2 \pi n y) \sum_{ab= n}  \left(\frac{a}{b} \right)^{s-1/2}
\end{equation}
It satisfies the functional equation $E^*(s,z) = E^*(1-s,z)$.
Moreover it is a meromorphic function of $s$ with poles precisely
at $s=0$ and $s=1$. In both cases the residue is the constant function.
\end{ex}

Motivated by this example, the reader may find
it helpful to keep in mind the ``dictionary'':
$f_{\Psi}(s,g)$ corresponds to $\pi^{-s} \Gamma(s) \zeta(2s)y^s= \xi(2s) y^s$,
and $E_{\Psi}(s,g)$ to $E^*(s,z)$ as defined in (\ref{eisid}).

\begin{rem}
Suppose $\Psi$ is invariant by $K_{\infty} \times \Kmax$.
Then $f_{\Psi}$ is a multiple of $g \mapsto \height(g)^s$,
as follows from the uniqueness of spherical functions
satisfying (\ref{eq:ftrans}).
Thus, for $\Re(s) > 1$,
$E_{\Psi}(s,g) = c(s) \sum_{\gamma \in B(F) \backslash \GL_2(F)}
\height(\gamma g)^{s}$.
\end{rem}

We now proceed to establish the ``standard''
properties of the Eisenstein series for $E_{\Psi}$.
It is convenient to first recall an explicit bound for archimedean Mellin transforms; the first
part is Tate's thesis, and the second will only be needed much later. 

\begin{lem} \label{lem:archimedean}
Let $v$ be archimedean and let $\Psi_v$ be a Schwarz function on $F_v$.
The integral $G(s) := \int_{x \in F_v^{\times}} \Psi_v(x) |x|^{s} d^{\times}x$ extends to a meromorphic
function and:

\begin{enumerate}
\item 
$\frac{G(s)}{\zeta_{F,v}(s)}$ is holomorphic, where $\zeta_{F,v}(s)$ is the local
factor of the Dedekind $\zeta$-function of $F$ at $v$. 

\item For  any $N \geq 0$, the function $G_N(s) := 
\prod_{i=0}^{ N } (s +i) G(s)$
is holomorphic in $\Re(s) \geq -N$, and the absolute value of $(1+|s|)^M G_N(s)$
in any strip $-N \leq \Re(s) \leq A$ is bounded by some Schwarz norm (depending on $A,N, M$;
see \eqref{schwarznormdef} for the definition) of $\Psi_v$.
\end{enumerate}
\end{lem}
\proof
The first assertion
is Tate's thesis, and we leave the second to the reader (if any). \qed

\begin{lem} \label{lem:fbound}
The function $s \mapsto s (1/2-s) f_{\Psi}(s,g)$
extends to a holomorphic function of $s$. It decays rapidly along vertical lines:
\begin{equation}
\label{eq:fbound}|(1+|\Im(s)|)^N s (1/2-s)f_{\Psi}(s,g)|
\ll_{\Psi} \height(g)^{\Re(s)},
\end{equation}
where the implicit
constant is uniform for $\Re(s)$ in a compact set.
\end{lem}
\proof
By (\ref{eq:ftrans}) and the Iwasawa decomposition, it
will suffice to prove the assertions in the special case $g \in K_{\infty} \times \Kmax$. 
So we write $g = k \in K_{\infty} \times \Kmax$ and denote by $k_v$ the component of $k$
in $\PGL_2(F_v)$. 
Moreover, without loss of generality, we may assume $\Psi$ is a product
of Schwarz functions at each place, i.e. $\Psi = \prod_{v} \Psi_v$. 
Then
 $$f_{\Psi}(s,g) =  \prod_{v \, \mathrm{infinite}} \int_{F_v^{\times}}
\Psi_v((0,t) k) |t|^{2s} d^{\times}t 
\prod_{v \, \mathrm{finite}} \int_{F_v^{\times}} \Psi_v((0,t) k_v) |t|^{2s} d^{\times}t$$
By Tate's thesis, it follows that that the product over finite places is of the form
$\zeta_F(2s) h(s)$, 
where $\zeta_F(\cdot)$ is the (finite part of the) Dedekind $\zeta$-function of the number field
$F$ and $h(s)$ is a holomorphic function with at most polynomial growth in
vertical strips (indeed, a
polynomial in $q^{\pm s}$ for various $q$).  All the assertions of the Lemma
now follow from Lem. \ref{lem:archimedean}, and standard
facts about the analytic properties of $\zeta_F$.

In fact, if the $\Psi_v$ for $v$ finite are regarded as fixed, then the implicit
constant in (\ref{eq:fbound}) is bounded by an appropriate
Schwarz norm, depending on $N$ and the compact set to which $\Re(s) $ is constrained, of $\prod_{v \, \mathrm{infinite}} \Psi_v$. This follows
from the second assertion of Lem. \ref{lem:archimedean}. 
%The claim that $s \mapsto s (1/2-s) f_{\Psi}(s,g)$
%is holomorphic, when $g$ is held fixed,  follows from Lem. \ref{lem:archimedean}. 
%
%Next, let $\Re(s) = \sigma$.
%It follows from the transformation property
%(\ref{eq:ftrans}) that, for $\Re(s) = \sigma$,
%$$|s(1/2-s) f_{\Psi}(s,g)| \leq \height(g)^{\sigma}
%\sup_{g \in K_{\infty} \times \Kmax } |s(1/2-s) f_{\Psi}(s,g)|.$$
%But, for fixed $g$, the function
%$s (1/2-s) f_{\Psi}(s,g)$ decays rapidly in vertical strips.
%From the fact $\Psi$ is $K_{\infty}\times \Kmax$-finite,
%we easily deduce the conclusion. 6COMPLETE; not valid]]
\qed

\begin{lem} \label{lem:eisconstterm}
The constant term $E^N_{\Psi}(s,g) := \int_{x \in F \backslash \adele_F}
E_{\Psi}(s, n(x) g) dx$ equals $f_{\Psi}(s,g) + f_{\widehat{\Psi}}(1-s, g)$. 
\end{lem}

\proof  (Sketch). A double coset decomposition shows that, for $s \gg 1$,  $E^N_{\Psi}(s,g) = f_{\Psi}(s, g) + \int_{n \in N(\adele_F)}
f_{\Psi}(s, wn g) dn$. So it will suffice to show that $\int_{n \in N(\adele_F)} f_{\Psi}(s, wng) = 
f_{\widehat{\Psi}}(1-s, g)$. The left-hand side may be expressed
as \begin{multline} |\det(g)|^s \int_{t \in \adele_F^{\times}} \int_{x \in \adele_F}
\Psi((t, tx) g) |t|^{2s} d^{\times} t \, dx \\ = |\det(g)|^s \int_{t \in \adele_F^{\times}/F^{\times}}
\int_{x \in \adele_F}  \sum_{\delta \in F^{\times}} \Psi(t (\delta, x) g) |t|^{2s} d^{\times} t \, dx \end{multline}

For any Schwarz function $\Psi$ on $\adele_F^2$, one has
$\sum_{\alpha \in F} \int_{y \in \adele_F} \Psi(\alpha, y) = 
\sum_{\beta \in F} \widehat{\Psi}(0, \beta)$. The result
follows from routine manipulation and use of Tate's functional equation. 
\qed

We set \begin{equation} \label{eq:baredef} \bar{E}_{\Psi}(s,g) = E_{\Psi}(s,g) -
f_{\Psi}(s,g) - f_{\widehat{\Psi}}(1-s,g), \end{equation}
so $\bar{E}_{\Psi}$ defines a function
on $B(F) \backslash PGL_2(\adele_F)$. It is a ``truncated'' Eisenstein series where we have removed
the constant term. Moreover, $\bar{E}_{\Psi}(s,g)$ is holomorphic in $s$ (this follows,
for example, by computing residues at each of the points $s=0,1/2,1$ and seeing they are all zero).  By definition, for $g \in \GL_2(\adele_F)$ we have an equality
\begin{equation} \label{eq:trunceis} 
E_{\Psi}(s,g)  = \bar{E}_{\Psi}(s,g) + f_{\Psi}(s,g) + f_{\widehat{\Psi}}(1-s,g).
\end{equation}

\begin{lem} \label{lem:eisbound}
Let $T, N > 0$ and let $\Re(s)$ lie in a fixed compact
subset of $\R$.
Then
\begin{equation} \label{eq:jb}
(1+|s|)^4\bar{E}_{\Psi}(s,g) \ll_{\Psi,N,T} \height(g)^{-N},\end{equation}
for $g \in \Sieg(T)$.
In particular, if $\Omega \subset \quot$ is compact,
then $s(1-s)E_{\Psi}(s,g)$ is uniformly bounded in
$|\Re(s)| \leq 2, g \in \Omega$.
\end{lem}
\proof

We first claim that, for $t \in \mathbb{R}$, we have $|(1+t^4) \bar{E}_{\Psi}(N+1+it,g)| \ll_{\Psi}
\height(g)^{-N+\epsilon}$.
Indeed, by definition, $$\bar{E}_{\Psi}(s,g)  = \sum_{\gamma \in B(F) \backslash \PGL_2(F),
\gamma \notin B(F)}f_{\Psi}(s,\gamma g) - f_{\widehat{\Psi}}(1-s,g).$$
In view of Lem. \ref{lem:fbound}, it will suffice to show that
\begin{equation}\label{eq:semmangudi}\sum_{\gamma \in B(F) \backslash \PGL_2(F), \gamma \notin B(F)}
\height(\gamma g)^{\sigma} \ll_{\epsilon} \height(g)^{1-\sigma+\epsilon},
\end{equation}
which follows from (\ref{eq:nottoo}) and (\ref{eq:bigone}).

Now (\ref{eq:jb}) follows
at once from the functional equation
$\bar{E}_{\Psi}(s,g) = \bar{E}_{\widehat{\Psi}}(1-s,g)$,
the maximal modulus principle
in the strip $|\Re(s)| \leq N+1$, and the previous Lemma.\footnote{To apply the maximal modulus principle in this context, one needs some {\em a priori} decay of $\bar{E}_{\Psi}$,
which follows easily from the corresponding properties of $E_{\Psi}$
and $f_{\Psi}$.}

The second assertion (involving $\Omega$) follows from (\ref{eq:fbound}) and (\ref{eq:jb}).
\qed

We now compute the Fourier coefficients of the Eisenstein
series in general. Recall that $e_F$ is a fixed additive character of $\adele_F/F$.
\begin{lem} \label{lem:foureis}
Set $W_{\Psi}(s,g) = \int_{x \in F \backslash \adele_F}
E_{\Psi}(s,n(x)g) e_F(x) dx$. Then, for $\Re(s) > 1$,
\begin{equation} \label{eq:foureis}W_{\Psi}(a(y)) = |y|^{1-s} \int_{t \in \adele_F^{\times}, x\in\adele_F}
\Psi(t, t x) e_F(xy) |t|^{2s} dx d^{\times}t. \end{equation}
In particular, if $\Psi = \otimes_v \Psi_v$,
then $W_{\Psi} = \prod_{v} W_{\Psi_v}$, where for $\Re(s) > 1$,
$$W_{\Psi_v}(a(y))= |y|_v^{1-s} \int_{t \in F_v^{\times}, x \in F_v}
\Psi_v(t, tx) e_F(xy)|t|_v^{2s} dx d^{\times}t,$$ for $y \in F_v$.
%In particular, if $v$ is nonarchimedean,
%$\Psi_v$ the characteristic function of $\order_{F_v}^2$,
%then $W_{\Psi_v,s}(1) = 1$.
 Finally, if $\Psi_v (x,y) = \varphi_1(x) \varphi_2(y)$,
$\omega_v$ a character of $F_v$, and
$\Re(s') + |\Re(s)| \gg 1$,
\begin{multline} \label{eq:melltrans}\int_{y \in F_v^{\times}}W_{\Psi_v}(a(y)) |y|^{s'} \omega_v(y) d^{\times} y
\\
= \int_{y \in F_v^{\times}} \varphi_1(y) |y|^{s'+s}  \omega_v(y) d^{\times}y
\int_{y \in F_v^{\times}} \widehat{\varphi_2}(y) |y|^{1+s'-s} \omega_v(y)
d^{\times}y,
\end{multline}
where $\widehat{\varphi_2}$ is the Fourier transform, defined
by $\widehat{\varphi_2}(y) = \int_{F_v} \varphi_2(y) e_{F_v}(yt) dt$. 

\end{lem}

\proof
By the Bruhat decomposition,
\begin{multline} W_{\Psi}(s,g) =\int_{F \backslash \adele_F} e_F(x)
\sum_{\gamma \in B(F) \backslash \PGL_2(F)} f_{\Psi}(\gamma n(x) g)
\\ = \int_{\adele_F} f_{\Psi}(w n(x) g)  e_F(x) dx.\end{multline}
Thus
\begin{multline}W_{\Psi}(s,g) = |\det(g)|^s\int_{x \in \adele_F} \int_{t \in \adele_F^{\times}}
\Psi((t,0) n(x) g) d^{\times}t |t|^{2s} e_F(x)
\\ = |\det(g)|^s \int_{t\in\adele_F, x \in \adele_F^{\times}}
\Psi((t,tx) g) |t|^{2s}e_F(x) dx d^{\times}t.
\end{multline}
The claimed conclusion follows upon substituting
$g = a(y)$, together with some routine computations.
\qed

\begin{rem} \label{rem:whiteis}
Remark that $W_{\Psi}(s,g)$ belongs to the
Whittaker model of a certain induced representation
of $\PGL_2(\adele_F)$, namely
that representation $\pi(s)$ induced
from the character $a(y) \mapsto |y|^{s-1/2}$ of the maximal torus (unitary
induction, so $\pi(s)$ is tempered for $\Re(s) =1/2$).  This representation
is the tensor product of local representations $\pi_v(s)$, analogously defined;
these local representations are irreducible and generic for all $s$.

Thus (\ref{eq:foureis}) determines
$W_{\Psi}$ uniquely (the theory of the Kirillov model).
Similarly the condition $W_{\Psi_v, s}(1)=1$ uniquely
determines the (spherical) vector $W_{\Psi_v,s}$.

We finally remark that $W_{\Psi_v, s}$, as $\Psi_v$ ranges
over all Schwarz-Bruhat functions on $F_v^2$ if $v$ is nonarchimedean,
or over all standard functions if $v$ is archimedean, exhausts the Whittaker model
of $\pi(s)$.  Indeed, the set of such functions $W_{\Psi_v, s}$ 
is a subspace of the Whittaker model of $\pi(s)$ that is stable under the action of the Hecke algebra
of $\PGL_2(F_v)$; this action is irreducible, whence the result.
\end{rem}

We recall that $\diff$ denotes the different (Sec. \ref{subsec:NF}) and we denote by $\zeta_{F,v}(s)$ or simply  $\zeta_v(s)$ the local 
factor of the Dedekind $\zeta$-function of $F$ at the place $v$. 

\begin{cor} \label{cor:unramcase}
Suppose $v$ is nonarchimedean, and $\Psi_v$ the characteristic function
of $\order_v^2$. Then  $W_v(a(y))$ satisfies
 \begin{equation}\label{wvf}\int_{F_v^{\times}} W_v(a(y))
|y|^{s'} d^{\times}y = q_v^{d_v(1+s'-s)} \zeta_v(s+ s')
\zeta_v(1-s+s'),\end{equation} with $d_v = v(\diff)$. Note that this specifies $W_v$,
because it is $K_v$-invariant. 
\end{cor}

In particular, for
each finite $v$ with $v(\diff) = 0$, the function $W_v(g)$ is the
unique spherical Whittaker function on $\GL_2(F_v)$ with Hecke
eigenvalue $q_v^{s} + q_v^{1-s}$, and with $W_v(1) = 1$. 

As is evident from (\ref{eq:four}),
the Eisenstein series themselves are not bounded.
They belong to $L^{2-\varepsilon}$, but not $L^2$.
To avoid some difficulties with growth,
we shall use wave-packets of Eisenstein series.
We now turn to their analysis.

\subsection{Regularization of Eisenstein series on $\PGL_2$.}
\label{sec:eisreg2}
Our aim in this section is to show that an appropriate ``wave packet''
of the Eisenstein series $E_{\Psi}(g,s)$ constructed in the previous section
lies in $L^{\infty}$.

Note that, in Example \ref{qexam} above
$E^*(s,z)$ differs from the usual unitary Eisenstein series
by a factor $\xi(2s)$. This factor ensures that $E^*(s,z)$ is holomorphic,
but this causes an inconvenience at $s=1/2$, which
will manifest itself in our construction of bounded wave-packets. Recall that this pole
can be interpreted rather naturally: see footnote on p. \pageref{extrasing}.

Let $\kappa > 0$, and let $\mathcal{H}(\kappa)$
be the family of functions
 holomorphic in an open neighbourhood of the strip $- \kappa \leq \Re(s) \leq
1+ \kappa$, with
rapid polynomial decay in vertical strips (i.e. $\sup_{t\in \mathbb{R}} (1+|t|)^N |h(\sigma+it)|$ is bounded, for each $N$, by a continuous function of $\sigma$)
and satisfying $h(0) = h(\frac{1}{2}) = h(1) = 0$. 
For each $N \in \mathbb{Z}$
we have a norm $\|\cdot\|_{N}$
on $\mathcal{H}(\kappa)$ defined via:
\begin{equation} \label{eq:normdef} \|h\|_N = \int_{-\infty}^{\infty} \left(|h(1+\kappa+it)|+ |h(-\kappa +
it)|\right) (1+|t|)^{N} dt.\end{equation}
\begin{lem} \label{eislinfty}
Let $h \in \mathcal{H}(\kappa)$, and
set $E_{h,\Psi}(g) = \int_{\Re(s) = 1+\kappa} h(s) E_{\Psi}(g,s) ds$.
Then:
$$\|E_{h,\Psi}(g)\|_{L^{\infty}} \ll_{\Psi, \kappa,F}
\|h\|_{0}$$
\end{lem}
\proof
In the notation of (\ref{eq:trunceis})
\begin{multline}E_{h,\Psi}(g) = \int_{\Re(s) = 1+\kappa}
\bar{E}_{\Psi}(s,g) h(s) ds+ \int_{\Re(s) = 1+\kappa}
h(s) f_{\Psi}(s,g)ds \\ + \int_{\Re(s) = 1+\kappa} h(s)
f_{\widehat{\Psi}}(1-s,g)ds.\end{multline}

Fix $T > 0$ so that $\Sieg(T)$ surjects onto $\quot$ (see Section \ref{sec:redtheory} for definitions).  We will bound each
term on the right-hand side of the above equation for $g \in \Sieg(T)$. 

By Lem. \ref{lem:eisbound}, the first term
on the right-hand side is $O_{\Psi,\kappa}(\|h\|_0)$.
By Lem. \ref{lem:fbound}, the function $f_{\widehat{\Psi}}(1-s,g)$
is uniformly bounded above in the region
$\Re(s) = 1+\kappa, g \in \Sieg(T)$; thus
the third term on the right-hand side
is also $O_{\Psi,\kappa}(\|h\|_0)$.

As for the second term, we shift contours
to the line $\Re(s) = -\kappa$. The
shift of contours is justified by the rapid
decay of $h(s)$ along vertical lines
and Lem. \ref{lem:fbound}. Moreover,
since $h(0) = h(1/2) = 0$,
the function  $s \mapsto h(s) f_{\Psi}(s,g)$
has no poles in between the contours.

Applying Lem. \ref{lem:fbound} one more
time to control the contour integral
along $\Re(s) = -\kappa$, we conclude.
\qed

\begin{rem} \label{rem:schwarz}
Suppose $\Psi = \prod_{v} \Psi_v$, and the $\Psi_v$ are regarded
as fixed for $v$ finite. Put $\Psi_f = \prod_{v \, \mathrm{finite}} \Psi_v$,
a Schwarz function on $\adele_{F,f}^2$, and $\Psi_{\infty} = \prod_{v \, \mathrm{infinite}} \Psi_v$.
Then the above argument gives the slightly more explicit bound
\begin{equation} \label{sahana} \|E_{h,\Psi}\|_{L^{\infty}} \ll_{\kappa, \Psi_f} \|h\|_{0} \mathcal{S}(\Psi_{\infty})
\end{equation} where $\mathcal{S}$ is a Schwarz norm on $F_{\infty}^2$. 
This follows by explicating the above argument, taking into 
account the last sentence of the proof of Lem. \ref{lem:fbound}. 
Indeed, one obtains even the corresponding bound for Sobolev norms, namely
\begin{equation} \label{eq:sobeisenstein} S_{\infty,d,\beta} (E_{h,\Psi}) \ll_{\kappa, \Psi_f} \|h\|_0 \mathcal{S}(\Psi_{\infty})\end{equation}
for an appropriate Schwarz norm of $F_{\infty}^2$.  One deduces this from (\ref{sahana})
upon noting that, if $\mathcal{D}$ belongs to the universal enveloping algebra of $\SL_2(F_{\infty})$,
then, by (\ref{eq:twistedequivariant}), $\mathcal{D} E_{\Psi}(s,g) = E_{\mathcal{D} \Psi}(s,g)$,
so also $\mathcal{D} E_{h, \Psi} = E_{h, \mathcal{D} \Psi}$. It is then easy to check
that a Schwarz norm of $\mathcal{D} \Psi_{\infty}$ is bounded by a Schwarz norm of $\Psi_{\infty}$. 
\end{rem}

\subsection{Regularization of Eisenstein series on $\PGL_2 \times
\PGL_2$.} \label{theartoffugue} In this section we carry out the
analogue of Lem. \ref{eislinfty} in the context of $\PGL_2 \times
\PGL_2$ (this amounts to regularizing the rank $2$ Eisenstein
series on $\PGL_2 \times \PGL_2$). 
%In this case we will not be able
%to construct a wave-packet that is bounded in $L^{\infty}$, but we will
%be able to construct one that is in $L^{4+\varepsilon}$. 

To ease the reader's path, we briefly mention what the point of this section
is in classical notation: Suppose $h(s_1, s_2)$ is holomorphic in two
variables inside the square $|\Re(s_1)| + |\Re(s_2)| \leq 1/2+\kappa$, and, moreover,
$h(s_1, s_2)$ has zeroes along the six planes defined by any of the linear constraints $s_1 = 0, s_1=1/2, s_1=-1/2, s_2 =0, s_2 = -1/2,s_2=1/2$.  

Define the wave-packet $E_h(z_1, z_2)$ on $
\SL_2(\Z) \backslash \mathbb{H} \times \SL_2(\Z) \backslash \mathbb{H}$
via $$E_h(z_1,z_2) = \int_{t,t'\in \mathbb{R}}
h(i t_1, it_2)  E^*(1/2+ i t, z_1)  E^*(1/2+i t', z_2)
dt dt'.$$ 
Here $E^*$ is as in Example \ref{qexam}. 
We shall show -- under mild decay conditions on $h$ --that  $E_h(z_1, z_2)$ is majorized,
on the product of two fundamental regions, by $A(y_1, y_2) :=\frac{\sqrt{y_1
y_2}}{y_1^{1/2+\kappa} + y_2^{1/2+\kappa}}$. Since
$\int_{y_1 \geq 1, y_2 \geq 1} A(y_1, y_2)^4 \frac{dy_1 dy_2}{y_1^2y_2^2}$
is finite, $E_h$ lies in $L^4$, and even in $L^{4+\epsilon}$ for $\epsilon$
small.

As the reader may verify at this point,
the majorization is little more than an exercise in complex integration,
using the
fact that the large contribution to the Eisenstein series
comes from the constant term.

We will need to repeatedly shift contours in
the setting of a function of two complex variables. To clarify matters,
we state the following Lemma, which we will use repeatedly {\em without
explicitly invoking it.}

\begin{lem} \label{lem:contourshift}
Suppose $U \subset \mathbb{R}^2$ is an open domain
and $f(z_1, z_2)$ a holomorphic function on the complex domain
$\{(z_1,z_2) \in \C^2: (\Re(z_1), \Re(z_2)) \in U\}$.
Suppose moreover that there is,
 a continuous function
$M: U \rightarrow \mathbb{R}$ such that
\begin{equation} \label{eq:decaycond}
\sup_{(t_1, t_2) \in \mathbb{R}^2}
|f(\sigma_1 + i t_1, \sigma_2 + i t_2)| (1+|t_1|+|t_2|)^3
\leq M(\sigma_1, \sigma_2).\end{equation}

Then the function
\begin{equation} \label{eq:square}
(\sigma_1, \sigma_2) \mapsto \int_{\Re(z_1) = \sigma_1, \Re(z_2) =
\sigma_2} f(z_1, z_2) dz_1 dz_2\end{equation}
is locally constant on $U$.
\end{lem}

We omit the easy proof.
%\proof
%Let $(\sigma_1, \sigma_2) \in U$
%and let $\epsilon > 0$ be such that
%$\sigma_1 \times [\sigma_2-\epsilon, \sigma_2+\epsilon] \subset U$.
%By (\ref{eq:decaycond}),
%the integral $\int_{\Re(z_1) = \sigma_1} f(z_1, z_2)$
%is absolutely convergent, uniformly for $z_2$ in compacta.
%In particular,
%$g(z_2) = \int_{\Re(z_1) = \sigma_1} f(z_1, z_2) dz_1$
%defines a holomorphic function
%of $z_2$ in the strip $|\Re(z_2) - \sigma_2| \leq \epsilon$.
%Moreover, the estimate (\ref{eq:decaycond}) implies that,
%for each $N > 0$, there is a constant $c(N)$ such that
%$|g(z_2)| \leq
%c(N) (1+|\Im(z_2)|)^{-N}$, whenever
%$|\Re(z) - \sigma_2| \leq \epsilon$.
%
%Now, by shifting of contours, we conclude
%that $\int_{\Re(z_2) = \sigma_2'} g(z) dz_2$
%is constant in the range $|\sigma_2' - \sigma_2| < \epsilon$.
%
%We conclude (speaking somewhat imprecisely) that the function defined by (\ref{eq:square})
%is constant when $\Re(z_1)$ is fixed and
%$\Re(z_2)$ is ``perturbed slightly.''
%From this we deduce the claim.
%\qed

We will now introduce a family of normed spaces
$\mathcal{H}^{(2)}(\kappa)$. In fact, the spaces themselves
are independent of $\kappa$, but the norm depends on $\kappa$. 
These are spaces of holomorphic functions in two variables $z_1, z_2$;
and the norm, roughly speaking, controls the behavior of $h$ when
the real parts of $(z_1, z_2)$ lie in the square $|\Re(z_1)| + |\Re(z_2)| \leq 1/2+\kappa$.

\begin{defn} \label{def:hdef} Let $0 < \kappa < 1$. 
Let $\mathcal{H}^{(2)}(\kappa)$
be the family of functions $h(z_1, z_2)$ in two complex variables,
holomorphic in a neighbourhood of $(0,0)$, and satisfying:
\begin{enumerate}
\item \label{cond:zero} Write $h' = \frac{h(z_1, z_2)}{z_1 z_2 (1/4-z_1^2)
(1/4-z_2^2)}$.
Then $h'$, originally a meromorphic function
in a neigbourhood of $0$, extends to a holomorphic function in the strip
$\{z_1: |\Re(z_1)| \leq 2\} \times \{z_2: |\Re(z_2)| \leq 2\}$.

\item Growth condition: for every $N \geq 0$,
$$ \sup_{(\sigma, \sigma') \in [-2,2]^2}
\sup_{t, t' \in \mathbb{R}^2} (1+|t|+|t'|)^N h(\sigma + it, \sigma' + it)  < \infty$$
\end{enumerate}
For each $N \in \mathbb{Z}$ we introduce a norm
on $H^{(2)}(\kappa)$ via:
\begin{equation}\begin{aligned}\|h\|_N = \int_{(t,t') \in \mathbb{R}^2}
\sum_{\epsilon_1, \epsilon_2 \in \{\pm 1 \}}
(1+|t|+|t'|)^N \\
\left(|h'(\epsilon_1(1/2 +\kappa) + it, it')|  +
|h'(it, \epsilon_2 (1/2+\kappa)+it')|\right) dtdt'.\end{aligned}\end{equation}
\end{defn}

%Note that, by a simple maximal modulus argument,
%the norm $\|h\|_N$ controls the integral
%of $h$ over any contour $\Re(t) = \sigma, \Re(t') = \sigma'$
%such that $|\sigma| + |\sigma'| \leq 1/2+\kappa$.
%\sum_{(\epsilon_1,\epsilon_2) \in \{\pm1\}^2}
%\left(
%|\tilde{h}'(\epsilon_1 (1/2 +\kappa) + it,
%\epsilon_2(1/2 +\kappa) + it')| \right) (1+|t|+|t'|)^N.$$

\begin{lem} \label{eislinfty2}
For $h \in H^{(2)}(\kappa)$, put
$$E_{h,\Psi,\Psi'}(g_1, g_2) = \int_{\Re(t) = 0}
\int_{\Re(t') = 0} h(t,t') E_{\Psi}(g_1,
1/2+t) E_{\Psi'}(g_2, 1/2 + t').$$
Then
$$E_{h,\Psi,\Psi'}(x_1, x_2) \ll_{\Psi, \Psi'} \|h\|_0 \frac{\height(x_1)^{1/2}
\height(x_2)^{1/2}}{\height(x_1)^{1/2+\kappa} +
\height(x_2)^{1/2+\kappa}}.$$
\end{lem}

\proof
We may assume that $\|h\|_0 = 1$.
Let notations be as established prior to Lem. \ref{lem:eisbound}. We will proceed
as in Lem. \ref{eislinfty}, expanding $E_{\Psi}$ via (\ref{eq:trunceis}).

It will suffice to give an upper
bound, in absolute value, for each of:
\begin{eqnarray}
I_0(g_1, g_2) = \int_{t,t'} h(t,t') \bar{E}_{\Psi_1}(g_1, 1/2 + t)
\bar{E}_{\Psi_2}(g_2, 1/2 + t') dt dt' \\
I_1(g_1, g_2) = \int_{t,t'} h(t,t')
\bar{E}_{\Psi_1}(g_1, 1/2 + t) f_{\Psi_2}(g_2, 1/2 \pm t') dt dt' \\
I_2(g_1, g_2) = \int_{t,t'} h(t,t')
f_{\Psi_1}(g_1, 1/2 \pm t) \bar{E}_{\Psi_2}(g_2, 1/2 + t') dt dt' \\
I_3(g_1, g_2) = \int_{t,t'} h(t,t')
 f_{\Psi_1}(g_1, 1/2 \pm t) f_{\Psi_2}(g_2, 1/2\pm t') dt dt',
\end{eqnarray}
whenever $\Psi_1, \Psi_2$ are Schwarz functions on $\adele_F^2$,
and in each case the contour of integration
is the surface $\Re(t) = \Re(t') = 0$. Moreover, in view of
condition (\ref{cond:zero}) in Def. \ref{def:hdef}, each integrand
extends to a holomorphic function of $(t,t')$ in the region
$|\Re(t)| \leq 1/2+\kappa, |\Re(t')| \leq 1/2+\kappa$. 

The bound $|I_0| \ll \height(g_1)^{-N} \height(g_2)^{-N}$ follows
from Lem. \ref{lem:eisbound}, whereas
the bounds $|I_1| \ll \height(g_1)^{-N} \height(g_2)^{-\kappa}$
and $|I_2 |\ll \height(g_2)^{-N} \height(g_1)^{-\kappa}$
follow from moving the $t'$ (in the case of $I_1$)
integral to the contour $\Re(t') = \pm (1/2+\kappa)$,
applying Lem. \ref{lem:eisbound} and Lem. \ref{lem:fbound}.

We now turn to $I_3$.
We shall consider the case where both signs are $+$,
the other cases being similar with appropriate interchanges
of sign.
Thus set \begin{equation}\begin{aligned} Z(t,t') =
h(t,t') f_{\Psi_1}(g_1, 1/2 + t) f_{\Psi_2}(g_2, 1/2 + t') \\
 = h'(t,t') t (1/4-t^2) f_{\Psi_1}(g_1, 1/2 + t)
t' (1/4 -t'^2) f_{\Psi_2}(g_2, 1/2 +  t') \end{aligned}\end{equation}
In view of Lem. \ref{lem:fbound}, the function $Z(t,t')$ satisfies the conditions
for $f$ in Lem. \ref{lem:contourshift}.
We apply Lem. \ref{lem:contourshift} to shift the contour to $\Re(t) = -1/2-\kappa, \Re(t') = 0$.
Now Lem. \ref{lem:fbound} implies that
$|\int_{\Re(t) = -1/2-\kappa, \Re(t') = 0}
Z(t,t')| \ll \height(g_2)^{1/2} \height(g_1)^{-\kappa}$.
A similar bound holds with $(g_1, g_2)$ interchanged,
so in fact we have the stronger bound
$|Z(t,t')| \ll \min(\height(g_2)^{1/2} \height(g_1)^{-\kappa},
\height(g_1)^{1/2} \height(g_2)^{-\kappa})$. This
may also be written $|Z(t,t')| \ll \frac{\height(g_1)^{1/2} \height(g_2)^{1/2}}{\height(g_1)^{1/2+\kappa}+
\height(g_2)^{1/2+\kappa}}$. 

Similar considerations apply to the terms in $I_3$ corresponding to other choices of sign, so we conclude that 
$|I_3| \ll \frac{\height(g_1)^{1/2}\height(g_2)^{1/2}}
{\height(g_1)^{1/2+\kappa} + \height(g_2)^{1/2 +\kappa}}$.
\qed
\begin{lem} \label{eislinfty3}
Let notations be as in the previous Lemma.
For any $p<\frac{4}{1-2\kappa}$, any $d,\beta > 0$, there exists $N$ such that
$$S_{p,d,\beta}(E_{h,\Psi}) \ll_{\Psi, \kappa,p,\beta}
\|h\|_{N}$$
\end{lem}
\proof
Indeed, we note that
$$\int_{y_1, y_2 \geq 1} \left(\frac{\sqrt{y_1y_2}}{y_1^{1/2+\kappa} +
y_2^{1/2+\kappa}}\right)^p \frac{dy_1 dy_2}{y_1^2 y_2^2}  <\infty$$
whenever $p < \frac{4}{1-2 \kappa}$.
 We apply the previous Lemma
and reduction theory to conclude.
%\footnote{Note that, in view
%of the relation $h . E_{\Psi}(s,\cdot) = E_{h . \Psi}(s, \cdot) |\det(h)|^s$, if $\mathcal{D}$ is any monomial
%in $\Lie(\PGL_2(F_{\infty}))$, we may express $ \mathcal{D} \cdot E_{\Psi}(s,g)$ 
%as a linear combination of $E_{\Psi_j}(s,g)$, for various $\Psi_j$ that are obtained
%by applying monomials in $\Lie(\PGL_2(F_{\infty}))$. )}
\qed

\section{Background on integral representations of $L$-functions.} \label{sec:rs}

The purpose of this section is as follows. The geometric method we have explained
in the text yields upper bounds for certain periods; to obtain subconvexity, we need to know that
$L$-functions can be expressed in terms of these periods.  This is the whole point
of the theory of integral representations of $L$-functions; however, we cannot quite
simply quote from that theory, as we often need e.g. some analytic control
on the choice of test vector for which there is no readily available reference. 

On occasion we have only sketched proofs in this section, as they
amount to simple explications of standard techniques such as the Rankin-Selberg method,
and moreover they are in some sense irrelevant to the main point of this paper
(which is to bound {\em periods}, not $L$-functions!) 

%{\em In this section we shall equip $\quot$ with the measure that arises
%as the quotient of the measure on $\PGL_2(\adele_F)$.} This makes unfolding \marginpar{Revised, unfortunately and annoyingly.} much more convenient. 
%This differs from the probability measure on $\quot$ by a constant depending only on $F$. 
%The only results from this section that are used in the text concern inequalities
%with implicit constants that depend on $F$ in any case. Therefore, the switch of measure
%should not cause any problem. For clarity, we have also commented on this switch of measure
%at those points in the text where we use results from the present section. 

\subsection{ Cuspidal triple product $L$-functions.}
\begin{hyp} \label{prop:intrep}
Let $\pi_2$ and $\pi_3$ be {\em fixed} automorphic cuspidal
representations of $\PGL_2$ over
$F$. Let $\pi_1$ be an automorphic cuspidal representation, 
whose finite conductor is a prime ideal, prime to the finite conductors
of $\pi_2$ and $\pi_3$. Suppose that $\pi_{1,\infty}$ (the representation of $\PGL_2(F_{\infty})$ underlying $\pi_1$) is restricted to a bounded set; let $\varphi_1$ be the new vector in $\pi_1$.

Then there exists finite collections of vectors $\mathcal{F}_2 \subset \pi_2, \mathcal{F}_3 \subset \pi_3$
so that, for any such $\pi_1$, there exist $\varphi_j \in \mathcal{F}_j \ (j=2,3)$ with \begin{equation}\label{eq:scream} \frac{L(\frac{1}{2}, \pi_1 \otimes \pi_2 \otimes \pi_3)}{
\left|\int_{\quot} \varphi_2(x \f) \varphi_3(x) \varphi_1(x) dx\right|^2}
\ll_{\epsilon, F, \pi_{1,\infty} }  N(\cond)^{1+\epsilon} \end{equation}
\end{hyp}

Note that no claim is made about the dependence of the constants in \eqref{eq:scream} on $\pi_2, \pi_3$ or the bounded set containing $\pi_{1,\infty}$; 
presumably with enough effort one could obtain polynomial dependence on the conductors. 

The proof of Hypothesis \ref{prop:intrep} should be, we believe, 
an elaborate but routine computation of certain $p$-adic integrals;
this has not carried out, but we expect it to be valid.

In the case when $F=\Q$ and $\pi_1, \pi_2, \pi_3$ holomorphic
Hypothesis \ref{prop:intrep} may follow
(in a slightly modified form, replacing $\PGL_2$ by a division algebra)
from the work of B\"ocherer and Schulze-Pillot. 
In any case, there exist good heuristic reasons to believe the Hypothesis:
based on a computation of the size of the relevant family, or alternately it is true
if one of the $\pi_j$ is Eisenstein. 

%Unfortunately their
%paper is notationally very complicated, and their result
%was later corrected in their joint paper with Sarnak;
%nevertheless, based on a brief study, I believe they establish
%essentially Hypothesis \ref{prop:intrep}.  I cannot guarantee this, however.

%
%There at least two very good reasons to believe Hypothesis \ref{prop:intrep}
%in the general case.
%There exists a fairly convincing
%heuristic argument,
%and it is true if one of the $\pi_j$ is Eisenstein, as is shown in the following section.
\subsection{Rankin-Selberg convolutions.}

\subsubsection{The Rankin-Selberg integral representation.}
Let $\pi_1, \pi_2$ be two automorphic representations,
with $\pi_2$ cuspidal. 
%Throughout this section we will equip
%$N(\adele_F)$ and $\adele_F$with the Haar measure that assigns mass $1$
%to the quotient $N(\adele_F) / N(F)$ (this choice is not very important, however, as
%any other choice of measure would lead to an absolute constant that is irrelevant to our results). 

Let $\Psi_v$ be a Schwarz-Bruhat function on $F_v^2$
such that, for almost all $v$, $\Psi_v$ is the
characteristic function of $\order_{F_v}^2$.
Put $\Psi = \prod_v \Psi_v$, a Schwarz function on $\adele_F^2$. 
Let $\varphi_j $ belong to the space of $\pi_j$ for $j=1,2$ and put
\begin{equation} I(\varphi_1, \varphi_2, \Psi,s) = \int_{\quot}
\varphi_1(g) \varphi_2(g) E_{\Psi}(s, g)  dg
\end{equation}

Unwinding, we see that for $\Re(s) > 1$:
\begin{multline}
I(\varphi_1, \varphi_2, \Psi, s) 
 =  c_F \int_{B(F) \backslash
\PGL_2(\adele_F)} f_{\Psi}(s,g) \varphi_1(g) \varphi_2(g)  \\ =
c_F \int_{B(F) \backslash \PGL_2(\adele_F)} f_{\Psi}(s,g)
\left(\int_{n \in N(F) \backslash N(\adele_F)} \varphi_1(ng)
\varphi_2(ng) dn \right) dg
\end{multline} 
Here the constant $c_F$ arises from change of measure: the measure on $\quot$
is the $\PGL_2(\adele_F)$-invariant probability measure, which is not the same
as the quotient measure from $\PGL_2(\adele_F)$. Note that $c_F$ will be unimportant in our arguments, as it depends only on $F$ and we are only interested in bounds.

Put $W_{1}(g) = \int_{F \backslash \adele_F}
 \varphi_1(n(x)g) e_F(x) dx$,
and define $W_2$ similarly but with $e_F$ replaced by $\overline{e_F}$.
Recall that our normalizations are so that the volume of $\adele_F/F$ is $1$. 
Fourier inversion shows that $\varphi_i(g) = \sum_{\alpha \in F^{\times}}
W_{i}(a(\alpha) g)$ if $\varphi_i$ is cuspidal. Thus, as long as one of $\varphi_1, \varphi_2$
is cuspidal, we see that:
\begin{multline}
I(\varphi_1, \varphi_2, \Psi,s) =  c_F
\int_{B(F) \backslash \PGL_2(\adele_F)} 
f_{\Psi}(s,g) \left(  \sum_{\alpha \in F^{\times}} W_1(a(\alpha) g)
W_2(a(\alpha) g) \right) \\ =  c_F
\int_{N(\adele_F) \backslash
\PGL_2(\adele_F) } W_1(g) W_2(g)  f_{\Psi}(s,g) dg \end{multline}

If $\varphi_1, \varphi_2$ are pure tensors,
then there is a corresponding product decomposition
$W_1 = \prod_{v} W_{1,v}, W_2 = \prod_{v} W_{2,v}$,
where $W_{j,v}$ belongs to the local Whittaker model
of $\pi_{j,v}$, a representation of $\PGL_2(F_v)$.
In that case,
\begin{equation} \label{eq:productformel} I(\varphi_1, \varphi_2, \Psi,s) =    c_F \prod_{v}
I_v(W_{1,v}, W_{2,v}, \Psi_v,s),\end{equation}
%where $c$ is a constant that depends only on the normalization of measures [[EQUALS ONE??]] and
%\footnote{Recall that we have normalized the measure on $N(\adele_F)$ so that the measure of the quotient
%$N(\adele_F)/N(F)$ is one; therefore a constant is introduced when we decompose
%this measure into measures on $N(F_v)$, which we still take to be the standard measures
%prescribed in Sec. \ref{subsec:measures}.}, and
where
\begin{multline}\label{eq:localdef}I_v(W_{1,v}, W_{2,v}, \Psi_v,s) = 
\\  \int_{N(F_v) \backslash \PGL_2(F_v)}
W_1(g_v) W_2(g_v) \left( |\det(g_v)|_v^{s} \int_{t \in F_v^{\times}}
\Psi((0,t) g_v) |t|^{2s} d^{\times} t\right) 
\end{multline}
We note that the bracketed quantity, defined {\em a priori} for 
$g_v \in \GL_2(F_v)$, descends to $\PGL_2(F_v)$.

\begin{lem} \label{lem:dooku}
Suppose $W_{1,v}, W_{2,v}$ the new vectors associated to spherical 
representations $\pi_{1,v}, \pi_{2,v}$,  $\Psi_v$ is the characteristic
function of $\order_v^2$, and $e_{F_v}$ is unramified.  Then $I_v(W_{1,v}, W_{2,v}, \Psi_v) = L_v(s,\pi_{1,v} \otimes \pi_{2,v})$.

If $W_{1,v}, W_{2,v}$ are nonzero and $\PGL_2(\order_v)$-invariant, $\Psi_v$ as above,
but $e_{F_v}$ is possibly ramified, then $I_v(W_{1,v}, W_{2,v}, \Psi_v)  = a q_v^{ks} L_v(s, \pi_{1,v} \times \pi_{2,v})$ where $k \in \Z$ is so that $e_{F_v}$ is trivial on $\varpi_v^{-k} \order_v$
but not on $\varpi_v^{-k-1} \order_v$. Moreover, $a=1$ if $W_{1,v}(\varpi_v^{-k}) = W_{2,v}(\varpi_v^{-k})=1$.

Suppose $W_{1,v}, W_{2,v}$ are the new vectors associated to $\pi_{1,v}$ spherical
and $\pi_{2,v}$ a Steinberg representation, and that $e_{F_v}$ is unramified. 
Then, with $\Psi_v$ the characteristic
function of $\order_v^2$, we have:
 $$I_v(\pi_{1,v}\left( \begin{array}{cc} 1 & 0 \\ 0 & \varpi_v \end{array}\right)  W_{1,v}, W_{2,v}, \Psi_v) = \pm\frac{q_v^s}{q_v+1} L(s, \pi_{1,v} \otimes \pi_{2,v}).$$
\end{lem}
\proof
See \cite[Thm 15.9]{jacquet} for the first assertion. The second assertion is an easy consequence. See Sec. \ref{subsec:localrankinselberg} for the final assertion. \qed 

Applying the Iwasawa decomposition to (\ref{eq:localdef}) yields the equivalent
\begin{multline} \label{eq:localdefiwasawa}I_v(W_{1,v}, W_{2,v}, \Psi_v,s) = \int_{y \in F_v^{\times}, k \in K_v}
W_1(a(y) k) W_2(a(y) k) |y|^{s-1} d^{\times}y  \\ \cdot  \left(
\int_{t \in F_v^{\times}} \Psi((0,t) k) |t|^{2s} d^{\times} t \right)\end{multline}

\subsubsection{Topologizing the space of local representations.} \label{topologizing}
The results in \cite{jacquet} provide ``good'' test vectors for the functionals $I_v$
when the local representations $\pi_{1,v}, \pi_{2,v}$ are fixed. On the other hand,
we will need such results with some mild uniformity in $\pi_{1,v}$. 
One can certainly extract the stronger results from the proofs in \cite{jacquet}.
For now, we will proceed by deducing the results ``by continuity'';
to do this, we will need to define the topology on the space of possible $\pi_{1,v}$. 
The considerations that follow are not really very crucial; it would be better
simply to explicate the implicit dependences in \cite{jacquet}.

Let $\mathscr{S}$ be a finite set of irreducible (continuous) representations
of $K_v = \PGL_2(\order_v)$. For any representation $W$ of $K_v$, we denote by
$W^{\mathscr{S}}$ that subspace of $W$ consisting of vectors
whose $K_v$-span contains only irreducibles that belong to $\mathscr{S}$. 
We shall say that elements of $W^{\mathscr{S}}$ are {\em of type $\mathscr{S}$}. 

Let $\mathscr{G}_v$ be the set
of isomorphism classes of generic irreducible representations of $\PGL_2(F_v)$. 
If $\pi$ is a generic irreducible representation which 
is a discrete series or supercuspidal, we shall define it to be isolated. 
Otherwise,
$\pi$ is induced from two quasicharacters $\mu, \nu: F_v^{\times} \rightarrow \C$. For $s \in \C$, let $(\pi(s), V(s))
$ be the representation induced
from the quasicharacters $\mu |\cdot|_v^s, \nu |\cdot|_v^{-s}$.
Then $\pi(s)$ is generic for all $s \in \C$ and irreducible in a neighbourhood of $0$.
We shall topologize $\mathscr{G}_v$ in such a way that
sets of the form $[\pi(s)]$, for $|s| < \varepsilon$ form a basis. 

If $E \subset \mathscr{G}_v$ is a closed subset that is bounded (when considered
as a subset of the set of isomorphism classes of irreducible admissible representations, 
and bounded in the sense of Sec. \ref{subsec:cuspformsbdd}),
then $E$ is compact, as one checks by direct verification. 

For each $\pi \in \mathscr{G}$, 
we have a Whittaker model $\mathscr{W}(\pi)$.
 Consider a function $\pi \mapsto W_{\pi}$, that assigns to each $\pi \in \mathscr{G}_v$
an element $W_{\pi}$ of its Whittaker model.  
We shall say that such an assignment $\pi \mapsto W_{\pi}$ is continuous
if there exists a neighbourhood of each $\pi$,
which we may assume to be of the form, 
$\{\pi(s): |s| < \varepsilon\}$,
and a set $\mathscr{S}$ of irreducible representations of $K_v$ 
so that:
\begin{enumerate}
\item \label{ktype} $W_{\pi(s)}$ is of type $\mathscr{S}$, for each $|s| < \varepsilon$
\item \label{cty} The assignment $s \mapsto W_{\pi(s)}(g)$
is continuous for each $g \in \PGL_2(F_v)$, uniformly for $g$ in any fixed compact. 
\end{enumerate}

It can be verified that if $W_{0}$ is an element of the Whittaker model of $\pi_0$, there exists a continuous assignment
$\pi \mapsto W_{\pi}$ in a neighbourhood of $\pi_0$ which has the value $W_0$ at $\pi_0$. 

The requirement (\ref{cty}) is not very strong, as it does not impose
any uniformity on all of $\PGL_2(F_v)$. 
However, in every context we shall consider, the necessary uniformity in $g$
is automatic. Let us sketch how one can prove such results. Assume that $v$ is finite;
the infinite case is similar although more technically involved.  One first observes that
if $\pi \mapsto W_{\pi}$ is a continuous assignment on some open set,  then,
for a fixed character $\chi_v$ of $F_v^{\times}$, the quotient
$\frac{\int_{y \in F_v^{\times}} W_{\pi}(a(y)) \chi_v(y) |y|^{s-1/2} d^{\times}y }{L_v(s, \pi_v \otimes \chi_v)}$
is a polynomial of the form $\sum_{k=-N}^{N} c_k q_v^{ks}$; moreover, the degree $N$
is locally bounded as $\pi$ varies, and all the coefficients $c_k$ can be taken to depend
continuously on $\pi$. To verify the local boundedness of the degree -- which requires only
property (\ref{ktype}) above -- one just notes
that there is (locally) a fixed $M$ such that $W_{\pi}(a(y))$ vanishes for $|y|_v > M$;
this, together with the functional equation, gives the local boundedness. To see that the
coefficients vary continuously, it suffices to check that, for any fixed integer $t$,
the integral $\int_{v(y) = t}  W_{\pi}(a(y)) \chi_v(y) |y|^{s-1/2} d^{\times} y$ varies continuously,
which follows from the definition of continuity for the assignment $\pi \mapsto W_{\pi}$. 
The archimedean case proceeds similarly, but one replaces the role of polynomials in $q_v^{\pm s}$
by functions of the form $c^s P(s)$, where $P$ is a polynomial and $c \in \mathbb{R}$. 

In the next few pages, we will make certain claims regarding
the continuity of various integrals involving $W_{\pi}$, if $\pi \mapsto W_{\pi}$
is a continuous assignment. One can  reduce all the claimed continuity statements (by standard ``Mellin transform'' arguments)
to the result just discussed. We will omit the details. 

%
%why this is true in the case
%where $v$ is finite.  Indeed, it suffices (by (\ref{ktype})) to consider
%the case when $g=a(y)$. Then, there is a fixed $M$ --depending only on the
%set $\mathscr{S}$ in (\ref{ktype}) -- so that $W_{\pi}(a(y))$ vanishes for $|y|_v > M$.
%As for the case where $|y|_v$ is small, one may apply the local functional
%equation for Jacquet-Langlands to control the behavior of $W_{\pi(s)}(a(y))$ uniformly
%in $\pi$ and $y$.   Alternately, one can check all of this by {\em direct computation},
%since we are dealing only with principal series for which the Whittaker functions may be explicitly computed. 

%Let us give one explicit example of this type of reasoning. Suppose that $v$ is a finite place, 
%and $\pi \mapsto W_{\pi}$ is a continuous assignment in some neighbourhood
%$\{\pi(s): |s| < \varepsilon\}$. 
%We claim that $\int_{F_v^{\times}} |W_{\pi}(a(y))|^2 d^{\times}y$ 
%varies continuously.  First, remark that $\int_{y} W_{\pi}(a(y)) \chi_v(y) |y|^{s-1//2}$
%is of the form $P_{\chi_v}(q_v^{s}) L_v(s, \pi_v \times \chi_v)$, for a certain polynomial
%$P_{\chi_v} \in \C[t, t^{-1}]$. Consequently, 
%$$\int_{y} |W_v(a(y))|^2 d^{\times}y = \sum_{\chi_v} \int_{y} |P_{\chi_v}(q_v^{1/2+it})|^2
%|L_v(1/2+it, \pi_v \otimes\chi_v)|^2 dt$$

\subsubsection{Choice of test vectors} 

\begin{lem} \label{lem:localchoice}
Let notation be as above.  
\begin{enumerate}
\item  \label{onee}
The quotient
$$\Xi_v (W_{1,v}, W_{2,v}, \Psi_v, s) :=
\frac{I_v(W_{1,v}, W_{2,v}, \Psi_v, s)}{L_v(s,\pi_{1,v} \otimes \pi_{2,v})}$$
is holomorphic in $s$. If $v$ is nonarchimedean, $\Xi_v$ is a polynomial in $q_v^{\pm s}$;
if $v$ is archimedean and $\Psi_v$ is standard, 
then $\Xi_v |a_v|_v^{2s}$ is a polynomial in $s$. 

\item  \label{twoe} For any fixed $s_0 \in \C$ we may choose
data $(W_{1,v}, W_{2,v}, \Psi_v)$ of the type described in (\ref{onee}) with
$\Xi_v(W_{1,v}, W_{2,v}, \Psi_v, s_0) \neq 0$.  

\item \label{threep}
If $\pi_{2,v}$ is regarded as fixed, and
$\pi_{1,v}$ remains within a fixed compact subset of $\mathscr{G}_v$ consisting
entirely of unitarizable representations,
then there exists a constant $C$ depending on the compact set so that
one may choose data as in (\ref{twoe}) in such a way that:
\begin{enumerate} \item \label{a}  $\Psi_v$ and $W_{2,v}$ may both be chosen from a finite list of size $\leq C$;
\item  \label{b} $\int_{F_v^{\times}} |W_{1,v}(a(y))|^2 d^{\times}y \leq C$; 
\item \label{c} $|\Xi_v(s_0)| \geq 1$ and, for all $s \in \mathbb{C}$, we have $|\Xi_v(s)| \ll C^{|\Re(s)|} (1+|s|)^{C}$.
\end{enumerate}
\end{enumerate}
\end{lem}

\proof The first two assertions are in \cite{jacquet}.

We will only sketch the last assertion. It can be also be proved directly by 
{\em exhibiting} such data by explicating the arguments of \cite{jacquet}.  In any case, we start by noting:
Given any continuous assignment $\pi_{1,v} \mapsto W_{1,v}, 
\pi_{2,v} \mapsto W_{2,v}$, the function $\Xi_v(W_{1,v}, W_{2,v}, \Psi_v, s)$
and $\int |W_{i,v}(a(y)|^2 d^{\times}y$
all varies continously in $\pi_{1,v}, \pi_{2,v}$.  
This assertion can be deduced by the methods explained in Section \ref{topologizing}. 
Here, when we speak of $\Xi_v(W_{1,v}, W_{2,v}, \Psi_v, s)$ varying continuously, we mean this in the ``strong sense'', i.e the statement
that $\Xi_v$ can be expressed as a polynomial in $q_v^s$ (nonarchimedean case)
or $b^s P(s)$ where $P$ is polynomial (archimedean case), so that all coefficients vary continuously with $\pi_{1,v}, \pi_{2,v}$. 

Now, given fix momentarily $\pi_{1,v}$ and $\pi_{2,v}$ and suppose we have chosen
data $(W_{1,v}, W_{2,v}, \Psi_v)$ as in (\ref{twoe}). Extend $W_{1,v}$ to a continuous assignment $\pi \mapsto W_{\pi}$
in a neighbourhood of $\pi_{1,v}$. By the remarks above, $(W_{\pi}, W_{2,v}, \Psi_v)$
will satisfy (\ref{b}) and (\ref{c}), for a suitable constant $C$, whenever $\pi$ belongs to a sufficiently small neighbourhood of $\pi_{1,v}$. 
Now a compactness argument demonstrates (\ref{threep}). 
%The fact that the degree of $\Xi_v$ in $s$ or $q_v^s$
%is locally bounded may be easily checked directly in the nonarchimedean case,
%and is implicit in the proofs in \cite[\S 17, \S 18]{jacquet} in the archimedean case. 
%The assertion about continuity we now leave to the reader, but note that -- given the assertion on degree just made -- it suffices
%to show it when $\Re(s) \gg 1$ so the integrals are absolutely convergent. 
%\footnote{It is still not quite formal because, as discussed previously, our
%definition of continuity for $\pi \mapsto W_{\pi}$ only provides uniformity for $g$ in a compact subset of $\PGL_2(F_v)$.)}
%
%We leave the verification of this in the general case to the reader, but point
%out that we only apply it, in the present paper in the case where $\pi_{1,v}, \pi_{2,v}$ are constituents
%of a cuspidal representation. In that case, 
%by the functional equation we may assume that $s_0 \geq 1/2$,
%and the Kim-Shahidi bounds
%show that the integral defining $I_v$ is absolutely convergent, in which case
%the result is much easier%In particular, if we extend $W_{1,v}$ and $W_{2,v}$ to continous
%sections $W_{1,v}, W_{2,v}$ in an arbitrary fashion, then $\Psi_v(W_{1,v}(s), W_{2,v}(w), \Psi_v, s_0)$ will be nonvanishing in a neighbourhood of $(0,0)$. This  is precisely the final assertion
%of the Lemma.
 \qed

We emphasize again that (\ref{eq:localdef})
is valid so long as {\em one} of $\pi_1, \pi_2$ is cuspidal.

\begin{lem} \label{lem:l2norm}
Let $\pi$ be an automorphic cuspidal representation
in $L^2(\quot)$, and let $\varphi \in \pi$ be so that
$W_{\varphi} := \int_{F \backslash \adele_F} e_F(x) \varphi(n(x) g)$
factorizes as a product $\prod_{v} W_v(g)$.  Then, for a certain
constant absolute constant $c$
\begin{equation} \int_{\quot} |\varphi(g)|^2 dg = c \mathrm{Res}_{s=1} \Lambda(s, \pi \otimes \tilde{\pi})
\prod_{v} \frac{  \int_{F_v^{\times}} |W_v(a(y))|^2 d^{\times}y}{L_v(s, \pi_v \otimes \tilde{\pi}_v)}
\end{equation}
\end{lem}
\proof This follows by taking the residue of $I(\varphi, \bar{\varphi}, \Psi,s)$ 
at $s=1$. Indeed, this residue equals, up to a constant depending only on the measure
normalization, $\left( \int_{\quot} |\varphi(g)|^2 dg \right)
\left( \int_{\adele_F^2} \Psi(x,y) dx dy\right)$ (see discussion of properties of $E_{\Psi}$
after (\ref{fourierpsidef}). 

On the other hand, by (\ref{eq:productformel}) and (\ref{eq:localdef})
$I(\varphi, \bar{\varphi}, \Psi, s)$ may be written as a product $c_F \prod_{v} I_v(W_{v}, \overline{W_{v}}, \Psi_v, s)$, where each $I_v$ is given by (\ref{eq:localdefiwasawa}).
%Applying the Iwasawa decomposition to (\ref{eq:localdef}), 
%we see that each local factor $I_v(s)$ may be expressed as the integral
%$\int_{y \in F_v^{\times}} \int_{K_v} |W_v(a(y)k)|^2 d^{\times}y |y|^{s-1}
%\int_{t \in F_v^{\times}} \Psi((0,t) k) |t|^{2s} d^{\times}t.$
The integral $\int_{F_v^{\times}} |W_v(a(y) k )|^2 d^{\times}y$
is independent of $k \in K_v$,  so $I_v(W_v, \overline{W_v}, \Psi_v, 1)$
factors as the product of
$\int_{y \in F_v^{\times}} |W_v(a(y))|^2 d^{\times}y $ and
$\int_{t \in F_v^{\times} , k \in K_v} \Psi_v((0
,t) k) |t|^2 dt$.
The latter integral differs from $\int_{F_v^2} \Psi_v(x,y) dxdy$
by a factor that depends only on the normalizations of measure;
moreover, this factor equals $(1-q_v^{-2})^{-1}$ for almost all $v$,
so the product of these factors is convergent. The conclusion easily follows. \qed

We now specialize to the cases of interest.
Fix $\pi_1$. We vary $\pi_2 := \pi$ through
a sequence  of automorphic cuspidal representations with
prime conductor $\cond$, prime to the conductor of $\pi_1$.
In particular, the local constituent of $\pi$ at $\cond$ is a special
representation. We denote by $\pi_{\infty}$ the representation
of $\GL_2(F_{\infty})$ underlying the representation $\pi$. 

\begin{lem} \label{lem:rsone}
 Suppose the archimedean constituent $\pi_{\infty}$ belongs to a bounded subset of $\widehat{\PGL_2(F_{\infty})}$ (in what follows the implicit
constants may depend on this subset) and regard $\pi_1$ as being fixed. 

Let $s_0 \in \C$. 
There exists a fixed finite set $\mathcal{F}$
of Schwarz Bruhat functions 
and a real number\footnote{depending on $\pi_1$ and the choice of bounded subset of $\widehat{\PGL_2(F_{\infty})}$} $C > 0$ so that, for any such $\pi$, 

There exist vectors $\varphi_1 \in \pi_1, \varphi \in \pi$ and $\Psi \in \mathcal{F}$ so that
$$\Phi(s) := \Norm(\cond)^{1-s}
\frac{ I(\f \cdot \varphi_1, \varphi, \Psi, s)}{ \Lambda(s, \pi_1 \otimes \pi)}$$
is holomorphic and satisfies:
\begin{enumerate}
\item  $|\Phi(s_0)| \gg 1$ and $|\Phi(s)| \ll C^{|\Re(s)|} (1+|s|)^C$;
\item 
At any nonarchimedean place $v$ such that
$\pi_{1}$ and $\pi$ are both unramified,
both $\varphi$ and $\varphi_1$ are 
 invariant by $\PGL_2(\order_{F_v})$;
 \item $\|\varphi_1\|_{L^{\infty}}$ is $O(1)$. 
 \item   $ \|\varphi\|_{L^2(\quot)} \ll_{\epsilon}
\Norm(\cond)^{\epsilon}$.
\end{enumerate}
\end{lem}
\proof
We first choose local data. 
For each place where $e_{F,v}$ and $\pi_1$ are not ramified, we take $W_v$ (resp.
$W_{v,1}$) to be the  new vector in the Whittaker model of $\pi_v$ (resp
$\pi_{1,v}$). We put $\Psi_v$
to be the characteristic function of $\order_v^2$. 

Let $\mathcal{B}$ be the set of remaining $v$. 
For $v \in \mathcal{B}$, the assumptions show that $\pi_v$ is restricted to a bounded set. 
We choose $W_v$, $W_{v,1}, \Psi_v$ for $v \in \mathcal{B}$
according to Lem. \ref{lem:localchoice}.  Finally we choose $\varphi$ so that $\int_{x \in F \backslash \adele_F} e_F(x) \varphi(n(x) g) = \prod_{v} W_v(g)$, and similarly for $\varphi_1$,
and take $\Psi = \prod_{v} \Psi_v$. 
The first two assertions of the Lemma are immediate.

To bound the $L^2$ norm of $\varphi$, use  Lem. \ref{lem:localchoice} (\ref{b}), 
Lem. \ref{lem:l2norm}, and Iwaniec's bounds 
on $L$-functions near $1$. As for $\varphi_1$, it in fact belongs to a fixed finite set of cusp forms,
so the third assertion is immediate. 
\qed

We continue to keep $\pi$ an automorphic cuspidal representation
of $\PGL_2(\adele_F)$ with prime conductor. 
\begin{lem} \label{lem:rsoneeis}
 Suppose $\pi_{\infty}$ belongs to a bounded subset of $\widehat{\PGL_2(F_\infty)}$ (in what follows
the implicit constants may depend on this bounded subset). 

Let $t_0, t_0' \in \C$.  There exists a fixed finite set $\mathcal{F}$
of Schwarz Bruhat functions 
and a real number $C > 0$ so that:

There exist vectors $\varphi \in \pi, \Psi_1, \Psi_2 \in \mathcal{F}$ so that:
\begin{equation}\label{eq:losing} \Phi(t,t') = \Norm(\cond)^{1/2-t} 
\frac{
\int_{\quot} \varphi(g) E_{\Psi_1}(g, \frac{1}{2}+t) E_{\Psi_2}(g \f, \frac{1}{2}+t') dg}
{ \Lambda(\frac{1}{2} + t + t', \pi)
\Lambda(\frac{1}{2} + t -t', \pi)}
\end{equation}
is holomorphic and satisfies:
\begin{enumerate}
\item $|\Phi(t_0, t_0')| \gg 1$ and $|\Phi(t, t')| \ll (1+|t|+|t'|)^C C^{|\Re(t)| + |\Re(t')|}$. 
\item For any nonarchimedean place $v$, each $\Psi_1$ and $\Psi_2$
is invariant by $\PGL_2(\order_{F_v})$. 
\item$ \|\varphi\|_{L^2(\quot)} \ll_{ \epsilon}
\Norm(\cond)^{\epsilon}$.
\end{enumerate}
\end{lem}

\proof
The proof is similar to that of the previous Lemma; recall that (\ref{eq:localdef}) was valid as long as one of $\pi_1, \pi_2$ were cuspidal.

Let $\Psi_2'$ be the translate of the Schwarz function $\Psi_2$ by $\f$.
Then by (\ref{eq:twistedequivariant}), 
$$   E_{\Psi_2'} (s,g) = \Norm(\p)^{-s} E_{\Psi_2}^{a([\mathfrak{p}])} (s,g) $$
%E_{h . \Psi}(s,g) = |\det(h)|^{-s} (h \cdot E_{\Psi}(s,g)),  $$ 

Suppose $\Psi_1, \Psi_2$ factorize as $\prod_{v} \Psi_{1,v}, \prod_{v} \Psi_{2,v}$, and
define $W_{\Psi_{1,v}}(s, g)$ and $W_{\Psi_{2,v}}(s,g)$ as in Lem. \ref{lem:foureis}.
Suppose moreover that $\int_{x \in F \backslash \adele_F} e_F(x) \varphi(n(x) g)$
factorizes as $\prod_{v} W_v(g)$. Then we can express
the global integral of (\ref{eq:losing}) as a product in two different ways, depending
on whether we let $E_{\Psi_1}$ or $E_{\Psi_2}$ play the role of $\pi_2$. Namely, 
as in (\ref{eq:localdefiwasawa}):

\begin{multline} \label{eq:lmr} \int_{\quot} \varphi(g) E_{\Psi_1}(g, \frac{1}{2}+t) E_{\Psi_2}(g \f, \frac{1}{2}+t') dg 
\\ = c_F \Norm(\p)^{1/2+t'} \prod_{v} I_v(W_{v}, W_{\Psi_{1,v}}(1/2+t, \cdot), \Psi_{2,v}', 1/2+t') \\ = 
c_F \prod_{v} I_v(W_{v}, W_{\Psi_{2,v}}(1/2+t', \cdot)^{a([\mathfrak{p}])_v}, \Psi_{1,v} ,1/2+t)\end{multline}

Here $W_{\Psi_{2,v}}(1/2+t', \cdot)^{a([\mathfrak{p}])_v}$ denotes the translate
of $W_{\Psi_{2,v}}$ by the $v$th component of $a([\mathfrak{p}])$. 
%Here $c$ is the constant of (\ref{eq:localdefiwasawa}), depending only on the normalization of measure.

For $v$ nonarchimedean (notations being similar
to that of the previous Lemma) we take $\Psi_{1,v}$ and $\Psi_{2,v}$ to be the characteristic
function of $\order_v^2$ for every finite $v$, and $W_v$ to be the new vector.

For $v$ archimedean we first apply Lem. \ref{lem:localchoice} 
with $s_0 = 1/2+t_0$, and $\pi_{2,v}$ the representation of $\PGL_2(F_v)$ spanned by $E_{\Psi_2}(1/2+t_0',g)$,
i.e. the representation unitarily induced from the character $a(y)
\mapsto |y|_v^{it_0'} $. Lem. \ref{lem:localchoice} provides
$W_v$ in the Whittaker model of $\pi_v$, $W_{2,v}$ in the Whittaker model
of $\pi_{2,v}$, and a Schwarz function $\Psi_{1,v}$ with $|I_v( W_v, W_{2,v}, \Psi_{1,v}, 1/2+t_0)| \geq 1$. 
The last comment of Rem. \ref{rem:whiteis} shows that there is a standard $\Psi_{2,v}$
so that $W_{\Psi_{2,v}}(1/2+t_0', g_v) = W_{2,v}(g_v)$ (notation of Lem. \ref{lem:foureis}). 
Moreover, Lemma \ref{lem:localchoice} also shows that $\Psi_{1,v}$ and $W_{2,v}$ (so also $\Psi_{2,v}$) may be chosen from a fixed finite
set of possibilities (depending, of course, on the original bounded set to which $\pi_{\infty}$ belongs,
as well as $t_0$ and $t_0'$). 

Again we put $\Psi_{i}= \prod_{v} \Psi_{i,v}$ for $i=1,2$ 
and take $\varphi$ with $\int_{x \in F\backslash \adele_F} \varphi(n(x) g) = \prod_{v} W_v(g)$.
From (\ref{eq:lmr}) we deduce that, with our choices, $|\Phi(t_0, t_0') | \gg 1$. 
The assertion about $\|\varphi\|_{L^2}$ follows as in the proof of the previous Lemma. 
The second assertion of the Lemma (concerning invariance of $\Psi_1, \Psi_2$) is immediate. 

It remains to prove that $\Phi$ is actually holomorphic in $(t,t')$ and that $|\Phi(t, t')| \ll (1+|t|+|t'|)^C e^{C |\Re(t)| + C|\Re(t')|}$. Put $\Xi_v = \frac{I_v}{L_v(\frac{1}{2} + t + t', \pi_v)
L_v(\frac{1}{2} + t -t', \pi_v)}$. It is simple to explicitly compute $\Xi_v$ for nonarchimedean $v$, 
using Cor. \ref{cor:unramcase} and Lem. \ref{lem:dooku}. One thereby sees that it will suffice 
to check, by similar arguments to those used in Lem. \ref{lem:localchoice}, the following statement for $v$ archimedean: $\Xi_v = c^s c'^{s'} P(s,s')$, where $P$ is a polynomial,
 and moreover $c,c',P$
vary continuously in $\pi_v$, if $\pi_v \mapsto W_{v}$ is a continuous assignment.  
We only sketch the proof of this. 
From (\ref{eq:localdefiwasawa}) and the fact that $W_v, \Psi_{1,v}, \Psi_{2,v}$ are all $K_v$-finite, 
it suffices to  prove the corresponding assertions for $\int_{F_v^{\times}} W_v(a(y)) W_{\Psi_{1,v}}(s, a(y)) |y|^{s'-1} d^{\times}y$.  For this we use Barnes' formula as in \cite{jacquet}. \qed

\subsection{Local Rankin-Selberg convolutions.} \label{subsec:localrankinselberg}
Let $v$ be a nonarchimedean place of $F$ with residue characteristic $q_v$.  Let $\pi_1,\pi_2$ be generic irreducible admissible representations
of $\mathrm{GL}(2,F_v)$ with trivial central character. (Since we shall work purely locally over $F_v$
throughout the present subsection, we shall use the notation $\pi_1$ rather than e.g. $\pi_{1,v}$).

Then $\pi_1, \pi_2$ are self-dual.
We assume that $\pi_2$ is unramified and $\pi_1$
has conductor $q_v$, and denote by $L(s,\pi_j)$ the local $L$-factors.

Fix once and for all an additive unramified
character $\psi$ of $F_v$.
Let $v : F_v^{\times} \rightarrow \mathbb{Z}$ be the valuation,
put $\order_{F_v} = \{x \in F_v^{\times}: v(x) \geq 0\}$,
and choose a uniformizer $\unift \in F_v^{\times}$.
Let $\order_{F_v}^{\times}$ be the multiplicative group of units in
$\order_{F_v}$.
Let $d^{\times}x, dx$ be Haar measures on $F_v^{\times}, F_v$ respectively,
assigning mass $1$ to $\order_{F_v}^{\times}$ and $\order_{F_v}$ respectively.
For $x \in F_v$, put $n(x) = \left(\begin{array}{cc} 1 & x \\0 & 1
\end{array}\right)$.
Also let $w = \left(\begin{array}{cc}  0 & 1 \\ -1 & 0  \end{array}\right).$
We choose a Whittaker model for $\pi_1$ transforming
by the character $n(x) \mapsto \psi(x)$, and a Whittaker model
for $\pi_2$ transforming by the character $n(x) \mapsto \overline{\psi(x)}$.

Let $\Psi_v$ be the characteristic function of $\mathfrak{o}_{F_v}^2$.
Set $W_1$ to be the new vector in the Kirillov model of $\pi_1$,
let $W_2^{*}$ be the new vector in the Kirillov model of $\pi_2$,
and set $W_2 =\pi_2(\left(\begin{array}{cc} 1 & 0 \\ 0 & \unift
\end{array} \right))  W_2^{*}$.
Then both $W_1, W_2$ are invariant by the subgroup

\begin{equation} \label{eqk} K_0 = \{\left(\begin{array}{cc} a & b \\c & d \end{array} \right):
a,b,d \in \order_{F_v}, c \in \unift \order_{F_v}\}.\end{equation}
Moreover $W_2$ is invariant by $n(\unift^{-1} \order_{F_v})$.

\begin{lem} \label{lem:rs}
Notations being as in (\ref{eq:localdef}),
let $L(s,\pi_1 \times \pi_2)$ be the local $L$-factor. Then:
$$\frac{I_v(W_{1,v}, W_{2,v}, \Psi_v)}{L(s,\pi_1 \times \pi_2)} = \pm \frac{q_v^{s}}{q_v+1}$$
\end{lem}
\proof We shall often use the following shorthand:
if $W$ is a function in the Whittaker model
of $\pi \in \{\pi_1, \pi_2\}$
and for $z \in F_v^{\times}$, we write $W(z)$ for
$W( \left(\begin{array}{cc} z & 0 \\ 0 & 1 \end{array} \right))$.
Thus, for instance, $W_2(z) = W_2^{*}(z \unift^{-1})$.
The function  $z \mapsto W(z)$ belongs to the Kirillov
model of $\pi$.

Let $\epsilon \in \{-1,1\}$ be the local root number of $\pi_1$
(it lies in $\{-1,1\}$ since $\pi_1$ is self-dual).  Then:
\begin{equation} \begin{aligned}
\label{eqn:fe1} \int_{a \in F_v^{\times}} W_1(a) |a|^{s-1/2} d^{\times} a
= L(s, \pi_1) , \\ \int_{a \in F_v^{\times}} W_2(a) |a|^{s-1/2}
d^{\times} a  = q_v^{-(s-1/2)} L(s, \pi_2),
\end{aligned}\end{equation}
as follows from defining properties of newforms and the fact $W_2(z) = W_2^{*}(z \unift^{-1})$;
 moreover
\begin{equation} \begin{aligned}\label{eqn:fe2}
\int_{a \in F_v^{\times}} \pi_1(w) W_1(a) |a|^{s-1/2} d^{\times} a
= \epsilon q_v^{ (s-1/2)} L(s, \pi_1), \\
\int_{a \in F_v^{\times}} \pi_2(w) W_2(a) |a|^{s-1/2}
d^{\times} a = q_v^{ (s-1/2)} L(s, \pi_2),\end{aligned}\end{equation}
as follows from local functional equation for
the standard $L$-function on $\mathrm{GL}(2)$:
see \cite[2.18]{jacquetlanglands}.
\footnote{
That is,  $\int_{a \in F^{\times}} W(a)
|a|^{s-1/2} d^{\times}a
 = \frac{L(s,\pi)}{\epsilon(s,\pi) L(1-s, \tilde{\pi})}
\int_{F^{\times}} W(aw) |a|^{1/2-s} \omega^{-1}(a)
d^{\times}a$.
In particular, if $\pi$ is a representation with
trivial central character, and $\chi$ a character of $F^{\times}$,
$\int_{a \in F^{\times}} W(a)
\chi(a) d^{\times}a
 = \frac{L(\frac{1}{2},\pi \otimes \chi)}{\epsilon(\frac{1}{2},\pi \otimes
\chi) L(\frac{1}{2}, \tilde{\pi}
\otimes \bar{\chi})}
\int_{F^{\times}} W(aw) \chi^{-1}(a) d^{\times}a.$}

Note moreover that $W_1$, $W_2$, $\pi_1(w) W_1$, and $\pi_2(w) W_2$
are all invariant by the maximal compact subgroup of the diagonal
torus of $\GL_2$.  Thus (\ref{eqn:fe1})
and (\ref{eqn:fe2}) completely determine
their restriction to the diagonal torus; we now
explicate this.

Choose $\alpha \in \C$ so that $L(s,\pi_1) = (1 - \alpha q_v^{-s})^{-1}$.
In fact,
$\alpha = -\epsilon q_v^{-1/2}$,
by \cite[Prop. 3.6]{jacquetlanglands}.
Choose $\gamma_1, \gamma_2$
so that $L(s,\pi_2) = ((1-\gamma_1 q_v^{-s}) (1-\gamma_2 q_v^{-s}))^{-1}$.
Recalling the notational convention established
in the paragraph prior to (\ref{eqn:fe1}), we see:
\begin{equation}\begin{aligned}\label{eqn:familyvalues}
W_1(\unift^r) = \begin{cases}  \alpha^r q_v^{-r/2}, r \geq 0 \\
0, r < 0 \end{cases}, \ \ \pi_1(w) W_1(\unift^r) = \begin{cases}
\epsilon \alpha^{r+1}
q_v^{-\frac{r+1}{2}}, r \geq -1\\ 0, r < -1 \end{cases} \\
 W_2(\unift^r) =
\begin{cases} 0, r \leq 0 \\ 1, r = 1\\
(\gamma_1^{r-1} + \gamma_1^{r-2} \gamma_2 + \dots + \gamma_2^{r-1})
q_v^{-\frac{{r-1}}{2}}, r \geq 2
\end{cases} \\
\pi_2(w) W_2(\unift^r) =
\begin{cases} 0, r < -1 \\ 1, r=-1 \\
(\gamma_1^{r+1} + \gamma_1^{r} \gamma_2 + \dots + \gamma_2^{r+1})
q_v^{-\frac{r+1}{2}},
r \geq 0
\end{cases}
\end{aligned}
\end{equation}
The local integral we wish to evaluate is the
right hand side of (\ref{eq:localdef}). In the case at hand,
with $N,G, Z$ denoting the $F_v$-points of the respective groups,
we have:
\begin{equation} I(s) :=
%\int_{N \backslash G} W_1(g)  W_2(g) \Psi((0,1) \cdot g)
%|\mathrm{det}(g)|^{s}
%dg  \\ =
\int_{Z N \backslash G}
W_1(g) W_2(g) |\det(g)|^s \left(
\int_{t} \Psi_v((0,t) \cdot g) |t|^{2s} d^{\times} t\right) dg\end{equation}

Using the Iwasawa decomposition, and recalling $\Psi_v$
was the characteristic function of $\order_v^2$ one finds:
$$I(s) = (1-q_v^{-2s})^{-1}
\int_{A \times K_v} \pi_1(k) W_1(a)\pi_2(k) W_2(a) |a|^{s-1} d^{\times}
a dk,  $$
where the measure $dk$ is the Haar measure of total mass $1$,
and $d^{\times}a$ assigns mass $1$ to $A \cap K_v$.

The function $k \mapsto \pi_1(k) W_1(a) \pi_2(k) W_2(a)$
is right invariant by $K_0$ (see \eqref{eqk} for definition) and left invariant by $N \cap K_v$.
There are two $( N \cap K_v, K_0)$ double cosets in $K_v$,
and we may therefore express $I(s)$ as a sum:
\begin{equation} \begin{aligned}
(1-q_v^{-2s}) I(s) = \frac{1}{q_v+1}
\int_{F_v^{\times}} W_1(a) W_2(a) |a|^{s-1} d^{\times} a
\\ +
\frac{q_v}{q_v+1}
 \int_{a \in F_v^{\times}}
\pi_1(w) W_1(a)
\pi_2(w) W_2(a) |a|^{s-1} d^{\times} a \\
\end{aligned}\end{equation}
To evaluate $I(s)$, we use (\ref{eqn:familyvalues}).
Noting that $L(s, \pi_1 \times \pi_2) =\frac{1}{
(1-\alpha \gamma_1 q_v^{-s})
(1-\alpha \gamma_2 q_v^{-s})}$,
an easy computation shows
$$I(s) = \frac{L(s,\pi_1 \times \pi_2)}{(q_v+1)(1-q_v^{-2s})}
\left(\alpha q_v^{-1/2} q_v^{-(s-1)}
+ \epsilon q_v q_v^{s-1}\right)$$
from where we obtain $I(s) = \epsilon \frac{q_v^s}{q_v+1} L(s, \pi_1 \times \pi_2)$.
Note also that $I(s)$ satisfies the necessary functional equation.
\qed

\subsection{Hecke-Jacquet-Langlands integral representations for
standard $L$-functions} \label{subsec:hjl}
Our goal here is to prove Prop. \ref{prop:hjlcusp} and \ref{prop:hjleis}, used in the text. 
This amounts to explicit computations connected to Hecke-Jacquet-Langlands
integral representations. Since, in the main text, we obtain subconvexity for $\GL(1)$ twists
of $\GL(2)$ $L$-functions, {\em with polynomial dependence in all parameters},
we will have to be somewhat more precise than in the case of Rankin-Selberg $L$-functions.

Let $\pi$ be a cuspidal representation of $\GL_2$ over $\adele_F$. Let $\chi$ be a unitary
character of $\adele_F^{\times}/F^{\times}$ of finite conductor
$\fcond$. Put $L_{unr}(s,\pi \times\chi)$ to be the unramified part of the (finite) standard $L$-function:
$$L_{unr}(s, \pi \times \chi):= \prod_{v \, \mathrm{finite}, \chi_v \, \mathrm{unramified}}L_v(s, \pi_v \times \chi_v).$$
Define $\mu_z$ as in (\ref{eq:nudef}), i.e. the measure on $\quotg$ defined as
$$\mu_z(f) = \int_{|y| =z} f(a(y)
n([\fcond])) \chi(y) d^{\times} y.$$
We refer to Sec. \ref{subsec:NF} and Sec. \ref{sec:gl2pgl2} for notation, 
as well as the start of Section \ref{sec:torus1} for a discussion of the meaning of $\mu_z$ in classical terms. 

\begin{lem} \label{lem:hjllocal}
Let $v$ be a nonarchimedean place of $F$ with residue
characteristic $q_v$, and $\pi_{v}$ an irreducible generic
representation of $\GL_2(F_v)$. Let $\psi_v$ be an unramified additive
character of $F_v$.
Let $\chi_v: F_v^{\times} \rightarrow \C$ a
multiplicative character of conductor $r$, $W_v$ be the new vector
in the $\psi_v$-Whittaker model of $\pi_v$. Then
$$\int_{y \in F_v^{\times}} W_v(a(y) n(\unif_v^{-r})) \chi_v(y)
|y|^{s-1/2} d^{\times}y =   \begin{cases} L_v(s, \pi_v
\times \chi_v), \, r= 0.\\
\theta, \, r \geq 1,\end{cases}$$ where $\theta$ is a scalar of
absolute value $q_v^{-r/2} (1-q_v^{-1})^{-1}$.
\end{lem}

\proof
If $r=0$, then $\chi_v$ is unramified, the result follows immediately
from the definition of the new vector. Otherwise, $\chi_v$ is ramified, and we
 rewrite the integral under consideration as
\begin{equation} \label{eq:kaligiyunte}
\int_{y \in F_v^{\times}} W_v(a(y)) \psi_v(\unif_v^{-r} y) \chi_v(y) |y|^{s-1/2}
d^{\times}y.\end{equation} Now $W_v(a(y))$ vanishes when $v(y) <
0$ and it is $\order_{F_v}^{\times}$-invariant. The integral
$\int_{v(y) = k} \chi_v(y) \psi_v(\unif_v^{-r} y) d^{\times}y$ is
nonvanishing only when $k=0$. In that case, it is a Gauss sum with
absolute value $\frac{q_v^{-r/2}}{ (1-q_v^{-1})}$, 
where the factor $(1-q_v^{-1})^{-1}$ arises from the measure normalization (cf.
Sec. \ref{subsec:measures}) namely
$\int_{v(y) = 0} d^{\times}y = 1$. The result follows. \qed

\begin{lem} \label{lem:callas}
Let $d,\beta \geq 0$. Then there exists $\varphi \in \pi$ such that, with
\begin{equation} \label{eq:tomyknees} \Phi(s) = \Norm(\fcond)^{1/2} \frac{\int_z \mu_z(\varphi) |z|^{s-1/2} d^{\times}z}{L_{unr}(s, \pi \times \chi)}\end{equation}
then $\Phi(s)$ is holomorphic and satisfies:
\begin{enumerate}
\item $ |\Phi(s)| \ll_{\Re(s),\epsilon} \Norm(\fcond)^{\epsilon}$ and $|\Phi(\frac{1}{2})| \gg_{\epsilon} \Norm(\fcond)^{-\epsilon}$. 
%\begin{equation} \left|\int_z \mu_z(\varphi) |z|^{s-1/2} d^{\times}z \right| \ll_{\Re(s)} 
%\Norm(\fcond)^{-1/2+\epsilon}  L_{unr}(s) 
%\end{equation}
%\item 
%\begin{equation} \label{eq:john2}
% \left|\int_{z} \mu_z(\varphi) d^{\times}z \right| \gg_{\epsilon}
%\Norm(\fcond)^{-1/2-\epsilon} |L(\frac{1}{2}, \pi \otimes \chi)|\end{equation}
\item
$\varphi$ is new at every finite place (i.e.,
for each finite prime $\q$ it is invariant by $K_0[\q^{s_{\q}}]$, 
where $s_{\q}$ is the local conductor of the local constituent $\pi_{\q}$). 
\item The Sobolev norms of $\varphi$ satisfy the bounds
(conductor notation as in Sec. \ref{iscond})
\begin{equation} \label{eq:john17}
S_{2,d,\beta}(\varphi) \ll_{\epsilon} \Cond_{\infty}(\pi)^{2d+\epsilon}
\Cond_{f}(\pi)^{\beta+\epsilon}
 \Cond_{\infty}(\chi)^{1/2 + 2d}
\end{equation}
\end{enumerate}
\end{lem}

\proof
For each infinite place $w$ of $F$,
denote by $\Cond_w(\chi)$ the contribution from $w$
to the Iwaniec-Sarnak analytic conductor of $\chi$
(see Sec. \ref{iscond}.)

The map $\varphi \mapsto W_{\varphi} =\int_{F \backslash \adele_F}
e_F(x) \varphi(n(x) g)$ is an isomorphism between the space of
$\pi$ and the Whittaker model of $\pi$. For each finite $v$, take
$W_v$ to be a new vector in the Whittaker model of $\pi_v$. A
point of caution is that $e_F$ may not be unramified on $F_v$; to
be absolutely concrete, we set $W_v(g) =
W_{v,\mathrm{new}}(a(\unif^{d_v}) g)$, where $W_{v, \mathrm{new}}$
is the new vector in the Whittaker model of $\pi_v$ taken w.r.t an
unramified additive character of $F_v$, and $d_v = v(\diff)$ is the local valuation of the different. 

Let us now choose $W_v$ at the infinite places.  Let $g_1$ be a smooth positive function of compact support on $F_v$. Let $\deg(v) = 2$ if $v$ is complex and $\deg(v)= 1$ if $v$ is real.
 For $\infty|v$, define 
 $$W_v(y) = \Cond_v(\chi) g_1(\Cond_v(\chi)^{1/\deg(v)} (y-1)).$$
  This is possible by the theory of the Kirillov model;
thus $W_v$ is a smooth (but not $K_v$-finite) vector.
In words, if $v$ is real, the function $W_v$ is supported in a 
neighbourhood of the identity of size $\Cond_v(\chi)^{-1}$ and takes values of
size $|\Cond_v(\chi)|$ there; if $v$ is complex, a similar
statement holds but now $W_v$ is supported in a disc around the identity
with area $\Cond_v(\chi)^{-1}$. 

Then there exists $\varphi \in \pi$ with $W_{\varphi} = \prod_{v} W_{v}$.
 By unfolding,  it follows that for $\Re(s) \gg
 1$
\begin{equation} \label{eq:com} \int_{z} \mu_z(\varphi) |z|^{s-1/2} d^{\times} z\\  = c_F \prod_{v}
\int_{y \in F_v^{\times}} W_v(a(y) n([\fcond]))
|y|^{s-1/2} \chi_v(y) d^{\times}y
\end{equation}
Here $c_F$ is a constant depending only on $F$, arising from change of measure; it is 
entirely unimportant as we will be only interested in bounds.\footnote{
The measure $\mu_z$ is normalized as a probability measure, whereas to unfold from $\adele_F^{\times}$
to $\prod_{v} F_v^{\times}$ we use the measures previously set up there (see Section \ref{subsec:measures}). } 

By Lem. \ref{lem:hjllocal}, with $\Phi(s)$ as in the statement of the Lemma, 
\begin{equation} \label{eq:phi} \Phi(s) = c_F \theta'  \cdot \Norm(\diff)^{s-1/2}
\prod_{\mathrm{infinite}\, v} \int_{F_v^{\times}} W_v(a(y)) |y|^{s-1/2} \chi_v(y)
\end{equation}
where $|\theta'| = \prod_{\q|\fcond} (1-\Norm(\q)^{-1})^{-1}$. For this choice of $\varphi$,
the second assertion if the Lemma is clear, and, if we
choose the support of $g_1$ to be small enough, the first  assertion also follows easily. 
\footnote{For the assertion concerning the lower bound for $|\Phi(1/2)|$, the 
point, in words, is that our choices are so that $\chi_v$ does not oscillate over the support
of $W_v$, cf. Remark \ref{rem:explication}. Note how convenient it is, here and elsewhere, to use 
smooth vectors rather than $K_{\infty}$-finite vectors; one could not e.g. achieve $W_v$
of compact support with $K_{\infty}$-finite vectors.}

 (\ref{eq:john17}) follows from
Lem. \ref{lem:sobnormk}, together with Lem. \ref{lem:l2norm}
and the upper bound for $L$-functions near $1$ due
to Iwaniec. See \cite[Chapter 8]{iwaniec} for this bound. \qed

The previous Lemma shows that $L(1/2, \pi \times \chi)$ may be ``well-approximated''
by an appropriate period integral. Unfortunately, this period integral is against
a measure of infinite mass, since $\adele_F^{\times}/F^{\times}$ is of infinite volume.
It is, therefore, convenient to know that the $\mu_z$-integral of (\ref{eq:tomyknees}) can be truncated
to a compact range without affecting the answer too much. This is, roughly speaking,
the geometric equivalent of the approximate functional equation in the classical theory,
and is provided by the next Lemma. It says, roughly speaking, that the integral of (\ref{eq:tomyknees})
can be truncated to the range where $z$ is around $\Norm(\fcond)^{-1}$. 

\begin{lem} \label{lem:callas2}
Let notation be as in Lem. \ref{lem:callas}.
Let $g_{+}, g_{-}$ be positive smooth functions on $\mathbb{R}_{\geq 0}$
such that $g_{+} + g_{-} = 1$, $g_{+} (t)= 1$ for $t \geq 2$
and $g_{-}(t) = 1$ for all $t \leq 1/2$.  
Then
$$I_{+} := \int_{z} \mu_z(\varphi) g_{+}(z/T)d^{\times}z  \ll_{g_{+},\epsilon} (\Norm(\fcond) T)^{-1/2}
(T \Cond(\pi) \Cond(\chi))^{\epsilon}$$
$$  I_{-} := \int_{z} \mu_z(\varphi) g_{-}(z/T) d^{\times} z \ll_{g_{-},\epsilon} (\Norm(\fcond) T)^{1/2}
(T \Cond(\chi))^{\epsilon}(\Cond_{\infty}(\chi) \Cond(\pi))^{1+\epsilon} $$
\end{lem}
\proof
Recall the definition of $\mu_z$ from (\ref{eq:nudef}). 
Put $\widehat{g_\pm}(s) = \int g_{\pm}(x) x^{s-1}  dx$, the Mellin transform of $g_{\pm}$;
then $g_{\pm}$ is holomorphic in $\pm \Re(s) < 0$ and for any $M \geq 0, \pm \sigma <0$
the integral $\int_{\Re(s) =  \sigma} |\widehat{g_\pm}(s)| (1+|s|)^M ds$ is convergent. 
Then, for any $\pm \sigma> 0$, we have, by the Plancherel formula on $\mathbb{R}^{\times}$, that: $$\int_{z} \mu_z(\varphi) g_{\pm}(z/T) d^{\times} z = \frac{1}{2 \pi i} \int_{\Re(s) = -\sigma} \left(\int \mu_z(\varphi) |z|^{-s} d^{\times} z \right) T^s \widehat{g_\pm}(s) ds.$$
So for any $M >0$, 
 $$|I_{\pm}| \ll_{\sigma, g_{\pm}, M} T^{-\sigma} \Norm(\fcond)^{-1/2} \sup_{\Re(s) = 1/2+\sigma}
\frac{|L_{unr}(s) \Phi(s)|}{(1+|s|)^M}$$ where $\Phi$ is as in the previous Lemma.Take $\sigma = 1/2 + \varepsilon$
in the $+$ case, $-1/2-\varepsilon$ in the $-$ case. 

Using Iwaniec's bounds on $L$-functions near $1$ \cite[Chapter 8]{iwaniec} and the functional equation,
we see that for sufficiently large $M$:
\begin{multline} \label{eq:vasu}
\sup_{\Re(s) = 1 + \varepsilon} |L_{unr}(s, \pi \times\chi)| \ll \Cond(\pi \otimes \chi)^{\varepsilon} \\
\frac{\sup_{\Re(s) = - \varepsilon }| L_{unr}(s, \pi \times \chi) |}{(1+|s|)^M}  \ll_M
\Cond(\pi \otimes \chi)^{1/2+\varepsilon}
\prod_{\chi_v \, \mathrm{ramified}\, \mathrm{finite}} \sup_{\Re(s) = - \varepsilon} | L_v(s, \pi_v \times \chi_v)|^{-1}
\end{multline}
For each $v$ where $\chi_v$ is ramified and $L_v(s, \pi_v \times \chi_v)$ is not identically $1$, 
the representation $\pi_v$ must also be ramified (i.e., not spherical). 
So one can bound the product on the second line on (\ref{eq:vasu}), using
trivial bounds towards the Ramanujan conjecture,  by $\Cond(\pi)^{1/2+2 \varepsilon}$. 
The fact that $\Cond(\chi) = \Cond_{\infty}(\chi) \Norm(\fcond)$, the bound \cite{BK} $\Cond(\pi \otimes \chi) \ll \Cond(\pi) \Cond(\chi)^2$, and the (easily verified) analogue of this bound of \cite{BK} at archimedean places, allows one to conclude. 
\qed
%
%\begin{lem} \label{lem:hjleis}
%Let $h, E_{h,\Psi}$ be defined as in Lem. \ref{eislinfty}. Then
%there exists a choice of Schwarz function $\Psi$, depending only
%on $\chi_{\infty}$,  such that $E_{h,\Psi}$ is $\Kmax$-invariant,
%and moreover:
%$$
%\left|\int_{\Re(s)=\sigma} h(s) \Lambda(\chi,s) \Lambda(\chi,1-s)
%ds \right| \ll \Norm(\fcond)^{1/2} \left|\int_{y \in
%\adele_F^{\times}/F^{\times}} E_{h,\Psi}(a(y)n ([\fcond])) \chi(y)
%d^{\times}y\right| $$
%\end{lem}

We now address the analogue of the previous Lemmas when $\pi$ is noncuspidal. \footnote{The content of the following Lemma, in classical language, is related to the following observation. Let $\chi$
be an even Dirichlet character mod $q$, 
and let $E^*(s,z)$ to be the Eisenstein series of (\ref{eq:four}), and
$\bar{E}^*(s,z) := E^*(s,z) - \xi(2s) y^s - \xi(2(1-s)) y^{1-s}$, then 
$\frac{1}{q} \int_{0}^{\infty} \sum_{1 \leq x \leq q-1} \bar{E}^*(s, \frac{x}{q}+iy) y^{s'} d^{\times}y$
coincides, up to some harmless factor, with $q^{-1/2} \Lambda(\chi, s+s') \Lambda(\chi, 1-s+s'))$, 
where $\Lambda(\chi,s)$ is the usual Dirichlet $L$-function completed to include the $\Gamma$-factor at $\infty$.  This particular expression is actually not quite suitable for our needs, because of the rapid decay of the $\Gamma$-factor swamps information about the finite $L$-function, and in fact the Lemma uses (the equivalent of) a different test vector
belonging to the automorphic representation underlying $E^*(z,s)$. }

\begin{lem} \label{lem:hjleis} 
Let $s_0, s_0' \in \C$. 
There is an absolute $C > 0$ (i.e., depending only on $F$)
and a Schwarz function $\Psi$ (depending on $\chi$) so that if we put
$$ \Phi(s,s') := \Norm(\fcond)^{1/2} \frac{ \int_{y \in
\adele_F^{\times}/F^{\times}} \bar{E}_{\Psi}(s,a(y)n ([\fcond])) \chi(y) |y|^{s'}
d^{\times}y}{ L(\chi,s+s') L(\chi,1-s+s')}$$
where $\bar{E}$ is defined as in (\ref{eq:baredef}),  then the integral
defining $\Phi$ is absolutely convergent in a right half-plane $\Re(s) \gg 1$. Moreover, 
$\Phi$ extends from $\Re(s), \Re(s') \gg 1 $ to a holomorphic
function on $\C^2$,  satisfying
\begin{enumerate}
\item  $|\Phi(1/2,0)| \gg 1$ and $|\Phi(s,s')| \ll C^{1+|\Re(s)| + |\Re(s')|} (1+|s|+|s'|)^C$. 

Moreover, given $N > 0$ we have that \begin{equation} \label{spd2}  |\Phi(s,s')|
(1+|s|+|s'|)^{N} \ll_{\Re(s), \Re(s'), N} \Cond_{\infty}(\chi)^{N'} \end{equation}where $N'$ and the implicit constant may be taken
to depend continuously on $N, \Re(s), \Re(s')$.

\item $\Psi$, and so also
$E_{\Psi}(s,g)$ is invariant by $\Kmax$;
\item Let $h \in \mathcal{H}(\kappa)$ be as in (\ref{eq:normdef}), 
and put $E_h := \int_{\Re(s) \gg 1} h(s) E_{\Psi}(s, g) dg$. 
For each $d,\beta$ there is $N $ such that $S_{\infty,d,\beta}(E_h) \ll_{\kappa} \|h\|_0 \Cond_{\infty}(\chi)^{N}$.
% (the additional
%content of this over Lem. \ref{eislinfty} is in removing the dependence of the implicit constant on $\Psi$.)
\end{enumerate}
\end{lem}
\proof
We shall not explicitly address details of convergence.
The manipulations that follow may be justified
by similar reasoning to that of Lem. \ref{eislinfty}.

We now define a Schwarz function $\Psi_v$ on $F_v^2$ for each
place $v$. For each finite place $v$, let $\Psi_v$ be the characteristic
function of $\order_v^2$.

For infinite $v$, we will first define a Schwarz function $\rho_v$ on $F_v$, and then take
$\Psi_v(x,y) = \rho_v(x) \widehat{\rho_v}(y)$; here $\widehat{\rho_v}$ is the inverse Fourier transform
of $\rho_v$, satisfying $\int_{F_v} \widehat{\rho_v}(y) e_{F_v}(xy) dy =\rho_v(x)$. 

 Let $g_1$ be a smooth positive function of compact support on $F_v$. Let $\deg(v) = 2$ if $v$ is complex and $\deg(v)= 1$ if $v$ is real.
 For $\infty|v$, define 
 $$\rho_v(y) = \Cond_v(\chi) g_1(\Cond_v(\chi)^{1/\deg(v)} (y-1)).$$\
  In words:
in the real (resp. complex) case, $\rho_v$ is localized in a real (resp. complex)
interval (resp. disc) around $1$, of length (resp. area) $\Cond_v(\chi)^{-1}$.
Now put $\Psi_v(x,y) = \rho_v(x)
\widehat{\rho_v}(y)$.  
The function $\Psi_v$ is not compactly supported; however,
it is of rapid decay. 
Indeed for each Schwarz norm $\mathcal{S}$, there is $M >0$ such that
\begin{equation} 
\label{eq:psibound} \mathcal{S}(\Psi_v) \ll \Cond_v(\chi)^M 
\end{equation}

Define a Schwarz function on $\adele_F^2$ via $\Psi(x,y) =
\prod_{v} \Psi_v(x,y)$.
 Define $W_{\Psi}(s,g)$ as in
Lem. \ref{lem:foureis} to be the Fourier coefficient of $E_{\Psi}(s,g)$.  The choice of $\Psi$ and Lem. \ref{lem:foureis} shows that $W_{\Psi}(s,g) = \prod_{v} W_v(g)$, where, 
for each finite $v$,  $W_v$ is given by Cor. \ref{cor:unramcase}, and satisfies
 \begin{equation}\int_{F_v^{\times}} W_v(a(y))
|y|^{s'} d^{\times}y = q_v^{d_v(1+s'-s)} L_v(|\cdot|^{s}, s')
L_v(|\cdot|^{1-s}, s'),\end{equation}
 For
infinite $v$, 
$W_v$ satisfies (Lem. \ref{lem:foureis}) 
\begin{equation} \label{wvif} \int_{F_v^{\times}} W_v(a(y)) \omega(y) |y|^{s'} d^{\times}y 
= \int_{F_v^{\times}} \rho_v(x) \omega(x) |x|^{s+s'} d^{\times}x
\int_{F_v^{\times}} \rho_v(x) \omega(x) |x|^{1-s+s'}d^{\times} x\end{equation}

By Fourier analysis and Lem. \ref{lem:eisconstterm},  $\bar{E}_{\Psi}(s,g)=
\sum_{\alpha \in F^{\times}} W_{\Psi}(a(\alpha) g)$. Thus, for $\Re(s) \gg 1$,
\begin{multline}\int_{y \in \adele_F^{\times}/F^{\times}}
\bar{E}_{\Psi}(s,a(y) n([\fcond])) \chi(y) |y|^{s'} d^{\times} y \ =
% \sigma} h(s) E^N(s,
%a(y) n([\fcond])) \chi(y) d^{\times}y
%\\ + 
\int_{y \in \adele_F^{\times}} 
W_{\Psi}(s,a(y) n([\fcond])) \chi(y) |y|^{s'} d^{\times}y \\ = 
\prod_{v} \int_{y \in F_v^{\times}} W_v(a(y) n([\fcond])) \chi(y) |y|^{s'} d^{\times}y \end{multline}
%The first quantity on the right hand side is seen to be zero: for
%$\q | \fcond$, the function $E^N(s,a(y)n([\fcond]))$ is invariant
%on cosets $y \order_{F_{\q}}^{\times}$, whereas $\chi$ is ramified
%at $\q$. 

For $s \gg 1$,  we use (\ref{wvf}), (\ref{wvif}) and Lemma \ref{lem:hjllocal}  to evaluate the local factors, obtaining:
\begin{multline} \Phi(s,s')
 = \theta' \cdot \Norm(\diff)^{1+s'-s} \prod_{v \, \mathrm{infinite}} \int_{y \in F_v^{\times}} W_v(a(y)) \chi(y) |y|^{s'} d^{\times}y \\ =  \theta' \cdot \Norm(\diff)^{1+s'-s}
 \int_{F_v^{\times}} \rho_v(x) \chi(x) |x|^{s+s'} d^{\times}x
\int_{F_v^{\times}} \rho_v(x) \chi(x) |x|^{1-s+s'}d^{\times} x\end{multline}
 where $|\theta'| =  \prod_{\q|\fcond} (1-\Norm(\q)^{-1})^{-1}$. 
 Now, by choice of $\varphi_v$, the integral $I_v(s) : = \int_{y \in F_v^{\times}} \rho_v(y) \chi(y) |y|^{s} d^{\times} y$ satisfies $|I_v(1/2)| \gg 1$ and $|I_v(s)| \ll (1+|s|)^C C^{1+|\Re(s)|}$, at least
 when we choose the support of $g_1$ to be sufficiently small.   It also satisfies
 $ |I_v(s)| (1+|s|)^N \ll_{N, \Re(s)} \Cond_{\infty}(\chi)^{N'}$, where $N'$ and the implicit constant may be taken
 to depend continuously on $N, \Re(s)$. 
 
 The corresponding facts (i.e. the first assertion of the Lemma) about $\Phi$ follow immediately. 
 The second assertion of the Lemma is immediate from our choice of $\Psi$. 
 
As for the third and final assertion, 
it follows from Rem. \ref{rem:schwarz} 
and (\ref{eq:psibound}). 
%
%a study of the proofs shows that, to control
%the dependence of the implicit constant of Lem. \ref{eislinfty} on $\Psi$,
%it suffices to do the same in Lem. \ref{lem:fbound},
%\begin{equation} \label{eq:rad} s (1/2-s)  (1+ |\Im(s)|^N) f_{\Psi}(s,g) \height(g)^{-\Re(s)},\end{equation}
%and the similar expression with $\Psi$ replaced by $\widehat{\Psi}$. 
% As was remarked in the
%final sentence of the proof of that Lemma, if $\Psi = \prod_{v} \Psi_v$
%and the $\Psi_v$ are  held fixed ($v$ finite), then 
%(\ref{eq:rad}) is bounded by a Schwarz norm of $\prod_{\infty|v} \Psi$. 
%i.e. to obtain explicit bounds for the absolute value of
%when $|\Re(s) | \leq M$ for some explicit $M$.   It suffices to do this for $g \in K_{\infty}$, as one
%verifies from the Iwasawa decomposition of a general $g \in \GL_2(\adele_F)$. 
%This, and the definition of $f_{\Psi}$, show that it suffices to bound, for $k \in K_{\infty}$, the quantity
%$s (1/2-s) (1+|\Im(s)|^N) \zeta_F(2s)  \int_{F_{\infty}} \Psi((0,t) k)|t|^{2s} d^{\times}t$, 
%where $\zeta_F$ is the finite part of the Dedekind $L$-function attached to $F$. 
%Thus (since the Dedekind $L$-function has at most polynomial growth along vertical lines) 
%it suffices to bound the quantity
%\begin{equation} \label{eq:everyday}(1+|\Im(s)|)^{N'} \int_{t \in F_{\infty}^{\times}} \Psi((0,t) k ) |t|^{2s} d^{\times} t,  \ \ \ \ \ 
%|\Re(s)| \leq M\end{equation}
%

%Note that even the fact that the the right-hand side of (\ref{eq:everyday}) is holomorphic
%is not entirely formal: it implicitly uses the fact that $\Psi_v(x,y)$ vanishes
%in a neighbourhood of $(0,0)$. 
\qed

\begin{lem} \label{lem:callas3}
Let notations be as in the previous Lemma. 
Assume $\chi$ is ramified at at least one finite place. 
Let $g_{+}, g_{-}$ be positive smooth functions on $\mathbb{R}_{\geq 0}$
such that $g_{+} + g_{-} = 1$, $g_{+} (t)= 1$ for $t \geq 2$
and $g_{-}(t) = 1$ for all $t \leq 1/2$.  

Then \begin{equation}\label{stupidestimate} \mu_z(E_h) \ll_{K,\Psi,h} \min(z^K, z^{-K})\end{equation} for any $K \geq 1$. 

Moreover, there is an absolute $N >0$ such that $$I_{+} := \int_{z} \mu_z(E_h) g_{+}(z/T)d^{\times}z  \ll (\Norm(\fcond) T)^{-1/2}
(T \Cond(\chi) )^{\epsilon} \|h\|_N $$
$$  I_{-} := \int_{z} \mu_z(E_h) g_{-}(z/T) d^{\times} z \ll  (\Norm(\fcond) T)^{1/2}
(T \Cond(\chi) )^{\epsilon} \Cond_{\infty}(\chi)^{1+\epsilon}  \|h\|_N$$
(Here the norms $\|\cdot\|_N$ are as in (\ref{eq:normdef}).)
\end{lem}

\proof 
Again, we shall leave verification of convergence to the reader. 
Recall that, with the relevant measure on $\adele_F^{\times}/F^{\times}$ having mass $1$: \begin{multline} \label{bathsheba} \mu_z(E_h) = \int_{y \in \adele_F^{\times}/F^{\times}, |y|=z} E_h(a(y) n[\fcond]) \chi(y) d^{\times}y
\\ = \int_{y \in \adele_F^{\times}/F^{\times},|y|=z} \int_{\Re(s) \gg 1} h(s) E_{\Psi}(s, a(y) n([\fcond])) \chi(y) d^{\times}y \\ =
 \int_{y \in \adele_F^{\times}/F^{\times},|y|=z} \int_{\Re(s) \gg 1} h(s) \bar{E}_{\Psi}(s, a(y) n([\fcond])) \chi(y) d^{\times}y  \end{multline}
Here, the last equality is justified by the fact that $(E_{\Psi}-\bar{E}_{\Psi})(s, a(y) n([\fcond]))$
is invariant under $y \mapsto y y'$, for $y' \in \prod_{v} \order_v^{\times}$. On the other hand, 
$\chi$ is nontrivial on $\prod_v \order_v^{\times}$, by assumption. 

Combining \eqref{bathsheba} with Lem. \ref{lem:hjleis}, we have
\begin{equation} \label{wiw} \int_{z} \mu_z(E_h) |z|^{s'} d^{\times}z  =  c_F \Norm(\fcond)^{-1/2} \int_{\Re(s) \gg 1} h(s) L(\chi, 1-s+s') L(\chi, s+s') \Phi(s,s') ds.\end{equation}
Here $c_F$ is an (unimportant) constant arising from measure normalization, as in \eqref{eq:com}. 

The assertion \eqref{stupidestimate} follows immediately from this, inverse Mellin transform,
and analytic properties of the right-hand side. 

Now proceed as in Lem. \ref{lem:callas2}; it follows that (for any $M$)
$$|I_{\pm}| \ll T^{\mp (1/2+\varepsilon)} \Norm(\fcond)^{-1/2} \sup_{\Re(s') = \pm (1/2+\varepsilon)} (1+|s'|)^{-M} \int 
h(s) L(\chi, 1-s+s') L(\chi, s+s') \Phi(s, s') ds.$$
We deal with the case of $I_{-}$. In that case, we take
the inner integral to be over $\Re(s)= 1/2$, and put $s' = -1/2-\varepsilon - it'$, and it will suffice to bound
$\int h(1/2+it) L(\chi, -\varepsilon - it - it') L(\chi, -\varepsilon + it - it') (1+|t|+|t'|)^C$.
This is bounded, up to an implicit constant depending on $\varepsilon$, by $\Cond(\chi)^{1+ 2 \varepsilon}  \|h\|_{M'} (1+|t'|)^{C'}$
for sufficiently big $M', C'$, whence the result. \qed

\end{document}